\documentclass[oneside]{amsart}
\pdfoutput=1

\newcommand\thmcounter{section}
\PassOptionsToPackage{final}{graphicx}
\usepackage{kazhdan--lusztig}
\title[A generalized Deodhar setting for Kazhdan--Lusztig polynomials]
	{A generalization of Deodhar's setting for questions in
		Kazhdan--Lusztig theory}
\author{Rohit Agrawal}
\author{Vladimir Sotirov}
\date{\today}
\begin{document}\setcounter{page}{0}
\maketitle\thispagestyle{empty}
\setcounter{tocdepth}1
\begin{abstract}
	
We make progress on a question of Skandera by showing that a product of
	Kazhdan--Lusztig basis elements indexed by maximal elements of
	parabolic subgroups admits a Kazhdan--Lusztig basis element as a
	quotient arising from operations in the Schur algebroid if and only if
	the sequence of parabolic subgroups satisfy both a rigidity condition
	and a combinatorial criterion.
For Weyl groups, the rigidity condition specializes to a necessary
	condition for smallness of Gelfand--MacPherson resolutions.
For Schubert varieties indexed by $4231$-avoiding permutations, we
	derive a stronger necessary condition that, up to an appropriate
	equivalence, is satisfied by at most one Gelfand--MacPherson
	resolution, and exactly one if and only if $45312$ is also avoided.
Moreover, we apply the combinatorial criterion to prove the (essentially
	unique) resolution is small when $34512$ and $45123$ are likewise
	avoided.

We develop the combinatorial criterion as part of a generalization
	of Deodhar's combinatorial setting for questions in Kazhdan--Lusztig
	theory, which in the case of Weyl groups we show captures the Bya\l
	ynicki--Birula decompositions of Gelfand--MacPherson resolutions.
In particular, we obtain new combinatorial interpretations of
	Kazhdan--Lusztig polynomials, new algorithms for computing them, 
	and in the case of Weyl groups, an equivalence of the existence of
	small Gelfand--MacPherson resolutions with that of certain
	factorizations and generating-function interpretations of the
	Poincar\'e polynomial of intersection cohomology of the Schubert
	variety.

\end{abstract}
\tableofcontents

\section*{Introduction}

\subsection{Background}
The Kazhdan--Lusztig basis for the Hecke algebra of a Coxeter
	system was introduced in 1979 by Kazhdan and Lusztig
	\cite{Kazhdan--Lusztig} as a tool for constructing certain
	representations of the Hecke algebra.
It is the unique basis invariant under a natural involution of the Hecke
	algebra whose coefficients relative to the standard basis, now known
	as the Kazhdan--Lusztig polynomials, have lowest possible degrees.

In the case of Weyl groups and affine Weyl groups, the Kazhdan--Lusztig
	polynomials are of particular importance as their coefficients record
	the dimensions of local intersection cohomology groups of
	Schubert varieties \cite{Kazhdan--Lusztig80} and consequently the
	multiplicities of Verma modules
	\cite{Brylinski--Kashiwara, Beilinson--Bernstein}.
Kazhdan and Lusztig conjectured that the coefficients are non-negative
	for any Coxeter system which was ultimately proven by Elias and
	Williamson in 2012 \cite{Elias--Williamson}.
Specifically, they verified a conjecture due to Soergel \cite{Soergel}
	that the degree condition holds for an involution-invariant basis with
	non-negative coefficients arising from a realization of the Hecke
	algebra as the Grothendieck group of the algebraically defined
	category of what are now known as Soergel bimodules.
	
By definition, the involution-invariance and degree condition
	amount to the Kazhdan--Lusztig polynomials satisfying a complicated
	doubly-recursive formula, which has motivated a search for
	non-recursive formulas and algorithms for computing them.
Additionally,  the non-negativity conjecture has motivated a search for
	combinatorial interpretations of the Kazhdan--Lusztig polynomials.

Most generally, Brenti discovered in 1997 and 1998 two non-recursive
	formulas \cite{Brenti97,Brenti98} which work well for effectively
	computing any \emph{individual} Kazhdan--Lusztig polynomial, but
	do not lead to a combinatorial interpretation.

In order to explain and contextualize the contribution of the present
	work, let us recall five interrelated lines of inquiry.

First, in 1981 Lascoux and Schutzenberger \cite{Lascoux--Schutzenberger}
	observed that for a Kazhdan--Lusztig basis element indexed by a
	\emph{cograssmanian} permutation, i.e.\ by a permutation representing
	a grassmanian Schubert variety, the Kazhdan--Lusztig polynomials
	satisfy a significantly simpler recurrence relation, which easily
	leads to both a combinatorial interpretation and a non-recursive
	formula.
Their simpler recurrence relation was extended by Enright and Shelton in
	1987 \cite{Enright--Shelton87} to Coxeter elements representing certain
	symplectic and orthogonal grassmanian Schubert varieties, while the
	combinatorial interpretation was generalized by Boe in 1988
	\cite{Boe88}.
In 1995, Lascoux \cite{Lascoux95} extended further the simpler
	recurrence relation and combinatorial interpretation to
	\emph{covexillary} permutations.

Second, in 1983, Zelivinski\u{\i} \cite{Zelevinskii83,Zelevinskii83-Russian}
	showed that Lascoux and Schutzenberger's formula can be interpreted
	geometrically as arising from a small Gelfand--MacPherson resolution
	of singularities of the corresponding grassmanian Schubert cell.
He did not, however, give a geometric explanation of Lascoux and
	Schutzenberger's combinatorial interpretation.

In 1994 and 1995 Sankaran and Vanchinathan
	\cite{Sankaran--Vanchinathan94, Sankaran--Vanchinathan95} extended
	Zelivinski\u{\i}'s  construction of small resolutions and resulting formulas to certain
	Schubert varieties in symplectic and orthogonal grassmanians.
	In 2007 Perrin \cite{Perrin07} generalized both Zelivinski\u{\i}'s and
	Sankaran and Vanchinathan's results by constructing small
	\emph{quasi-resolutions} (relative minimal models, in fact) for all
	\emph{minuscule} Schubert varieties.
He also gave an easy to check combinatorial criterion for when
	his quasi-resolutions are small resolutions of singularities, but left
	as an open question the relationship with the calculation and
	combinatorial interpretations of Kazhdan--Lusztig polynomials due to
	Boe.

Third, in 1990 Deodhar \cite{Deodhar} gave a recursive algorithm for
	computing Kazhdan--Lusztig basis elements in terms of what he called
	a \emph{combinatorics of subexpressions}, which works well if one aims
	to compute \emph{all} Kazhdan--Lusztig polynomials.
Furthermore, if one aims to compute the Kazhdan--Lusztig polynomials of
	an individual basis element, Scott and Sprowl \cite{Scott--Sprowl}
	have shown that simply skipping most of the steps involved in
	Deodhar's recursion yields a feasible non-recursive algorithm.

Unlike Brenti's formulas, when combined with the non-negativity
	conjecture, Deodhar's algorithm implies the existence of generating
	functions realizing the Kazhdan--Lusztig polynomials as counting the
	elements of some set of subexpressions of a fixed reduced word for a
	Coxeter element (weighted by a certain \emph{defect} statistic).
Furthermore, Jones and Woo \cite{Jones-Woo} have shown that in
	the case of basis elements indexed by cograssmanian permutation, one
	can explicitly choose the set of subexpressions so as to recover
	the formulas due to Lascoux and Schutzenberger
	\cite{Lascoux--Schutzenberger} and Zelevinski\u{\i}
	\cite{Zelevinskii83-Russian,Zelevinskii83} (although Lascoux's more
	general formula for covexillary permutations \cite{Lascoux95} remains
	outside Deodhar's combinatorics of subexpressions).
However, Jones and Woo's choices of a set of subexpressions are not at
	all canonical, leading them to suspect that Deodhar's combinatorial
	model ought to be augmented with criteria on subexpressions to guarantee
	uniqueness.

Fourth, a more drastic method than Scott and Sprowl's of simplifying
	Deodhar's algorithm is to skip \emph{all} of the recursive steps,
	valid, of course, only for \emph{some} Kazhdan--Lusztig basis
	elements, and corresponding to the generating functions realizing
	Kazhdan--Lusztig polynomials counting elements of the canonical set of
	\emph{all} subexpressions of a reduced word.

In 1999, Billey and Warrington \cite{Billey--Warrington} extracted from
	Deodhar's combinatorial framework a criterion for when a
	Kazhdan--Lusztig basis element $C'_w$ indexed by a permutation
	$w$ admits such a simplified computation, equivalently, such a
	canonical combinatorial interpretation.
Combining the criterion with a subtle graph-theoretic argument, they
	obtained a characterization in terms of pattern-avoidance of such
	permutations $w$.
In 2007, Billey and Jones \cite{Billey--Jones} extended the
	characterization of such elements, which they called 
	\emph{Deodhar elements}, to all finite Coxeter systems.

Billey and Warrington also observed that Deodhar elements admit two
	other equivalent descriptions.
First, that in the case of (affine) Weyl groups the Poincar\'e
	polynomial of the full intersection cohomology of a Schubert cell
	indexed by a Coxeter element $w$ factors as $(1+q)^{\ell(w)}$ where
	$\ell(w)$ is the length of a reduced word for $w$.
Second, that the Kazhdan--Lusztig basis element $C'_w$ indexed by $w$
	factors as a product $C'_{s_{i_1}}\cdots C'_{s_{i_r}}$ of
	Kazhdan--Lusztig basis elements where $s_{i_1}\cdots s_{i_r}=w$ is
	any reduced word for $w$ in the Coxeter system.

Fifth, in 2007, Skandera \cite{Skandera} characterized permutations $w$
	with smooth Schubert cells in a way resembling the second description
	of the elements characterized by Billey and Warrington.
Explicitly, relying on a classification in terms of pattern-avoidance
	due to Lakshmibai and Sandhya \cite{Lakshmibai--Sandhya} of
	permutations $w$ indexing smooth Schubert varieties, he showed that
	that satisfaction of the pattern-avoidance condition is equivalent to
	the Kazhdan--Lusztig basis element $C'_w$ being a quotient by a
	particular Laurent polynomial of a monomial
	$C'_{\mxml{J_1}}\cdots C'_{\mxml{J_r}}$ of Kazhdan--Lusztig basis
	elements for $\mxml{J_1}\cdots \mxml{J_r}$ a particular factorization
	of $w$ into \emph{reversal permutations}, i.e.\ into maximal elements
	of parabolic subgroups.
Observing that neither his nor Billey and Warrington's result
	generalizes the other, but that the two agree on permutations
	satisfying both of their respective pattern-avoidance conditions,
	Skandera posed the question of determining a precise generalization of
	both, e.g.\ of characterizing which Kazhdan--Lusztig basis elements
	can be expressed as a quotient of a monomial $C'_{\mxml{J_1}}\cdots
	C'_{\mxml{J_r}}$.

The results of this work, which we describe below, contribute to each of
	these lines of inquiry and establish closer relations between them.

\subsection{A partial resolution to Skandera's question}
The hypothesis of the main result of this work, \cref{thm:2},
	arises from the following subtlety with Skandera's original question.

If a monomial $C'_{\mxml{J_1}}\cdots C'_{\mxml{J_r}}$ has a
	quotient by a Laurent polynomial that is Kazhdan--Lusztig basis
	element, then the Laurent polynomial must be the coefficient of the
	leading (relative to the Bruhat order on the Coxeter
	group) term $w(J)$ of $C'_{\mxml{J_1}}\cdots C'_{\mxml{J_r}}$.
Unfortunately, the quotient by this leading coefficient may fail to be a
	Kazhdan--Lusztig basis element not because of failing the degree
	condition on its coefficients, but because of its coefficients
	failing to be polynomials rather than Laurent polynomials. 

Fortunately, the composition laws of the Schur algebroid introduced by
	Williamson \cite{Williamson} provide an ample supply of the kind of
	quotients we need: given any two Hecke algebra elements $h_1$ and
	$h_2$ that are respectively invariant under the induced right and left
	action of a finite parabolic subgroup, their product $h_1\cdot h_2$
	has a quotient by the Laurent Poincar\'e polynomial of that parabolic
	subgroup such that the coefficients of that quotient are polynomials
	rather than Laurent polynomials as desired.

In particular, the quotients studied by Skandera, it turns out, arise by
	repeated application of the composition laws of the Schur algebroid
	and thereby suggest a narrower but more tractable version of his
	question: to determine which Kazhdan--Lusztig basis elements belong to
	the closure of the set of Kazhdan--Lusztig basis elements 
	$C'_{\mxml J}$ under the composition laws of the Schur algebroid.
Whether there exist quotients of monomials $C'_{\mxml{J_1}}\cdots
	C'_{\mxml{J_r}}$ whose coefficients are polynomials rather than
	Laurent polynomials outside the above closure is an open question,
	whose negative answer would completely subsume Skandera's question
	under our generalization of Deodhar's framework.

With this caveat in mind, our main result which we record as
	\cref{thm:2} asserts that if $C'_{\mxml{J_1}}\cdots C'_{\mxml{J_r}}$
	has a quotient arising from the composition laws of the Schur
	algebroid, then the \emph{factorization} $J_1,\dots,J_r$ of $w(J)$ has
	a \emph{contraction} that is \emph{admissible absolutely bidescent}
	and \emph{tight}.

To explain the conclusion, we need the following terminology:
	given a Coxeter system $(W,S)$, we say that two subsets of $S$
	\emph{commute} if each element of one is distinct from and commutes
	with each element of the other, and say that a non-empty $J'\subseteq
	J$ is a \emph{connected component} of $J$ if it commutes with with
	$J\setminus J'$ and is not itself a union of commuting subsets.
Then a \emph{contraction} is a simplification of $J_1,\dots,J_r$ by a
	process which removes a connected component of $J_i$ if it is
	contained in $J_k$ and commutes with $J_j$ for $j$ strictly between
	$i$ and $k$.
Being \emph{bidescent} means that every connected component of $J_i$ 
	that commutes with $J_j$ when $i<j$ (respectively when $j<i$) is also
	a connected component of $\rdes{w(J)}
		=\set{s\in S:\len{w(J)s}<\len{w(J)}}$ (respectively
	$\ldes{w(J)}=\set{s\in S:\len{sw(J)}<\len{w(J)}}$), where $S$ is the
	set of generators of the Coxeter system $(W,S)$ and $\len w$ -- the
	corresponding length of $w\in W$.
Being \emph{absolutely} bidescent means that every \emph{slice} $J_i,\dots,J_k$
	is also bidescent, as is every slice of a sequence \emph{equivalent}
	to $J_1,\dots,J_r$ under the relation of moving to $J_k$ a connected
	component of $J_i$ that commutes with each $J_j$ for $j$ between $i$
	and $k$, inclusive of $k$.

Finally, being \emph{admissible} refers to a notion arising in our
	generalization of Deodhar's framework which amounts to asserting that
	$1$ is the polynomial indexed by the leading term $w(J)$ in the
	quotient under consideration. 
Being \emph{tight} is the other key notion of our generalization, which
	amounts to the quotient \emph{being} the Kazhdan--Luszig basis element
	$C'_{w(J)}$ (we borrow the term from Lusztig \cite{Lusztig}).
	
The conclusion of \cref{thm:2} is quite a rigid condition since even
	though we show in \cref{sec:bidescent} that every $w\in W$ has an
	absolutely bidescent \emph{factorization}, i.e.\ absolutely bidescent
	sequence $J_1,\dots,J_r$ with $w=w(J)$, we also show that no simple
	approach would work to prove that every $w$ admits an
	\emph{admissible} such factorization, leaving the existence of such an
	open question (\cref{ques}).
The significance of the conclusion is also evidenced by the results we
	obtain in type $A$, which we describe below.

The proof of \cref{thm:2} amounts to a careful combination of
	Soergel's conjecture (specialized to the Schur algebroid by
	Williamson), \cref{cor:factors} which asserts a factorization is tight
	only if it is absolutely tight, the standard fact that
	Kazhdan--Lusztig polynomials $P_{x,w}$ satisfy $P_{xs,w}=P_{x,w}$ for
	all $s\in\rdes w$, a precise description of the composition laws of
	the Schur algebroid we give in \cref{cor:invariance-positive}, and a
	technical result for Coxeter systems we record as
	\cref{lem:wz-w-inequality}.

\subsection{Tight factorizations of permutations avoiding $4231$}
The importance of our main result, \cref{thm:2}, is evidenced by the
	results we obtain for type $A$, i.e.\ for the
	Coxeter system of the symmetric group generated by the set $S$ of
	adjacent transposition. 
In order to explain the results, recall that subsets $J\subseteq S$ can
	be identified with unions of intervals, so that $J_1,\dots,J_r$ being
	bidescent implies in particular that connected components of $J_i$
	commuting with $J_j$ when $i<j$ correspond to \emph{descent intervals}
	$[a,b]$ of $w$, i.e.\ $[a,b]$ such that
	$w(a-1)<w(a)>w(a+1)>\cdots>w(b)<w(b+1)$.

Our first result in type $A$, \cref{thm:4231-45312-avoid}, is that a
	\emph{minimal} (i.e.\ without contractions) absolutely bidescent
	factorization of a $4231$-avoiding permutation is tight only if it has
	the property of being \emph{monotone}, in which case the factorization
	is absolutely \emph{$4231$-$45312$-avoiding}.
Moreover, we show all monotone factorizations are equivalent and can be
	computed via the following algorithm (\cref{prop:monotone-algorithm}),
	which a priori can be applied to any permutation $w$ but succeeds if
	and only if $w$ is $4231$-$45312$-avoiding (hence also gives a
	criterion for a permutation to be $4231$-$45312$-avoiding):

\begin{enumerate}
	\item Let the last term of the factorization have connected components
		corresponding to the \emph{strong descent} intervals $[a,b]$ of $w$,
		i.e.\ on which $w$ is decreasing, and satisfies $w(p)<w(a)$ if $p<a$, and
		$w(b)<w(r)$ if $b<r$.
	The existence of such intervals for non-identity permutations is
		guaranteed by $4231$-avoidance.
	\item For each interval, let 
		$\ell=\max\set{q\in[a,b]:w(p)<w(q)$ if $p<a}$ and
		$r=\min\set{q\in[a,b]:w(q)>w(r)$ if $b<q}$.
	The relevant sets being non-empty is equivalent to $45312$-avoidance.
	\item Define $w'$ by sorting in increasing order $w(a),\dots,w(b)$,
		and then in decreasing order the $w(p)$ with 
		$p\in[a,a+b-\ell)\cup(b-r+a,b]$.
	\item If $w'$ is not the identity, apply steps 1-4 to it to obtain the
		term of the factorization preceding the one just computed.
\end{enumerate}
The key to the proof of \cref{thm:4231-45312-avoid} is
	\cref{lem:unbounded}, which is a technical result regarding type $A$
	Coxeter systems, similar to \cref{lem:wz-w-inequality} used in the
	proof of \cref{thm:2}.
The rest of the proof is a matter of carefully analyzing the notion of
	monotonicity and the consequences of \cref{lem:unbounded} for $4231$-
	and $45312$-avoiding permutations.

Two remarks are in order.
First, the permutation $4231$ has two non-equivalent tight
	factorizations obtained from the $i\mapsto 5-i$ symmetry, indicating
	that the combinatorial problem of determining tight factorizations for
	$4231$-realizing permutations is significantly richer than that for
	$4231$-avoiding ones.
Second, since a permutation $w$ indexes a smooth Schubert variety if and
	only if it avoids the patterns $3412$ and $4231$, our result applies
	to such permutations and recovers (contractions of) Skandera's
	factorizations.
A natural direction for future work is to determine if for other Weyl
	groups the elements indexing smooth Schubert varieties also have
	unique minimal admissible absolutely bidescent factorizations.

Our second result in type $A$ is \cref{thm:4}, which asserts that for
	$4231$-$45312$-$34512$-$45123$-avoiding permutations,
	the monotone factorizations are actually tight.

The hypotheses of avoiding the patterns $34512$ and $45123$ are
	explained by our \cref{lem:diamond-pattern}, which asserts that
	avoidance of the pattern $m_2+1\ m_2+2\ \dots\ m_2+m_1\ 1\ 2\ \dots\
	m_2$ by a $4231$-$45312$-avoiding permutation corresponds to a certain
	restriction on the complexity of the monotone factorization, 
In particular, $3412$-avoidance of Skandera's results correspond to 
	$m_1=m_2=2$, whereas $34512$-$45123$-avoidance correspond to the more
	general $\set{m_1,m_2}=\set{2,3}$.
Extending this result for higher $m_1,m_2$ is a natural direction for
	future work.

The proof of \cref{thm:4} is an application of the combinatorial
	criterion for a factorization being tight, which we develop as part of
	our generalization of Deodhar's framework.
Our generalization of Deodhar's framework also allows us to reinterpret
	our results in type $A$ as asserting that a Schubert variety indexed
	by a $4231$-avoiding permutation has a small Gelfand--MacPherson
	resolution only if it is also $45312$-avoiding, in which case the only
	Gelfand--MacPherson resolutions that could be small correspond to
	factorizations with monotone contractions, and moreover Schubert
	varieties indexed by $4231$-$45312$-$34512$-$45123$-avoiding
	permutations do have small Gelfand--MacPherson resolutions.

This interpretation of our results complements nicely results of
	Larson \cite{Larson}, who has shown that Schubert varieties indexed by
	$3412$-$52341$-$635241$-avoiding permutations also have small
	resolutions.
In particular, the intersection of ours and Larson's results are
	precisely the results of Skandera for $3412$-$4231$-avoiding
	permutations.  
Incorporating Larson's completely different approach into our framework
	is another avenue of future work.
\subsection{Our generalization of Deodhar's framework and the Schur algebroid}
Our generalization of Deodhar's framework can be geometrically motivated
	by a careful analysis of the Bya\l ynicki--Birula decomposition of
	Gelfand--MacPherson resolutions. 
We only sketch the analysis in \cref{sec:gm} because the connection to
	geometry is only helpful in translating our results on tight
	factorizations to results on smallness of Gelfand--MacPherson
	resolutions in the case of Weyl groups. 
Indeed, the development of our generalization requires only algebra and
	combinatorics.
Nevertheless, the connection is significant as it implies that every
	small Gelfand--MacPherson resolution gives both a combinatorial
	interpretation of and a formula for computing Kazhdan--Lusztig
	polynomials, clarifying the relationship between the two lines of
	inquiry following Lascoux and Schutzenbgerger, and Zelivinski\u{\i}.

Explicitly, our generalization of Deodhar's framework associates to any
	factorization $J_1,\dots,J_r$, a family of generating functions we
	call \emph{defect polynomials}.
These generating functions count \emph{equivalence classes} of sequences
	of elements $\sigma_k$ (which following Billey and Warrington
	\cite{Billey--Warrington} we call \emph{masks}) of the parabolic
	subgroups $W_{J_k}$ generated by the $J_k$, weighted by
	a generalization of Deodhar's defect statistic. 

The equivalence relation we impose on masks is the main point of
	subtlety in our generalization.
Without it, our defect polynomials would give a combinatorial
	interpretation of the monomial $C'_{\mxml{J_1}}\cdots C'_{\mxml{J_r}}$,
	whereas with it they ought to interpret a quotient by a Laurent
	polynomial arising from the composition laws of the Schur algebroid.
Specifically, the Laurent polynomial and equivalence relation we choose
	are determined by \emph{overlap subsets} $\overlap Jik\subseteq
	J_i\cap J_k$ that commute with $J_j$ for each $j$ between $i$ and $k$,
	with the Laurent polynomial given by the product of Laurent Poincar\'e
	polynomials of the parabolic subgroups $W_{\overlap Jik}$, and 
	the equivalence relation induced by the joint left- and right-
	actions of $W_{\overlap Jik}\leq W_{J_i}\cap W_{J_k}$ on $W_{J_i}$ and
	$W_{J_k}$.

The subtlety with this (or any other) choice of equivalence relation is
	the same subtle problem as Skandera's original question: we have no
	choice but to verify the quotient arises from compositions in the
	Schur algebroid, which we do in \cref{cor:quotient}.

In order to prove \cref{cor:quotient}, we first give a recursive formula
	for computing defect polynomials (\cref{cor:recursive-defect-polys})
	and then show it corresponds to compositions in the Schur algebroid
	(\cref{prop:invariant_hecke_mult}).
The former relies on a careful anaylsis in terms of the Coxeter
	system of the defects counted by the defect statistic, which we record
	as \cref{prop:defect-ideals,cor:unique-reps}.
The latter is complicated by the fact that the defect polynomials are
	the coefficients of the relevant quotient in the Hecke algebra only up
	to a uniform degree shift.
In particular, we again require the careful analysis of the composition
	laws of the Schur algebroid in \cref{sec:hecke}.

\subsection{The combinatorial criterion}
Having settled on an equivalence relation on masks, the next issue is
	with extending not only Deodhar's notion of \emph{defects}, but also
	Billey and Warrington's notions of \emph{bounces} and
	\emph{crossings}, in a way that is invariant under the equivalence.
The latter are necessary to formulate our generalization of Billey and
	Warrington's \emph{combinatorial criterion} for a factorization being
	tight: that the number of defect bounces exceeds the number of
	non-defect bounces for all but the ``maximal'' mask (\cref{thm:combi}).

The most interesting aspect of the proof is actually our definition of
	defects, bounces, and crossings for arbitrary Coxeter systems. 
Billey and Warrington's definitions in the case of Coxeter systems of
	type $A$ are of defects, bounces, and crossings \emph{of pairs of
	strings} in a \emph{string diagram} that arises from interpreting a
	mask using the natural action of the permutation group.
It is a pleasant surprise that even though arbitrary Coxeter groups lack
	an analogous natural action, we have nevertheless been able to 
	formulate a notion of a \emph{virtual pair of strings} of a mask, 
	which has allowed us to write our definitions and arguments so that
	they apply word for word both for arbitrary Coxeter systems and for
	the concrete special case of permutation groups.

We also give a conceptual explanation of the combinatorial criterion by
	showing that the generating functions counting the negative of the
	excess actually express (without a degree shift!) the
	coefficients of the quotient in the Hecke algebra in terms of a basis
	more suitable for identifying Kazhdan--Lusztig basis elements
	(\cref{lem:normalized-factorization}).
	
Finally, developing our generalization as above allows for Billey and
	Warrington's argument that the combinatorial criterion also
	characterizes whether the Poincar\'e polynomial of intersection
	cohomology of a Schubert cell factors in a way to go
	through word for word in our more general setting
	(\cref{cor:schubert}), although with the factorization of the
	polynomial depending on the sequence $J_1,\dots,J_r$.

\subsection{Computing Kazhdan--Lusztig polynomials}
Our final contribution is showing that Deodhar's theorem
	\cite[Theorem~4.12]{Deodhar} generalizes to our setting: given an
	admissible factorization $J_1,\dots,J_r$ of a Coxeter element $w$, the
	family of Kazhdan--Lusztig polynomials $(P_{x,w})_{x\in W}$ can be
	combinatorially interpreted as a family of generating functions that
	is a restriction of the family of defect polynomials.
Using ideas of \cite{Scott--Sprowl}, we actually obtain a
	generic variation of Deodhar's and Scott and Sprowl's algorithms by
	replacing the use of singleton factorizations induced by reduced words
	with arbitrary admissible factorizations
	(\cref{thm:1}).
Moreover, we suggest using admissible \emph{right-descent}
	factorizations to compute Kazhdan--Lusztig polynomials, as described
	in \cref{exmp:new-algo}, as a potentially more efficient algorithm.

The most interesting part of the proof is that whereas
	the final part of Deodhar's original argument relies on the
	non-negativity of the coefficients of the Kazhdan--Lusztig polynomials
	(established for Weyl and affine Weyl groups by Kazhdan and Lusztig
	\cite{Kazhdan--Lusztig80}, and for arbitrary Coxeter systems by Elias
	and Williamson \cite{Elias--Williamson}), the final part of our
	argument relies on Williamson's work on the Schur algebroid
	\cite{Williamson} and on the more general conjecture of Soergel that
	Elias and Williamson established in \cite{Elias--Williamson}.

Curiously, we obtain along the way a new derivation of the M\"obius
	function \cite[Section~3.7]{Stanley} of the Bruhat order using only
	Deodhar's original framework, as explained in \cref{rmrk:mobius}.

\subsection{Roadmap}
We have organized the exposition of our theory and results as follows.

In \cref{sec:gm}, we give a sketch of how in
	the case of Weyl groups our generalization of Deodhar's framework
	arises from a precise analysis of the Bya\l ynicki--Birula
	decomposition of Gelfand--MacPherson resolutions.
This result is for the reader versed in algebraic geometry and none of
  the arguments in the sequel depend on it (though the geometric
  interpretation of the results does).

The remaining material can be classified in two groups:
	\crefrange{sec:prelim}{sec:bidescent} comprise our results valid for
	all Coxeter systems, while \crefrange{sec:4231-45312}{sec:height2}
	comprise those for type $A$ Coxeter systems.
Their contents are as follows.

In \cref{sec:prelim}, we review the basic facts about Coxeter
	systems that we will use and also establish the conventions and
	notation that we will follow.
We also prove a few crucial results relevant to the notion of defects,
	including the key \cref{lem:wz-w-inequality}.

In \cref{sec:defects}, we introduce the term \emph{factorization}
	for a sequence of subsets of generators $J_k\subseteq S$ generating
	finite parabolic subgroups, and develop the basic properties of
	\emph{masks} and \emph{defects} associated to a factorization 
	that generalize Deodhar's combinatorics of subexpressions of a reduced
	word.
We also give useful illustrations of the concepts for type $A$ Coxeter
	systems.

In \cref{sec:masks-poly} we associate to each factorization a
	family of generating functions we call \emph{defect polynomials} and
	derive their basic properties including a recursive formula for their
	computation.
We also introduce the relation of \emph{equivalence} of factorizations,
	for which equivalent factorizations have the same defect polynomials,
	and the relation of one factorization being a \emph{factor} of
	another, for which the corresponding defect polynomials are quotients
	of one another. 
Crucial notions introduced here are those of a property of a
	factorization holding \emph{absolutely}, i.e.\ for all factors, and
	the notions of a factorization being \emph{admissible} and
	\emph{tight}.

In \cref{sec:combi} we define our generalization of Billey and
	Warrington's notions of \emph{bounces} and \emph{crossings} and prove
	the combinatorial criterion for when a factorization is tight.

In \cref{sec:hecke} we review the basic facts concerning the
	involution of the Hecke algebra and identify a uniform degree shift
	necessary to realize a family of polynomials as the coefficients of an
	involution-invariant Hecke algebra element, which we call a
	\emph{normalized Hecke algebra element} of the family.
We then review the \emph{Schur algebroid} introduced by
	Williamson \cite{Williamson} and use it to show that normalized Hecke
	algebra elements of defect polynomials are involution-invariant
	quotients of products of Kazhdan--Lusztig basis elements of the Hecke
	algebra indexed by maximal elements of parabolic subgroups.

Also in this part, we give our conceptual explanation of the
	combinatorial criterion: that the normalized Hecke algebra element of
	the family of defect polynomials, expressed in the alternative basis
	for the Hecke algebra, has coefficients given by the generating
	functions of half the negative excess of non-defect bounces over
	defect bounces.

In \cref{sec:deodhar} we complete our generalization of Deodhar's
	framework by showing how to use admissible factorizations to compute
	and give combinatorial interpretations of Kazhdan--Lusztig
	polynomials, and that tight factorizations correspond to expressing
	Kazhdan--Lusztig basis elements as quotients.
Along the way we show that the properties of being admissible,
	respectively tight, hold only if they hold absolutely.

We also show that factorization of a Weyl group element $w$ is tight if
	and only if the Poincar\'e polynomial of the full intersection
	cohomology of the Schubert variety $X_w$ can be factored in a
	particular way, in which case we give a combinatorial interpretation
	in terms of defects and bounces, and also obtain a new derivation of
	the M\"obius function of the Bruhat order.

In \cref{sec:schur} we give our partial resolution to
	Skandera's question by showing that any Kazhdan--Lusztig basis element
	in the closure of those indexed by maximal elements under the Schur
	composition laws arises from a tight factorization.
To do so we introduce the \emph{contraction} partial order on
	factorizations, and identify the key properties of a factorization
	being \emph{right-descent} and \emph{left-descent}.

In \cref{sec:bidescent} we show that any minimal in the contraction
	partial order tight factorization is \emph{bidescent}, i.e.\ both
	right-descent and left-descent.
Accordingly, we suggest that Kazhdan--Lusztig basis elements may be more
	efficiently computed using admissible right-descent factorizations,
	rather than arbitrary admissible factorizations

We also show that any Coxeter element has an absolutely bidescent
	factorization, and indicate why the proof cannot work to show it has
	an absolutely admissible bidescent factorization, which a tight
	factorization would have to be.
We end by introducing the property of a factorization being \emph{strong
	right-descent}, which, unlike the property of being right-descent,
	holds absolutely if it holds.

With \cref{sec:4231-45312} we begin our specific results valid for type
	$A$.
In this section, we show that any tight factorization of a
	$4231$-avoiding permutation has to be absolutely
	$4231$-$45312$-avoiding with a \emph{monotone} contraction.
We also introduce the notions of \emph{permutation-directedness}, which
	allow for a visual definition of the condition of being monotone,
	that realizes being monotone as a more restricted version of being
	strong right-descent.

In \cref{sec:45312} we introduce the \emph{heap of connected components}
	of a factorization and use it to characterize monotone factorizations
	as the bidescent factorizations of $45312$-avoiding permutations, and
	to show that every monotone factorization is absolutely monotone .

In \cref{sec:rigid} we continue using the heap of connected components
	to show that strong bidescent factorizations have a rigid structure
	allowing us to conclude that monotone factorizations are all
	absolutely $4231$-$45312$-avoiding, and that every
	$4231$-$45312$-avoiding permutation has a unique up to equivalence
	minimal monotone factorization.

Finally, in \cref{sec:height2}, we use the rigidity of bidescent
	factorizations to relate their complexity to pattern avoidance, and
	then apply the combinatorial criterion to conclude that
	$4231$-$45312$-$45123$-$34512$-avoiding permutations have exactly one
	(up to equivalence) minimal tight factorization.

\section{Deodhar's framework and Gelfand--MacPherson resolutions}

\label{sec:gm}
As described by Jones and Woo in \cite[Section~5]{Jones-Woo}, for the
	case of type $A$ Weyl groups and the general linear group, Deodhar's
	framework arises naturally from the Bia\l ynicki--Birula decomposition
	of Bott--Samelson resolutions of Schubert varieties.
Our generalization arises from the analogous decomposition of the more
	general Gelfand--MacPherson resolutions, as we sketch below.

\subsection{Cellular decompositions of flag varieties}
We begin by reviewing the cellular decomposition of flag varieties and
	giving a precise description of Schubert cells in quotients by any
	parabolic subgroup.
For more details, we direct the reader to \cite[Chapter~21]{Milne} and
	\cite[Section~7.1]{Kumar}.

A \define{split reductive group} is a pair $(G,T)$ of a \emph{reductive}
	group $G$ and a maximal torus subgroup $T\leq G$ that is \emph{split},
	in which case \cite[Theorem~21.11]{Milne} there is a set $\Phi$ of
	of non-zero characters called \define{roots} such that for every
	$\alpha\in\Phi$ there exists a (necessarily unique) subgroup
	$U_\alpha$ of $G$, called a \define{root group}, that is stable under
	conjugation by $T$ and admits a $T$-equivariant isomorphism to the
	$1$-dimensional (algebraic) representation determined by $\alpha$.

A \define{cocharacter}, i.e.\ a group homomorphism from the
	multiplicative group of scalars into the torus, is \define{regular} if
	the induced action by conjugation on each root
	subgroup is non-trivial, in which case we have a \define{system of
	positive roots} $\Phi^+\subseteq\Phi$ consisting of those roots on
	whose root subgroups the action by the multiplicative group of scalars
	induced by the cocharacter is to multiply by a positive power of the
	scalar.

To each system of positive roots corresponds a \define{Borel subgroup}
	$B\le G$, i.e.\ a solvable subgroup with the quotient $G/B$ complete,
	containing $T$, and explicitly given as a semi-direct product
	$B=T\ltimes B_u$ where $B_u$ has the property
	\cite[Theorem~21.68]{Milne} that for any subgroup $U\le B_u$
	normalized by $T$, and any total ordering on $\set{\alpha\in\Phi^+
	:U_\alpha\le U}$, the group operation of $G$ gives an isomorphism
	$\ds\prod_{\alpha\in\Phi^+:U_\alpha\le U}U_{\alpha}\cong U$.

The intersection of the normalizer $N$ of $T$ with $B$ is $T$ itself,
	and its quotient by the torus is an etal\'e group scheme called the
	\define{Weyl group}, which is in fact a constant group scheme, i.e.\
	a finite group. 
By definition the conjugation action of the normalizer permutes the
	root groups and so gives a left action on the roots in $\Phi$.

Any subgroup $P$ containing a Borel subgroup $B$ is a \define{parabolic
	subgroup}, i.e.\ such that $G/P$ is a complete (and in fact
	projective and connected) variety, and determines a subgroup 
	$N\cap P/B\cong W'\leq W$ of the Weyl group, such that the quotient
	$W/W'$ of the Weyl group embeds in the \define{flag variety} $G/P$,
	and is the constant scheme of fixed points for the left action of $T$
	on $G/P$ \cite[Lemma~7.1.22]{Kumar}.

The \define{Schubert cells} of the flag variety $G/P$ are the $B$-,
	equivalently, $B_u$-orbits of the $T$-fixed points $W/W'$.
If $n_w$ is a representative in $N\leq G$ of $w\in W$ and $w'$ is the
	image of $w$ in $W/W'$ , then the isotropy subgroup $U_{w'}^P$ of the
	$T$-fixed point $w'$ for the $B_u$-action on $G/P$ is given by
	$U_w^P=B_u\cap n_wPn_w^{-1}$.
Since $U_w^P$ is a $T$-normalized subgroup of $B_u$, it follows that
	the group operation induces an isomorphism $U_w^P
		\cong\ds\prod_{\alpha\in\Phi^+\setminus N_P^w}U_\alpha$ where
 $N_P^w=\set{\alpha\in\Phi^+:U_{w^{-1}(\alpha)}$ not contained in $P}$.
Since $B_u\cong\ds\prod_{\alpha\in\Phi^+}U_\alpha$, it follows that
	$B_uw'\cong B_u/U^P_{w'}\cong\ds\prod_{\alpha\in N_P^w}U_\alpha$
	(this argument extends the special case of $P=B$ proved in
	\cite[Proposition~21.79]{Milne}).
In particular, the group operation gives an isomorphism
	$\ds\prod_{\alpha\in N_P^w}U_\alpha\cong U^w_P$ for a subscheme
	$U^w_P$ of $B_u$ such that $B_uw'=U^w_Pw'\cong U^w_P$.
Note that $N^w_P$ and $U^P_w$ depend only on the left coset 
	$wW'=w'\in W/W'$ and that, for $Q$ containing $P$, the induced map 
	$U_P^w\to U_Q^w$ is the projection $\ds\prod_{\alpha\in N^w_P}
		U_\alpha\to\ds\prod_{\alpha\in N^w_Q}U_\alpha$ corresponding to
	the natural inclusion of indices $N^w_P\hookleftarrow N^w_Q$.

The \define{Bya\l ynicki--Birula decomposition} of the flag variety is
	given by \cite[Corollary~13.49]{Milne} $G/P=\bigudj_{w'\in
	W'}Y(w')$ (disjoint union of locally closed subvarieties) where
	$Y(w')$ is the concentrator scheme of the the fixed point $w'$ in
	$G/P$ for the action by scalars determined by a regular cocharacter
	for which the induced by conjugation actions on root groups $U_\alpha$
	with $\alpha\in\Phi^+$ are multiplication by distinct (positive)
	powers of scalars.
In fact, $Y(w')$ are precisely the Schubert cells, so that
	the Bya\l ynicki--Birula decomposition coincides with the
	\define{cellular decomposition} $G/P=\ds\bigudj_{w'\in W'}B_uw'P/P
		=\bigudj_{w'\in W'}U_P^ww'P/P$ \cite[Lemma~7.1.22]{Kumar}.
\subsection{Bya\l ynicki--Birula decompositions and smallness of 
	Gelfand--\linebreak{}MacPherson resolutions}
Given a sequence $P_1,\dots,P_m$ of parabolic subgroups containing $B$,
	consider the iterated base change 
	$P_1/\qy{P_1\cap P_2}\times_{G/P_2}G/\qy{P_2\cap P_3}\times_{G/P_3}\cdots\times G/\qy{P_{m-1}\cap P_m}\times_{G/P_m} G/B$
	\[\xymatrix@!0@C=3.5em{
	&&&&P_1/\qy{P_1\cap P_2}\times_{G/P_2}G/\qy{P_2\cap P_3}\times_{G/P_3}\cdots\times G/\qy{P_{m-1}\cap P_m}\times_{G/P_m} G/B
	\ar[dllll]\ar[dll]\ar[d]\ar[drr]\ar[drrrr]\\
	P_1/\qy{P_1\cap P_2}\ar[dr]&&G/\qy{P_2\cap P_3}\ar[dl]\ar[dr]&&\ar[dl]\cdots\ar[dr]&&G/\qy{P_{m-1}\cap P_m}\ar[dl]\ar[dr]&&G/B\ar[dl]\\
	&G/P_2&&G/P_3&&\cdots&&G/P_m}\]
	which is the vetex of the limiting cone above, and is smooth and
	connected \cite[Proposition~1.10]{Larson}.
In other words, it is the scheme of those $(g_1\qy{P_1\cap P_2},\dots,g_{m-1}\qy{P_{m-1}\cap P_m},g_mB)$ satisfying 
	$g_1\in P_1$ and $g_kP_{k+1}=g_{k+1}P_{k+1}$.

Since the diagram the iterated base change is the limit of consists of
	projective varieties equipped with a left actions of $B$ and
	$B$-equivariant morphisms between them, the iterated base change is
	also a projective scheme with a left $B$-action. 
In particular, it has a locally affine left $T$-action because each of
	the factors in the diagram do, whence the cocharacter giving the
	Bya\l ynicki--Birula decomposition of the flag varieties also gives a
	Bya\l ynicki--Birula decomposition of the iterated base change.

Taking the fixed points with respect to action of $T$ and
	concentrator schemes into those fixed points must preserve the
	limiting cone. 
In particular, if each $P_k$ determines the subgroup $W_k$ of the Weyl
	group, then the constant scheme of fixed points for the torus action
	is the iterated base change $W_1/\qy{W_1\cap
	W_2}\times_{W/P_2}G/\qy{W_2\cap W_3}\times_{W/P_3}\cdots\times
	W/\qy{W_{m-1}\cap W_m}\times_{W/P_m}W$, and for a $T$-fixed point
	$(w'_1,\dots,w'_m)$, its concentrator subscheme is given by the
	iterated base change $U^{w_1}_{P_1\cap P_2}
		\times_{U^{w_2}_{P_2}}U^{w_2}_{P_2\cap P_3}
		\times_{U^{w_2}_{P_3}}\cdots\times_{U^{w_m}_{P_m}}U^{w_m}_B$
	(which makes sense because $w_kW_{k+1}=w_{k+1}W_{k+1}$).
But as noted above, the diagram the limiting cone is over consists of
	projection maps from products, so the result is the product of
	$U_\alpha$'s indexed by the colimit of the diagram of subsets of
	$\Phi^+$:
\[\xymatrix@!0@C=2.5em{N^{w_1}_{P_1\cap P_2}&&N^{w_2}_{P_2\cap P_3}&&\cdots&&N^{w_{m-1}}_{P_{m-1}\cap
P_m}&&N^{w_m}_B\\
&N^{w_1}_{P_2}\ar@{_{(}->}[ul]\ar@{^{(}->}[ur]&&N^{w_2}_{P_3}\ar@{_{(}->}[ul]\ar@{^{(}->}[ur]&&\cdots\ar@{_{(}->}[ul]\ar@{^{(}->}[ur]&&N^{w_{m-1}}_{P_m}\ar@{_{(}->}[ul]\ar@{^{(}->}[ur]}\]
and so indexed by the disjoint union of sets (not subsets)
$\qy{N^{w_1}_{P_1\cap P_2}\setminus N^{w_1}_{P_2}}
	\udj\cdots\udj\qy{N^{w_{m-1}}_{P_{m-1}\cap P_m}\setminus
	N^{w_{m-1}}_{P_m}}\udj~N^{w_m}_B$.

The image of the projection morphism
	$P_1/\qy{P_1\cap P_2}\times_{G/P_2}G/\qy{P_2\cap P_3}
		\times_{G/P_3}\cdots\times G/\qy{P_{m-1}\cap P_m}\times_{G/P_m}G/B
		\to G/B$ is always a \define{Schubert variety}, i.e.\ the closure of
	a Schubert cell $U^w_Bw$ in $G/B$, and when birational, the projection morphism is a
	\define{resolution of singularities} \cite[Definition~4.1]{Larson},
	\cite[Section~2.11]{Gelfand--MacPherson}, called a
	\define{Gelfand--MacPherson resolution}.
The map is \define{small} \cite[Section~6.2]{Goresky--MacPherson} if the
	locus of points whose pre-images have dimension at least $d>0$ has
	codimension strictly more than $2d$.

Since the projection morphism 
	$P_1/\qy{P_1\cap P_2}\times_{G/P_2}G/\qy{P_2\cap P_3}
		\times_{G/P_3}\cdots\times G/\qy{P_{m-1}\cap P_m}\times_{G/P_m}G/B
		\to G/B$ 
	is $B$-equivariant, the pre-image of any point in $G/B$ is isomorphic
	to the pre-image of a $T$-fixed point $w_m$ in $G/B$.
In particular, the locus is $B$-invariant and hence a union of Schubert
	cells.
Consequently, the smallness condition is equivalent to the codimension
	of the Schubert variety indexed by $x$ in the Schubert variety indexed
	by $w$, i.e.\ the difference $\abs{N^w_B}-\abs{N^x_B}$ of the
	dimensions, being more than twice the dimension of the pre-image of
	$x$.

But the pre-image of $x$ is similarly an iterated base change
\[\xymatrix@!0@C=3.5em{
	&&&&P_1/\qy{P_1\cap P_2}\times_{G/P_2}G/\qy{P_2\cap P_3}\times_{G/P_3}\cdots\times G/\qy{P_{m-1}\cap P_m}\times_{G/P_m}n_xB/B
	\ar[dllll]\ar[dll]\ar[d]\ar[drr]\ar[drrrr]\\
	P_1/\qy{P_1\cap P_2}\ar[dr]&&G/\qy{P_2\cap P_3}\ar[dl]\ar[dr]&&\ar[dl]\cdots\ar[dr]&&G/\qy{P_{m-1}\cap P_m}\ar[dl]\ar[dr]&&n_xB/B\ar[dl]\\
	&G/P_2&&G/P_3&&\cdots&&G/P_m}\]
	and a similar analysis to \cite[Proposition~1.10]{Larson} shows it is
	a connected smooth projective variety with an action of $T$, $T$-fixed
	points $(w'_1,\dots,w'_{m-1})$ with $w_1=\sigma_1\in W_1$ and
	$w_kW_{k+1}=w_{k+1}W_{k+1}$ where $w_m=x$, and concentrator schemes
	given by iterated base change
	$U^{w_1}_{P_1\cap P_2}
		\times_{U^{w_2}_{P_2}}U^{w_2}_{P_2\cap P_3}
		\times_{U^{w_2}_{P_3}}\cdots\times_{U^{w_{m-1}}_{P_m}}n_xB/B$, which
	is the product of $U_\alpha$'s indexed by the disjoint union of sets 
	$\qy{N^{w_1}_{P_1\cap P_2}\setminus N^{w_1}_{P_2}}
		\udj\cdots\udj\qy{N^{w_{m-1}}_{P_{m-1}\cap P_m}\setminus N^{w_{m-1}}_{P_m}}$. 

Since in a Bya\l ynicki--Birula decomposition
	\cite[Corollary~13.49(c)]{Milne} exactly one of the concentrator
	subschemes is an open and dense subscheme, it follows that the
	dimension of the pre-image is the dimension of one of the above
	concentrator subschemes. 
Consequently, the projection morphism is small onto the Schubert cell
	indexed by $w$ if, for every $x$ and sequence
	$(w_1,\dots,w_{m-1},w_m)\in W^m$ with $w_1\in W_1$,
	$w_kW_{k+1}=w_{k+1}W_{k+1}$, and $w_m=x$, we
	have $\ds\abs{N^w_B}-\abs{N^x_B}
		>2\sum_{k=1}^{m-1}\abs{N^{w_k}_{P_k\cap P_{k+1}}
			\setminus N^{w_k}_{P_k}}$ whenever the right-hand side is
		non-zero.

If the right-hand side is zero, then the pre-image of $x$ is a finite
	subscheme of $T$-fixed points, and our
	\cref{cor:recursive-defect-polys} will imply it consists of a unique
	point.
In particular, the morphism will be birational, and hence a
	Gelfand--MacPherson resolution of singularities, if and only if the
	right-hand side is always zero for $x=w$, which our
	\cref{prop:bounded-monic} will imply already follows from smallness.
\subsection{A generalization of Deodhar's framework for Weyl groups}
\label{subsec:weyl-gm}
For each root $\alpha\in\Phi$, the Weyl group contains an involution
	$s_\alpha$, called a \define{reflection}, such that
	$s_\alpha(\alpha)=-\alpha\in\Phi$, $s_\alpha=s_{-\alpha}$
	\cite[Theorem~21.11]{Milne}, and $\set{\beta\in\Phi:U_\beta$ not
	contained in $P} =\set{\beta\in\Phi\setminus\Phi^+:s_\beta\not\in W'}$
	\cite[Theorem~21.91]{Milne}.
Consequently, $N_P^w=\qy{\Phi^+\setminus w(\Phi^+)}
			\cap w\qy*{\set{\beta\in\Phi\setminus\Phi^+:s_\beta\not\in W'}}
=N_B^w\cap w\qy*{\set{\beta\in\Phi\setminus\Phi^+:s_\beta\not\in W'}}$.
But $w\mapsto N_B^w$ satisfies $N_B^{xy}=\Phi^+\setminus xy(\Phi^+)
	=\qy{\Phi^+\setminus x(\Phi^+)}+x(\Phi^+\setminus y(\Phi^+))$ where 
	``$+$'' denotes symmetric difference of sets, i.e.\ $w\mapsto N^w_B$
	is a \define{crossed homorphism} to the left $W$-module of finite
	subsets of $\Phi^+$.
	
Moreover, the set $S=\set{w\in W:\abs{N^w_B}=1}=\set{w:N^w_B=\set w}$ of
	\define{simple reflections} is such that
	$\Delta=\set{\alpha\in\Phi^+:s_\alpha\in S}$ is the base for the
	system of positive roots $\Phi^+$, and hence is a generating set for
	the Weyl group \cite[Corollary~21.39]{Milne}.
Thus, $\lref w=\set{s_\alpha:\alpha\in N^w_B}$ is completely determined
	by being a crossed homomorphism such that $\lref s=\set s$ for $s\in
	S$, from which we may recover 
	$N^w_B=\set{\alpha\in\Phi^+:s_\alpha\in\lref w}$.

Also, for each subset $I\subseteq S$, there is a unique parabolic
	subgroup $P_I$ containing $B$ such that the associated subgroup of the
	Weyl group is $W_I=\left<I\right>$, and all parabolic subgroups arise
	this way \cite[Theorem~21.91]{Milne}. 
The subgroups $W_I\leq W$ are themselves called the \define{parabolic
	subgroups} of $W$. 

We can thus express the dimension formulas obtained in the previous
	section fully in terms of the Weyl group $W$, its set of generators
	$S\subseteq W$, and the crossed homomorphism $\lref[]{}$ to the
	$W$-module of finite subsets of reflections, as follows.

First, we have $N^w_P=\set{\alpha\in\Phi^+
	:s_\alpha\in\lref w\setminus wW_Iw^{-1}}$.
Second, if the parabolic subgroup $P_1,\dots,P_m$ correspond to
	subgroups $W_{J_1},\dots,W_{J_m}$ of the Weyl group, we have
	$N^{w_k}_{P_k\cap P_{k+1}}\setminus N^{w_k}_{P_k}=\set{\alpha\in\Phi^+
:s_\alpha\in\lref{w_k}\cap w_k\qy{W_{J_{k+1}}\setminus W_{J_k\cap J_{k+1}}}w_k^{-1}}$.

The notion of \emph{tight factorization} in our generalization of
	Deodhar's framework then arises as follows.
We first note that sequences $(w_1,\dots,w_m)\in W^m$ such that $w_1\in
	W_{J_1}$, $w_kW_{J_{k+1}}=w_{k+1}W_{J_{k+1}}$, and $w_m=x$ correspond
	bijectively to sequences $(\sigma_1,\dots,\sigma_m)\in
	W_{J_1}\times\cdots\times W_{J_m}$ with $w_k=\sigma_1\cdots\sigma_k$. 
Writing $\sigma\to_Jx$ for such a sequence, for each reflection
	$s_\alpha$ we define its sef of defects 
	$\rdef\sigma{s_\alpha}=\set{k:s_\alpha
		\in\lref{w_{k-1}}\cap w_{k-1}\qy{W_{J_k}\setminus W_{J_{k-1}\cap J_k}}w_{k-1}^{-1}}$,
		and $d_R(\sigma)=\ds\sum_{s_\alpha}\abs{\rdef\sigma{s_\alpha}}$.
Then the sequence of parabolic subgroups $P_{J_1},\dots,P_{J_r}$ determines
	a small Gelfand--MacPherson resolution if and only if 
	$\abs{\lref w}-\abs{\lref x}>2\ds d_R(\sigma)$ for every
	$\sigma\to_Jx$ with $d_R(\sigma)\neq0$.

Moreover, the sets of defects are invariant for equivalence classes
	under the quotient $W_{J_1}\times\cdots\times W_{J_m}\to
	W_{J_1}\times_{W_{J_1\cap J_2}}W_{J_2}\times_{W_{J_2\cap
	J_3}}\cdots\times_{W_{J_{m-1}}\cap J_m}W_{J_m}$,.
Thus, for each $x\in W$ we have a generating function
	$P^J_x=\ds\sum_{[\sigma]\to_Jx}q^{d_R(\sigma)}$
	where the sum is over equvalence classes $[\sigma]$ for the quotient
	map such that $\sigma_1\cdots\sigma_m=x$.
These are our generalization of Deodhar's \emph{defect polynomials}
	(\cref{defn:defect-polynomials}) and the formula determining smallness
	above will coincide with our definition of tightness
	(\cref{prop:bounded-monic}).

\section{Preliminaries on Coxeter systems and defects}

\label{sec:prelim}
We review the notation and basic facts concerning Coxeter systems that
	we will use.
We also prove five elementary but non-standard technical results
	(\cref{lem:connected}, \cref{lem:length-coset-rep},
	\cref{lem:conjugated-lref-coset-rep}, \cref{prop:T_L},
	\cref{lem:wz-w-inequality}) relevant to our generalization of
	Deodhar's notion of defects.

For the most part we follow the notation of the standard reference 
	\cite{Bjorner--Brenti}, but emphasize the set of associated
	reflections of a Coxeter system because of their central role in our
	generalization of Deodhar's framework.
In particular, we \emph{define} Coxeter systems using Dyer's \cite{Dyer}
	characterization in terms of those sets and also sketch proofs
	in those terms of the standard results (which may be found in
	\cite{Bjorner--Brenti}) that we review.

\begin{defn}[\cite{Dyer}]
\label{defn:dyer}
	A \define{Coxeter system} is a group $W$
		equipped with a generating set $S$ of order $2$ elements called
		\define{simple reflections} for which there exists a (necessarily
		unique) function sending elements $w\in W$ to their 
		\define{sets of right associated reflections} $\rref w$ satisfying 
	\begin{enumerate}
		\item $\rref s=\set{s}$ for each $s\in S$, and
		\item $\rref{xy}=y^{-1}\rref xy+\rref y$, where:
	\end{enumerate}
		``$+$'' denotes \define{symmetric difference} of sets, each
		$\rref w $ is a subset of the \define{set of reflections
		$T=\displaystyle\bigcup_{w\in W}w^{-1}Sw$} of the system, and
		$y^{-1}\rref{x}y=\set{y^{-1}ty:t\in\rref x}$.
\end{defn}

\begin{exmp}
\label{exmp:type-A-tr}
	The group $S_n$ of permutations of the integers $\set{1,\dots,n}$ is
		generated by the set $S$ of \define{adjacent transpositions}
		$S=\set{(i,i+1)\in S_n:1\leq i\leq n-1}$ where
		$(i,i+1)$ sends $i\mapsto i+1$, $i+1\mapsto i$, and $j\mapsto j$
		when $j\neq i,i+1$.
	The permutations $(i,j)$ (defined when $i\neq j)$ sending 
		$i\mapsto j$, $j\mapsto i$, and $k\mapsto k$ when $k\neq i,j$ are
		the \define{transpositions} of $S_n$ and comprise the union
		$\displaystyle T=\bigcup_{w\in S_n}w^{-1}Sw$ of all conjugates of
		adjacent transpositions.

	The pairs $(S_n,S)$ are the Coxeter systems of \define{type $A$};
		the sets of right associated reflections to a permutation $w$ are
		given by $\rref w=\set{(i,j)\in S_n:i<j,w(i)>w(j)}$ and can be
		naturally identified with the set of \define{position-based
		inversions} of $w$, i.e.\ with the set of unordered pairs
		$\set{i,j}$ of positions at which the outputs of $w$ are put
		out-of-order.

	Indeed, a pair $\set{i,j}$ with $i<j$ can be a position-based
		inversion of $xy$ in exactly two ways: either $y(i)>y(j)$ and
		$xy(i)>xy(j)$, or $y(i)<y(j)$ and $xy(i)>xy(j)$, i.e.\ either
		$\set{i,j}$ is a position-based inversion of $y$, or
		$\set{y(i),y(j)}$ is a position-based inversion of $x$, but not
		both.
	Since $(i,j)=y^{-1}(y(i),y(j))y$, applying the proposed identification
		of position-based inversions with right-associated reflections to
		the above disjunction yields the defining identity
		$\rref{xy}=y^{-1}\rref{x}y+\rref{y}$ of right associated
		reflections.
\end{exmp}

Dyer's characterization of Coxeter systems is self-dual in the sense of
	\cref{lem:self-dual} below.
We will freely use the duals of the properties we prove and will
	generally prove only one result of any dual pair.
\begin{lem}
\label{lem:self-dual}
	Given a Coxeter system $(W,S)$, defining the \define{set of left
		associated reflections by $\lref w=\rref{w^{-1}}$}
		gives us an \define{opposite Coxeter system $(W^{op},S)$} satisfying
	\begin{enumerate}
		\item $\lref s=\set s$ for $s\in S$, and
		\item $\lref{xy}=\lref x+x\lref yx^{-1}$.
	\end{enumerate}
	Furthermore, $\rref e=\emptyset=\lref e$ and
		$w\rref ww^{-1}=\lref w$ so that 
		$\abs{\rref w}=\abs{\lref w}$ and
		$\abs{\rref{xy}}=\abs{\rref x+\lref y}=\abs{\lref{xy}}$ (where $e$
		is the identity of $W$).
	\begin{proof}
	The main claim follows from the definitions; the supplemental
		claims follow from applying the main claim to $ee=e=ww^{-1}$.
	\end{proof}
\end{lem}

\begin{exmp}
	In the case of Coxeter systems of type $A$, the set of left associated
		reflections $\lref w$ can be identified with the set of
		\define{output-based inversions} of $w$, i.e.\ is given by
		$\lref w=\set{(w(i),w(j))\in S_n:i<j,w(i)>w(j)}$, which is naturally
		in bijection with the set of unordered pairs $\set{w(i),w(j)}$ of
		outputs of $w$ that are out-of-order.
	The bijection $\lref w\cong\rref w$ given by $\lref w=w\rref ww^{-1}$
		then corresponds to the evident bijection between output-based and
		position-based inversions.
\end{exmp}

The first non-standard technical result we need
	(\cref{lem:connected}) concerns parabolic subgroups and its proof
	relies on properties of reduced words which we review below.
We also review the Bruhat order on a Coxeter system, which we will also
	use.
\begin{defn}
\label{defn:coxeter_features}
	Given a Coxeter system $(W,S)$, its \define{length function}
		$\ell\colon W\to\NN$ is given by $\ell(w)=\abs{\rref w}$ ($\rref w$
		are finite by \cref{lem:subset} below).
	A \define{word} $s_{i_1}\cdots s_{i_r}$ of generators $s_{i_k}\in S$
		is \define{reduced} if $\ell(s_{i_1}\cdots s_{i_r})=r$.
\end{defn}

\begin{lem}
\label{lem:subset}
	For each word $s_{i_r}\cdots s_{i_r}=y$, we have $\rref y
		=\set{t_1}+\cdots+\set{t_r}\subseteq\set{t_k\in T:
			t_k=s_{i_r}\cdots s_{i_{k+1}}s_{i_k}s_{i_{k+1}}\cdots s_{i_r}}$. 
	Consequently $\ell(s_{i_1}\cdots s_{i_r})\leq r$, with the inequality
		strict if and only it there exist $j$ and $k$ such that $1\leq
		j<k\leq r$ and $t_j=t_k$, i.e.\ $s_{i_1}\cdots s_{j-1}s_{j+1}\cdots
		s_{k-1}s_{k+1}\cdots s_{i_r}=y$.
	In particular every word for an element has a reduced subword for that
		same element, and a word $s_{i_1}\cdots s_{i_r}$ is reduced if and
		only if $r=\min\set{k:s_{i_1}\cdots s_{i_r}=s_{j_1}\cdots s_{j_k}}
			=\ell(s_{i_1}\cdots s_{i_r})$.
\begin{proof}
	Immediate by induction on $r$ using \cref{defn:dyer}.
\end{proof}
\end{lem}

\begin{lem}
\label{lem:injective-cocycle}
	$S\cap \rref w=\emptyset$ if and only if $w=e$; moreover 
		$T_R\colon W\to 2^T$ injective.
	\begin{proof}
		One direction of the main claim follows from \cref{lem:subset}
			applied to a reduced word and the other from the supplemental
			claim of \cref{lem:self-dual}. 
		The supplemental claim follows by induction on $\ell(x)=\ell(y)$
			using the main claim to pick $s\in\rref x=\rref y$ yielding 
			$\rref{xs}=s\rref xs\setminus\set s=s\rref ys\setminus\set s
				=\rref{ys}$.
	\end{proof}
\end{lem}

\begin{rmrk}
\label{rmrk:tr=tl}
	\cref{lem:injective-cocycle} also implies that $w=w^{-1}$ if and
		only if $\rref w=\rref{w^{-1}}$, i.e.\ if and only if 
		$\rref w=\lref w$.
\end{rmrk}

\begin{defn}
\label{defn:bruhat}
	The \define{Bruhat order} on $W$ is the reflexive and transitive
		relation generated by $x<y$ if $x=yt$ for some $t\in\rref y$.
\end{defn}

\begin{lem}
\label{lem:t-in-rref-t}
	For any $t \in T$ we have that $t\in \rref t$, and in general
		for any $x\in W$ we have $t\in\rref{xt}=t\rref xt+\rref t$
		if and only if $t\not\in\rref x$, i.e.\ $x$ and $xt$ are
		comparable in the Bruhat order for any $t\in T$, and $x<xt$ if and
		only if $t\not\in\rref x$. 
	\begin{proof} 
		Using \cref{lem:self-dual} we get
			$\rref t=\rref{w^{-1}sw}=t\rref wt+\set t+\rref w$. 
		Since $t^3=t$ implies that $t\in\rref w$ if and only if 
		$t\in t\rref wt$, it follows that $t\not\in\rref w+t\rref wt$,
			so that $t\in\rref t$, as desired.
	\end{proof}
\end{lem}

\begin{lem}
\label{lem:bruhat}
	Given a reduced word $s_{i_1}\cdots s_{i_r}=y$, the Bruhat relation
		$x\leq y$ holds only if $x$ is realizable as a (reduced) subword of
		the chosen reduced word for $y$.
	Consequently, the Bruhat order is antisymmetric (hence indeed a
		partial order), and $\ell\colon W\to\NN$  is order-preserving ($x<y$
		implies $\ell(x)<\ell(y)$).
	In particular, the Bruhat order is well-founded, i.e.\ every subset of
		$W$ has a minimal element and furthermore the identity $e\in W$ is
		the minimum element of $W$.
	\begin{proof}
		Immediate from combining \cref{defn:bruhat},
			\cref{lem:subset}, and \cref{lem:t-in-rref-t}.
	\end{proof}
\end{lem}

\begin{rmrk}
	A well-founded order has a minimum element if and only if it has a
		unique minimal element.
\end{rmrk}

\begin{rmrk}
	The standard result that the converse of \cref{lem:bruhat} also holds,
		i.e.\ that every element $x$ realizable as a subword of a chosen
		reduced word for $y$ is below $y$ in the Bruhat order follows from
		our \cref{prop:defect-ideals}.
\end{rmrk}

\begin{defn}
	Given a subset $J\subseteq S$, the subgroup $W_J\leq W$ generated by
		$J$ is called a \define{parabolic subgroup}.
\end{defn}

\begin{lem}
\label{lem:subsystem}
	For any subset $J\subseteq S$, we have that $x\in W_J$ if and only if
		$\rref x\subseteq W_J$, i.e.\ if and only if 
		$\rref x=\rref x\cap W_J$.
	In particular, the pair $(W_J,J)$ is a \define{Coxeter
		subsystem} of $(W,S)$, i.e.\ its right associated
		reflections, and hence length function and Bruhat order, are
		restrictions of those of the Coxeter system $(W,J)$.
	In particular, $S\cap W_J=J$ since the generators $J\subseteq W_J$ are
			characterized by $\ell(s)=1$, and hence $W_{J_1}\cap
			W_{J_2}=W_{J_1\cap J_2}$ for any pair of subsets $J_1,J_2
				\subseteq S$.
	\begin{proof}
		$x\in W_J$ implies $\rref x\subseteq W_J$ by \cref{lem:subset}.
		The converse holds by induction on $\ell(x)$ as we may pick
			$t\in\rref x\subseteq W_J$, whence $\rref t\subseteq W_J$ by
			Lemma~\ref{lem:subset} and $\rref{xt}=t\rref xt+\rref t\subseteq
			tW_Jt\cup W_J=W_J$ with $\ell(xt)<\ell(x)$ by
			\cref{lem:t-in-rref-t}.
	\end{proof}
\end{lem}

\begin{defn}
	We define $\nbd J=J\cup\set{s\in S:s$ does not commute with
		some element of $J}$.
	We define the \define{connected components of $J$} to be the
		equivalence classes of the equivalence relation on elements of $J$
		generated by the relation $s\in\nbd{\set{s'}}$, i.e.\ by the
		(symmetric) relation $s\sim s'$ if $s=s'$ or $ss'\neq s's$.

	We say two subsets $J_1,J_2\subseteq S$ 
		\define{commute} if $J_1\cap\nbd{J_2}=\emptyset$, i.e.\  
		if they are disjoint and each element of $J_1$ commutes with each
		element of $J_2$, or equivalently if each is a union of connected
		components of $J_1\cup J_2$.
	We say that a subset $J\subseteq S$ is \define{connected} if it has
		exactly one connected component, i.e.\ if it is neither empty nor
		the union of two commuting non-empty subsets of $S$.

	For every element $w\in W$, the \define{set $\rdes w$ of
		right-descents of $w$} is given by $\rdes w=S\cap \rref w$.
	Similarly, the \define{set $\ldes w$ of left-descents
		of $w$} is given by $\ldes w = S\cap \lref w$.
\end{defn}

\begin{exmp}
\label{exmp:interval_decomposition}
	In the case of a Coxeter system of type $A$, each connected subset
		$J\subseteq S$ corresponds to an interval $[i,j]$ with $i<j$ where
		$s_k=(k,k+1)\in J$ if and only if $i\leq k<k+1\leq j$.
	Furthermore, two intervals are disjoint if and only if their
		corresponding connected subsets of $S$ commute.

	We will sometimes identify a not necessarily connected subset
		$J\subseteq S$ with its \define{interval decomposition}
		$[a_1,b_1],[a_2,b_2],\cdots,[a_m,b_m]$ of intervals corresponding to
		the connected components of $J$, where
		$a_1<b_1<a_2<b_2<\cdots<a_m<b_m$.

	In particular, connected components of the set $\rdes{w}$ of
		right-descents of a permutation $w$ correspond to
		\define{right-descent intervals}, i.e.\ to maximal intervals $[i,j]$
		with $i<j$ such that $w(i)>w(i+1)>\dots>w(j)$.
	The set $\rdes w$ of right-descents is then in bijection with the set
		of adjacent position-based inversions of $w$ contained in a
		right-descent interval of $w$.

	Similarly, the connected components of the set $\ldes w$ of
		left-descents of a permutation $w$ correspond to
		\define{left-descent intervals}, i.e.\ to maximal intervals $[i,j]$
		with $i<j$ such that $w^{-1}(i)>w^{-1}(i+1)>\cdots>w^{-1}(j)$.
	The set $\ldes w$ is then in bijection with all adjacent output-based
		inversions contained in a left-descent interval of $w$.
\end{exmp}

\cref{lem:des-noncommuting-refl} below is a technical result we use in
	the proof of \cref{lem:connected}, which is itself the key technical
	result used in the proofs of \cref{lem:wz-w-inequality} and
	\cref{prop:rigidity}.
\begin{lem}
\label{lem:des-noncommuting-refl}
	Given a reduced word $s_{i_1}\cdots s_{i_r}=x$ and $s\in S$, 
		the word $ss_{i_1}\cdots s_{i_r}s$ is reduced if and only
		if $sx\neq xs$ and $s\not\in\rref x\cup\lref x$, in which case
		$s\in\rref{sxs}\cap\lref{sxs}$.
	In particular, if $s\in\nbd{\set{s_{i_1},\dotsc,s_{i_r}}}
			\setminus\set{s_{i_1},\dotsc,s_{i_r}}$ 
		then $ss_{i_1}\cdots s_{i_r}s$ is a reduced word with
		$s\in\rref{ss_{i_1}\cdots s_{i_r}s}
			\cap\lref{ss_{i_1}\cdots s_{i_r}s}$.
	\begin{proof}
		By \cref{lem:subset}, the word $ss_{i_1}\cdots
			s_{i_r}s=sxs$ is reduced if and only if
			$\abs{\rref x}+2=\abs{\rref{sxs}}=\abs{s\rref{sx}s+\set s}
				=\abs{s\qy*{x^{-1}\set sx+\rref x}s+\set s}
				=\abs{sx^{-1}\set sxs+s\rref xs+\set s}$. 
		Since $\abs{\rref x}=\abs{s\rref xs}$, the equality is equivalent to
			having both $s\not\in\rref x\cup x\rref xx^{-1}
					=\rref x\cup\lref x$ and $s\neq sx^{-1}sxs$, i.e.\ and 
				$sx\neq xs$, in which case $\rref{sxs}
					=x^{-1}\set sxs\udj s\rref xs\udj\set s\ni s$. 

		For the supplemental claim, by \cref{lem:subsystem}
			$s\in\nbd{\set{s_{i_1},\dotsc,s_{i_r}}}
				\setminus\set{s_{i_1},\dotsc,s_{i_r}}$ implies that
			$s\not\in W_{\set{s_{i_1},\dots,s_{i_r}}}
				\supseteq\lref x\cup\rref x$. 
		To verify it also implies $sx\neq xs$, it will suffice to show
			it implies $s\rref xs\neq \rref x$ since $sx=xs$ implies $\rref x
				=\rref{sxs}=\set{sx^{-1}sxs}+s\rref xs+\set s=s\rref xs$.

		To that end, let $k$ be maximal such that
			$s\in\nbd{\set{s_{i_k}}}$, i.e.\ such that 
			$s\in\nbd{\set{s_{i_k}}}
				\setminus\nbd{\set{s_{i_{k+1}},\dots,s_{i_r}}}$.
		Then $t_k=s_{i_r}\cdots s_{i_k}\cdots s_{i_r}\in\rref x$
			satisfies $\rref{t_k}\cup\lref{t_k}
				\subseteq W_{\set{s_{i_k},\dots,s_{i_r}}}\not\ni s$ by
				\cref{lem:subsystem}, and
			$st_k=s_{i_r}\cdots ss_{i_k}\cdots s_{i_r}
				\neq s_{i_r}\cdots s_{i_k}s\cdots s_{i_r}=t_ks$ by definition. 
		Thus the main claim implies $W_{\set{s_{i_k},\dots,s_{i_r}}}
				\not\ni s\in\rref{st_ks}$.
		Since by \cref{lem:subsystem} 
			$st_ks\in W_{\set{s_{i_1},\dots,s_{i_r}}}$ would imply the
			contradictory $s\in\rref{st_ks}
				\subseteq W_{\set{s_{i_1},\dots,s_{i_r}}}\not\ni s$,
			it follows that $s\rref xs\ni st_ks\not\in
			W_{\set{s_{i_1},\dots,s_{i_r}}}\supseteq\rref x$, and hence that
			$s\rref xs\neq\rref x$, as desired. 
	\end{proof}
\end{lem}

\begin{lem}
\label{lem:connected}
	Given a subset $J\subseteq S$ and a connected $J'\subseteq S$ not
		contained in $J$, for each $s\in J\cap\nbd{J'}\setminus J'$ there
		exist $s'\in J'\setminus J$ and $z\in W_{(J'\cap J)\cup\set s}$ such
		that $s\in\rref{z^{-1}s'z}$.
	\begin{proof}
		Let $s_{i_0}=s\in J\cap\nbd{J'}\setminus J'$, and observe that by
			definition there exists $s_{i_1}\in J'\cap \nbd{\set{s_{i_0}}}$.
		Since $J'$ is connected and $J'\setminus J$ is non-empty while
			$s_{i_1}\in J'$, it follows that the sequence $s=s_{i_0},s_{i_1}$
			can be extended to a sequence $s=s_{i_0},s_{i_1},\dots,s_{i_r}$
			such that $s_{i_k}\in J'\cap\nbd{\set{s_{i_{k-1}}}}$ for every $k$
			such that $1\leq k\leq r$ and such that
			$s_{i_r}\in J'\setminus J$.

		Note that if $s_{i_k}\not\in J'\cap J$ for some $k$ such that
			$1\leq k<r$, then $s_{i_k}\in J'\setminus J$ and hence
			$s_{i_0},\dots,s_{i_{k-1}},s_{i_k}$ would be a shorter sequence
			with the same property.
		Similarly, if $s_{i_j}\in\set{s_{i_{j+1}},\dots,s_{i_r}}$, e.g.\ if
			$s_{i_j}=s_{i_k}$ where $j<k$, then 
			$s_{i_0},\dots, s_{i_j}, s_{i_{k+1}},\dots, s_{i_r}$
			would also be a shorter sequence with the same property.
		Thus, any sequence of minimal length satisfying the property
			satisfies $s_{i_k}\in J'\cap J$ for $1\leq k<r$ and
			$s_{i_j}\in\nbd{\set{s_{i_{j+1}}}}
				\setminus\set{s_{i_{j+1}},\dots,s_{i_r}}$ for $0\leq j<r$.
		Consequently, taking $s'=s_{i_r}\in J'\setminus J$ for such a
			minimal sequence of length $r$, we have
			$z=s_{i_{r-1}}\cdots s_{i_0}\in W_{\qy{J\cap J'}\cup\set{s}}$, and
			by induction using \cref{lem:des-noncommuting-refl} that
			$z^{-1}s'z=s_{i_0}\cdots s_{i_r}\cdots s_{i_0}$ is a reduced word
			with $s=s_{i_0}\in\rref{s_{i_0}\cdots s_{i_r}\cdots s_{i_0}}
				=\rref{z^{-1}s'z}$.
	\end{proof}
\end{lem}

The next pair of technical results we need
	(\cref{lem:length-coset-rep} and
	\cref{lem:conjugated-lref-coset-rep}) have to do with minimal
	coset decompositions in a Coxeter system.
\begin{defn}
	Given a subset $J\subseteq S$, a \define{minimal left coset
		representative} of $x\in W$ is an element $\rcos Jx\in xW_J$ such
		that $\rref{\rcos Jx}\cap W_J=\emptyset$, in which case we denote by
		$x_J\in W_J$ the element such that $x=\rcos Jxx_J$.
	The \define{set of minimal left coset representatives} is denoted by
		$\rcos JW$, i.e.\ $\rcos JW=\set{x\in W:\rref x\cap W_J=\emptyset}$.

	Dually, the \define{minimal right coset representative} of $x\in W$ is
		an element $\lcos Jx\in W_Jx$ such that $\lref{\lcos Jx}\cap
		W_J=\emptyset$ and the \define{set of
		minimal left coset representatives} is denoted by $\lcos JW$, i.e.\ 
		$\lcos JW=\set{x\in W:\lref x\cap W_J=\emptyset}$.
\end{defn}

\begin{exmp}
	For a Coxeter system of type $A$, the set $W^J$ of minimal left coset
		representatives consists of those permutations $x$ such that
		$x(a_i)<x(a_i+1)<\dots<x(b_i)$ for each interval $[a_i,b_i]$ of
		the interval decomposition describing $J$.

	Accordingly, the minimal left coset representative $\rcos Jw$ of a
		permutation $w$ is then the permutation whose one-line notation
		$w^J(1),\dots,w^J(n)$ is obtained from the one-line notation
		$w(1),\dots,w(n)$ of $w$ by sorting in increasing order the
		sequences of outputs $w(a),\dots,w(b)$ for each interval $[a,b]$
		corresponding to a connected component of $J$.
	In other words, the permutation $w_J\in W_J$ such that $w=\rcos Jww_J$
		arranges the increasing sequence
		$\rcos Jw(a)<\dots<\rcos Jw(b)$ into 
		the sequence of outputs $w(a),\dots,w(b)$ for each interval $[a,b]$
		corresponding to a connected component of $J$.
\end{exmp}

\begin{lem}
\label{lem:coset-properties}
	For any subset $J\subseteq S$ we have that:
	\begin{enumerate}
		\item\label{lem:coset-exist}
			$\rcos JW\cap xW_J$ is non-empty and consists of the minimal
				elements of $xW_J$ under the Bruhat order and in particular
				those of minimal length; 
		\item\label{lem:coset-unique}
			for $\rcos Jx\in \rcos JW\cap xW_J$ and $y\in W_J$ we have
			$\rref{\rcos Jxy}=y^{-1}\rref{\rcos Jx}y\bigudj\rref y$, i.e.\ 
			$\ell(\rcos Jxy)=\ell(\rcos Jx)+\ell(y)$; in particular, the
			element $x_J$ such that $x=\rcos Jxx_J$ is uniquely determined by
			$\rref{x_J}=\rref x\cap W_J$ whence $\rcos JW\cap xW_J$ consists
			only of the unique minimal (and hence minimum) element $\rcos Jx$
			of $xW_J$ under the (well-founded) Bruhat order; furthermore the
			group operation $W\times W\to W$ restricts to a bijection $\rcos
			JW\times W_J\cong W$ with inverse given by 
			$(\rcos Jx,x_J)\mathrel{\reflectbox{\ensuremath{\mapsto}}}x$ and
			satisfying $\ell(\rcos Jx)+\ell(x_J)=\ell(x)$;
		\item\label{lem:component-decompose}
		$W_{S\setminus J}\subseteq\rcos JW$ and when $J$ commutes with
			$S\setminus J$ the group operation restricts to a group
			isomorphism $W_{S\setminus J}\times W_J\to W$ and we have 
			$\rcos JW=W_{S\setminus J}$ and $\rcos{S\setminus J}W=W_J$;
			in particular the inverse of such a group isomorphism is given by
			$(\rcos Jx,\rcos{S\setminus J}x)=(x_{S\setminus J},x_J)
			\mathrel{\reflectbox{\ensuremath{\mapsto}}}x$, and furthermore
			$T\cap W=\qy*{T\cap W_{S\setminus J}}\udj\qy*{T\cap W_J}$, so that
			every parabolic subgroup $W_J$ decomposes as the (internal) direct
			sum $W_J=\bigoplus W_{J_i}$ for $J_i$ (unions of) connected
			components of $J$, in which case $T\cap W_J=\bigudj T\cap
			W_{J_i}$;
	\end{enumerate}
	\begin{proof}
(\ref{lem:coset-exist}) 
	$t\in\rref x\cap W_J$ implies $xt\in xW_J$ by definition, whence
		$xt<x$ by \cref{defn:bruhat} and $\ell(xt)<\ell(x)$ by
		\cref{lem:bruhat}.

(\ref{lem:coset-unique}) 
	$\rref{\rcos Jyx}=y^{-1}\rref{\rcos Jx}y\udj\rref y$ is immediate from
	$y\in W_J$ implying $\rref y\subseteq W_J$ by
	\cref{lem:subsystem} and $\rref{\rcos Jx}\cap W_J=\emptyset$,
	while uniqueness of $x_J$ follows from
	\cref{lem:injective-cocycle}.

(\ref{lem:component-decompose})
	$W_{S\setminus J}\subseteq\rcos JW$ because by
		\cref{lem:subsystem} $x\in W_{S\setminus J}$ implies
		$\rref x\cap W_J\subseteq T\cap W_{S\setminus J}\cap W_J
			=T\cap W_{(S\setminus J)\cap J}=T\cap W_\emptyset=\emptyset$. 
	When $J$ commutes with $S\setminus J$, it follows from 
		$W_{S\setminus J}\cap W_J\subseteq\rcos JW\cap W_J=\emptyset$ and
		consideration of reduced words that the group operation restricts to
		an isomorphism $W_{S\setminus J}\times W_J\cong W_J$. 
	It is then easy to check that $\rref x=\rref{x_{S\setminus J}x_J}
			=\qy{\rref x\cap W_{S\setminus J}}\udj\qy{\rref x\cap W_J}$
		because by \cref{lem:subsystem} $x_J\in W_J$ commutes with
		each element of $\rref{x_{S\setminus J}}\subseteq W_{S\setminus J}$.
	In particular, $t\in T\cap W$ implies by \cref{lem:t-in-rref-t}
		that $t\in\qy*{\rref t\cap W_{S\setminus J}}
			\udj\qy*{\rref t\cap W_J}$.
		\end{proof}
\end{lem}

\cref{lem:length-coset-rep,lem:conjugated-lref-coset-rep} are our
	second and third technical results establishing the connection between
	sets of associated reflections and the notion of defects used in the
	proof of \cref{lem:y-well} 
\begin{lem}
\label{lem:length-coset-rep}
	Given $x\in W$ and subsets $J, J_*\subseteq S$, we have 
		$(x_J)_{J_*}=x_{J\cap J_*}$ and that
			$\rcos[J]{J_*}x=\rcos{J_*}{(x_J)}=\rcos{J\cap J_*}{(x_J)}$
		satisfies $\rref{\rcos[J]{J_*}x}=x_{J\cap J_*}\qy{\rref x
			\cap W_J\setminus W_{J_*}}x_{J\cap J_*}^{-1}$ and so has length
		$\len{\rcos[J]{J_*}x}=\abs{\rref x\cap W_J\setminus W_{J_*}}$.
	\begin{proof}
		\cref{lem:coset-properties}(\ref{lem:coset-unique}) and
			\cref{lem:subsystem} yield 
			$\rref{(x_J)_{J_*}}=\rref{x_J}\cap W_{J_*}
				=\rref x\cap W_J\cap W_{J_*}
				=\rref x\cap W_{J\cap J_*}=\rref{x_{J\cap J_*}}$, so that
			$(x_J)_{J_*}=x_{J\cap J_*}$ by \cref{lem:injective-cocycle}.
		\cref{lem:coset-properties}(\ref{lem:coset-unique}) then gives
			$\rcos[J]{J_*}x=\rcos[J]{J\cap J_*}x$ and $\rref{x_J}
				=\rref{\rcos[J]{J_*}x\rcos[J\cap J_*]{}x}
				=(x_{J\cap J_*})^{-1}\rref{\rcos[J]{J_*}x}x_{J\cap J_*}
					\udj\rref{x_{J\cap J_*}}
				=(x_{J\cap J_*})^{-1}\rref{\rcos[J]{J_*}x}x_{J\cap J_*}
					\udj\qy*{\rref{x_J}\cap W_{J_*}}$
		so that $\rref{\rcos[J]{J_*}x}
				=x_{J\cap J_*}\qy*{\rref{x_J}\setminus W_{J_*}}
					x_{J\cap J_*}^{-1}
				=x_{J\cap J_*}\qy*{\rref x\cap W_J\setminus W_{J_*}}
						x_{J\cap J_*}^{-1}$.
	\end{proof}
\end{lem}

\begin{rmrk}
\label{rmrk:rcos-notation}
	It is \emph{not} the case that
		$\rcos{J_*}{\qy*{\rcos[J]{}x}}=\rcos[J]{}{\qy*{\rcos{J_*}x}}$,
		but we will always use the notation $\rcos[J]{J_*}x$ to refer to
		$\rcos{J_*}{\qy*{\rcos[J]{}x}}$.
\end{rmrk}

\begin{lem}
\label{lem:conjugated-lref-coset-rep}
	Given $x\in W$ and subsets $J,J_*,J'\subseteq S$, we have
		$\lref x\cap x\qy{W_J\setminus W_{J_*}}x^{-1}
			=x\qy*{\rref x\cap W_J\setminus W_{J_*}}x^{-1}
			=\rcos{J\udj\qy*{J'\setminus\nbd J}}x
				\lref{\rcos[J]{J_*}x}\qy*{\rcos{J\udj\qy*{J'\setminus\nbd J}}
					x}^{-1}$.
	\begin{proof}
		\cref{lem:length-coset-rep} implies $\rref{\rcos[J]{J_*}x}
					=x_{J\cap J_*}\qy*{\rref x\cap W_J\setminus W_{J_*}}
							x_{J\cap J_*}^{-1}$ and so by applying
			\cref{lem:coset-properties}(\ref{lem:coset-unique}) and
			\cref{lem:self-dual} we obtain
			\begin{multline*}
				x\qy*{\rref x\cap W_J\setminus W_{J_*}}x^{-1}
				=\rcos{J\cap J_*}x\rref{\rcos[J]{J_*}x}
						\qy*{\rcos{J\cap J_*}x}^{-1}\\
				=\rcos Jx\rcos[J]{J_*}x\rref{\rcos[J]{J_*}x}
					\qy*{\rcos[J]{J_*}x}^{-1}\qy*{\rcos Jx}^{-1}
				=\rcos Jx\lref{\rcos[J]{J_*}x}\qy*{\rcos Jx}^{-1}.
			\end{multline*}
			Since $\rcos[J]{J_*}x\in W_J$ implies
				$\lref{\rcos[J]{J_*}x}\subseteq W_J$ by
				\cref{lem:subsystem}, we have $\lref{\rcos[J]{J_*}x}
					=x_{J'\setminus\nbd J}
						\lref{\rcos[J]{J_*}x}x_{J'\setminus\nbd J}^{-1}$ 
				and hence $\rcos Jx\lref{\rcos[J]{J_*}x}\qy*{\rcos Jx}^{-1}
					=\rcos{J\udj J'\setminus\nbd J}x
						\lref{\rcos[J]{J_*}x}
						\qy*{\rcos{J\udj J'\setminus\nbd J}x}^{-1}$.
	\end{proof}
\end{lem}

We will also make use of the interaction of minimal coset decompositions
	with the Bruhat order.
\begin{prop}
\label{prop:bruhat-coset}
	Given a subset $J\subseteq S$, $x\in\rcos JW$, $y\in W_J$, and 
		$w\in W$, then $xy\leq w$ implies $x\leq\rcos Jw$, i.e.\ the
		minimal coset representative map $W\twoheadrightarrow\rcos JW$
		preserves the Bruhat order. 
	If $x=\rcos Jw$, then $xy\leq w$ holds if and only if $y\leq w_J$
		holds, i.e.\ $W_J\xhookrightarrow{x\cdot}W$ preserves and reflects
		the Bruhat order.
	Furthermore, $x\leq\rcos Jw$ and $y\leq w_J$ imply $xy\leq w$, e.g.\
		$\rcos JW\xhookrightarrow{\cdot y}W$ preserves and reflects the
		Bruhat order.
	\begin{proof}
		The main claim follows immediately from analyzing the generating
			for the Bruhat order relation $xy=wt$ for $t\in\rref w$ in
			combination with
			\cref{lem:coset-properties}(\ref{lem:coset-unique}).

		The supplemental claim follows easily in the singleton case
			$J=\set s$ using the fact that $xy=wt$ for $y\in W_J$ and
			$t\in\rref w$ is a generating relation $x\prec w$ for the Bruhat
			order on $\rcos JW$.
		The case of arbitrary $J$ follows from reducing to the case $y=w_J$,
			then inducting on $\ell(w_J)$ using
			\cref{lem:injective-cocycle} to pick $s\in
			S\cap\rref{w_J}\subseteq S\cap W_J=J$, and applying the singleton
			case $\set s$ to $xw_Js<\rcos Jww_Js$.
	\end{proof}
\end{prop}

\begin{rmrk}
\label{rmrk:lift}
	Unless $x=\rcos Jw$, it is not necessarily the case that $xy\leq w$
		with $x\in\rcos JW$ and $y\in W_J$ implies $y\leq w_J$, i.e.\ the
		map $W\twoheadrightarrow W_J$ given by $w\mapsto w_J$ does not
		necessarily preserve the Bruhat order.
	For example, $s_1s_2\neq s_2s_1$ if and only if
		$s_1\not\in\rref{s_1s_2}=s_2\set{s_1}s_2\udj\set{s_2}$, in which
		case $x=s_1<s_1s_2=w$ but $w_J=e<s_1=x_J$ for $J=\set{s_1}$.
\end{rmrk}

The fourth and fifth technical results we need
	(\cref{prop:T_L} and \cref{lem:wz-w-inequality})
	concern finite parabolic subgroups and their maximal elements.
The former will be used in a variety of results concerning right-descent
	factorizations, while the latter is the key result needed  for the
	proof of \cref{prop:invariant_mult_connected}, which builds up to our
	main result, \cref{thm:2}.

The supplemental claim of \cref{prop:finitary} below is a
	somewhat stronger formulation that will more useful to us than the
	usual one.
\begin{defn}
\label{defn:maximal}
	The \define{maximal element} of a parabolic subgroup $W_J$ is the
		(necessarily unique by \cref{lem:injective-cocycle}) element 
		$\mxml J\in W$ such that $\rref{\mxml J}=T\cap W_J$.
	Following Williamson \cite{Williamson}, we say that $J\subseteq S$ is
		\define{finitary} if $\mxml J$ exists.
\end{defn}

\begin{exmp}
	The maximal element $\mxml J$ of a parabolic subgroup $W_J$ in type $A$
		is the unique permutation reversing the sequences $a,\dots,b$
		for which the intervals $[a,b]$ correspond to connected components
		of $J$.
	Accordingly, $\rcos Jw\mxml J$ is the permutation obtained from $w$ by
		sorting in \emph{decreasing order} the outputs $w(a),\dots,w(b)$ for
		each interval $[a,b]$ corresponding to a connected component of $J$.
\end{exmp}

\begin{prop}
\label{prop:finitary}
	Given a subset $J\subseteq S$, the maximal element $\mxml J$ of the
		parabolic subgroup $W_J$ is equivalently characterized by maximality
		of $\mxml J$ in the Bruhat order on $W_J$, and also by maximality of
		$\ell(\mxml J)$.
	In particular, its existence, i.e.\ $J$ being finitary, is equivalent
		to $W_J$ being finite.

	Furthermore, given $w\in W$ we have that $\mxml J$ exists and
		coincides with $w_J$ if and only if $\rref w\supseteq J$, i.e.\ if
		and only if $J\subseteq\rdes w$.
	In particular, $w_{\rdes w}=\mxml{\rdes w}$ so that $\rdes w\subseteq
		S$ (and dually $\ldes w\subseteq S$) is finitary for all $w\in W$.
	\begin{proof}
		The main claim follows immediately from \cref{lem:t-in-rref-t,lem:bruhat}, and from $\ell(\mxml J)=\abs{T\cap W_J}$
			implying finiteness of $2^{T\cap W_J}$ and hence of $W_J$ by
			\cref{lem:injective-cocycle}.

		For the supplemental claim, since
			\cref{lem:coset-properties}(\ref{lem:coset-unique}) implies
			$\rref{w_J}=\rref w\cap W_J$, we have by definition that
			$w_J=\mxml J$ if and only if $T\cap W_J=\rref w\cap W_J$, i.e.\ if
			and only if $T\cap W_J\subseteq\rref w$, in which case 
			$J\subseteq\rref w$ by \cref{lem:subsystem}. 

		Conversely, it will suffice by induction on $\ell(v)$ to show
			that $J\subseteq\rdes w$ implies $v\leq w_J$ for all $v\in W_J$.
		The base case $\ell(v)=0$ holds because by
			\cref{lem:injective-cocycle} $\ell(v)=0$ implies $v=e$
			and \cref{lem:bruhat} implies $e\leq w_J$.
		For the inductive step, suppose that $s\in\rdes v\subseteq J
			\subseteq\rdes{w_J}$, where the first containment holds by
			\cref{lem:subsystem}.
		Then $v=\rcos{\set s}vs$ with $\len{\rcos{\set s}v}<\ell(v)$ by
			\cref{lem:bruhat}, so inductively we have 
			$\rcos{\set s}v\leq w_J$, which by
			\cref{prop:bruhat-coset}
			implies $v=\rcos{\set s}vs\leq\rcos{\set s}{w_J}s=w_J$ as desired.
	\end{proof}
\end{prop}

\begin{prop}
\label{prop:T_L}
	Given a Coxeter element $w\in W$, and a subset $J\subseteq S$, an
		element $v\in W$ satisfies $T\cap W_J
			=\qy{\rref w\cap W_J}\udj\lref v$
		if and only if $J$ is finitary and $wv=\rcos Jw\mxml J$, in which
		case $J=\qy{\rdes w\cap J}\udj\ldes v$. 
	In particular, $\rdes{\rcos Jw\mxml J}=\qy*{\rdes w\cap
		\rdes{\rcos Jw\mxml J}}\udj\ldes v\subseteq\rdes w\cup J$.
	\begin{proof}
		By definition of $\mxml J$ and by
			\cref{lem:coset-properties}(\ref{lem:coset-unique})
			we have that $J$ is finitary and $\rcos Jw\mxml J=wv$ if and only
			if $\rref{wv}\cap W_J=T\cap W_J$, in which case $v=w_J^{-1}\mxml J
				\in W_J$.
		Since for any $v\in W_J$ we have $T\cap W_J=v\qy{T\cap W_J}v^{-1}$
			and $v\rref vv^{-1}=\lref v\subseteq W_J$ by
			\cref{lem:self-dual} and \cref{lem:subsystem}, it
			follows that
			\begin{align*}
				v\qy{\rref{wv}\cap W_J}v^{-1}
					&=v\qy{\rref{wv}}v^{-1}\cap W_J
					=v\qy{v^{-1}\rref wv+\rref v}v^{-1}\cap W_J\\
				&=\qy{\rref w+v\rref vv^{-1}}\cap W_J
					=\qy{\rref w\cap W_J}+\lref v
			\end{align*}
		Consequently, $\rcos Jw\mxml J=wv$, i.e.\ $T\cap
			W_J=\rref{wv}\cap~W_J$, holds if and only if 
			$T\cap W_J=\qy{\rref w\cap W_J}+\lref v$, i.e.\ if and only if 
			$T\cap W_J=\qy{\rref w\cap W_J}\udj\lref v$.

		\cref{lem:subsystem} then also implies 
			$J=T\cap W_J\cap S=\qy{\rref w\cap W_J\cap S}
				\udj\qy{\lref v\cap S}=\qy{\rdes w\cap J}\udj\ldes v$.
		Applying the main claim to $w$ and the (finitary by
			\cref{prop:finitary}) $\rdes{\rcos Jw\mxml J}\subseteq S$
			in place of $J\subseteq S$ yields $\rdes{\rcos Jw\mxml J}
				=\big(\rdes w\cap\rdes{\rcos Jw\mxml J}\big)\udj\ldes
				v\subseteq\rdes w\cup J$.
	\end{proof}
\end{prop}

\begin{lem}
\label{lem:wz-w-inequality}
	Given $w\in W$, a subset $J_*\subseteq\rdes w$, and a finitary subset
		$J\subseteq S$, then for each $z\in W_{J_*\cap \rdes{w^J\mxml J}}$,
		we have that $\len{\rcos[J]{J_*}{(wz)}}-\len{\rcos[J]{J_*}w}
			=\abs{\qy*{T\cap W_J\setminus\rref w}\setminus zW_Jz^{-1}}$.
	In particular, if $J'$ is a connected component of $J$ not contained
		in $\rdes w$, then for each $s\in\rdes w\cap\nbd{J'}\setminus J'$
		there exists $z\in W_{\qy{J'\cap\rdes w}\cup\set s}
			\subseteq W_{\rdes w}$ such that
		$\len{\rcos[J]{\rdes w}{\qy{wz}}}>\len{\rcos[J]{\rdes w}w}$.
	\begin{proof}
		Let $v\in W$ be such that $wv=\rcos Jw\mxml J$ so that
			by \cref{prop:T_L} we have $v\in W_J$ and $T\cap W_J
				=\qy{\rref w\cap W_J}\udj\lref v$.
		Since by \cref{prop:finitary} and
			\cref{lem:coset-properties}(\ref{lem:coset-unique}) we have 
			$T\cap W_{J_*}\subseteq T\cap W_{\rdes w}=\rref{\mxml{\rdes w}}
				=\rref{w_{\rdes w}}=\rref w\cap W_{\rdes w}\subseteq\rref w$, it
			follows that $T\cap W_J\setminus W_{J_*}
				=\qy{\rref w\cap W_J\setminus W_{J_*}}\udj\lref v$.

		Thus to see that $\len{\rcos[J]{J_*}{(wz)}}-\len{\rcos[J]{J_*}w}
			=\abs{\qy*{T\cap W_J\setminus\rref w}\setminus zW_Jz^{-1}}$,  it
			will suffice to show that $T\cap W_J\setminus W_{J_*}
			=\qy{\rref{wz}\cap W_J\setminus W_{J_*}}
				\udj\qy*{z^{-1}\lref vz\cap W_J}$.
		Indeed, then \cref{lem:length-coset-rep} would show that
			$\len{\rcos[J]{J_*}{(wz)}}-\len{\rcos[J]{J_*}w}
			=\abs{\rref{wz}\cap W_J\setminus W_{J_*}}
				-\abs{\rref w\cap W_J\setminus W_{J_*}}
			=\abs{\lref v} -\abs{z^{-1}\lref vz \cap W_J}
			=\abs{\lref v} -\abs{\lref v\cap zW_Jz^{-1}}
			=\abs{\qy*{T\cap W_J\setminus\rref w}\setminus zW_Jz^{-1}}$, 
			the last equality holding by \cref{prop:T_L}.
	
		To complete the argument, note that \cref{prop:finitary} implies
			$wzz^{-1}v=wv
			=\rcos{\rdes{\rcos Jw\mxml J}}{\qy{wv}}\mxml{\rdes{\rcos Jww_0(J}})
			=\rcos{\rdes{\rcos Jw\mxml J}}{\qy{wz}}\mxml{\rdes{\rcos Jww_0(J}})$.
		Then \cref{prop:T_L} yields 
			$T\cap W_{\rdes{\rcos Jw\mxml J}}
				=\qy{\rref{wz}\cap W_{\rdes{\rcos Jw\mxml J}}}\udj\lref{z^{-1}v}$
			and  because $J\subseteq\rdes{\rcos Jw\mxml J}$ we obtain
			$T\cap W_J=\qy{\rref{wz}\cap W_J}
				\udj\qy*{\lref{z^{-1}v}\cap W_J}$.
		Since $z\in W_{J_*\cap\rdes{\rcos Jw\mxml J}}\subseteq W_{J_*}$, it
			follows by \cref{lem:self-dual} and 
			\cref{lem:subsystem} that $\lref{z^{-1}}
				=\rref z\subseteq W_{J_*}$ and hence 
		$\lref{z^{-1}v}\setminus W_{J_*}
			=\qy*{\lref{z^{-1}}+z^{-1}\lref vz}\setminus W_{J_*}
				=z^{-1}\lref vz\setminus W_{J_*}$. 
		But $z\in W_{J_*}$ implies $z^{-1}\lref vz\setminus W_{J_*}
			=z^{-1}\qy{\lref v\setminus W_{J_*}}z=z\lref vz^{-1}$.
		Thus 
		$T\cap W_J\setminus W_{J_*}=\qy{\rref{wz}\cap W_J\setminus W_{J_*}}
				\udj\qy*{z^{-1}\lref vz\cap W_J}$.
		
	For the supplemental claim, \cref{lem:connected} implies that
		corresponding to $s\in\rdes w\cap\nbd{J'}\setminus J'$ there exist
		$s'\in J'\setminus\rdes w$ and $z\in W_{\qy*{J'\cap\rdes w}
			\cup\set s}$ such that $s\in\rref{z^{-1}s'z}$. 
		Since $J'$ being a connected component of $J$ means that
			$\nbd{J'}\cap J=J'$, then on the one hand 
			$s\in\nbd{J'}\setminus J'=\nbd{J'}\setminus J$ implies 
			$z^{-1}s'z\not\in W_J$ by \cref{lem:subsystem}, i.e.\
			$s'\not\in zW_Jz^{-1}$.
		On the other hand, by \cref{lem:subsystem} we have the containment 
			$s'\in J'\setminus\rdes w\subseteq J\setminus\rdes w$, so that
			$s'\in\qy{T\cap W_J\setminus\rref w}\setminus zW_Jz^{-1}$. 
		Thus by the main claim we have 
			$\len{\rcos[J]{\rdes w}{(wz)}}-\len{\rcos[J]{\rdes w}w}\geq1$.
		\end{proof}
\end{lem}

\section{Factorizations, their masks, and their defects}

\label{sec:defects}

We introduce the term \emph{factorization} for a sequence of finitary
	subsets of generators $J_k\subseteq S$, the notion of a
	\emph{factorization with overlaps}, and associate to each
	factorization with overlaps notions of \emph{masks} and \emph{defects}
	generalizing Billey and Warrington's reformulation of Deodhar's
	combinatorics of subexpressions of a reduced word.
We give a detailed description of how to visualize our generalization in
	the case of type $A$, and also generalize to our setting several
	useful properties of Deodhar's original framework proved by Deodhar,
	and Jones and Woo.

\begin{defn}
	We define a \define{factorization} to be a sequence $J_1,\dots,J_r$ of
		finitary subsets $J_k\subseteq S$.
	In the special case of a factorization $J$ with $J_k=\{s_{i_k}\}$
		for every $k$, we say that $J$ is a \define{singleton
		factorization}, and if furthermore $s_{i_1}\cdots s_{i_r}$ is a
		reduced word, we say that $J$ is a factorization \define{induced by
		a reduced word}.
\end{defn}

In Deodhar's combinatorics of subexpressions \cite{Deodhar}, a
	\emph{subexpression} of a reduced word $s_{i_1}\cdots s_{i_r}$ is
	a sequence $(\sigma_1,\cdots,\sigma_r)$ of Coxeter elements such that
	$\sigma_{k-1}^{-1}\sigma_k\in\set{e,s_k}$, in which case indices $k$
	such that $s_{i_k}\in\rref{\sigma_{k-1}}$ are called
	\emph{defects} of the subexpression.
Deodhar then constructs a family of generating functions
	$\qy*{{\ds\sum_{\sigma:\sigma_r= x}
	q^{d(\sigma)}}}_{x\in W}$ where $d(\sigma)$ is the number of
	defects of a subexpression~$\sigma$.

An easy generalization of Deodhar's combinatorics of subexpressions
	of a reduced word $s_{i_1}\cdots s_{i_r}$ to one of ``subexpressions''
	of a factorization $J_1,\dots,J_r$ would be obtained by defining a
	\emph{subexpression} to be a sequence $(\sigma_1,\dots,\sigma_r)$ such
	that $\sigma_{k-1}^{-1}\sigma_k\in W_{J_k}$, and defining a
	\emph{defect} of a subexpression to be a pair of an index $k$ and a
	reflection $t\in\rref{\sigma_{k-1}}\cap W_{J_k}$.
The resulting family of generating functions
	$\qy*{{\ds\sum_{\sigma:\sigma_r= x}
		q^{d(\sigma)}}}_{x\in W}$ would then coincide (up to a uniform
		degree shift) with the family of coefficients of the monomial
		$C'_{\mxml{J_1}}\cdots C'_{\mxml{J_r}}$ relative to the standard basis
		of the Hecke algebra of the Coxeter system.
However, we want a generalization that can also produce quotients of
	those monomials by Laurent polynomials, such as the quotients
	associated to Skandera's zig-zag factorizations of
	$3412$-$4231$-avoiding permutations and more generally the quotients
	arising in the composition laws of Williamson's Schur algebroid.

We will obtain such generalizations by first equipping a factorization
	with what we call \emph{overlap subsets}, which we use to define an
	equivalence relation on the set of subexpressions, then restricting
	the notion of defects of subexpressions so that it depends only on
	the equivalence class, and finally replacing the above sum counting
	defects over subexpressions with a sum counting restricted defects
	over equivalence classes.
The equivalence relation determined by overlap subsets will have a
	cumbersome expression in terms of subexpressions, but a simple one in
	terms of Billey and Warrington's reformulations of Deodhar's framework
	\cite[Sections~4,5]{Billey--Warrington}.
Accordingly, we begin by generalizing their reformulations rather than
	directly generalizing Deodhar's combinatorics of subexpressions.

Billey and Warrington's first reformulation represents a subexpression
	of a reduced word $s_{i_1}\cdots s_{i_r}$ by a boolean \emph{mask}
	indicating whether $s_{i_k}=\sigma_{k-1}^{-1}\sigma_k$ or not.
For the purposes of our generalization, we consider the
	mask representing a subexpression of a reduced word not as a sequence
	of booleans, but as the sequence of elements
	$\sigma_{k-1}^{-1}\sigma_k\in\set{e,s_{i_k}}=W_{\set{s_{i_k}}}$.
From here onward, $(\sigma_1,\dots,\sigma_r)$ will always denote a mask
	in this sense, as made precise in \cref{defn:mask} below. 
\begin{defn}
\label{defn:mask}
	Given a factorization $J_1,\dots,J_r$, we define a \define{mask
		$\sigma$ on $J$ leading to $x$}, denoted $\sigma\to_J x$, to be a
		sequence $\sigma=(\sigma_1,\dots,\sigma_r)$ with $\sigma_k\in
		W_{J_k}$ for $1\leq k\leq r$ such that $\sigma_1\cdots\sigma_r=x$.
	The \define{right-subexpression} represented by a mask $\sigma$ is the
		sequence $(\sigma_1,\sigma_1\sigma_2,\dots,\sigma_1\cdots\sigma_r)$.

	We define the \define{slice $\slice Jkt$} of the factorization
			$J$ to be the factorization $J_k,\dots,J_t$.
	Given a mask $\sigma$ on $J_1,\dots,J_r$, we write $\slice \sigma kt$
		for the mask on $\slice Jkt$ given by $(\sigma_k,\dots,\sigma_t)$.
	For convenience, we write $\upslice Jt$ for the slice $\slice J1t$ and
		$\downslice Jk$ for the slice $\slice Jrk$, and similarly 
		$\upslice\sigma t$ and $\downslice\sigma k$ for the masks
		$\slice\sigma1t$ and $\slice\sigma kr$.

	We recursively define the \define{right-canonical mask
		$\omega_R(J)\to_Jw(J)$} to be the mask
		$\omega_R(J)=(\omega_R(J)_1,\dots,\omega_R(J)_r)$ such that
		$w(\upslice Jk)=\omega_R(J)_1\cdots\omega_R(J)_k
			=\rcos{J_k}{w(\upslice J{k-1})}\mxml{J_k}
			=\rcos{J_k}{(\omega_R(J)_1\cdots\omega_R(J)_{k-1})}\mxml{J_k}$ for
		each $k$ such that $1\leq k\leq r$, and we say that $J_1,\dots,J_r$
		is a \define{factorization of $w=w(J)$}.
\end{defn}

\begin{rmrk}
	We use the (a priori ambiguous) term ``factorization'' instead of the
		more specific ``right-factorization'' because
		\cref{cor:quotient} will show that the notion of
		a factorization of $w$ is self-dual. 
	Another way to understand the self-duality is that $J$ being a
		factorization of $w$ is equivalent to
		$\mxml{J_1}*\mxml{J_2}*\cdots*\mxml{J_r}=w$  being a product in the
		monoid used by Richardson and Springer in
		\cite{Richardson--Springer}, and which Knutson and Miller
		\cite{Knutson--Miller} called the \emph{Demazure product}. 
\end{rmrk}

Billey and Warrington's second reformulation, which they develop only
	for Coxeter systems of type $A$, visualizes a reduced word
	$s_{i_1}\cdots s_{i_r}$ as a sequence of lattice points $(k,i_k)$, and
	a boolean mask as a \emph{string diagram} of strings headed upward
	such that each pair of strings either bounce or cross whenever they
	\emph{meet}, that is, whenever they pass near the same lattice point. 

\begin{figure}[htbp]
	\centering
	\begin{tikzpicture}
	\node (1) at (0,0) {$1$};
	\node (2) at (1,0) {$2$};
	\node (3) at (2,0) {$3$};
	\node (4) at (3,0) {$4$};
	\node (5) at (4,0) {$5$};

	\node[anchor=east] (y1) at (-0.5, 1) {$[2,4]\approx\{s_2, s_3\} = J_1$};
	\node[anchor=east] (y2) at (-0.5, 2) {$[1,2],[3,5]\approx\{s_1, s_3, s_4\} = J_2$};
	\node[anchor=east] (y3) at (-0.5, 3) {$[2,4]\approx\{s_2, s_3\} = J_3$};

	\node[anchor=west] (s1) at (4.5, 1) {$\sigma_1 = (24)$};
	\node[anchor=west] (s2) at (4.5, 2) {$\sigma_2 = (345)$};
	\node[anchor=west] (s3) at (4.5, 3) {$\sigma_3 = (23)$};

	\draw[very thick] (1, 1) -- (3, 1);

	\draw[very thick] (0, 2) -- (1, 2);
	\draw[very thick] (2, 2) -- (4, 2);

	\draw[very thick] (1, 3) -- (3, 3);

	\node (1e) at (0, 4) {$1$};
	\node (2e) at (4, 4) {$2$};
	\node (3e) at (3, 4) {$3$};
	\node (4e) at (2, 4) {$4$};
	\node (5e) at (1, 4) {$5$};
	\node (wJ) at (2, 4.5) {$\sigma \to_J x$};

	\draw (1) -- (0, 1.7) -- (0.5, 2) -- (0, 2.3) -- (1e);
	\draw[densely dashdotted] (2) -- (1, 0.7) -- (2, 1) -- (3, 1.3) -- (3, 2) -- (4, 2.3)  -- (2e);
	\draw[dotted] (3) -- (2, 1.7) -- (3, 2) -- (3, 2.7)  -- (2, 3) -- (3, 3.3) --  (3e);
	\draw[dashed] (4) -- (3, 0.7) -- (2, 1) -- (1, 1.3) -- (1, 1.7) -- (0.5, 2) -- (1, 2.3) -- (1, 2.7) -- (2, 3)  -- (4e);
	\draw (5) -- (4,1.7) -- (3,2) -- (2,2.3) -- (2, 3) -- (1, 3.3) -- (5e);
\end{tikzpicture}
	\caption{Type $A$ visualization of a mask on a factorization of $45312$.}
	\label{fig:mask}
\end{figure}
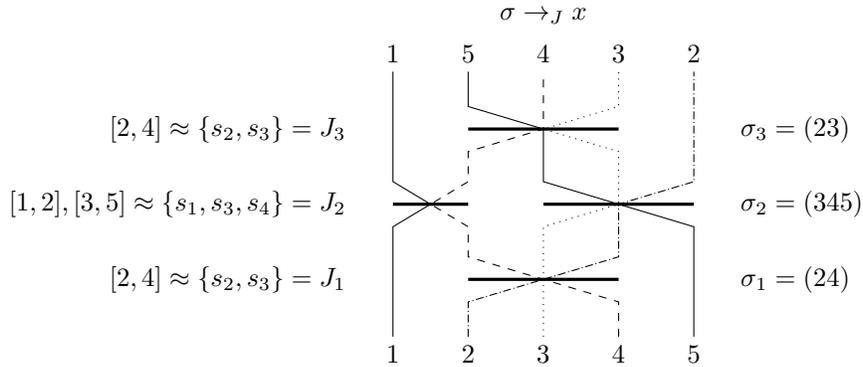

In the case of a Coxeter system of type $A$, our generalization admits a
	similar visualization (and gives essentially the same result as that
	of Billey and Warrington in the case of a factorization induced by a
	reduced word).
Specifically, we picture a factorization $J_1,\dots,J_r$ as a
	configuration of horizontal line segments with integer coefficients,
	with the line segments at height $k$ corresponding to the interval
	decomposition of $J_k$ (described in
	\cref{exmp:interval_decomposition}).
We then visualize a mask as upward-directed strings that are allowed to
	rearrange whenever they go through one of the horizontal line
	segments, in which case the permutation the mask leads to is given by
	the final arrangement of strings, as illustrated in \cref{fig:mask}.

\begin{figure}[htbp]
	\centering
	\begin{tikzpicture}
	\node (1) at (0,0) {$1$};
	\node (2) at (1,0) {$2$};
	\node (3) at (2,0) {$3$};
	\node (4) at (3,0) {$4$};
	\node (5) at (4,0) {$5$};

	\node[anchor=east] (y1) at (-0.5, 1) {$[2,4]\approx\{s_2, s_3\} = J_1$};
	\node[anchor=east] (y2) at (-0.5, 2) {$[1,2],[3,5]\approx\{s_1, s_3, s_4\} = J_2$};
	\node[anchor=east] (y3) at (-0.5, 3) {$[2,4]\approx\{s_2, s_3\} = J_3$};

	\draw[very thick] (1, 1) -- (3, 1);

	\draw[very thick] (0, 2) -- (1, 2);
	\draw[very thick] (2, 2) -- (4, 2);

	\draw[very thick] (1, 3) -- (3, 3);

	\node (1e) at (3, 4) {$1$};
	\node (2e) at (4, 4) {$2$};
	\node (3e) at (2, 4) {$3$};
	\node (4e) at (0, 4) {$4$};
	\node (5e) at (1, 4) {$5$};
	\node (wJ) at (2, 4.5) {$w(J)$};

	\draw[dashed] (1) -- (0, 1.7) -- (0.5, 2) -- (1, 2.3) -- (1, 2.7) -- (2, 3) -- (3, 3.3) -- (1e);
	\draw[dashed] (2) -- (1, 0.7) -- (2, 1) -- (3, 1.3) -- (3, 1.7) -- (3, 2) -- (4, 2.3) -- (2e);
	\draw[dotted] (3) -- (2, 1) -- (2, 1.7) -- (3, 2) -- (3, 2.7) -- (2, 3) -- (3e);
	\draw (4) -- (3, 0.7) -- (2, 1) -- (1, 1.3) -- (1, 1.7) -- (0.5, 2) -- (0, 2.3) -- (4e);
	\draw (5) -- (4,1.7) -- (3, 2) -- (2, 2.3) -- (2, 3) -- (1, 3.3) -- (5e);
\end{tikzpicture}
	\caption{String diagram of a right-canonical mask.}
	\label{fig:right-canonical}
\end{figure}
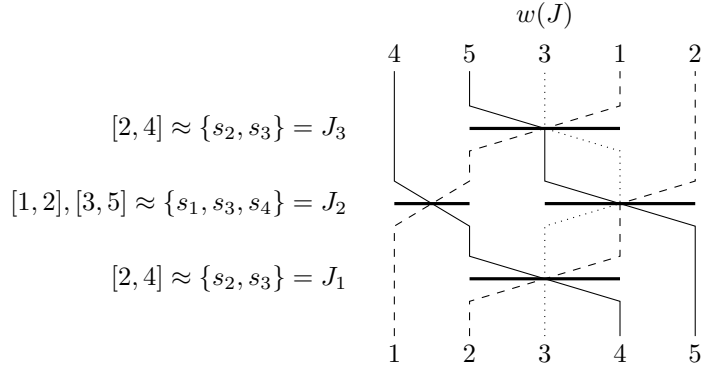

\begin{exmp}
	In the case of a Coxeter system of type $A$, the right-canonical mask
		corresponds to the unique string diagram such that each horizontal
		line segment sorts in decreasing order the strings passing through
		it, as illustrated in \cref{fig:right-canonical}.
\end{exmp}

\begin{exmp}
\label{exmp:reduced-w}
	Since a word $s_{i_1}\cdots s_{i_r}$ is reduced if and only if
		$s_{i_k}\not\in\rdes{s_{i_1}\cdots s_{i_{k-1}}}$ for all $k$,
		i.e.\ if and only if $s_{i_1}\cdots s_{i_k}
			=\rcos{\set{s_{i_k}}}{(s_{i_1}\cdots s_{i_{k-1}})}s_{i_k}$ for all
		$k$, it follows that the right-canonical mask
		$\omega_R(J)$ on a singleton factorization $J_1,\dots,J_r$ given by
		$J_k=\set{s_{i_k}}$ satisfies $\omega_R(J)_k=s_{i_k}$ for all $k$ if
		and only if $s_{i_1}\cdots s_{i_r}$ is a reduced word.
	In particular, a singleton factorization $J_1,\dots,J_r$ given by
		$J_k=\set{s_{i_k}}$ is induced by a reduced word $s_{i_1}\cdots
		s_{i_r}$ if and only if $s_{i_k}\not\in\rdes{w(\upslice J{k-1})}$
		for all $k$.
\end{exmp}

\begin{rmrk}
\label{rmrk:skandera}
	Skandera's \emph{zig-zag factorizations} of a ($3412$-$4231$-avoiding)
		permutation $w$ \cite[Section~3]{Skandera} are all factorizations 
		\emph{in the Coxeter group} of the form 
		$\mxml{J_1}\cdots \mxml{J_{2r+1}}=w$ with the property that 
		$\rcos[J_{2k}]{}{\qy*{\mxml{J_1}\cdots \mxml{J_{2k-1}}}}
			=\mxml{J_{2k}}$
		for all $k$ such that $1\leq k\leq r$. 
	It is straightforward to check the property implies that
		$J_1,\dots,J_r$ is a factorization of $w$ in our sense, i.e.\ that
		$w(J)=w$ and that the factorization $J_1,J_3,\dots,J_{2r+1}$
		obtained by removing the even-indexed factors is also a
		factorization of $w$ in our sense.
	In fact the simpler factorization $J_1,J_3,\dots,J_{2r+1}$ is what we
		call a \emph{contraction} of the factorization $J_1,\dots,J_{2r+1}$
		(see \cref{exmp:skandera}).
	Our \cref{thm:2} asserts that to answer Skandera's question (for the
		kind of quotients arising from the Schur algebroid), it suffices to
		analyze such simpler factorizations.
\end{rmrk}

In the case of Coxeter systems of type $A$, the equivalence relation on
	masks determined by overlap subsets can be formulated in terms of
	strings ``travelling together''.
Informally, we say that a pair of strings ``travel together'' at a
	certain height if they are in the same ``overlap interval'' at that
	height, and our relation will identify two string diagrams if they put
	the same pairs of strings in the same ``overlap intervals,'' i.e.\ if
	they have exactly the same pairs of strings ``travelling together''.

We generalize and make precise this equivalence relation to arbitrary
	Coxeter systems (for which the notion of strings has no immediate
	analogue) by abstracting the fact that two masks are equivalent if and
	only if one can be obtained from the other by permuting the strings
	within each overlap interval.
Since this amounts to acting in a particular way on the elements of the
	mask by elements of the parabolic subgroups generated by the
	overlap intervals, it suffices to generalize the notion of
	overlap intervals.
\begin{defn}
\label{defn:slice-descents-equivalence}
	Fix a factorization $J=J_1,\dots,J_r$.
	We define the \define{set $\rdeslev jJ$ of level $j$ right-descents
		of $J$} by $\rdeslev jJ
			=\ds J_j\setminus\qy{\bigcup_{k=j+1}^r\nbd{J_k}}
			=\set{s\in J_j:s$ is distinct from and commutes with each
				element of $J_k$ when $j<k\leq r}$,
		and similarly the \define{set $\ldeslev jJ$ of level $j$
		left-descents of $J$} by $\ldeslev jJ
			=\ds J_j\setminus\qy{\bigcup_{i=1}^{j-1}\nbd{J_i}}
			=\set{s\in J_j:s$ commutes with each element of $J_i$ when $1\leq
				i<j}$.
	Then the \define{set $\rdes J$ of right-descents of $J$} and the
		\define{set $\ldes J$ of left-descents of $J$} are given by
		$\rdes J=\ds\bigudj_{j=1}^r\rdeslev jJ$
		and $\ldes J=\ds\bigudj_{j=1}^r\ldeslev jJ$
		respectively, and satisfy $\rdes{\upslice Jk}
			=J_k\udj\qy{\rdes{\upslice J{k-1}}\setminus\nbd{J_k}}$ and
		$\ldes{\downslice Ji}=J_i\udj\qy{\ldes{\downslice J{i+1}}
			\setminus\nbd{J_i}}$.

	We define a \define{factorization with overlaps} $J$ to be a
		factorization $J=J_1,\dots,J_r$ equipped with, for each pair $i,k$
		such that $1\leq i<k\leq r$, an \define{overlap subset} $\overlap Jik
			\subseteq
			\ds J_i\cap J_k\setminus\qy{\bigcup_{j=i+1}^{k-1}\nbd{J_j}}
			=\set{s\in J_i\cap J_k:s$ is distinct from and commutes with each
				element of $J_j$ when $i<j<k}
			=\rdeslev i{\upslice Jj}\cap\ldeslev k{\downslice J{j+1}}$ for
			any $j$ such that $i\leq j<k$.
	We denote by $\overlap J*k$ the union $\overlap J*k
		=\ds\bigudj_{i=1}^{k-1}\overlap Jik$, and by
		$\overlap Ji*$ the union $\overlap Ji*
		=\ds\bigudj_{k=i+1}^r\overlap Jik$.
	We will implicitly consider a factorization $J$ with unspecified
		overlaps as the factorization with \define{maximal overlaps} given
		by $\overlap Jik\ds 
			=J_i\cap J_k\setminus\qy*{\bigcup_{j=i+1}^{k-1}\nbd{J_j}}
			=\set{s\in J_i\cap J_k:s$ is distinct from and commutes with each
				element of $J_j$ when $i<j<k}
			=\rdeslev i{\upslice Jj}\cap\ldeslev k{\downslice J{j+1}}$ for
			any $j$ such that $i\leq j<k$.

	Finally, we define the \define{equivalence relation on masks} to
		be generated by the symmetric relations
		$(\sigma_1,\dots,\sigma_i,\dots,\sigma_k,\dots,\sigma_r)\sim(
		\sigma_1,\dots,\sigma_is,\dots,s\sigma_k,\dots,\sigma_r)$ for
		some $s\in \overlap{J}{i}{k}$. 
	Then given a mask $\sigma\to_J x$ we denote by $[\sigma]\to_Jx$ the
		\define{equivalence class of masks leading to $x$} that contains
		$\sigma$ (well-defined by \cref{lem:action}).
\end{defn}

\begin{rmrk}
\label{rmrk:factorization-descents}
	We refer to the elements of $\rdes J$ as ``right-descents'' of the
		factorization $J$ because \cref{cor:recursive-defect-polys}
		will show that $\rdes J\subseteq\rdes{w(J)}$. 
	Moreover, \emph{right-descent factorizations}
		(\cref{defn:descent}), which play a key role in the proof of our
		\cref{thm:2} and will in fact by characterized by
		the condition that $\rdups Jk=\rdes{w(\upslice Jk)}$ for all $k$
		(\cref{prop:rdescent-inductive-step}).
\end{rmrk}

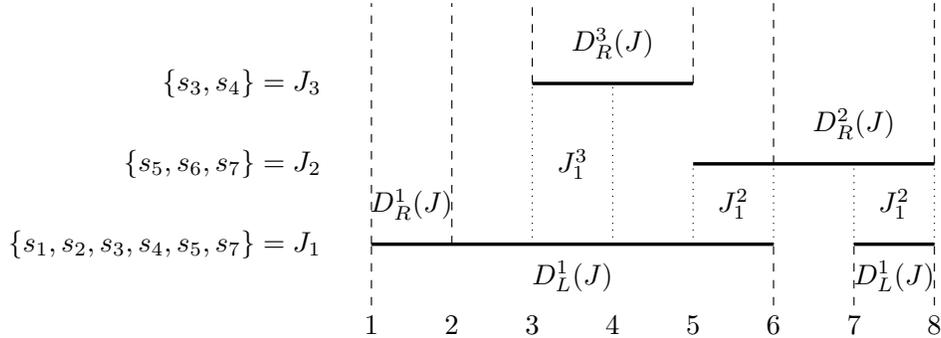
\begin{figure}[htbp]
	\centering
	\begin{tikzpicture}
	\node (1) at (0,0) {$1$};
	\node (2) at (1,0) {$2$};
	\node (3) at (2,0) {$3$};
	\node (4) at (3,0) {$4$};
	\node (5) at (4,0) {$5$};
	\node (6) at (5,0) {$6$};
	\node (7) at (6,0) {$7$};
	\node (8) at (7,0) {$8$};

	\node[anchor=east] (y1) at (-0.5, 1) {$\{s_1, s_2, s_3, s_4, s_5, s_7\} = J_1$};
	\node[anchor=east] (y2) at (-0.5, 2) {$\{s_5, s_6, s_7\} = J_2$};
	\node[anchor=east] (y3) at (-0.5, 3) {$\{s_3, s_4\} = J_3$};

	\draw[very thick] (0, 1) -- (5, 1);
	\draw[very thick] (6, 1) -- (7, 1);

	\draw[very thick] (4, 2) -- (7, 2);

	\draw[very thick] (2, 3) -- (4, 3);

	\draw[dashed] (0,1) -- (0, 4);
	\draw[dashed] (1,1) -- (1, 4);
	\node (rdes1) at (0.5, 1.5) {$D_R^1(J)$};

	\draw[dotted] (2,1) -- (2, 3);
	\draw[dotted] (3,1) -- (3, 3);
	\node (J13) at (2.5, 2) {$J_1^3$};

	\draw[dashed] (2, 3) -- (2, 4);
	\draw[dashed] (4, 3) -- (4, 4);
	\node (rdes3) at (3, 3.5) {$D_R^3(J)$};

	\draw[dotted] (4,1) -- (4, 2);
	\draw[dotted] (5,1) -- (5, 2);
	\node (J12) at (4.5, 1.5) {$J_1^2$};

	\draw[dashed] (7) -- (6, 1);
	\draw[dashed] (5, 2) -- (5, 4);
	\draw[dashed] (8) -- (7, 1);
	\draw[dashed] (7,2) -- (7, 4);
	\node (rdes2) at (6, 2.5) {$D_R^2(J)$};

	\draw[dotted] (7,1) -- (7,2);
	\draw[dotted] (6,1) -- (6,2);
	\node (J12) at (6.5, 1.5) {$J_1^2$};

	\draw[dashed] (1) -- (0, 1);
	\draw[dashed] (6) -- (5, 1);
	\node (ldes11) at (2.5, 0.6) {$D_L^1(J)$};
	\node (ldes12) at (6.5, 0.6) {$D_L^1(J)$};
\end{tikzpicture}
	\caption{Factorization with labelled overlap and descent intervals}
	\label{fig:overlap}
\end{figure}

\begin{exmp}
\label{exmp:descents_overlaps}
	In the case of a Coxeter system of type $A$, the connected components
		of $\overlap Jik$ are what we informally called \define{overlap
		intervals}.
	\cref{fig:overlap} illustrates the connected components of maximal
		overlap subsets and of sets of left- and right-descents of a
		factorization in a Coxeter system of type $A$.
\end{exmp}

\begin{exmp}
\label{exmp:trivial-overlaps}
	Given a singleton factorization $J$ given by $J_k=\set{s_{i_k}}$ 
		we have by definition that 
		$\emptyset\neq\overlap Jjk$ only if $s_{i_j}=s_{i_k}$ and 
		$s_{i_1}\cdots s_{i_j}\cdots s_{i_k}
			=s_{i_1}\cdots s_{i_{j-1}}s_{i_{j+1}}\cdots s_{i_{k-1}}$, so that
		in particular $J$ is not induced by a reduced word.
	Thus for singleton factorizations induced by reduced words, each
		overlap subset is trivial, as is consequently our equivalence
		relation on masks.
\end{exmp}

\begin{exmp}
\label{exmp:gm-factorization}
	Say that a factorizations $J'$ has \emph{Gelfand--MacPherson} overlaps
		if $\overlap{J'}ik=\begin{cases}J'_i\cap J'_k&k-i=1\\
			\emptyset&\text{otherwise}\end{cases}$.
	In the case of Weyl groups, such factorizations correspond precisely
		to the iterated base changes as described in \cref{subsec:weyl-gm}.
	More generally, to any factorization with overlaps $J_1,\dots,J_r$,
		we can associate a factorization $J'_1,\dots,J'_r$ whose maximal
		overlaps are Gelfand--MacPherson by setting 
		$J'_j =\ds J_j\udj\bigudj_{i<j<k}\overlap Jik$ and
		$\overlap{J'}j{j+1}=\ds\bigudj_{i\leq j<k}\overlap Jik$, and for
		which $J$  is a \emph{contraction} of $J'$ in the sense of
		\cref{defn:contraction}.
	The reason we develop the theory of general factorizations instead of
		only of ones with Gelfand--MacPherson overlaps is that the
		conditions on factorizations imposed by our major result 
		\cref{thm:2} are in terms of their \emph{contractions}, whose
		maximal overlaps need not be Gelfand--MacPherson overlaps.
\end{exmp}

\cref{lem:action} below gives the more precise description of the
	equivalence classes of masks as orbits under a group action on masks,
	while \cref{lem:action-subexprs-vpstrings} will describe the
	induced action on right-subexpressions
	$(\sigma_1,\sigma_1\sigma_2,\dots,\sigma_1\cdots\sigma_r)$.
\begin{lem}
\label{lem:action}
	The actions on masks on a factorization with overlaps $J_1,\dots,J_r$ 
		given by 
		$(\sigma_1,\dots,\sigma_r)\mapsto\begin{cases}
			(\sigma_1,\cdots,x\sigma_j,\cdots,\sigma_r)&x\in W_{\ldeslev jJ}\\
			(\sigma_1,\dots,\sigma_iy^{-1},\dots, y\sigma_k,\dots,\sigma_r)
				&y\in W_{\overlap Jik}\\
			(\sigma_1,\cdots,\sigma_jz^{-1},\cdots,\sigma_r)
				&z\in W_{\rdeslev jJ}\end{cases}$ commute and determine an
		action of the (external) direct product $W_{\ldes J}
			\times\ds\qy*{\prod_{1\leq i<k\leq r} W_{\overlap Jik}}\times W_{\rdes J}$ on the set of masks of $J$,
		with the property that $(x, y, z)\in W_{\ldes J}\times\qy*{
			\ds\prod_{1\leq i<k\leq r} W_{\overlap Jik}}
				\times W_{\rdes J}$ sends $\sigma\to_J w$ to $\sigma'\to_J xwz$.

	Furthermore, the equivalence classes of masks are the orbits of
		masks under the induced action of the (external) direct product
		$\ds\prod_{1\leq i<k\leq r}W_{\overlap Jik}$.
	In particular,  we have an induced action of $W_{\ldes J}\times
		W_{\rdes J}$ on equivalence classes of masks and this action has the
		property that each $(x,z)\in W_{\ldes J}\times W_{\rdes J}$ induces
		a bijection $\set{[\sigma]\to_Jy}\cong\set{[\tau]\to_J xyz^{-1}}$.
	\begin{proof}
		By definition, the actions of $W_{\rdeslev jJ}$ and 
			$W_{\ldeslev jJ}$ are respectively induced by the right and left
			actions of $W_{J_j}$ on the $j^\text{th}$ component of the mask,
			and hence commute with one another.
		Since actions on different components of the mask also commute with
			one another, and since by definition $\ldes J
			=\ds\bigudj_{j=1}^r\ldeslev jJ$ where the 
			$\ldeslev jJ$ commute with each other and hence
			$W_{\ldes J}\cong\ds\prod_{j=1}^rW_{\ldeslev jJ}$ and
			$W_{\rdes J}\cong\ds\prod_{j=1}^rW_{\rdeslev jJ}$ by
			\cref{lem:coset-properties}\ref{lem:component-decompose}),
			it follows that we do have an induced action of the (external)
			direct product $W_{\ldes J}\times W_{\rdes J}$.

		Similarly, the action of $W_{\overlap Jik}$ is induced by the left
			action of $W_{J_k}$ on the $k^\text{th}$ component and the right
			action of $W_{J_i}$ on the $i^\text{th}$ component, so certainly
			commutes with actions on components except possibly other left
			actions on the $i^\text{th}$ component and other right actions on
			the $k^\text{th}$ component.
		But because by definition $\overlap Jik$ commutes with $\ldeslev kJ$,
			and with $\overlap Jjk \subseteq J_j$ when $i\neq j<k$, it follows
			that the right actions on the $k^\text{th}$ component of 
			$W_{\overlap Jik}$, $W_{\ldeslev kJ}$, and $W_{\overlap Jjk}$ with
			$i\neq j<k$ commute. 
		Similarly, the left actions on the $i^\text{th}$ component of
			$W_{\ldeslev iJ}$, $W_{\overlap Jik}$, and $W_{\overlap Jij}$ with
			$i<j\neq k$ commute as well.
		It follows that we have a well-defined action on masks by the
			(external) direct product $W_{\ldes J}\times\ds
				\prod_{1\leq i<k\leq r}W_{\overlap Jik}\times W_{\rdes J}$.

		Now, since $s\in\overlap Jij$ commutes with each element of $J_k$
			and hence with $\sigma_k\in W_{J_k}$ when $i<k<j$, it follows that
			for each $y\in W_{\overlap Jij}$ we have
			$\sigma_1\cdots\sigma_iy^{-1}\cdots y\sigma_j\cdots\sigma_r
				=\sigma_1\cdots\sigma_{j-1}y^{-1}y\sigma_j\cdots\sigma_r
				=\sigma_1\cdots\sigma_r$.
		Similarly, since $s\in\ldeslev jJ$ commutes with each
			element of $J_i$ when $i<j$, it follows that for $x\in W_{\ldeslev
			jJ}$ we have $\sigma_1\cdots
			x\sigma_j\cdots\sigma_r=x\sigma_1\cdots\sigma_r$ and dually for
			$W_{\rdeslev jJ}$.
		Thus, each $(x,y,z)\in\ds W_{\ldes{J}}
			\times\prod_{1\leq i<k\leq r}W_{\overlap Jik}\times W_{\rdes{J}}$
			sends $\sigma\to_Jw$ to $\tau\to_Jxwz^{-1}$.
	\end{proof}
\end{lem}
Next, we come to they key concept in our generalization of Deodhar's
	combinatorics of subexpressions: the notion of \emph{defects}.
In Billey and Warrington's second reformulation of Deodhar's
	combinatorics of subexpressions, defects occur whenever a pair of
	strings meet having crossed an odd number of times, i.e.\ whenever
	they meet ``out-of-order''.
Unfortunately, the order in which a pair of string meets (and hence
	whether or not they realize a defect) is not preserved under our
	notion of equivalence unless the factorization's overlaps are empty.
Fortunately, it \emph{is} preserved if we restrict the notion of meeting
	to refer only to those pairs of strings that meet without having
	traveled together up to their meeting.
\begin{exmp}
	In \cref{fig:mask} the pair of strings labeled by $2$ and $3$
		meet (in the unrestricted sense) out-of-order at the middle-right
		segment, but this does not count as a meeting in the restricted
		sense (and thus also not as a defect) because they travel
		together up to this meeting along the same overlap interval of the
		bottom segment and the middle-right segment. 
	In fact they had their meeting (in the restricted sense) at the bottom
		segment, where they were not out-of-order.
\end{exmp}
To express our generalization of defects to factorizations with overlaps
	on arbitrary Coxeter systems, we need to remove its dependence on
	individual strings that the permutation group acts on (since such
	strings lack an immediate analogue for arbitrary Coxeter systems).
This can be done because even though individual strings have no analogue
	in arbitrary Coxeter systems, \emph{pairs of strings} do.
Specifically, we can identify a pair of strings in type $A$
	with the sequence of unordered pairs of positions that the strings are
	put in, and then further identify that sequence of unordered pairs
	with a sequence of transpositions in the permutation group.

Then a pair of strings meets out-of-order at a line segment at height
	$k$ precisely when the $k-1^\text{st}$ transposition in the
	identification corresponds to a position-based inversion of
	$\sigma_1\cdots\sigma_{k-1}$ contained in an interval corresponding
	to the line segment at height $k$.
Furthermore, a pair of strings ``travel together'' up to
	height $k$ precisely when the $k-1^{\text{st}}$ transposition in
	the identification is contained in an overlap interval
	up to $k$.

Since \cref{exmp:trivial-overlaps} implies these restrictions are
	vacuous in the case of singleton factorizations induced by reduced
	words, and since restricted meetings and restricted defects play a more
	central role than the unrestricted generalization, we will use the
	terms ``meetings'' and ``defects'' from here onward to refer only to
	what we have informally called restricted meetings and restricted
	defects above.

\begin{defn}
\label{defn:defects}
	Given a reflection $t\in T$ and a mask
		$\sigma=(\sigma_1,\dots,\sigma_r)$ on a factorization with overlaps
		$J_1,\dots,J_r$, we define the \define{virtual pair of strings
		$t^\sigma$} of $\sigma$ right-labeled by $t$ to be the sequence
		$(t^\sigma_0,\dots,t^\sigma_r)$ given by $t^\sigma_k
			=(\sigma_1\cdots\sigma_k)^{-1}t(\sigma_1\cdots\sigma_k)$.

	Given $k$ such that $1\leq k\leq r$, we say that $k$ is a
		\define{right-meeting} of $t^\sigma$ if $t^\sigma_{k-1}\in
			W_{J_k}\setminus W_{\overlap J*k}$ and denote by $\rmeet\sigma t$
		the \define{set of right-meetings} of $t^\sigma$.
	We say that $k$ is furthermore a \define{right-defect} of $t^\sigma$
		if $k\in\rmeet\sigma t$ and $t^\sigma_{k-1}
			\in\rref{\sigma_1\cdots\sigma_{k-1}}$, i.e.\ if $t^\sigma_{k-1}
			\in\rref{\sigma_1\cdots\sigma_{k-1}}\cap W_{J_k}
			\setminus W_{\overlap J*k}$, and
	we denote by $\rdef\sigma t$ the \define{set of right-defects} of
		the virtual pair of strings $t^\sigma$.

	We define the \define{right-defect statistic
		$d_R(\sigma)$} of a mask $\sigma$ to be given by
		$d_R(\sigma)=\ds\sum_{t\in T}\abs{\rdef\sigma t}
			=\sum_{k=1}^r\abs{\rref{\sigma_1\cdots\sigma_{k-1}}\cap
			W_{J_k}\setminus W_{\overlap J*k}}$: the
		total number of right-defects over all virtual pair of strings
		$t^\sigma$ for a given mask $\sigma$.
	We also denote by $m_R(J)$ the total number of right-meetings of any
		mask given by $m_R(J)=\ds\sum_{k=1}^r\abs{T\cap W_{J_k}
			\setminus W_{\overlap J*k}}$.
\end{defn}

\begin{rmrk}
\label{rmrk:gm-defects}
	Since $t\in\lref{w_{k-1}}\cap w_{k-1}\qy{W_{J_k}
	  \setminus W_{J_{k-1}\cap J_k}}w_{k-1}^{-1}$ if and only if 
	$w_{k-1}^{-1}tw_{k-1}\in w_{k-1}^{-1}\lref{w_{k-1}}w_k\cap W_{J_k}
		\setminus W_{J_{k-1}\cap J_k}
		=\rref{w_{k-1}}\cap W_{J_k}\setminus W_{J_{k-1}\cap J_k}$,
	in the case of Weyl groups the notion of defect for a factorization
		$J$ with overlaps $\overlap Jik=\begin{cases}J_i\cap J_k&k-i=1
		\\\emptyset&\text{otherwise}\end{cases}$ matches the notion
	described in \cref{subsec:weyl-gm}.
\end{rmrk}

For Coxeter systems of type $A$, we initially described the equivalence
	relation on masks in terms of pairs of strings, rather than as the
	orbits induced by an action.
\cref{lem:action-subexprs-vpstrings} describes
	the action induced by \cref{lem:action} on virtual pairs of
	strings for arbitrary Coxeter systems, while
	\cref{prop:defects-are-meetings} derives from this action
	the desired invariance of defects under our equivalence relation on
	masks.
\begin{lem}
\label{lem:action-subexprs-vpstrings}
	Given a mask $\sigma=(\sigma_1,\dots,\sigma_r)$ on a factorization
		with overlaps $J_1,\dots,J_r$ suppose that $\tau$ is the image of
		$\sigma$ under the action specified in \cref{lem:action} of
		$s\in\rdeslev Ji$ (or $s\in\overlap Jik$).
	Then for each $m$ such that $1\leq m\leq r$ we have
		$(\tau_1\cdots\tau_m,t^\tau_m)=\begin{cases}
			(\sigma_1\cdots\sigma_ms,st^\sigma_ms)&i\leq m(<k)\\
			(\sigma_1\cdots\sigma_m,t^\sigma_m)&\text{otherwise}\end{cases}$. 
	\begin{proof}
		Suppose that $\tau$ is the image of $\sigma$ under $s\in\rdeslev iJ$
			or $s\in\overlap Jik$, i.e.\ that
			$\tau=(\sigma_1,\cdots,\sigma_is,\cdots,\sigma_r)$ or
			$\tau=(\sigma_1,\cdots,\sigma_is,\cdots,s\sigma_k,\dots,\sigma_r)$.
		If $\sigma_1\cdots\sigma_m=\tau_1\cdots\tau_m$, then
			$t^\sigma_m=t^\tau_m$ by definition.
		Otherwise, by definition of $\tau$ we have
			$\sigma_1\cdots\sigma_m\neq\tau_1\cdots\tau_m$ only if
			$i\leq m$.
		Since $s\in \rdeslev iJ$ (respectively $s\in\overlap Jik$) commutes
			with $\sigma_j\in W_{J_j}$ for $j$ such that $i<j$ 
			(respectively such that $i<j<k$), it follows that
			$\tau_1\cdots\tau_m
				=\sigma_1\cdots\sigma_is\cdots\sigma_m
				=\sigma_1\cdots\sigma_m s$ if $i\leq m$ (resp.\ if $i\leq m<k$).
		Furthermore, it follows that if $k\leq m$, then
			$\tau_1\cdots\tau_m=\sigma_1\cdots\sigma_{k-1}ss\sigma_k\cdots\sigma_m
				=\sigma_1\cdots\sigma_m$
		Thus $\sigma_1\cdots\sigma_m\neq\tau_1\cdots\tau_m$ if and only
			if $i\leq m$ (respectively $i\leq m<k)$, in which case
			$t^\tau_m=(\tau_1\cdots\tau_m)^{-1}t(\tau_1\cdots\tau_m)
				=s(\sigma_1\cdots\sigma_m)^{-1}t(\sigma_1\cdots\sigma_m)s
				=st^\sigma_ms$, as desired.
	\end{proof}	
\end{lem}

For Coxeter systems of type $A$ one can see from the visualization that
	for each transposition $t$ whether the pair of strings labeled by
	$t$ has a right-defect at $k$, i.e.\ whether it meets out-of-order at
	a horizontal line segment at height $k$, depends only on the
	equivalence class of the mask $\sigma$. 
This is also true for arbitrary Coxeter systems and will allows us to
	define the generating functions associated to a factorization with
	overlaps.

\begin{prop}
\label{prop:defects-are-meetings}
	Given a mask $\sigma=(\sigma_1,\dots,\sigma_r)$ on a factorization
		with overlaps $J_1,\dots,J_r$ and $t\in T$, whether
		$k\in\rmeet\sigma t$ and whether $k\in\rdef\sigma t$ depend only on
		the $\ds\qy*{\prod_{1\leq i<k\leq r}
			W_{\overlap Jik}}\times W_{\rdes J}$-orbit of the mask $\sigma$.
	\begin{proof}
		Suppose $\tau$ is the image of $\sigma$ under the action of
			$s\in\ds\qy*{\bigudj_{1\leq i<k\leq r}
				\overlap Jik}\udj\qy*{\bigudj_{i=1}^r\rdeslev iJ}$,
			and suppose $m$ is such that $1\leq m\leq r$.
		If $\tau_1\cdots\tau_{m-1}=\sigma_1\cdots\sigma_{m-1}$, then 
			by definition $t^\tau_{m-1}=t^\sigma_{m-1}$ for any $t\in T$, in
			which case $m\in\rmeet\tau t$ is equivalent
			to $m\in\rmeet\sigma t$, and $m\in\rdef\tau t$ is equivalent
			to $m\in\rdef\sigma t$ by definition.
		If $\tau_1\cdots\tau_{m-1}\neq\sigma_1\cdots\sigma_{m-1}$, then
			by \cref{lem:action-subexprs-vpstrings} we have
			$\tau_m\cdots\tau_{m-1}=\sigma_m\cdots\sigma_{m-1}s$ and
			$t^\tau_{m-1}=st^\sigma_{m-1}s$ where either:
			$s\in\rdeslev iJ$ for $i$ such that $i\leq m-1$, or 
			$s\in\overlap Jik$ for $i,k$ such that $i\leq m-1<k$.

		If $k=m$, then $s\in\overlap Jim\subseteq\overlap J*m\subseteq J_m$.
		If otherwise $s\in\rdeslev iJ$ or $s\in\overlap Jik$ with $m<k$,
			then in particular $s\not\in\nbd{J_m}$  so that in both cases we
			have $sW_{J_m}s=W_{J_m}$ and $sW_{\overlap J*m}s
				=W_{\overlap J*m}$.
		Consequently, $m\in\rmeet\sigma t$, i.e.\ 
			$t^\sigma_{m-1}\in W_{J_m}\setminus W_{\overlap J*m}$, if and only
			if $t^\tau_{m-1}=st^\sigma_{m-1}s
					\in s\qy*{W_{J_m}\setminus W_{\overlap J*m}}s
				=W_{J_m}\setminus W_{\overlap J*m}$, i.e.\ if and only if
			$m\in\rmeet\tau t$.

		Next, since $\rdeslev iJ\cap J_m=\emptyset$ when $i\leq m-1$ and
			$\overlap Jik\cap J_m=\emptyset$ if $i<m<k$, it follows that 
			$s\in J_m$ holds only if $s\in\overlap Jim\subseteq\overlap J*m$,
			i.e.\ we have $s\not\in W_{J_m}\setminus W_{\overlap J*m}$.
		Thus, $m\in\rdef\tau t$, i.e.\ $t^\tau_{m-1}=st^\sigma_{m-1}s
			\in\rref{\tau_1\cdots\tau_{m-1}}
				\cap W_{J_m}\setminus W_{\overlap J*m}
			=\rref{\sigma_1\cdots\sigma_{m-1}s}
				\cap W_{J_m}\setminus W_{\overlap J*m}
			=\qy*{s\rref{\sigma_1\cdots\sigma_{m-1}}s+\set s}
				\cap W_{J_m}\setminus W_{\overlap J*m}
			=s\rref{\sigma_1\cdots\sigma_{m-1}}s
				\cap W_{J_m}\setminus W_{\overlap J*m}
			=s\qy*{\rref{\sigma_1\cdots\sigma_{m-1}}
				\cap W_{J_m}\setminus W_{\overlap J*m}}s$
			if and only if $t^\sigma_{m-1}
				\in\rref{\sigma_1\cdots\sigma_{m-1}}
					\cap W_{J_m}\setminus W_{\overlap J*m}
				=\rref{\sigma_1\cdots\sigma_{m-1}}
					\cap W_{J_m}\setminus W_{\overlap J*m}$
			i.e.\ $m\in\rdef\sigma t$.
	\end{proof}
\end{prop}	

\begin{rmrk}
	In the special case of a singleton factorization induced by a reduced
		word, the invariance of $\rdef\sigma t$ under the action of
		$W_{\rdes J}$ asserted by
		\cref{prop:defects-are-meetings} is a generalization of
		\cite[Proposition~2.3(ii)]{Deodhar}.
\end{rmrk}

\begin{exmp}
\label{exmp:singleton-defects}
	For a singleton factorization $J_1,\dots,J_r$ given by
		$J_k=\set{s_{i_k}}$, we have that $k\in\rmeet\sigma t$ if and only
		if $t^\sigma_{k-1}\in T\cap W_{J_k}
				\setminus W_{\overlap J*k}
			=\set{s_{i_k}}\setminus W_{\overlap J*k}
			=\set{s_{i_k}}\setminus\qy*{S\cap W_{\overlap J*k}}
			=\set{s_{i_k}}\setminus\overlap J*k$, i.e.\ if and only
			if $t^\sigma_{k-1}=s_{i_k}\not\in\overlap J*k$.
	Thus, $k\in\rdef\sigma t$ if and only if
		$t_{k-1}^\sigma=s_{i_k}\in\rref{\sigma_1\cdots\sigma_{k-1}}
			\setminus\overlap J*k$, i.e.\ if and only if 
		$t_{k-1}^\sigma=s_{i_k}\in\rdes{\sigma_1\cdots\sigma_{k-1}}
			\setminus\overlap J*k$, so that we have
			$m_R(J)=\abs{\set{k:s_{i_k}\not\in\overlap J*k}}$ and
			$d_R(\sigma)=\abs{\set{k:s_{i_k}
				\in\rdes{\sigma_1\cdots\sigma_{k-1}}
					\setminus\overlap J*k}}$.

	In the special case of a singleton factorization induced by a reduced
		word, \cref{exmp:trivial-overlaps} shows all overlap
		subsets are trivial, and hence implies that 
		$m_R(J)=\abs{\set{k:s_{i_k}\not\in\overlap J*k}}=r=\ell(w(J))$ and
		$d_R(\sigma)=\abs{\set{k:s_{i_k}\in\rdes{\sigma_1\cdots\sigma_{k-1}}}}$,
		thereby recovering Deodhar's original notion of defect.
	Furthermore, we have in that case that $d_R(\omega_R(J))=0$ since 
		$\omega_R(J)_k=s_{i_k}\not\in\rdes{s_{i_1}\cdots s_{i_{k-1}}}
			=\rdes{\omega_R(J)_1\cdots\omega_R(J)_{k-1}}$ for all
		$k$ such that $1\leq k\leq r$, as explained in 
		\cref{exmp:reduced-w}.
\end{exmp}

We conclude this section by describing the connection between the
	equivalence relation on masks and their associated subexpressions.
\begin{lem}
\label{lem:y-well}
	Given a mask $\sigma\to_Jx$ on a factorization with overlaps 
		$J_1,\dots,J_r$, set 
		$x^\sigma_k=\rcos{\rdups Jk}{\qy{\sigma_1\cdots\sigma_k}}$ and
		$y^\sigma_k=\rcos[J_k]{\overlap J*k}
			{(\sigma_1\cdots\sigma_{k-1})}$.
	Then $\set{t\in T:k\in\rdef\sigma t}=x^\sigma_k\lref{y^\sigma_k}\qy{
		x^\sigma_k}^{-1}$ so that $\ell(y^\sigma_k)=
			\abs{\set{t\in T:k\in\rdef\sigma t}}$ and
			$d_R(\sigma)=\ds\sum_{k=1}^r\ell(y^\sigma_k)$.
	Furthermore $x^\sigma_k$ and $y^\sigma_k$ are constant on
		$\ds\qy*{\prod_{1\leq i<k\leq r}
			W_{\overlap Jik}}\times W_{\rdes J}$-orbits of masks.
	\begin{proof}
		We have by the definition of right-defects that
			\begin{align*}
				\set{t\in T&:k\in\rdef\sigma t}\\
					&=\set{t\in T:t^\sigma_{k-1}\in\rref{\sigma_1\cdots\sigma_{k-1}}
						\cap W_{J_k}\setminus W_{\overlap J*k}}\\
					&=\qy{\sigma_1\cdots\sigma_{k-1}}
						\qy*{\rref{\sigma_1\cdots\sigma_{k-1}}
							\cap W_{J_k}\setminus W_{\overlap J*k}}
						\qy{\sigma_1\cdots\sigma_{k-1}}^{-1}
			\end{align*}
			Applying \cref{lem:conjugated-lref-coset-rep} with 
				$J=J_k$, $J_*=\overlap J*k$, and $J'=\rdups J{k-1}$ yields that
				this latter term is $uVu^{-1}$ for $u=
				\rcos{J_k\udj\qy*{\rdups J{k-1}\setminus\nbd{J_k}}}{
					\qy{\sigma_1\cdots\sigma_{k-1}}}$
			and $V=
				\lref{\rcos[J_k]{\overlap J*k}{\qy{\sigma_1\cdots\sigma_{k-1}}}}
				=\lref{y^\sigma_k}$.
			Since $J_k\udj\qy*{\rdups J{k-1}\setminus\nbd{J_k}}=\rdups Jk$ and
			$\sigma_k\in W_{J_k}$ it follows that also $u=x^\sigma_k$, as
			desired.

		For the supplemental claim, since \cref{lem:injective-cocycle}
			implies $y_k^\sigma$ is uniquely determined by $\lref{y_k^\sigma}$
			and since \cref{prop:defects-are-meetings} implies the
			set $\set{t\in T:k\in\rdef\sigma t}$ is invariant under the action
			of $\ds\qy*{\prod_{1\leq i<k\leq r}W_{\overlap Jik}}
				\times W_{\rdes J}$, by the main claim it will suffice to
			show that $x^\sigma_k$ is also invariant under this action.
		Suppose $\tau$ is the image of $\sigma$
			under the action of $s \in\ds
				\qy*{\bigudj_{1\leq i<k\leq r}
					\overlap Jik}\udj\qy*{\bigudj_{i=1}^r\rdeslev iJ}$.
		If $\tau_1\cdots\tau_k=\sigma_1\cdots\sigma_k$, then there
			is nothing to prove.
		If $\tau_1\cdots\tau_k\neq\sigma_1\cdots\sigma_k$, then by
			\cref{lem:action-subexprs-vpstrings} we have
			$\tau_1\cdots\tau_k=\sigma_1\cdots\sigma_ks$ where either 
			$s\in\rdeslev iJ\subseteq\rdeslev i{\upslice Jk}
				\subseteq\rdups Jk$ with $i\leq k$
			or $s\in\overlap Jij\subseteq\rdeslev i{\upslice Jk}
				\subseteq\rdups Jk$ with $i\leq k<j$.
		Thus,	$x^\tau_k=\rcos{\rdes{\upslice Jk}}{(\tau_1\cdots\tau_k)}
			=\rcos{\rdes{\upslice Jk}}{(\sigma_1\cdots\sigma_ks)}
			=\rcos{\rdes{\upslice Jk}}{(\sigma_1\cdots\sigma_k)}
			=x^\sigma_k$, as desired.
	\end{proof}
\end{lem}

\begin{prop}
\label{prop:defect-ideals}
	Given a mask $\sigma\to_Jx$ on a factorization with overlaps 
		$J_1,\dots,J_r$, let $y^\sigma=(y^\sigma_1,\dots,y^\sigma_r)$ and
		$z^\sigma=(z^\sigma_1,\dots,z^\sigma_r)$ be given by
		$y^\sigma_k=\rcos[J_k]{\overlap J*k}
			{(\sigma_1\cdots\sigma_{k-1})}$ and
		$z^\sigma_k=\rcos[\overlap J*k]{}
			{\qy*{\sigma_1\cdots\sigma_{k-1}}}$.
	Then for each $x'\in W$ such that
		$\rcos{\rdes J}{x'}\leq \rcos{\rdes J}x$, each $(y_1,\dots,y_r)
			\in\ds\prod_{k=1}^r\rcos[J_k]{\overlap J*k}W$
			such that $y_k\leq y^\sigma_k$ for all $k$, and each
			$(z_1,\dots,z_r)
				\in\ds\prod_{k=1}^r W_{\overlap J*k}$,
		there is a unique mask $\tau \to_J x'$ such that
			$y^\tau=(y_1,\dots,y_r)$ and $z^\tau=(z_1,\dots,z_r)$.
	Furthermore, for any $\sigma\to_Jx$ we have $x\leq w(J)$. 
	\begin{proof}
		We argue by induction on the length $r$ of the factorization.
		For the base case $r=1$, note that every mask
			$\sigma\to_Jx$ is given by $\sigma_1=x\in W_{J_1}=W_{\rdes J}$ and
			satisfies $y^\sigma_1=e$, $W_{\overlap J*1}
				=\set e \ni z^\sigma_1=e$, and $W_{J_1}\ni x\leq \mxml{J_1}$
			by \cref{prop:finitary}.
		But then $\rcos{\rdes J}{x'}
			\leq\rcos{\rdes J}x=e$ if and only if $\rcos{\rdes J}{x'}=e$ by
			\cref{lem:bruhat}, i.e.\ if and only if $x'\in W_{\rdes J}
				=W_{J_1}$
			by \cref{lem:coset-properties}(\ref{lem:coset-unique}).

		For the inductive step, suppose the claims hold for factorizations
			of length $r-1$ and consider a factorization $J_1,\dots,J_r$.
		For any mask $\sigma\to_Jx$ we have
			$\upslice\sigma{r-1}\to_{\upslice J{r-1}}x\sigma_r^{-1}$.
		The inductive hypothesis then gives us $x\sigma_r^{-1}\leq
			w\qy*{\upslice J{r-1}}$, and hence
			\cref{prop:bruhat-coset}
			implies that $\rcos{J_r}x=\rcos{J_r}{(x\sigma_r^{-1})}
				\leq\rcos{J_r}{\qy*{w\qy*{\upslice J{r-1}}}}$,
			so since \cref{prop:finitary} implies $\rcos[J_r]{}x
				\leq \mxml{J_r}$,
			\cref{prop:bruhat-coset} gives 
				$x\leq\rcos{J_r}{\qy*{w\qy*{\upslice J{r-1}}}}\mxml{J_r}=w(J)$
				and the supplemental claim.

		For the main claim, note that we have $\upslice\sigma{r-1}\to_{\upslice J{r-1}}
			x\sigma_r^{-1}
			=\rcos{J_r}{\qy*{x\sigma_r^{-1}}}
				\rcos[J_r]{}{\qy*{x\sigma_r^{-1}}}
			=\rcos{J_r}{\qy*{x\sigma_r^{-1}}}
			\rcos[J_r]{\overlap J*r}{\qy*{x\sigma_r^{-1}}}
				\rcos[\overlap J*r]{}{\qy*{x\sigma_r^{-1}}}
			=\rcos{J_r}xy^\sigma_rz^\sigma_r$ by
			\cref{lem:length-coset-rep}.
		Thus, given $x'\in W$ and $z_r
				\in W_{\overlap J*r}$, it follows that a mask
			$\tau\to_J x'$ satisfies $y^\tau_r=y_r$ and $z^\tau_r=z_r$ if and
			only if it satisfies $\upslice\tau{r-1}
				\to_{\upslice J{r-1}}x'^{J_r}y_rz_r$, 
			and thus if and only if $\tau_r=\qy*{y_rz_r}^{-1}(x')_{J_r}
				\in W_{J_r}$.
		By \cref{prop:bruhat-coset}, we have that 
			$\rcos{\rdes J}x\geq\rcos{\rdes J}{(x')}$ together with
			$y^\sigma_r\geq y_r$ implies 
			$x^{\rdes J}y^\sigma_r\geq x'^{\rdes J}y_r$, and hence
			$\rcos{\rdups J{r-1}}{\qy*{x^{\rdes J}y^\sigma_r}}
				\geq \rcos{\rdups J{r-1}}{\qy*{x'^{\rdes J}y_r}}$.
		Since $z^\sigma_r,z_r\in W_{\overlap J*r}\subseteq 
			W_{J_r\cap\rdups J{r-1}}\subseteq
				W_{\rdups J{r-1}}$, $y^\sigma_r,y_r\in W_{J_r}$, and 
				$\rdes J=J_r\udj\qy*{\rdups J{r-1}\setminus\nbd{J_r}}$, it follows that 
			\begin{align*}
			\rcos{\rdups J{r-1}}{\qy*{x\sigma_r^{-1}}}
			&=\rcos{\rdups J{r-1}}{\qy*{x^{\rdes J}y^\sigma_rz^\sigma_r
				x_{\rdups J{r-1}\setminus\nbd{J_r}}}}
				\\
			&=\rcos{\rdups J{r-1}}{\qy*{x^{\rdes J}y^\sigma_r}}\\
			&\geq \rcos{\rdups J{r-1}}{\qy*{x'^{\rdes J}y_r}}\\
			&=\rcos{\rdups J{r-1}}{\qy*{{x'}^{\rdes J}y_rz_r
				x'_{\rdups J{r-1}\setminus\nbd{J_r}}}}\\
			&=\rcos{\rdups J{r-1}}{\qy*{{x'}^{J_r}y_rz_r}}\\
			&=\rcos{\rdups J{r-1}}{\qy*{x'\tau_r^{-1}}}.
			\end{align*}
		The inductive hypothesis then implies the existence of a unique mask
			$\upslice\tau{r-1} \to_{\upslice J{r-1}}\rcos{J_r}{x'}y_rz_r$ with
			$y^{\upslice\tau{r-1}}=(y_1,\dots,y_{r-1})$ and
			$z^{\upslice\tau{r-1}}=(z_1,\dots,z_{r-1})$.
	\end{proof}
\end{prop}

\begin{exmp}
	In the case of a singleton factorization induced by a reduced word
		$s_{i_1}\cdots s_{i_r}=x$, \cref{prop:defect-ideals} implies that
		every $x'\leq x$ can be realized as a subword of that reduced word,
		i.e.\ that the converse to \cref{lem:bruhat} holds.
\end{exmp}

\begin{cor}
\label{cor:unique-reps}
	The map $\sigma\mapsto(y^\sigma,\rcos{\rdes J}x)$ (where $y^\sigma$ is
		defined as in \cref{prop:defect-ideals}) is injective on
		the set of $\ds\qy*{\prod_{1\leq i<k\leq r}W_{\overlap Jik}}
			\times W_{\rdes J}$-orbits of masks, and for each $x\in W$ and
		each family of subsets $(D_t)_{t\in T}\subseteq[1,r]$ there is at
		most one equivalence class of masks $\sigma$ such that
		$[\sigma]\to_Jx$ and $D_t=\rdef\sigma t$ for all $t\in T$.
	\begin{proof}
		\cref{prop:defect-ideals} asserts that the map
			$\sigma\to_Jx\mapsto(\rcos{\rdes J}x, y^\sigma,x_{\rdes J},
				z^\sigma)$ (using the definition of $z^\sigma$ from that
			\namecref{prop:defect-ideals})
			is injective on masks, and that for each fixed
			$(y^\sigma,\rcos{\rdes J}x)$ its restriction
			$\sigma\to(z^\sigma,\rcos[\rdes J]{}x)$ is bijective onto
			$\ds\prod_{k=1}^r\qy*{W_{\overlap J*k}
				\times W_{\rdeslev kJ}}
				\cong\qy*{\prod_{1\leq i<k\leq r}W_{\overlap Jik}}
					\times W_{\rdes J}$.
		By \cref{lem:y-well} $\rcos{\rdes J}x$ and $y^\sigma$ are
			constant on $\ds\prod_{1\leq i<k\leq r}\qy*{W_{\overlap Jik}}
				\times W_{\rdes J}$-orbits of the action, and thus
			$\sigma\to_Jx\mapsto (y^\sigma,\rcos{\rdes J}x)$ is injective on
			the set of $\ds\qy*{\prod_{1\leq i<k\leq r}W_{\overlap Jik}}
				\times W_{\rdes J}$-orbits.
			
		Given $x\in W$ and $(D_t)_{t\in T}\subseteq[1,r]$, by
			\cref{prop:defect-ideals}, to show that there is at
			most one equivalence class of masks such that $D_t=\rdef\sigma t$
			for all $t$, it will suffice to show there is exactly one mask
			$\sigma\to_Jx$ such that $z^\sigma_k=e$ for all $k$ and
			$D_t=\rdef\sigma t$ for all $t\in T$.
		We argue by induction on the length $r$.

		The base case $r=1$ is immediate since there is only one mask
			$\sigma\to_Jx$ and it satisfies 
			$z^\sigma_1={\qy*{\sigma_1}}_{\overlap J*1}=e$ and
			$\rdef\sigma t=\emptyset$ for all $t\in T$.
		For the inductive step, consider two equivalent masks $\sigma\to_Jx$
			and $\tau\to_Jx$ such that $z^\sigma_k=e=z^\tau_k$ for all $k$. 
		Then $\rcos{\rdups Jr}{(\sigma_1\cdots\sigma_r)}
			=\rcos{\rdups Jr}x=\rcos{\rdups Jr}{(\tau_1\cdots\tau_r)}$
			and $\set{t\in T:r\in\rdef\sigma t}=\set{t\in T:r\in\rdef\tau t}$
			imply by \cref{lem:y-well} that $y^\sigma_r=y^\tau_r$.
		Consequently, $\upslice\sigma{r-1}\to_{\upslice J{r-1}}x\sigma_r^{-1}
			=\rcos{J_r}xy^\sigma_rz^\sigma_r=\rcos{J_r}xy^\tau_rz^\tau_r
				=x\tau_{r^{-1}}$, whence $\sigma_r=\tau_r$.
		But since by definition $z^{\upslice\sigma{r-1}}_k
			=z^{\upslice\tau{r-1}}_k$ and 
			$\rdef{\upslice\sigma{r-1}}t=\rdef\sigma t\setminus\set r
				=\rdef\tau t\setminus\set r=\rdef{\upslice\tau{r-1}}t$, the
			inductive hypothesis implies $\upslice\sigma{r-1}
				=\upslice\tau{r-1}$, and hence $\sigma=\tau$, as desired.
	\end{proof}
\end{cor}

\begin{rmrk}
	In the case of Weyl groups and factorizations with Gelfand--MacPherson
	 overlaps, the argument of \cref{cor:unique-reps} can be used to show
	 that the pre-image of $x$ is isomorphic to a closed subscheme of
	 $P_{J_1}/P_{J_2}\times P_{J_2}/P_{J_3}\times\cdots\times
	 P_{J_{m-1}}/P_{J_m}$.
\end{rmrk}

\begin{rmrk}
	Without the $W_{\rdes J}$-orbits claim, the special case of
		\cref{cor:unique-reps} for singleton factorizations associated to
		reduced words is \cite[Lemma~3.4]{Jones-Woo}, while the
		even more special case with $D=\emptyset$ is
		\cite[Proposition~2.3(iii)]{Deodhar}.
	The special case for singleton factorizations associated to reduced
		words of the $W_{\rdes J}$-orbits claim is a slight generalization
		of \cite[Proposition~2.3(iv)]{Deodhar}.
\end{rmrk}

\section{Defect polynomials and their basic properties}

\label{sec:masks-poly}
We continue our generalization of Deodhar's combinatorics of
	subexpressions and Billey and Warrington's reformulations thereof by
	associating to each factorization with overlaps a family of generating
	functions we call \emph{defect polynomials} and examining their
	properties, including a recursive formula for computing them.

We then use the Cartier--Foata partially commutative monoid generated by
	connected finitary subsets to define a notion of \emph{equivalence} of
	factorizations for which defect polynomials are invariant. 
We also generalize Billey and Jones' notion of \emph{factors}, for which
	we show that one factorization being a factor of another implies the
	defect polynomials of one are quotients by the same factor of the
	defect polynomials of the other. 

Finally, we generalize Billey and Jones' notion of \emph{Deodhar
	elements} by introducing the notion of \emph{tight} factorizations,
	which we show necessarily satisfy the property of being
	\emph{admissible}, and which will ensure their families of defect
	polynomials coincide with Kazhdan--Lusztig polynomials.
These properties will in fact both be \emph{absolute}, in the sense that
	they hold for a factorization only if they hold for each of its
	factors.

\begin{defn}
\label{defn:defect-polynomials}
	Given a factorization with overlaps $J_1,\dots,J_r$, we define its
		associated \define{family $(P^J_x)_{x\in W}\in\NN[q]$ of defect
		polynomials} to be given by
		$P^J_x=\ds\sum_{[\sigma]\to_Jx}q^{d_R(\sigma)}$ where
		$[\sigma]\to_Jx$ ranges over all equivalence classes of masks on $J$
		leading to $x$.
\end{defn}

\begin{rmrk}
	We use the a priori ambiguous term ``defect polynomials'' rather than
		the unambiguous ``right-defect polynomials'' because
		\cref{cor:quotient} implies that defect polynomials are a
		self-dual notion, even though right-defects themselves are not a
		self-dual notion.
\end{rmrk}

\begin{exmp}
\label{exmp:admissible-singleton}
	In the case of a singleton factorization $J$ induced by a reduced
		word, masks are in one-to-one correspondence with subwords, and
		hence $\sigma\to_Jw(J)$ if and only if $\sigma=\omega_R(J)$.
	In particular, $P^J_{w(J)}(q)=1$ by
		\cref{exmp:singleton-defects}.
\end{exmp}

\begin{exmp}
\label{exmp:r=1f}
	In the case of a factorization $J_1$ of length $1$, masks are in
		one-to-one correspondence with elements of $W_{J_1}$ while 
		overlap subsets (and hence the equivalence relation) is necessarily
		trivial.
	Thus, $P^J_x(q)=\begin{cases}1&x\in W_{J_1}\\
		0&x\not\in W_{J_1}\end{cases}=\abs{\set x\cap W_{J_1}}$.
\end{exmp}

\begin{exmp}
	In the case of a factorization $J$ with Gelfand--MacPherson overlaps
		(see \cref{exmp:gm-factorization}) and a Weyl group, the defect
		polynomials defined above match those described in
		\cref{subsec:weyl-gm}.
\end{exmp}

\begin{defn}
\label{defn:right-parabolic-product}
	Given a family of Laurent polynomials $(P_x)_{x\in W}\in\ZZ[q^{\frac12}
		,q^{-\frac12}]$, we define its \define{set $\rdes P$ of associated
		right-descents} to be $\rdes P=\set{s\in S:P_{xs}=P_x$ for all $x\in
		W}$.
	Similarly, we define its \define{set $\ldes{P}$ of associated
		left-descents} to be $\ldes{P}=\set{s\in S:P_{sx}=P_x$ for all $x\in
		W}$.
	Furthermore, given also a finitary subset $J\subseteq S$, for each subset 
		$J_*\subseteq \rdes P\cap J$, we denote by 
		$\qy{\rparp P{J_*}J_x}_{x\in W}$ the family given by 
		$\rparp P{J_*}J_x
			=\ds\sum_{y\in\rcos[J]{J_*}W}q^{\ell(y)}P_{\rcos Jxy}$.
\end{defn}

\begin{cor}
\label{cor:recursive-defect-polys}
	The defect polynomials of a factorization with overlaps 
		$J_1,\dots,J_r$ satisfy $\qy{P^{\upslice Jk}_x}_{x\in W}
		=\rparp{P^{\upslice J{k-1}}}{\overlap J*k}{J_k}$, i.e.\
		\[
			P^{\upslice Jk}_x
				= \sum_{y\in \rcos[J_k]{\overlap J*k}W}
				q^{\ell(y)}P^{\upslice J{k-1}}_{\qy*{\rcos{J_k}xy}}
		\]
		for every $1\leq k\leq r$ and $x\in W$.
	Furthermore, $P^J_{x'}(q)-P^J_x(q)\in\NN[q]$ whenever $x'\leq x$ in
		the Bruhat order, and $P^J_x\neq0$ if and only if $x\leq
		w(J)$ in which case $P^J_x(0)=1$.
	\begin{proof}
		Expanding the definition of the defect polynomials, this
			amounts to verifying
			\[
				\sum_{[\sigma]\to_{\upslice Jk} x} q^{d_R(\sigma)}
				=
				\sum_{y\in\rcos[J_k]{\overlap J*k}W}q^{\ell(y)}
				\sum_{[\tau] \to_{\upslice J{k-1}} x^{J_k}y}
				q^{d_R(\tau)}.
			\]
		By \cref{prop:defect-ideals} and
			\cref{cor:unique-reps} we have that
			each $[\sigma]\to_{\upslice Jk}x$ is defined uniquely by the choice
			of $y^{\upslice\sigma k}\in \prod_{j=1}^k \rcos[J_j]{\overlap J*j}W$,
			 and that choosing a representative of
			this class with $z^{\upslice\sigma k}_k = e$ uniquely defines an
			equivalence class of masks on $\upslice J{k-1}$ by
				$[\tau]=[\upslice\sigma{k-1}]
					\to_{\upslice J{k-1}}\rcos{J_k}xy^\sigma_k$
			with $\ell(y^\sigma_k)
				=\abs{\set{t\in T:k\in \rdef{\upslice\sigma k}t}}
				=d_R(\upslice\sigma k)-d_R(\tau)$,
			as desired.
		The supplemental claims are generating-function formulations
			of \cref{prop:defect-ideals}.
	\end{proof}
\end{cor}

\begin{lem}
\label{lem:pre-descent}
	For any family $(P_x)_{x\in W}\in\ZZ[q^{\frac12},q^{-\frac12}]$ we
		have that $\rdes P\subseteq \rdes{w}$ and $\ldes P\subseteq
		\ldes w$ for each $w\in W$ maximal in the Bruhat order among $x$
		such that $P_x\neq 0$.
	In particular, for any factorization with overlaps $J_1,\dots,J_r$ its
		family of defect polynomials $(P^J_x)_{x\leq w(J)}$ satisfies 
		$\rdes J\subseteq\rdes{P^J}\subseteq\rdes{w(J)}$.
	\begin{proof}
		Let $w$ be some maximal element such that $P_w\neq0$.
		If $s\in\rdes P$ then by definition we
			have that $P_{ws}=P_w\neq 0$ and thus $ws$ is not greater than
			$w$ in the Bruhat order, i.e.\ $ws<w$ and 
			$s\in S\cap\rref w=\rdes w$ by \cref{lem:t-in-rref-t}.
		Similarly, if $s\in\ldes P$, then $sw=ww^{-1}sw$ is not greater
			than $w$ in the Bruhat order, so $w^{-1}sw\not\in\rref w$, i.e.\
			$s\in S\cap w\rref ww^{-1}=S\cap\lref w=\ldes w$.

		For the supplemental claim, the invariance of defects under
			the action of $W_{\rdes J}$ by \cref{prop:defects-are-meetings}
			implies that $\rdes J \subseteq \rdes{P^J}$,
			and $\rdes{P^J}\subseteq\rdes{w(J)}$ follows from
			the main claim and the fact that $P^J_x\neq 0$ if and only
			if $x\leq w(J)$ by \cref{cor:recursive-defect-polys}.
	\end{proof}
\end{lem}

\begin{rmrk}
	\label{rmrk:polynomial-descents}
	We call the set $\rdes P$ the ``right-descents'' of the family $P$
	for an analogous reason to that discussed in
	\cref{rmrk:factorization-descents} in the case of
	factorizations: because \cref{lem:pre-descent} above shows
	that $\rdes P\subseteq\rdes w$ for families
	$(P_x)_{x\leq w}$ with $P_w\neq 0$. 
\end{rmrk}

\begin{rmrk}
	Although the inequalities $\rdes J\subseteq \rdes{P^J}\subseteq
	\rdes{w(J)}$ from \cref{lem:pre-descent} need not be equalities
	in general, we will show in \cref{prop:rdescent-inductive-step}
	that their equality characterizes the important class of
	right-descent factorizations. 
\end{rmrk}

The defect polynomials and most of our theory can be developed
	in a more general setting than factorizations with overlaps by dropping the
	requirement of a linear ordering on the finitary subsets, and imposing
	instead the structure of a \emph{Viennot heap} \cite{Viennot}.
Using this language exclusively, however, is cumbersome, and our
	arguments only require the following equivalence relation on
	factorizations in terms of the equivalent 
	\emph{Cartier--Foata monoid}.

\begin{defn}
\label{defn:cartier--foata}
	The \define{Cartier--Foata partially commutative
		monoid of finitary subsets $J\subseteq S$} (see
		\cite{Cartier--Foata}) is the semigroup freely generated by the
		finitary subsets $J\subseteq S$ subject to the relations
		$J_1J_2=J_1\sqcup J_2=J_2J_1$ if $J_1$ and $J_2$ commute (i.e.\ if
		each element of $J_1$ commutes with and is distinct from each
		element of $J_2$).
	We say that two factorizations are \define{equivalent} if they induce
		the same element of the Cartier--Foata monoid.
\end{defn}

\begin{rmrk}
	$\emptyset\udj J=J=J\udj\emptyset$ ensures that $\emptyset$ is the
		unit of the Cartier--Foata monoid.
\end{rmrk}

\begin{exmp}
	The factorizations of $45312$ given by
		$\set{s_2,s_3},\set{s_1,s_3,s_4},\set{s_2,s_3}$ and
		$\set{s_2,s_3},\set{s_3,s_4},\set{s_1},\set{s_2,s_3}$ determine the
		same element of the Cartier--Foata monoid, which is distinct from
		the element determined by the factorization
		$\set{s_2,s_3},\set{s_1,s_2,s_4},\set{s_2,s_3}$.
	The first is illustrated in \cref{fig:right-canonical}, while the
		second is its mirror image.
\end{exmp}

\begin{exmp}
\label{exmp:cartier--foata-visualization}
	In terms of the visualization of factorizations as in
		\cref{exmp:interval_decomposition} and \cref{fig:mask},
		the generating relations of the Cartier--Foata monoid consist of
		partitioning a level of disjoint segments into two levels, or
		joining two levels of disjoint segments into a single level.
\end{exmp}

\begin{prop}
\label{prop:expansions}
	If $J=J_1\sqcup J_2$ where each element of $J_1$ commutes with each
		element of $J_2$, then for any finite
		$(P_x)_{x\in W}\in\ZZ[q^{\frac12},q^{-\frac12}]$ and set
		$J_*\subseteq J\cap \rdes P$ we have
		$\rparp P {J_*} J=\rparp{\qy{\rparp P
			{J_*\cap J_1} {J_1}}} {J_*\cap J_2} {J_2}$.
	\begin{proof}
		By \cref{lem:coset-properties}\ref{lem:coset-unique}) 
			we have that the group operation of $W$ induces isomorphisms
			$W^{J_2}\times W_{J_2}\cong W\cong W^{J_1}\times W_{J_1}$ and
			$W\cong W^J\times W_J\cong W^J\times W_{J_1}\times W_{J_2}\cong
			W^J\times W_{J_2}\times W_{J_1}$.
		It follows that $x^J=(x^{J_1})^{J_2}$ and that $(xz)^{J_1}=x^{J_1}z$
			for any $z\in W_{J_2}$.
		Thus, for any $y\in W_{J_1}$ we have
			$(xz)^{J_1}y=x^{J_1}zy=x^{J_1}yz$.

		Since $J_*\subseteq\rdes P$, the expression $\rparp P{J_*\cap
			J_1}{J_1}$ is well-defined.
		Since for any $s\in J_*\cap J_2$, we have not only $s\in\rdes P$ but
			also $(xs)^{J_1}y=x^{J_1}ys$ for any $x\in W$ and any $y\in
			W_{J_1}$ such that $\rcos Jxy\leq w$, it follows that 
			$\qy*{\rparp P{J_*\cap J_1}{J_1}}_{xs}
				=\ds\sum_{y\in W_{J_1}^{J_1\cap J_*}}
					q^{\ell(y)}P_{\rcos{J_1}{\qy{xs}}y}
				=\sum_{y\in W_{J_1}^{J_1\cap J_*}}
					q^{\ell(y)}P_{\rcos{J_1}xys}
				=\sum_{y\in W_{J_1}^{J_1\cap J_*}}q^{\ell(y)}P_{x^{J_1}y}
				=\qy{\rparp P{J_1\cap J_*}{J_1}}_x$,
			and hence that $s\in\rdes{\rparp P{J_1\cap J_*}{J_1}}$.
		Thus, $J_2\cap J_*\subseteq \rdes{P*_{J_*\cap J_1}P^{J_1}}$ so that
			$\rparp{\qy{\rparp P {J_*\cap J_1}{J_1}}}{J_*\cap J_2}{J_2}$ is
			well-defined.

		Next, the group operation also induces isomorphisms 
			$W_J^{J_*}\times W_{J_*}\cong W_J\cong W_{J_1}\times W_{J_2}
				\cong W_{J_1}^{J_1\cap J_*}\times W_{J_1\cap J_*}
					\times W_{J_2}^{J_2\cap J_*}\times W_{J_2\cap J_*}
			\cong W_{J_1}^{J_1\cap J_*}\times W_{J_2}^{J_2\cap J_*}
				\times W_{J_1\cap J_*}\times W_{J_2\cap J_*}
			\cong W_{J_1}^{J_1\cap J_*}\times W_{J_2}^{J_2\cap J_*}
				\times W_{J_*}$ and hence an isomorphism 
			$\rcos[J]{J_*}W
				\cong\rcos[J_2]{J_2\cap J_*}W\times\rcos[J_1]{J_1\cap J_*}W$.
		Thus
			$\qy{\rparp{\qy{\rparp P {J_*\cap J_1}{J_1}}}{J_2\cap J_*}{J_2}}_x
			=\ds\sum_{z\in\rcos[J_2]{J_2\cap J_*}W}q^{\ell(z)}
				\qy{\rparp P {J_*\cap J_1}{J_1}}_{\rcos{J_2}xz}
			=\sum_{z\in\rcos[J_2]{J_2\cap J_*}W}
				q^{\ell(z)}\sum_{y\in\rcos[J_1]{J_1\cap J_*}W}
					q^{\ell(y)}P_{\rcos{J_1}{\qy{\rcos{J_2}xz}}y}
			=\sum_{z\in\rcos[J_2]{J_2\cap J_*}W}
				\sum_{y\in\rcos[J_1]{J_1\cap J_*}W}
					q^{\ell(zy)}P_{\rcos{J_1}{\qy{\rcos{J_2}x}}zy}
			=\break\sum_{z\in\rcos[J_2]{J_2\cap J_*}W}
				\sum_{y\in\rcos[J_1]{J_1\cap J_*}W}q^{\ell(zy)}P_{\rcos Jxzy}
			=\sum_{y\in\rcos[J]{J_*}W}q^{\ell(y)}P_{\rcos Jxy}
			=\qy{\rparp P {J_*} J}_x$.
	\end{proof}
\end{prop}

\begin{cor}
\label{cor:spread_defect}
	The map $J\mapsto(P^J_x)_{x\leq w(J)}$ induces a well-defined
		function from the Cartier--Foata monoid to positive families
		$(P_x)_{x\leq w(J)}\in\NN[q]$.
	In particular, the factorization obtained from $J_1,\dots,J_r$ by
		replacing each $J_k$ with a sequence of all of its connected
		components has the same defect polynomials.
	\begin{proof}
		Let $J_k=J_{k,1}\sqcup J_{k,2}$ with $J_{k,1}$ commuting with
			$J_{k,2}$.
		If we denote by $\upslice J{k-1},J_{k,1}$ and $\upslice J{k-1},J_{k,1},J_{k,2}$ the factorizations
			$J_1,\dots,J_{k-1},J_{k,1}$ and
			$J_1,\dots,J_{k-1},J_{k,1},J_{k,2}$ respectively, then
			$\rdes{\upslice J{k-1}}\cap J_{k,1}
				=\rdes{\upslice J{k-1}}\cap J_k \cap J_{k,1}$ because
			$J_k\supseteq J_{k,1}$, and $\rdes{\upslice J{k-1}}\cap J_{k,2}
				=\qy*{\rdes{\upslice J{k-1}}\setminus \nbd{J_{k,1}}}\cap J_{k,2}
				=\qy*{J_{k,1}\udj \qy*{\rdes{\upslice J{k-1}}
					\setminus \nbd{J_{k,1}}}}\cap J_{k,2}
				=\rdes{\upslice J{k-1},J_{k,1}}\cap J_{k,2}$, where the first
				equality is because $J_{k,1}$ commutes with $J_{k,2}$.

		\cref{cor:recursive-defect-polys} and
			\cref{prop:expansions} then yield $P^{\upslice Jk}
				=\rparp{P^{\upslice J{k-1}}}{\rdes{\upslice J{k-1}}\cap J_k}{J_k}
				=\rparp{\qy*{\rparp{P^{\upslice J{k-1}}}{\rdes{\upslice J{k-1}}
					\cap J_{k,1}}{J_{k,1}}}}{\rdes{\upslice J{k-1}}\cap J_{k,2}}{J_{k,2}}
				=\rparp{P^{\upslice J{k-1},J_{k,1}}}{\rdes{\upslice J{k-1},J_{k,1}}\cap J_{k,2}}{J_{k,2}}
				=P^{\upslice J{k-1},J_{k-1},J_{k-2}}$. 
		Thus, any two factorizations $J^1$ and $J^2$ related by
			the generating relations of the Cartier--Foata monoid satisfy
			$P^{J^1}=P^{J^2}$, and consequently $J\mapsto P^J$ is well-defined
			on the monoid as desired.
	\end{proof}
\end{cor}

\begin{defn}
	We say that a factorization $J'_1,\dots,J'_t$ is a \define{factor} of
		a factorization $J_1,\dots,J_r$ if $J$ induces an element of the
		Cartier--Foata monoid contained in the two-sided principal ideal
		generated by the element induced by $J'$, i.e.\ if $J'_1,\dots,J'_t$
		it is a slice of a factorization $J''$ inducing the same element of
		the Cartier--Foata monoid as $J$.
	We say that a factorization $J_1,\dots,J_r$ is \define{absolutely $P$}
		for a property $P$ if every factor $J'$ of $J_1\cdots J_r$ satisfies
		property $P$.
	We say a property $P$ is itself \define{absolute} if any factorization 
	  satisfying $P$ is absolutely $P$.
\end{defn}

\begin{exmp}
\label{exmp:short-braid-factor}
	Every reduced word of a Coxeter element determines an element of the
		Cartier--Foata monoid.
	Iwahori's theorem \cite{Matsumoto} asserts that one can pass from any
		reduced word of a Coxeter element $w$ to any other reduced word by
		successive replacements of contiguous subwords
		$\underbrace{s_is_js_is_js_i\cdots}_{m(i,j)}
			\mapsto\underbrace{s_js_is_js_i\cdots}_{m(i,j)}$,
		where $m(i,j)$ is the order of $s_is_j$.
	By definition $s_is_js_i$ is a \emph{short braid} if and only if
		$m(i,j)\geq3$, and it follows that all reduced words for a Coxeter
		element $w$ induce the same element of the Cartier--Foata monoid if
		and only if the element of the Cartier--Foata monoid induced by some
		reduced word for $w$ has no factor that is induced by a short braid.
	Elements $w$ with this property are exactly the \emph{fully
		commutative elements} introduced in \cite{Stembridge96}.
\end{exmp}	

\begin{rmrk}
	Billey and Jones \cite[Section~5]{Billey--Jones} say that \emph{$y$ is
		a factor of $w$} if there exist elements $x$ and $z$ such that
		$w=xyz$ with $\ell(w)=\ell(x)+\ell(y)+\ell(z)$.
	This is equivalent to asserting that \emph{some} factorization
		associated to a reduced word for $y$ is a factor in our sense of a
		factorization associated to a reduced word for $w$.
	Furthermore, since reduced words of a short braid avoiding element
		induce a unique element of the Cartier--Foata monoid, in that case
		their notion is completely recovered from ours.
\end{rmrk}

As we show later (\cref{cor:factors}) the following properties of
	factorizations with overlaps are absolute.
\begin{defn}
\label{defn:tight}
	We say a factorization with overlaps $J$ is \define{admissible} if
		$P^J_{w(J)}=1$, and \define{tight} if $\deg P^J_x(q)
			-\frac12\qy*{\deg P^J_{w(J)}(q)+\len{w(J)}-\len x}<0$ when
		$x<w(J)$.
\end{defn}

\begin{lem}
\label{lem:admissible-double-max}
	Given a family $(P_x)_x\in\ZZ[q^{\frac12},q^{-\frac12}]$ such that
		$P_x\neq0$ only if $x\leq w$ and subsets $J_1,J_2\subseteq S$ such
		that $J_1\cup J_2\neq\emptyset$, $P_{sx}=P_x$ for $s\in J_1$ and
		$P_{xs}=P_x$ for $s\in J_2$, we have 
		$\deg P_x(q)-\frac12\qy{\deg P_w(q)-\deg P_w(q^{-1})+\ell(w)
			-\ell(x)}<0$ for $x\in W_{J_1}wW_{J_2}\setminus\set w$ if and only
			if $\deg P_w(q)+\deg P_w(q^{-1})\leq 0$, i.e.\ if and only if
			$P_w$ is a monomial.
	\begin{proof}
		By definition of $J_1$ and $J_2$, for each $x\in W_{J_1}wW_{J_2}$ we
			have $P_x=P_w$, so that 
			$\deg P_x(q)-\frac12\qy{\deg P_w(q)-\deg P_w(q^{-1})+\ell(w)
				-\ell(x)}<0$ is equivalent to 
			$\deg P_w(q)+\deg P_w(q^{-1})<\ell(w)-\ell(x)$.
		Note that $(P_x)_{x\leq w}$ implies that for each 
			$x\in W_{J_1}wW_{J_2}\setminus\set w$ either $x<w$ or $0=P_x$.
		In the latter case, $\deg P_w(q)+\deg P_w(q^{-1})
			=-\infty<\ell(w)-\ell(x)$ (see
			\cref{rmrk:bounded-has-pw-nonzero}), and in the former case
			$\emptyset\neq J_1\cup J_2$ implies $\ell(w)-1
				=\max\set{y\in W_{J_1}wW_{J_2}$ and $\ell(y):y<w}$.
		Thus, $\deg P_w(q)+\deg P_w(q^{-1})<\ell(w)-\ell(y)$ for all such
			$y$ is equivalent to $\deg P_w(q)+\deg P_w(q^{-1})<1$, i.e.\ to 
			$\deg P_w(q)+\deg P_w(q^{-1})\leq 0$, as desired.
	\end{proof}
\end{lem}

\begin{prop}
\label{prop:bounded-monic}
	A factorization $J_1,\dots,J_r$ of $w$ is tight only if it is
		admissible, in which case it satisfies $\deg
		P^J_x(q)-\frac12\qy*{\len{w(J)}-\len x}<0$ when $x<w(J)$.
	\begin{proof}
		If $\rdes J=\emptyset$, the fact that 
			$\rdups Jk =J_k\udj\qy*{\rdups J{k-1}\setminus\nbd{J_k}}$ implies
			recursively that $J_k=\emptyset$ for all $k$, in which case
			$w=w(J)=e$ so that $P^J_w=P^J_e=1$.

		If $\rdes J\neq\emptyset$, then since $P^J_x(0)=1$ when $x<w$ and
			$P^J_x=0$ for when $w<x$ by \cref{cor:recursive-defect-polys}, it
			follows that $J$ is tight only if $\deg P_x(q)-\frac12\qy{\deg
			P^J_w(q)-\deg P^J_w(q^{-1})+\ell(w)-\ell(x)}<0$ for all 
			$x\in wW_{\rdes J}\setminus\set w$.
		Since $s\in\rdes J$ implies $P_{xs}=P_x$ by \cref{lem:action}, 
			\cref{lem:admissible-double-max} implies $P^J_w$ is a monomial
			with $P^J_w(0)=1$, i.e.\ $P^J_w=1$.
	\end{proof}
\end{prop}

\begin{rmrk}
	A reduced word for an element $w$ inducing a tight factorization
		coincides with Billey and Jones's \cite{Billey--Jones} notion of $w$
		being a \emph{Deodhar element}.
	Because \cref{exmp:short-braid} will show that factorizations induced
		by short braids are not tight, \cref{cor:factors} will imply that
		short braid avoiding, i.e.\ fully commutative elements, include all
		Deodhar elements.
	In particular tight factorizations being absolute (\cref{cor:factors})
		generalizes the fact (\cite[Proposition~5.2]{Billey--Jones}) that
		every factor in the sense of Billey and Jones of a Deodhar element
		$w$ is itself a Deodhar element.
\end{rmrk}

\section{The combinatorial criterion}

\label{sec:combi}
We complete our generalization of Deodhar's combinatorics of
	subexpressions and Billey and Warrington's reformulations thereof by
	generalizing Billey and Warrington's combinatorial criterion for an
	element of a Coxeter group being Deodhar to a criterion for when a
	factorization is tight.
As an example application, we recover the facts proved by Fan and Green
	\cite[Proposition~2.1.3]{Fan--Green} that reduced words of Deodhar
	elements are \emph{short braid avoiding} and that a permutation is
	Deodhar only if it is \emph{hexagon-avoiding} (reproved in \cite[Proof
	of Theorem~1.~2$\Rightarrow$3]{Billey--Warrington}).
Along the way, we generalize more of the basic properties of defect
	polynomials in Deodhar's original framework proven in \cite{Deodhar},
	and in the end summarize the extent to which Billey and Warrington's
	defect graph argument would translate to our setting.

We begin by generalizing a second statistic on masks used by Deodhar,
	together with its complementary statistic identified by
	Billey and Warrington in their reformulation of the combinatorics of
	subexpressions.

Recall from \cref{defn:defects} that, given a mask
	$\sigma=(\sigma_1,\dots,\sigma_r)$ on a factorization $J_1,\dots,J_r$
	and a reflection $t\in T$, we have a right-meeting 
	$i\in\rmeet\sigma t$ if and only if 
	$t^\sigma_{i-1}\in W_{J_i}\setminus W_{\overlap J*i}$ where
	$(t^\sigma_1,\dots,t^\sigma_r)$ is the virtual pair of strings of
	$\sigma$ right-labeled by $t$ given by
	$t^\sigma_m=(\sigma_1\cdots\sigma_m)^{-1}t(\sigma_1\cdots\sigma_m)$.

\begin{defn}
\label{defn:meetings}
	Given a mask $\sigma$ on a factorization with overlaps $J_1,\dots,J_r$
		we say that a virtual pair of strings $t^\sigma$
		\define{right-travel together after} $j$ for $j$ such that
		$1\leq j\leq r$ if either $t^\sigma_j\in W_{\overlap Jik}$ for $i,k$
		such that $i\leq j<k$, or $t^\sigma_j\in W_{\rdeslev iJ}$ for $i$
		such that $i\leq j$.

	Given a factorization with overlaps $J_1,\dots,J_r$, we say that a
		right-meeting $m\in\rmeet\sigma t$ has \define{right-separation
		after $k$} if $k$ is minimal such that $m\leq k$ and $t^\sigma$ do
		not right-travel together after $k$ (where the minimum of
		the empty subset of $\set{1,\dots,r}$ is defined to be $r$).

	We say that a right-meeting $i\in\rmeet\sigma t$ with
		right-separation after $k$ is a \define{right-crossing} if
		$t^\sigma_{i-1}\in\lref{\sigma_i\cdots\sigma_k}$ and a
		\define{right-bounce} otherwise.

	We denote by $\rcross\sigma t$ the \define{set of right-crossings} 
		of the virtual pair of strings $(t^\sigma_0,\dots,t^\sigma_r)$ of
		$\sigma$ right-labeled by $t\in T$, and by $\rbounce\sigma t$ its
		\define{set of right-bounces}.
	We define the \define{right-crossing statistic $c_R(\sigma)$} of a
		mask $\sigma$ to be given by $c_R(\sigma)=\ds\sum_{t\in
		T}\abs{\rcross\sigma t}$ and the \define{right-bounce statistic
		$b_R(\sigma)$} to be given by $b_R(\sigma)=\ds\sum_{t\in
		T}\abs{\rbounce\sigma t}$.
\end{defn}

\begin{exmp}
	In the case of Coxeter systems of type $A$, right-separation of a
		right-meeting at height $i$ of pair of strings is the lowest height
		$k$ at least $i$ after which the pair of strings do not go into the
		same overlap interval.
	If the order of the pair of strings is reversed from right before
		their right-meeting to right after the right-separation, then the
		right-meeting is a right-crossing, and otherwise it is a
		right-bounce. 
	
	For example, in \cref{fig:mask}, the right-meeting at height $1$
		of the pair of strings labeled by $(2,3)$ has a right-separation
		at height $2$ since it is there they first enter distinct overlap
		intervals above height $1$.
	Furthermore, this right-meeting is
		a right-crossing since they are in order right before height $1$ but
		are out-of-order right after height $2$, where the fact that they
		are out of order right after height $1$ is irrelevant since they
		are still right-travelling together at that point.

	On the other hand, the pair of strings labeled by $(3, 4)$ 
		has a right-meeting at height $1$ with a right-separation at height
		$1$ that is a right-crossing, and a right-meeting at height $3$ with
		a right-separation at height $3$ that is a right-bounce.
\end{exmp}

\begin{exmp}
\label{exmp:singleton-meet}
	For a factorization $J$ given by $J_j=\{s_{i_j}\}$ induced by a
		reduced word $s_{i_1}\cdots s_{i_r}=w$,
		\cref{exmp:singleton-defects} asserts that 
		$m\in\rmeet{\sigma}{t}$ if and only if $t^\sigma_{m-1}=s_{i_m}$,
		i.e.\ if and only if
		$t^\sigma_m=s_{i_m}t^\sigma_{i_{m-1}}s_{i_m}=s_{i_m}s_{i_m}s_{i_m}=s_{i_m}$.
	Since \cref{exmp:trivial-overlaps} also asserts that
		overlap subsets are empty for such factorizations, we also have that
		the right-separation of $i=m$ is either $k=r$ if $s_{i_m}$ commutes
		with and is distinct from $s_{i_j}$ for all $j$ such that $m<j$, or
		$k=m$, i.e.\ itself.
	In the first case that $k=r$, we have that $m\in\rcross\sigma t$
		if and only if $t^\sigma_m=s_{i_m}\in\lref{\sigma_m\cdots\sigma_r}
		\cap W_{J_m}=\qy*{\lref{\sigma_m}
		\udj\lref{\sigma_{m+1}\cdots\sigma_r}}\cap W_{J_m}
		=\lref{\sigma_m}$, and in the second case that $k=m$
		also if and only if $t^\sigma_m=s_{i_m}\in\lref{\sigma_m}$,
		i.e.\ for $J$ induced by a reduced word we have
		$m\in\rcross\sigma t$ if and only if $t^\sigma_m=s_{i_m}=\sigma_m$.
	In particular, the total number of right-crossings is given by
		$\ds c_R(\sigma)=\sum_{t\in T}\abs{\rcross\sigma t}
		=\abs{\set{k:s_{i_k}=\sigma_k}}=\sum_{k=1}^r\ell(\sigma_k)$, i.e.\
		$c_R(\sigma)$ coincides with the statistic Deodhar denotes by
		$\ell(\sigma)$ \cite[Section~2]{Deodhar}.
	Similarly, $b_R(\sigma)=\ds\sum_{t\in T}
		\abs{\rbounce\sigma t}=\abs{\set{k:s_{i_k}\neq\sigma_k}}$.
\end{exmp}

The following proposition ensures that right-crossings and
	right-bounces are invariant under the equivalence relation on masks.
\begin{prop}
\label{prop:crossing-bounce-equivalence}
	Given a factorization with overlaps $J_1,\dots,J_r$, 
		a reflection $t\in T$, and a mask
		$\sigma$ on $J$, whether $m\in\rmeet\sigma t$ has right-separation
		$k$ depends only on the $\ds\qy*{\prod_{1\leq i<k\leq r}
			W_{\overlap Jik}}\times W_{\rdes J}$-orbit of the mask $\sigma$,
		and whether $m\in\rbounce\sigma t$ or $m\in\rcross\sigma t$ depends
		only on the $\ds\prod_{1\leq i<k\leq r} W_{\overlap Jik}$-orbit of
		$\sigma$ (and if $k<r$, only on the $\ds\qy*{
			\prod_{1\leq i<k\leq r} W_{\overlap Jik}}\times W_{\rdes J}$-orbit).
	\begin{proof}
		\cref{lem:action} implies it will suffice
			to consider the case where $\tau$ is the image of $\sigma$ under
			the action of $s\in\ds\qy*{\bigudj_{1\leq i<k\leq r}
			\overlap Jik}\udj\qy*{\bigudj_{i=1}^r\rdeslev iJ}$.
		Then \cref{lem:action-subexprs-vpstrings} asserts that for any
			$m$ we have either $t^\tau_m=t^\sigma_m$, or
			$t^\tau_m=st^\sigma_ms$ where either $s\in\rdeslev iJ$ with
			$i\leq m$ or $s\in\overlap Jik$ with $i\leq m<k$.
		In particular, $s$ commutes with each element of $\rdeslev{i'} J$ so
			that $sW_{\rdeslev{i'}J}s=W_{\rdeslev{i'}J}$ for any $i'\leq m$,
			and given $i'$, $k'$ such that $i'\leq m<k'$, if
			$i'=i$ and $k'=k$ then $s\in\overlap J{i'}{k'}$, and
			otherwise either $i<i'<k$ or $i<j'<k$ so that by definition
			$s$ commutes with each element of $\overlap J{i'}{k'}\subseteq
			J_{i'}\cap J_{k'}$.
		Either way, $sW_{\overlap J{i'}{k'}}s=W_{\overlap J{i'}{k'}}$.

		In particular $t^\sigma$ right-travels after $m$, i.e.\ 
			$t^\sigma_m\in W_{\rdeslev{i'}J}$ for $i'$ such that $i'\leq m$ or
			$t^\sigma_m\in W_{\overlap J{i'}{k'}}$ for $i'$,$k'$ such that
			$i'\leq m<k'$ if and only if $t^\tau_m\in W_{\rdeslev{i'}J}$ or
			$t^\tau_m\in W_{\overlap J{i'}{k'}}$, i.e.\ if and only if
			$t^\tau$ right-travels after $m$. 
		Since $\rmeet\sigma t=\rmeet\tau t$ by
			\cref{prop:defects-are-meetings} we have
			$m\in\rmeet\sigma t$ with right-separation after $k$ if and only
			if $m\in\rmeet\tau t$ with right-separation after $k$.

		Next, $\rbounce\sigma t\udj\rcross\sigma t=\rmeet\sigma
			t=\rmeet\tau t=\rbounce\tau t\udj\rcross\tau t$ and thus 
			$\rbounce\sigma t=\rbounce\sigma t$ and $\rcross\sigma t
				=\rcross\tau t$ hold if and only if $m\not\in\rbounce\sigma
				t+\rbounce\tau t =\rcross\sigma t+\rcross\tau t$ for all $m$.
		Since \cref{prop:defects-are-meetings} also gives us
			$\rdef\sigma t=\rdef\tau t$, it follows that $\rbounce\sigma
			t=\rbounce\tau t$ and $\rcross\sigma t=\rcross\tau t$ hold
			if and only if $m\not\in\qy*{\rcross\sigma t+\rdef\sigma t}
				+\qy*{\rcross\tau t+\rdef\tau t}$ for any $m\in\rmeet\sigma
				t=\rmeet\tau t$.
		
		If $m\in\rmeet\sigma t=\rmeet\tau t$ with right-separation after
			$k$, then %
			$m\not\in\rdef\sigma t+\rcross\sigma t$ if and only if 
			$t^\sigma_{m-1}\not\in\rref{\sigma_1\cdots\sigma_{m-1}}
						+\lref{\sigma_m\cdots\sigma_k}
					=\rref{\sigma_1\cdots\sigma_{m-1}}
						+\qy*{\sigma_m\cdots\sigma_k}\rref{\sigma_m\cdots\sigma_k}
							\qy*{\sigma_m\cdots\sigma_k}^{-1}
					=(\sigma_m\cdots\sigma_k)\rref{\sigma_1\cdots\sigma_k}
						(\sigma_m\cdots\sigma_k)^{-1}$.

		Since $t^\sigma_{m-1}=(\sigma_m\cdots\sigma_k)t^\sigma_k
			(\sigma_m\cdots\sigma_k)^{-1}$, we thus have that
			$m\not\in\qy*{\rcross\sigma t+\rdef\sigma t}
				+\qy*{\rcross\tau t+\rdef\tau t}$ is equivalent to
			$\abs{t^\sigma_k\cap\rref{\sigma_1\cdots\sigma_k}}
			=\abs{t^\tau_k\cap\rref{\tau_1\cdots\tau_k}}$.

		To show that $\abs{t^\sigma_k\cap\rref{\sigma_1\cdots\sigma_k}}
			=\abs{t^\tau_k\cap\rref{\tau_1\cdots\tau_k}}$, note that
			\cref{lem:action-subexprs-vpstrings} yields either that
			$\tau_1\cdots\tau_k=\sigma_1\cdots\sigma_k$ and
			$t^\sigma_k=t^\tau_k$, in which case we are done, or otherwise
			$\tau_1\cdots\tau_k=\sigma_1\cdots\sigma_ks$, $t^\tau_k =
			st^\sigma_ks$ where either $s\in\rdeslev{i'}J$ for $i'$ such that
			$i'\leq k$ or $s\in\overlap J{i'}{k'}$ for $i'$, $k'$ such that
			$i'\leq k<k'$. 
		Then the desired $\abs{t^\sigma_k\cap\rref{\sigma_1\cdots\sigma_k}}
				=\abs{t^\tau_k\cap\rref{\tau_1\cdots\tau_k}}$
			is equivalent to $t^\tau_k\not\in\rref{\tau_1\cdots\tau_k}
					+s\rref{\sigma_1\cdots\sigma_k}s
				=\rref{\sigma_1\cdots\sigma_ks}+s\rref{\sigma_1\cdots\sigma_k}s
				=\set s$, i.e.\ to $t^\tau_k\neq s$.

		If $k=r$, then the hypothesis that $s\in\ds\bigudj_{1\leq
			i<j\leq r}\overlap Jij$ but $s\not\in\rdes J$ implies by
			\cref{lem:action} that
			$\tau_1\cdots\tau_k=\sigma_1\cdots\sigma_k$, which is the already
			handled trivial case.
			
		Otherwise, $m\in\rmeet\tau t$ with right-separation after $k<r$
			means that $t^\tau$ is not right-traveling after $k$, which in
			particular requires that $t^\tau_k\not\in W_{\overlap J{i'}{k'}}
				\ni s$ and $t^\tau_k\not\in W_{\rdeslev{i'}J}\ni s$, as desired.
	\end{proof}
\end{prop}

Similar to our derivation of \cref{cor:recursive-defect-polys}, we
	derive our combinatorial criterion for tightness of a factorization by
	picking specific representatives of equivalence classes of masks.
This time we pick representatives that generalize the property of the
	right-canonical mask $\omega_R(J)$ that restricting to an initial
	slice $\upslice Jk$ does not affect whether a right-meeting at height
	less than $k$ is a right-crossing or not.
\begin{prop}
\label{prop:equiv-induct}
	Each equivalence class of masks $[\sigma]\to_Jx$ on a factorization
		with overlaps $J_1,\dots,J_r$ has a unique representative $\sigma$
		such that $k\in\rcross\sigma t$ if and only if
		$t^\sigma_{k-1}\in\lref{\sigma_k}$, in which case
		$\rcross{\upslice\sigma{r-1}}t=\rcross\sigma t\setminus\set r$.
	\begin{proof}
		We argue by induction on the length $r$ of the factorization.
		For the base case, $r=1$, the result is trivial because the
			equivalence relation is trivial so that $[\sigma]\to_{J_1}x$
			implies $\sigma_1=x$, $\rcross{\upslice\sigma0}t=\emptyset
				=\rbounce{\upslice\sigma0}t$,
				and $1\in\rcross\sigma t$ if and only if
				$t=t^\sigma_0\in\lref{\sigma_1}$ by definition.
		
		For the inductive step, suppose such a representative exists for
			equivalence classes of masks on factorizations of length $r-1$,
			and suppose that $\sigma\to_Jx$ for a factorization
			$J_1,\dots,J_r$.
		Since by definition $r\in\rcross\tau t$ if and only if
			$t^\tau_{r-1}\in\lref{\tau_r}\setminus W_{\overlap J*r}$, it
			follows that $r\in\rcross\tau t$ is equivalent to
			$t^\tau_{r-1}\in\lref{\tau_r}$ if and only if $W_{\overlap
				J*r}\cap\lref{\tau_r}=\emptyset$, i.e.\ if and only if
			$\tau_r\in\lcos[J_r]{\overlap J*r}W$.
		But since \cref{lem:action} implies that the map
			$\set{\tau\to_Jx:\tau\sim\sigma}\to W_{J_r}$ given by
			$\tau\mapsto\tau_r$ has range the right coset
			$W_{\bigudj_{i=1}^{r-1}\overlap Jir}\sigma_r=
				W_{\overlap J*r}\sigma_r$, it follows that for $\tau$
			equivalent to $\sigma$ we have $r\in\rcross\tau t$ equivalent
			to $t^\tau_{r-1}\in\lref{\tau_r}$ if and only if
			$\tau_r=\lcos{\overlap J*r}\sigma_r$.
		
		Furthermore, \cref{lem:action} also implies that two masks
			$\tau$ and $\tau'$ such that $\tau_r=\tau'_r$ are equivalent if
			and only if $\upslice\tau{r-1}$ and $\upslice{\tau'}{r-1}$ are
			equivalent.
		Consequently, by the inductive hypothesis we have that $\sigma$ is
			equivalent to a unique mask $\tau$ such that on the one hand
			$r\in\rcross\tau t$ if and only if $t^\tau_{r-1}\in\lref{\tau_r}$,
			and on the other hand $k\in\rcross{\upslice\tau{r-1}}t$ if and
			only if $t^\tau_{k-1}\in\lref{\tau_k}$ for $k$ such that 
			$1\leq k\leq r-1$ .
		
		Thus, it will suffice to show that $r\in\rcross\tau t$ being
			equivalent to $t^\tau_{r-1}\in\lref{\tau_r}$ is equivalent to
			$\rcross{\upslice\tau{r-1}}t=\rcross\tau t\setminus\set r$.
		To that end, note that by definition
			$\rmeet{\upslice\tau{r-1}}t=\rmeet\tau t\setminus\set r$, and suppose
			that $m\in\rmeet{\upslice\tau{r-1}}t$ has right-separation after
			$k$, while $m\in\rmeet\tau t$ has right-separation after $k'$.
		Then our goal is to show that $m\not\in\rcross{\upslice\tau{r-1}}t
			+\rcross\tau t$, that is, $t^\tau_{m-1}\not\in
			\lref{\sigma_m\cdots\sigma_k}+\lref{\sigma_m\cdots\sigma_{k'}}$.

		If $k'=r$, then for all $j$ such that $m\leq j<r$ we have that
			there either exists $i\leq j$ such that
			$t^\tau_j\in W_{\rdeslev iJ}\udj W_{\overlap Jir}
				\subseteq W_{\rdeslev i{\upslice J{r-1}}}$,
			or there exists $i\leq j<n<r$ such that
			$t^\tau_j\in W_{\overlap Jin}$,
			so that in either case $t^{\upslice\tau{r-1}}$
			is right-traveling together after $j$ and $k=r-1$.
		Then $t^\tau_{m-1}\in\lref{\tau_m\cdots\tau_{r-1}}
				+\lref{\tau_m\cdots\tau_r}
				=(\tau_m\cdots\tau_{r-1})\lref{\tau_r}
					(\tau_m\cdots\tau_{r-1})^{-1}$
				if and only if $t^\tau_{r-1}
					=(\tau_m\cdots\tau_{r-1})^{-1}t^\tau_{m-1}
						(\tau_m\cdots\tau_{r-1})\in\lref{\tau_r}$.
		But $m\in\rmeet\tau t\setminus\set r$ having right-separation after
			$r$ implies that either $t^\tau_{r-1}\in W_{\rdeslev iJ}$ for
			$i\leq r-1$, so that $t^\tau_{r-1}\not\in
			W_{J_r}\supseteq\lref{\tau_r}$,
			or that $t^\tau_{r-1}\in W_{\overlap J*r}$ 
			which since $\lref{\tau_r}\cap W_{\overlap J*r}=\emptyset$ and
			$\lref{\tau_r}\subseteq W_{J_r}$ implies
			$t^\tau_{r-1}\not\in\lref{\tau_r}$. 
		Either way, $m\not\in\rcross{\upslice\tau{r-1}}t+\rcross\tau t$, as
			desired.

		If $k'<r$, we first claim that $k\geq k'$, since $t^\tau$
			right-traveling together after $k$ implies that either
			$t^\tau_k\in W_{\rdeslev iJ} \udj W_{\overlap Jir}
				\subseteq W_{\rdeslev i{\upslice J{r-1}}}$
			for some $i\leq k$, or $t^\tau_k\in\overlap Jin$ for some
			$i\leq k<n<r$, and for both of these we
			have that $t^{\upslice\tau{r-1}}$ right-travels together
			after $k$.
		If $k=k'$, the result is immediate.
		Finally, if $k>k'$, then
			$t^\tau_{m-1}\in\lref{\sigma_m\cdots\sigma_{k'}}
			+\lref{\sigma_m\cdots\sigma_k}$ if and only if
			$t^\tau_{k'}\in\lref{\sigma_{k'+1}\cdots\sigma_k}$.
		But in this case by definition $t^{\upslice\tau{r-1}}$
			right-travels together after $k'$ but $t^\tau$ does
			not, so we must have $t^\tau_{k'}\in\rdeslev i{\upslice J{r-1}}
			\setminus \rdeslev iJ$ for some $i\leq k'$, and in particular
			$t^\tau_{k'}$ commutes with every element of $J_j$ for $k'<j\leq k$,
			so that $t^\tau_{k'}\not\in\lref{\sigma_{k'+1}\cdots\sigma_k}$.
	\end{proof}
\end{prop}

\begin{prop}
\label{prop:length-meetings}
	Given a mask $\sigma\to_Jx$ on a factorization $J_1,\dots,J_r$, for
		each $t\in T$ we have that
		$\abs{\set{t^\sigma_r}\cap\rref{\sigma_1\cdots\sigma_r}}
			=\abs{\set t\cap\lref x}
			=\abs{\rcross\sigma t\setminus\rdef\sigma t}
			-\abs{\rcross\sigma t\cap\rdef\sigma t}$.
	In particular, $\ell(x)=\abs{\lref x}
			=\ds\sum_{t\in T}
				\abs{\rcross\sigma t\setminus\rdef\sigma t}
				-\abs{\rcross\sigma t\cap\rdef\sigma t}$.
	\begin{proof}
		Since $x=\sigma_1\cdots\sigma_r$, $\rcross\sigma t$, and
			$\rdef\sigma t$ depend only on the equivalence class of $\sigma$
			by \cref{lem:action} and
			\cref{prop:crossing-bounce-equivalence}, it follows
			from \cref{prop:equiv-induct} that it will suffice to
			show the claim for masks $\sigma$ such that for every $t\in T$
			and $k$ such that $1\leq k\leq r$,
			$k\in\rcross{\sigma}t$ if and only if
			$t^{\sigma}_{k-1}\in\lref{\sigma_k}$.
		We do so by induction on the length $r$ of the factorization, with
			the base case $r=0$ holding since for every $t\in T$,
			$t^\sigma_0=t\not\in\emptyset=\rref e$ and $\rcross\sigma t
				=\rdef\sigma t \subseteq\rmeet\sigma t=\emptyset$.

		For the inductive step, suppose the claim holds for factorizations
			of length $r-1$.
		Because
			$t^{\sigma}_r=(\sigma_r)^{-1}t^{\sigma}_{r-1}\sigma_r$
			while $\rref{\sigma_1\cdots\sigma_r}
				=(\sigma_r)^{-1}\rref{\sigma_1\cdots\sigma_{r-1}}
					\sigma_r+\rref{\sigma_r}$, it follows that
			$t^{\sigma}_r\in\rref{\sigma_1\cdots\sigma_r}$ if and only
			if $t^{\sigma}_{r-1}\in\rref{\sigma_1\cdots\sigma_{r-1}}
				+\lref{\sigma_r}$.
		In particular, this gives us
		\begin{multline*}
		\abs{\set{t^\sigma_r}\cap\rref{\sigma_1\cdots\sigma_r}}
			=\abs{\set{t^\sigma_{r-1}}\cap\rref{\sigma_1\cdots\sigma_{r-1}}}
			+\abs{\set{t^\sigma_{r-1}}\cap\lref{\sigma_r}}\\
			-2\abs{\set{t^\sigma_{r-1}}\cap\lref{\sigma_r}
				\cap\rref{\sigma_1\cdots\sigma_{r-1}}}
		\end{multline*}
		By our choice of $\sigma$, we have that $t^\sigma_{r-1}
			\in\lref{\sigma_r}$ if and only if $r\in\rcross\sigma t
			\subseteq\rmeet\sigma t$, so that $t^\sigma_{r-1}\in\lref{\sigma_r}
				\cap\rref{\sigma_1\cdots\sigma_{r-1}}$
			if and only if $r\in\rcross\sigma t\cap\rdef\sigma t$, and thus
			\begin{multline*}
				\abs{\set{t^\sigma_{r-1}}\cap\lref{\sigma_r}}
					-2\abs{\set{t^\sigma_{r-1}}\cap\lref{\sigma_r}
					\cap\rref{\sigma_1\cdots\sigma_{r-1}}}\\
				=\abs{\set r\cap\rcross\sigma t\setminus\rdef\sigma t}
					-\abs{\set r\cap \rcross\sigma t\cap\rdef\sigma t}
			\end{multline*}
		Furthermore, since $\rcross{\upslice{\sigma}{r-1}}t
				=\rcross{\sigma}t\setminus\set r$ by
				\cref{prop:equiv-induct} and
			$\rdef{\upslice{\sigma}{r-1}}t=\rdef{\sigma}t\setminus\set r$ by
			definition, the inductive hypothesis implies
			$\abs{\set{t^\sigma_{r-1}}\cap\rref{\sigma_1\cdots\sigma_{r-1}}}
				=\abs{\rcross\sigma t\setminus\rdef\sigma t\setminus\set r}
					-\abs{\rcross\sigma t\cap\rdef\sigma t\setminus\set r}$.
		Together, this gives us
			\begin{align*}
				\abs{\set{t^{\sigma}_r}
						\cap\rref{\sigma_1\cdots\sigma_r}}
				&=\abs{\rcross\sigma t\setminus\rdef\sigma t\setminus\set r}
					-\abs{\rcross\sigma t\cap\rdef\sigma t\setminus\set r}\\
				&\qquad+
				\abs{\set r\cap\rcross\sigma t\setminus\rdef\sigma t}
					-\abs{\set r\cap \rcross\sigma t\cap\rdef\sigma t}\\
				&=\abs*{\rcross\sigma t\setminus\rdef\sigma t}
					-\abs{\rcross\sigma t\cap\rdef\sigma t}\qedhere
		\end{align*}
	\end{proof}
\end{prop}
\begin{rmrk}
	\cref{prop:length-meetings} is equivalent to 
		$\abs{\rcross\sigma t}
			-\abs{\set{t^\sigma_r\cap\rref{\sigma_1\cdots\sigma_r}}}
			=2\cdot\abs{\rcross\sigma t\cap\rdef\sigma t}$ 
		and hence implies 
			$\ds c_R(\sigma)-\ell(x)
				=2\sum_{t\in T}\abs{\rcross\sigma t\cap\rdef\sigma t}
				\leq 2\sum_{t\in T}\abs{\rdef\sigma t}=2d_R(\sigma)$, which in
		the special case of a singleton factorization induced by a reduced
		word is the inequality of claim (i) of
		\cite[Proposition~2.3]{Deodhar}.
\end{rmrk}

Finally, claim (v) of \cite[Proposition~2.3]{Deodhar} is the special
	case of \cref{cor:length+defects<=meetings} below for a
	singleton factorization induced by a reduced word.
\begin{cor}
\label{cor:length+defects<=meetings}
	For any mask $\sigma \to_J x$ on a factorization $J_1,\dots, J_r$,
		we have $\ell(x)+d_R(\sigma)=\ds\sum_{t\in
		T}\abs{\rcross\sigma t+\rdef\sigma t}
			\leq\ds\sum_{t\in T}\abs{\rmeet\sigma t}=m_R(J)$
		with equality if and only if $\rdef\sigma t=\rbounce\sigma t$ for
		each $t\in T$, or equivalently, if and only if
		$\sigma$ is equivalent to the right-canonical mask $\omega_R(J)$, in
		which case $x=w(J)$.
	In particular, $\deg P^J_{w(J)}(q)=d_R(\omega_R(J))$ and 
		$\deg P^J_x(q)<\deg P^J_{w(J)}(q)+\ell(w(J))-\ell(x)$.
	\begin{proof}
		By \cref{prop:length-meetings}, the definition of
			$d_R(\sigma)$, and since $\rcross\sigma t\udj\rbounce\sigma t
				=\rmeet\sigma t\supseteq\rdef\sigma t$, we have 
			$\ell(x)+d_R(\sigma)
				=\ds\sum_{t\in T}\abs{\rcross\sigma t
					\setminus\rdef\sigma t}-\abs{\rcross\sigma t\cap\rdef\sigma t}
						+\abs{\rdef\sigma t}
				=\sum_{t\in T}\abs{\rcross\sigma t+\rdef\sigma t}
				\leq\sum_{t\in T}\abs{\rmeet\sigma t}=m_R(J)$, with equality
			if and only if $\rdef\sigma t=\rbounce\sigma t$ for all $t$.

		To see that $\rdef\sigma t=\rbounce\sigma t$ for every $t\in T$
			if and only if $\sigma$ is equivalent to $\omega_R(J)$, we first
			show that $\sigma=\omega_R(J)$ implies
			$\rdef{\sigma}t=\rbounce{\sigma}t$ for all $t$.
		To that end, fix $t\in T$ and suppose $m\in\rmeet\sigma t$ has
			right-separation $k$.
		Then by definition $t^\sigma_k\in W_{\rdeslev iJ}
			\subseteq W_{\rdeslev i{\upslice Jk}}$ or 
		$t^\sigma_k\in W_{\overlap Ji*}
			\subseteq W_{\rdeslev i{\upslice Jk}}$ for some $i$ such that
			$i\leq k$. 
		Since $\rdups Jk\subseteq \rdes{w(\upslice Jk)}$ by
			\cref{cor:recursive-defect-polys}, we have that 
			\cref{prop:finitary,prop:T_L} imply
			$t^\sigma_k\in T\cap W_{\rdups Jk}\subseteq \rref{w(\upslice Jk)}
				=\rref{\sigma_1\cdots\sigma_k}$, where the last equality
			holds because by definition $\omega_R(\upslice Jk)
				\to_{\upslice Jk}w(\upslice Jk)$. 
		Thus, we have $t^\sigma_{m-1}=(\sigma_m\cdots\sigma_k)t^\sigma_k
			(\sigma_m\cdots\sigma_k)^{-1}
				\in (\sigma_m\cdots\sigma_k)\rref{\sigma_1\cdots\sigma_k}
					(\sigma_m\cdots\sigma_k)^{-1}
				=\rref{\sigma_1\cdots\sigma_{m-1}}+\lref{\sigma_m\cdots\sigma_k}$,
			so that by definition $m\in\rdef\sigma t+\rcross\sigma t$, that is,
			$m\in\rdef\sigma t$ if and only if $m\in\rbounce\sigma t$ as
			desired.

		Next, we will show that any pair of masks $\sigma\to_Jx$ and
			$\tau\to_Jx'$ such that $\rdef\sigma t=\rbounce\sigma t$ and
			$\rdef\tau t=\rbounce\tau t$ for all $t$ are necessarily
			equivalent, and hence equivalent to $\omega_R(J)$.
		Because the equalities $\rdef\sigma t=\rbounce\sigma t$, and
			$\rdef\tau t=\rbounce\tau t$  depend only on the
			equivalence classes of $\sigma$ and $\tau$,
			\cref{prop:equiv-induct} implies we may assume 
			that $\sigma$ and $\tau$ additionally satisfy the
			properties that $k\in\rcross\sigma t$ if and only if
			$t^\sigma_{k-1}\in\lref{\sigma_k}$ and $k\in\rcross\tau t$ if and
			only if $t^\tau_{k-1}\in\lref{\tau_k}$. 

		Assuming these additional hypotheses, we will show that
			$\sigma=\tau$, or equivalently that
			$\sigma_1\cdots\sigma_k=\tau_1\cdots\tau_k$ for all $k$.
		We proceed using induction on $k$ with trivial base case $k=0$.
		For the inductive step, we assume that $\sigma_1\cdots\sigma_{k-1}
			=\tau_1\cdots\tau_{k-1}$ and wish to show that $\sigma_k=\tau_k$. 
		For every $t\in T$ we have by definition that 
			$t^\sigma_{k-1}=t^\tau_{k-1}$ and hence that
			$k\in\rmeet\sigma t$ if and only if $k\in\rmeet\tau t$, and
			that $k\in\rdef\sigma t=\rbounce\sigma t$ is equivalent to
			$k\in\rdef\tau t=\rbounce\tau t$.
		Hence, for all $t\in T$ we have $k\in\rcross\sigma t
			=\rmeet\sigma t\setminus\rbounce\sigma t$ if and only if
			$k\in\rcross\tau t=\rmeet\tau t\setminus\rbounce\tau t$.

		But our assumptions on $\sigma$ and $\tau$ imply
			$\set{t\in T:t^\sigma_{k-1}\in\lref{\sigma_k}}
				=\set{t\in T:k\in\rcross\sigma t}
				=\set{t\in T: k\in\rcross\tau t}
				=\set{t\in T:t^\tau_{k-1}\in\lref{\tau_k}}$
			so that $\lref{\sigma_k}=\lref{\tau_k}$ and hence 
			$\sigma_k=\tau_k$ by \cref{lem:injective-cocycle}.

		Finally, the supplemental claims follow from $\deg P^J_x(q)
			=\max\set{d_R(\sigma):\sigma\to_J x}\leq m_R(J)-\ell(x)
			=d_R(\omega_R(J))+\ell(w(J))-\ell(x)$ being a strict inequality
			if and only if no mask $\sigma\to_Jx$ satisfies 
			$\rdef\sigma t=\rbounce\sigma t$, i.e.\ if and only if no mask
			$\sigma\to_Jx$ is equivalent to $\omega_R(J)$, which amounts to
			$x\neq w(J)$.
	\end{proof}
\end{cor}

\begin{cor}
\label{thm:combi}
	Given a factorization with overlaps $J_1,\dots,J_r$ of $w$, for every
		mask $\sigma\to_Jx$ we have $\deg P_x(q)
				-\frac12\qy{\deg P_w(q)+\ell(w)-\ell(x)}
			=\max\set{\sigma\to_Jx:d_R(\sigma)-\frac12\qy{m_R(J)-\ell(x)}}$
			and $d_R(\sigma)-\frac12\qy{m_R(J)-\ell(x)}
			=\ds\frac12\sum_{t\in T}\abs{\rbounce\sigma t
			\cap\rdef\sigma t}-\abs{\rbounce\sigma t\setminus\rdef\sigma t}$.
	In particular, $J$ is tight if and only if for every mask
	$\sigma\to_Jx<w(J)$ the non-right defect right-bounces \emph{exceed}
	the right-defect right-bounces in the sense that 
		$\ds\sum_{t\in T}\abs{\rbounce\sigma t \cap\rdef\sigma t}
			<\sum_{t\in T}\abs{\rbounce\sigma t\setminus\rdef\sigma t}$.
	\begin{proof}
		$\ell(x)=\ds\sum_{t\in T}
			\abs{\rcross\sigma t\setminus\rdef\sigma t}
				-\abs{\rcross\sigma t\cap\rdef\sigma t}$ by 
				\cref{prop:length-meetings} and hence:
		\begin{align*}
			d&_R(\sigma)-\frac12\qy{m_R(J)-\ell(x)}\\
			&=\sum_{t\in T}\abs{\rdef\sigma t}
				-\frac12\qy[\big]{
					\abs{\rmeet\sigma t}
					-\abs{\rcross\sigma t\setminus\rdef\sigma t}
					+\abs{\rcross\sigma t\cap\rdef\sigma t}}\\
			&=\sum_{t\in T}\abs{\rdef\sigma t}
				-\frac12\qy[\big]{
					\abs{\rbounce\sigma t}+\abs{\rcross\sigma t}
					-\abs{\rcross\sigma t\setminus\rdef\sigma t}
					+\abs{\rcross\sigma t\cap\rdef\sigma t}}\\
			&=\sum_{t\in T}\abs{\rdef\sigma t}
				-\frac12\qy[\big]{
					\abs{\rbounce\sigma t}+2\abs{\rcross\sigma t\cap\rdef\sigma t}}\\
			&=\frac12\sum_{t\in T}2\abs{\rdef\sigma t}
				-2\abs{\rcross\sigma t\cap\rdef\sigma t}
				-\abs{\rbounce\sigma t}\\
			&=\frac12\sum_{t\in T}2\abs{\rbounce\sigma t\cap\rdef\sigma t}
				-\abs{\rbounce\sigma t}
			=\frac12\sum_{t\in T}\abs{\rbounce\sigma t\cap\rdef\sigma t}
				-\abs{\rbounce\sigma t\setminus\rdef\sigma t}\qedhere
		\end{align*}
	\end{proof}
\end{cor}

\begin{exmp}\cite[Lemma~2.1.2]{Fan--Green} 
\label{exmp:short-braid}
	Suppose we have a singleton factorization $J$ given by
		$J_j=\{s_{i_j}\}$ for a reduced word $s_{i_1}\cdots s_{i_r}=w$ for
		which we have a mask
		$\sigma\to_Jx$ given by
		$\sigma=(s_{i_1},\dots,s_{i_{k-1}},e,s_{i_{k+1}},\dots,s_{i_r})$
		such that $\ell(x)<\ell(w)-1$, i.e.\ such that $s_{i_1}\cdots
		s_{i_{k-1}}s_{i_{k+1}}\cdots s_{i_r}$ is not a reduced word.
	Then there exists $j$ such that $k<j$,
		$s_{i_j}\in\rref{\sigma_1\cdots\sigma_{j-1}}$ and $\upslice Jj$ is
		induced by a reduced word $s_{i_1}\cdots s_{i_j}$.

	From \cref{exmp:singleton-meet} the mask
		$\tau=(s_{i_1},\dots,s_{i_{k-1}},e,s_{i_{k+1}},\dots,s_{i_{j-1}},e)$
		on $\upslice Jj$ satisfies 
		$m\in\rmeet\tau t$ if and only if $t^\tau_m=s_{i_m}$,
		and that $m\in\rbounce\tau t$ if and only if
		$t^\tau_m=s_{i_m}$ and $\tau_m=e$. 
	Thus, $\ds\sum_{t\in T}\abs{\rbounce\tau t}=2$, and
		furthermore $t\in T$ such that $t^\sigma_{k-1}=s_{i_k}
			\not\in\rref{s_{i_1}\cdots s_{i_{k-1}}}
			=\rref{\tau_1\cdots\tau_{k-1}}$ has $k\in\rbounce\tau t\setminus
				\rdef\tau t$,
		and $t'\in T$ such that $(t')^\tau_j=s_{i_j}\in\rref{
			\tau_1\cdots\tau_{j-1}}$ has $j\in\rbounce\tau t\cap
				\rdef\tau t$.
	In particular, we have
		$\ds\sum_{t\in T}\abs{\rbounce\tau t\cap\rdef\tau t}=1
			=\sum_{t\in T}\abs{\rbounce\tau t\setminus\rdef\tau t}$, implying
			$\upslice Jj$ is not tight. 
		By \cref{cor:factors} it will follow that $J$ itself is not
			tight.
	In particular, factorizations induced by \emph{short braids}, i.e.\ by
		reduced words $s_is_js_i$ with $s_is_j\neq s_js_i$ are never tight.
\end{exmp}
\begin{rmrk}
	As we explain in \cref{rmrk:simply-laced}, in the case of simply-laced
		Coxeter systems, short braid avoidance is a special case of one of
		our main results (\cref{thm:2}) that tight factorizations
		have minimal contractions that are absolutely admissible bidescent.
\end{rmrk}

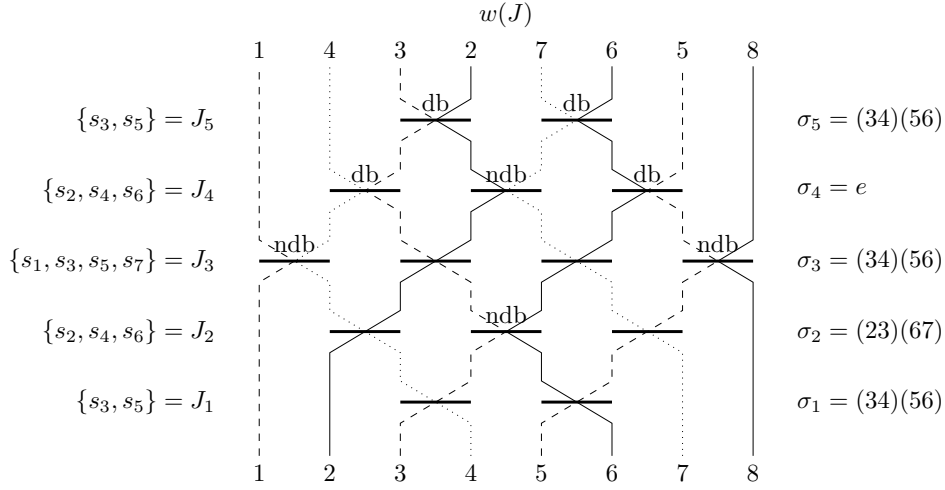
\begin{figure}[htbp]
	\centering
	\begin{tikzpicture}
	\node (1) at (0,0) {$1$};
	\node (2) at (1,0) {$2$};
	\node (3) at (2,0) {$3$};
	\node (4) at (3,0) {$4$};
	\node (5) at (4,0) {$5$};
	\node (6) at (5,0) {$6$};
	\node (7) at (6,0) {$7$};
	\node (8) at (7,0) {$8$};

	\node[anchor=east] (y1) at (-0.5, 1) {$\{s_3, s_5\} = J_1$};
	\node[anchor=east] (y2) at (-0.5, 2) {$\{s_2, s_4, s_6\} = J_2$};
	\node[anchor=east] (y3) at (-0.5, 3) {$\{s_1, s_3, s_5, s_7\} = J_3$};
	\node[anchor=east] (y4) at (-0.5, 4) {$\{s_2, s_4, s_6\} = J_4$};
	\node[anchor=east] (y5) at (-0.5, 5) {$\{s_3, s_5\} = J_5$};

	\node[anchor=west] (s1) at (7.5, 1) {$\sigma_1 = (34)(56)$};
	\node[anchor=west] (s2) at (7.5, 2) {$\sigma_2 = (23)(67)$};
	\node[anchor=west] (s3) at (7.5, 3) {$\sigma_3 = (34)(56)$};
	\node[anchor=west] (s4) at (7.5, 4) {$\sigma_4 = e$};
	\node[anchor=west] (s5) at (7.5, 5) {$\sigma_5 = (34)(56)$};
	
	\node[anchor=south] (b36) at (3.5, 2) {ndb};
	\node[anchor=south] (b14) at (0.5, 3) {ndb};
	\node[anchor=south] (b58) at (6.5, 3) {ndb};
	\node[anchor=south] (d34) at (1.5, 4) {db};
	\node[anchor=south] (d27) at (3.5, 4) {ndb};
	\node[anchor=south] (d56) at (5.5, 4) {db};
	\node[anchor=south] (d23) at (2.5, 5) {db};
	\node[anchor=south] (d67) at (4.5, 5) {db};

	\draw[very thick] (2, 1) -- (3, 1);
	\draw[very thick] (4, 1) -- (5, 1);

	\draw[very thick] (1, 2) -- (2, 2);
	\draw[very thick] (3, 2) -- (4, 2);
	\draw[very thick] (5, 2) -- (6, 2);

	\draw[very thick] (0, 3) -- (1, 3);
	\draw[very thick] (2, 3) -- (3, 3);
	\draw[very thick] (4, 3) -- (5, 3);
	\draw[very thick] (6, 3) -- (7, 3);

	\draw[very thick] (1, 4) -- (2, 4);
	\draw[very thick] (3, 4) -- (4, 4);
	\draw[very thick] (5, 4) -- (6, 4);

	\draw[very thick] (2, 5) -- (3, 5);
	\draw[very thick] (4, 5) -- (5, 5);

	\node (1e) at (0, 6) {$1$};
	\node (2e) at (3, 6) {$2$};
	\node (3e) at (2, 6) {$3$};
	\node (4e) at (1, 6) {$4$};
	\node (5e) at (6, 6) {$5$};
	\node (6e) at (5, 6) {$6$};
	\node (7e) at (4, 6) {$7$};
	\node (8e) at (7, 6) {$8$};
	\node (wJ) at (3.5, 6.5) {$w(J)$};

	\draw[dashed] (1)--(0, 2.7)--(0.5, 3)--(0, 3.3)--(1e);
	\draw (2)--(1, 1.7)--(1.5, 2)--(2, 2.3)--(2, 2.7)--(2.5, 3)--(3, 3.3)
		--(3, 3.7)--(3.5, 4)--(3, 4.3)--(3, 4.7)--(2.5, 5)--(3, 5.3)--(2e);
	\draw[dashed] (3)--(2, 0.7)--(2.5, 1)--(3, 1.3)--(3, 1.7)--(3.5, 2)
		--(3, 2.3)--(3, 2.7)--(2.5, 3)--(2, 3.3)--(2, 3.7)--(1.5, 4)
		--(2, 4.3)--(2, 4.7)--(2.5, 5)--(2, 5.3)--(3e);
	\draw[dotted] (4)--(3, 0.7)--(2.5, 1)--(2, 1.3)--(2, 1.7)--(1.5, 2)
		--(1, 2.3)--(1, 2.7)--(0.5, 3)--(1, 3.3)--(1, 3.7)--(1.5, 4)
		--(1, 4.3)--(4e);
	\draw[dashed] (5)--(4, 0.7)--(4.5, 1)--(5, 1.3)--(5, 1.7)--(5.5, 2)
		--(6, 2.3)--(6, 2.7)--(6.5, 3)--(6, 3.3)--(6, 3.7)--(5.5, 4)
		--(6, 4.3)--(5e);
	\draw (6)--(5, 0.7)--(4.5, 1)--(4, 1.3)--(4, 1.7)--(3.5, 2)--(4, 2.3)
		--(4, 2.7)-- (4.5, 3)--(5, 3.3)--(5, 3.7)--(5.5, 4)--(5, 4.3)
		--(5, 4.7)--(4.5, 5)--(5, 5.3)--(6e);
	\draw[dotted] (7)--(6, 1.7)--(5.5, 2)--(5, 2.3)--(5, 2.7)--(4.5, 3)
		--(4, 3.3)--(4, 3.7)--(3.5, 4)--(4, 4.3)--(4, 4.7)--(4.5, 5)
		--(4, 5.3)--(7e);
	\draw (8)--(7, 2.7)--(6.5, 3)--(7, 3.3)--(8e);
\end{tikzpicture}
	\caption{Hexagon factorization illustrating combinatorial criterion}
	\label{fig:hexagon}
\end{figure}

\begin{exmp}[{\cite[Proof of Proposition~2.2.1]{Fan--Green}]}]
	For Coxeter systems of type $A$, the mask $\sigma$ on the
		\emph{hexagon} in \cref{fig:hexagon} has four right-defect
		right-bounces: $4\in\rbounce\sigma{(3,4)}\cap\rdef\sigma{(3,4)}$,
		$4\in\rbounce\sigma{(5,6)}\cap\rdef\sigma{(5,6)}$, 
		$5\in\rbounce\sigma{(2,3)}\cap\rdef\sigma{(2,3)}$
		and $5\in\rbounce\sigma{(6,7)}\cap\rdef\sigma{(6,7)}$, and
		four non-defect right-bounces:
		$2\in\rbounce\sigma{(3,6)}\setminus\rdef\sigma{(3,6)}$,
		$3\in\rbounce\sigma{(1,4)}\setminus\rdef\sigma{(1,4)}$, $3\in\rbounce\sigma{(5,8)}\setminus\rdef\sigma{(5,8)}$,
		$4\in\rbounce\sigma{(2,7)}\setminus\rdef\sigma{(2,7)}$. 

	Thus, $\ds\sum_{t\in T}\abs{\rbounce\tau t\cap\rdef\tau t}=4
			=\sum_{t\in T}\abs{\rbounce\tau t\setminus\rdef\tau t}$ so that
			$J$ is unbounded.
	Consequently, hexagon factorizations are not tight (Fan and Green's
		proof of Proposition~2.2.1 considers the case of a hexagon with an
		added $\set{s_4}$ on top and bottom, which fails to be tight by
		\cref{cor:quotient} because it contains a hexagon as a
		slice.)
\end{exmp}

\begin{rmrk}
\label{rmrk:defect-graph}
	Billey and Warrington's defect graph criterion
		\cite[Sections 5,6]{Billey--Warrington} arises from the
		observation that if we associate to each pair $(j,t)$ with
		$j\in\rbounce{\sigma}{t}\cap\rdef{\sigma}{t}$ a set
		$b(t,j)\subseteq\ds\bigsqcup_{t\in T}\rbounce{\sigma}{t}$
		of ``critical bounces'' such that
		$\abs{b(t,j)}=3$, and $b(t,j)\cap b(t',j')\cap
		b(t'',j'')\neq\emptyset$ only if
		$\abs{\set{(t,j),(t',j'),(t'',j'')}}\leq2$, then applying the
		inclusion-exclusion principle to the quantity
		$\ds\abs{\bigcup_{t\in
			T}\bigcup_{j\in\rbounce{\sigma}{t}}b(t,j)}$
		yields the difference $3\abs{V}-\abs{E}$ where $(V,E)$ is the graph
		whose vertices are
		$V=\set{(t,j):t\in T,j\in\rbounce{\sigma}{t}\cap\rdef{\sigma}{t}}$
		and whose number of edges between distinct vertices is
		$\abs{b(t,j)\cap b(t',j')}$.
	In particular, we would have that $3\abs{V}-\abs{E}
		\leq\ds\sum_{t\in T}\abs{\rbounce{\sigma}{t}}$ and then 
		$\abs{V}=\ds\sum_{t\in
		T}\abs{\rbounce{\sigma}{t}\cap\rdef{\sigma}{t}}$
		implies that $2\abs{V}-\abs{E}\leq\ds\sum_{t\in
		T}\abs{\rbounce{\sigma}{t}\setminus\rdef{\sigma}{t}}$.
	Thus, the factorization is tight if $\abs{V}<
		2\abs{V}-\abs{E}$, i.e.\ if $\abs{E}<\abs{V}$, so in particular if
		the graph $(V,E)$ is a forest.

	Their defect graph argument is that for factorizations arising from
		reduced words associated to $321$-avoiding permutations, the
		``critical bounces'' can be chosen in such a way that the
		corresponding defect graph failing to be a forest guarantees the
		factorization contains a hexagon and is in particular not tight. 
	Unfortunately, we have not found a way to make this argument work for
		the monotone factorizations we analyze in
		\crefrange{sec:4231-45312}{sec:height2}, and have to
		rely on a different method to prove \cref{thm:4}.
\end{rmrk}

\section{Families of defect polynomials and the Schur composition laws}

\label{sec:hecke}
We review the definitions of the Hecke algebra and its natural
	involution and determine that a uniform degree shift of $-m_R(J)/2$ is
	necessary to realize a family of defect polynomials of a factorization
	$J$ as the coefficients of an \emph{involution-invariant} element of
	the Hecke algebra.
We then prove the shift is sufficient by showing that the recursive
	formula for computing families of defect polynomials corresponds to
	repeatedly combining Kazhdan--Lusztig basis elements 
	$C'_{\mxml J}$ for finitary $J\subseteq S$ using the composition laws
	of Williamson's \emph{Schur algebroid} \cite{Williamson}.
Along the way we give a more precise description of the composition laws
	and their interaction with involution-invariance.
\begin{defn}
	Given a Coxeter system $(W,S)$, the \define{Hecke algebra} $\mathcal
		H$ (cf.\ \cite[Section~1]{Kazhdan--Lusztig}) is the unital
		associative $\ZZ[q^{\frac12},q^{-\frac12}]$-algebra with basis
		$\set{T_w:w\in W}$ and bilinear multiplication satisfying
	\begin{enumerate}
		\item $T_xT_s=T_{xs}$ if $s\not\in\rdes{x}$;
		\item $(T_s)^2=(q-1)T_s+qT_e$ (where $e$ is the identity of $W$).
	\end{enumerate}
	Equivalently (cf.\ \cite{Lusztig83}), the Hecke algebra is the unital
		associative $\ZZ[q^{\frac12},q^{-\frac12}]$-algebra with basis
		$\set{H_w:w\in W}$ and bilinear multiplication satisfying
	\begin{enumerate}
		\item $H_xH_s=H_{xs}$ if $s\not\in\rdes{x}$;
		\item $(H_s)^2=(q^{\frac12}-q^{-\frac12})H_s+H_e$;
	\end{enumerate}
	where the equivalence of the two definitions is obtained from the
		relations $H_x=q^{-\ell(x)/2}T_x$ and $T_x=q^{\ell(x)/2}H_x$.
\end{defn}

\begin{rmrk}
	Both \cite{Elias--Williamson} and \cite{Williamson} use the latter
		presentation of the Hecke algebra as a $\ZZ[\nu,\nu^{-1}]$-algebra
		where $\nu=q^{-\frac12}$. 
	We use both presentations because families of defect polynomials are
		more naturally interpreted in terms of the basis $\set{T_x}$, while
		algorithms for computing Kazhdan--Lusztig basis elements and (our
		generalization of) Billey and Warrington's combinatorial criterion
		are more naturally expressed in terms of the basis $\set{H_x}$.
\end{rmrk}

\begin{defn}[{\cite[2.0a]{Kazhdan--Lusztig}}]
	The Hecke algebra admits an \define{involution
		$\inv\cdot\colon\mathcal H\to\mathcal H$} given by
		$q^{\frac12}\mapsto q^{-\frac12}$ and $T_s\mapsto
		(T_s)^{-1}=(-q)^{-1}((q-1)T_e-T_s)$.
	We say that an element $h\in\mathcal H$ of the Hecke algebra is
		\define{invariant} if $\inv h=h$.
	The \define{$R$-polynomials}
		$R_{x,y}$ are the unique elements of $\ZZ[q^{\frac12},q^{-\frac12}]$
		such that $\inv{T_y}=\qy*{T_{y^{-1}}}^{-1}
			=\ds(-q)^{-\ell(y)}\sum_x (-1)^{\ell(x)}R_{x,y}T_x$.
\end{defn}

\begin{exmp}
	$(T_s)^{-1}=(-q)^{-1}((q-1)T_e-T_s)$ asserts that $R_{e,s}=q-1$ and
		$R_{s,s}=1$.
	Furthermore, note that
		$\set{p(q^{\frac12})\in\ZZ[q^{\frac12},q^{-\frac12}]:\inv p=p}
			=\ZZ[q^{-\frac12}+q^{\frac12}]$ since by strong induction on
		degree, if $\inv p=p$, then there exists $a_n\in\ZZ$ so that
		$a_nq^{\pm\frac n2}$ are the highest and lowest degree terms of $p$,
		hence that $\inv{p-a_n(q^{-\frac12}+q^{\frac12})^n}
			=p-a_n(q^{-\frac12}+q^{\frac12})^n$ is a polynomial whose highest
		and lowest degrees are respectively lower and higher than those of
		$\inv p=p$.
\end{exmp}

\begin{prop}[{\cite[Equations~2.0b-c, 2.2a, Lemmas~2.1(i), 2.6(vi)]{Kazhdan--Lusztig}}]
\label{prop:R-polynomials}
	The $R$-polynomials satisfy:
	\begin{enumerate}
		\item $R_{x,ys}=(q-1)R_{x,y}+qR_{xs,y}$ and $R_{xs,ys}=R_{x,y}$ 
			whenever $s\not\in\rdes x\cup\rdes y$, i.e.\ whenever
			$x,y\in\rcos{\set s}W$.
		\item\label{item:R-constant}
			$R_{x,y}\in\ZZ[q]$, $R_{x,y}=0$ unless $x\leq y$, in which
			case $R_{x,y}(0)=(-1)^{\ell(y)-\ell(x)}$ and $\deg
			R_{x,y}(q)=\ell(y)-\ell(x)$, so that $R_{x,x}=1$;
		\item $\inv{R_{x,y}}=(-q)^{\ell(x)-\ell(y)}R_{x,y}$;
		\item\label{item:invariance}
			$\ds\sum_{x\in W}P_xT_x=\sum_{x\in W}p_xH_x$ is
			involution-invariant if and only if $(P_x)_{x\in
			W}\in\ZZ[q^{\frac12},q^{-\frac12}]$ satisfies
			$\inv{P_x}=\ds\sum_y q^{\ell(x)}R_{x,y}P_y$ for all $x$,
			or equivalently, if and only if 
			$(p_x)_{x\in W}\in\ZZ[q^{\frac12},q^{-\frac12}]$ satisfies
			$\inv{p_x}=\ds\sum_yq^{(\ell(x)-\ell(y))/2}R_{x,y}p_y$
			for all~$x$.
		\item\label{item:finitary-invariant}
			If $J\subseteq S$ is finitary then
			$\ds\sum_{x\leq y\leq \mxml J}R_{x,y}(q)
				=q^{\ell(\mxml J)-\ell(x)}$ for all $x\in W_J$, or equivalently,
				then $q^{-\ell(\mxml J)/2}\ds\sum_{x\in W_J}T_x$ is
				involution-invariant.
	\end{enumerate}
	\begin{proof}
		(1) Follows from expanding $\inv{T_{ys}}=\inv{T_y}\cdot\inv{T_s}$.

		(2) and (3)  follow by induction on $\ell(y)$ using (1).

		(4) follows from the definition of $R$-polynomials and (3):
		\begin{align*}
				\inv{\sum_y P_y T_y}
				&=\sum_y\inv{P_y}(-q)^{-\ell(y)}\sum_x(-1)^{\ell(x)}R_{x,y}T_x\\
				&=\sum_x\qy*{\sum_y\inv{P_y}(-q)^{-\ell(y)}(-1)^{\ell(x)}R_{x,y}}T_x\\
				\inv{\sum_y\inv{P_y}(-q)^{-\ell(y)}(-1)^{\ell(x)}R_{x,y}}
				&=\sum_y P_y(-q)^{\ell(y)}(-1)^{\ell(x)}\inv{R_{x,y}}\\
				&=\sum_y P_y(-q)^{\ell(y)}(-1)^{\ell(x)}(-q)^{\ell(x)-\ell(y)}R_{x,y}\\
				&=\sum_y q^{\ell(x)}R_{x,y}P_y
		\end{align*}
		Thus, it is necessary and sufficient that
			$\inv{P_x}=\ds\sum_yq^{\ell(x)}R_{x,y}P_y$.
		Since $P_x=q^{-\ell(x)/2}p_x$ and $P_y=q^{-\ell(y)/2}p_y$, we have
			that the condition is equivalent to
			$\inv{q^{-\ell(x)/2}p_x}=\ds\sum_y q^{\ell(x)}R_{x,y}P_y
				=\sum_y q^{\ell(x)-\ell(y)/2}R_{x,y}p_y$ and hence to
			$\inv{p_x}=\ds\sum_y q^{(\ell(x)-\ell(y))/2}R_{x,y}p_y$.

		The main claim of (5) follows from (1) by induction on
			$\ell(\mxml J)-\ell(x)$ and the equivalent formulation then follows
			from (4).
	\end{proof}
\end{prop}

\begin{rmrk}
	Deodhar \cite{Deodhar1985} showed that property (1) can be
		interpreted combinatorially.
	Given a factorization $J_1,\dots,J_r$ let us write
		$\sigma\twoheadrightarrow_J x$ for a mask $\sigma\to_Jx$ such that
		$\rref{\sigma_1\cdots\sigma_{k-1}}\subseteq\lref{\sigma_k}$ for all
		$k$.
	In the special case of a factorization whose overlap subsets are
		trivial, we have $\sigma\twoheadrightarrow x$ if and only if
		$\rbounce\sigma t\subseteq\rdef\sigma t$ for all $t\in T$, and in
		the special case of a factorization $J_1,\dots,J_r$ induced by a
		reduced word, \cref{exmp:singleton-meet} implies
		$\sigma\twoheadrightarrow x$ if and only if
		$s_{i_k}\in\rref{\sigma_1\cdots\sigma_{k-1}}$ implies
		$\sigma_k=s_{i_k}$.
		
	Deodhar observed \cite[Lemma~5.1]{Deodhar1985} that in the special
		case of a singleton factorization induced by a reduced word, every
		mask $\sigma\twoheadrightarrow_Jx$ is uniquely determined by its
		restriction $\upslice\sigma{r-1}=\sigma'\twoheadrightarrow_{\upslice
		J{r-1}}x'$, in which case the resulting bijection has by definition
		the property that 
		$d_R(\sigma)=1+d_R(\sigma')$ if $\sigma_r\neq s_{i_r}$, that
		$b_R(\sigma)=1+b_R(\sigma')$ if $\sigma_r=s_{i_r}\not\in\rdes x$,
		and that $d_R(\sigma)=d_R(\sigma')$ and $b_R(\sigma)=b_R(\sigma')$
		if $\sigma_r=s_{i_r}\in\rdes x$.
	He concluded the $R$-polynomials have a combinatorial interpretation
		given by \cite[Theorem~1.3, Lemma~5.1]{Deodhar1985}:
		$\ds R_{x,y}=\sum_{\sigma\twoheadrightarrow_J x}
			q^{d_R(\sigma)}\qy*{q-1}^{b_R(\sigma)}$.
\end{rmrk}

\cref{prop:R-polynomials}(\ref{item:invariance}) shows that
	for any family $(P_x)_{x\in W}\in\ZZ[q^{\frac12},q^{-\frac12}]$ with
	a unique maximal $w$ such that $P_w\neq0$, e.g.\ for any family
	$(P^J_x)_{x\in W}$ of defect polynomials associated to a factorization
	of $w$, the Hecke algebra element $\ds\sum_xP_xT_x$ is
	involution-invariant only if $\inv{P_w}=q^{\ell(w)}P_w$, or
	equivalently $\inv{q^{\ell(w)/2}P_w} =q^{\ell(w)/2}P_w$.
Consequently, in order for our association of Hecke algebra elements to
	families as above to result in \emph{involution-invariant} Hecke algebra
	elements, it is necessary to introduce a uniform degree shift as in
	\cref{defn:normalized} below, transforming the invariance
	requirement of
	\cref{prop:R-polynomials}(\ref{item:invariance}) into that
	of \cref{exmp:normalized-invariance-equation}.

\begin{defn}
\label{defn:normalized}
	We denote by $(P_x)_{x\leq w}$ a family $(P_x)_{x\in W}$ such that
	 $0\neq P_x$ only if $x\leq w$ in the Bruhat order.
	We say that a family 
		$(P_x)_{x\leq w}\in\ZZ[q^{\frac12},q^{-\frac12}]$ is a
		\define{positive family} if $P_x\in\NN[q]$ with
		$P_x(0)=1$ whenever $x\leq w$.
	Given a family $(P_x)_{x\leq w}\in\ZZ[q^{\frac12},q^{-\frac12}]$ with
		$0\neq P_w$, we define its \define{normalized Hecke algebra element
		$\underline h^P$} by $\underline h^P\ds
			=q^{-(\deg P_w(q)-\deg P_w(q^{-1})+\ell(w))/2}\sum_xP_xT_x$,
		i.e.\ by $\underline h^P=\ds\sum_x p_xH_x$ where $p_x
			=q^{-(\deg P_w(q)-\deg P_w(q^{-1})+\ell(w)-\ell(x))/2}P_x$.
\end{defn}

\begin{rmrk}
\label{rmrk:bounded-has-pw-nonzero}
	We use the convention that $P=0\in\ZZ[q^{\frac12},q^{-\frac12}]$ is
		unique such that $\deg P(q)=-\infty$ and unique such that $\deg
		P(q^{-1})=-\infty$, in which case $\deg P(q)+\deg
		P(q^{-1})=-\infty$.
	This makes it so that $P(q)\in\ZZ[q^{\frac12},q^{-\frac12}]$ satisfies
		$P(q)\in\ZZ[q^{\frac12}]$ with $0\neq P(0)$ if and only if
		$0=\deg P(q^{-1})$, i.e.\ so that a family $(P_x)_{x\leq
		w}\in\ZZ[q^{\frac12},q^{-\frac12}]$ has a normalized Hecke algebra
		element given by $\underline h^P=\ds q^{-(\deg
		P_w(q)+\ell(w))/2}\sum_xP_xT_x$ if and only if
		$P_w\in\ZZ[q^{\frac12}]$ with $0\neq P_w(0)$.
\end{rmrk}

\begin{exmp}
\label{exmp:normalized-invariance-equation}
	Given a family $(P_x)_{x\leq w}\in\ZZ[q^{\frac12},q^{-\frac12}]$,
		with $P_w\neq 0$,
		\cref{prop:R-polynomials}(\ref{item:invariance})
		implies that $\underline h^P$ is involution-invariant if
		and only if it satisfies the \emph{invariance equation}
		$\ds
			q^{\deg P_w(q)-\deg P_w(q^{-1})+\ell(w)-\ell(x)}\inv{P_x}
			=\sum_yR_{x,y}P_y$ for every $x\leq w$.
\end{exmp}

\begin{lem}
\label{lem:normalized-factorization}
	Given a factorization with overlaps $J_1,\dots,J_r$, the family of
		defect polynomials $(P^J_x)_{x\in W}$ satisfies 
		$(P^J_x)_{x\leq w(J)}$ and is positive.
	Furthermore, its normalized Hecke algebra element is given by
		 $\underline h^{P^J}
			=q^{-m_R(J)/2}\ds\sum_x
				\sum_{[\sigma]\to_Jx}q^{d_R(\sigma)}T_x=\sum_xp^J_xH_x$,
		where $m_R(J)=\ds\sum_{k=1}^r\abs{T\cap W_{J_k}
				\setminus W_{\overlap J*k}}$ 
		(as in \cref{defn:defects}) and
		$p^J_x=\ds\sum_{[\sigma]\to_Jx}
			q^{d_R(\sigma)-\frac12\qy{m_R(J)-\ell(x)}}$ are generating
			functions with
		$d_R(\sigma)-\frac12\qy{m_R(J)-\ell(x)}
			=\ds\frac12\sum_{t\in T}\abs{\rbounce\sigma t
			\cap\rdef\sigma t}-\abs{\rbounce\sigma t\setminus\rdef\sigma t}$
			($\rbounce\sigma t$ and $\rcross\sigma t$ as in
			\cref{defn:meetings}).
	\begin{proof}
		$(P^J_x)_{x\leq w(J)}$ is positive by
			\cref{cor:recursive-defect-polys}, 
			satisfies $\deg P^J_{w(J)}(q^{-1})=0$ by
			\cref{rmrk:bounded-has-pw-nonzero} and
			$\deg P^J_{w(J)}(q)+\ell(w{(J)})=d_R(\omega_R(J))+\ell(w(J))
				=m_R(J)$ by \cref{cor:length+defects<=meetings}, and
			the definition of the normalized Hecke algebra element
			yields $\underline h^{P^J}=q^{-m_R(J)/2}
					\ds\sum_x\sum_{[\sigma]\to_Jx}q^{d_R(\sigma)}T_x
				=\sum_xp^J_xH_x$ 
			where $p^J_x=\ds\sum_{[\sigma]\to_Jx}q^{d_R(\sigma)
				+\frac12(m_R(J)-\ell(x))}$.
		The combinatorial interpretation then follows from \cref{thm:combi}.
\end{proof}	

\end{lem}
\begin{exmp}
\label{exmp:r=1}
	In the case of a factorization of length $1$, i.e.\ of a finitary
		subset $J\subseteq S$, we have that $\underline h^{P^J}
			=q^{-\ell(\mxml J)/2}\ds\sum_{x\in W_J}T_x$ because 
		$P^J_x=\abs{\set x\cap W_J}$ by \cref{exmp:r=1f}, and that
		$\underline h^{P^J}$ is involution-invariant by
		\cref{prop:R-polynomials}(\ref{item:finitary-invariant}).
	In fact, $\underline h^{P^J}=C'_{\mxml J}$ is by definition the
		Kazhdan--Lusztig basis element indexed by $\mxml J$ (see
		\cref{exmp:finitary-C'w}).
\end{exmp}

The remainder of this section is dedicated to showing that, given a
	factorization with overlaps $J_1,\dots,J_r$, the normalized
	Hecke algebra element $\underline h^{P^J}$ of the family of defect
	polynomials $(P^J_x)_{x\in W}$ is involution-invariant, and more
	precisely, an involution-invariant quotient of the monomial of
	Kazhdan--Lusztig basis elements $C'_{\mxml{J_1}}\cdots C'_{\mxml{J_r}}$.
We do so by induction on the length of the factorization with overlaps
	$J_1,\dots,J_r$.

The base case is simply \cref{exmp:r=1}.
For the inductive step, we will show that the recursive formula of
	\cref{cor:recursive-defect-polys} for computing families
	of defect polynomials corresponds to combining the
	involution-invariant normalized Hecke algebra elements $\underline
	h^{P^{J_k}}=C'_{\mxml{J_k}}$ using the composition laws of Williamson's
	\emph{Schur algebroid} \cite{Williamson}, which we review below.

\begin{defn}
	Given a finite subset $A\subseteq W$, its \define{Poincar\'e
		polynomial $A(q)$} is given by $A(q)=\ds
			\sum_{x\in A}q^{\ell(x)}$.
\end{defn}

\begin{lem}
\label{lem:invariant_Poincare}
	Given a finitary $J\subseteq S$, $q^{-\ell(\mxml J)/2}W_J(q)$ is
		involution-invariant.
	\begin{proof}
			$\ds\inv{q^{-\ell(\mxml J)/2}W_J(q)}
				=q^{\ell(\mxml J)/2}\sum_{x\in W_J}\!q^{-\ell(x)}
				=q^{-\ell(\mxml J)/2}\sum_{x\in W_J}\!q^{\ell(\mxml J)-\ell(x)}
				\break=q^{-\ell(\mxml J)/2}\sum_{x\in W_J}q^{\ell(x^{-1}\mxml J)}
				=q^{-\ell(\mxml J)/2}\sum_{v\in W_J}\!q^{\ell(v)}
				=q^{-\ell(\mxml J)/2}W_J(q)$
	\end{proof}	
\end{lem}

\begin{defn}
\label{defn:schur}
	Given a finitary $J_*\subseteq S$, we define its associated
		\define{Schur composition law} to be the partially-defined binary
		operation on the Hecke algebra given by
		$h_1,h_2\mapsto h_2\circ_{J_*}h_1
			=\qy{q^{-\ell(\mxml{J_*})/2}W_{J_*}(q)}^{-1}h_1\cdot h_2$ and
		defined when $h_1=\ds\sum_xP_xT_x$ satisfies 
			$J_*\subseteq \ldes P$ and $h_2=\ds\sum_xQ_xT_x$
		satisfies $J_*\subseteq\rdes Q$.
		
\end{defn}

\begin{exmp}
	The Schur composition law $\circ_\emptyset$ associated to $\emptyset$
		is precisely the multiplication law of the Hecke algebra. 
\end{exmp}

\begin{rmrk}
	The composition laws of Williamson's Schur algebroid
		\cite[Definition~2.3.1]{Williamson} are given by the same formula
		but require that $h_1$ and $h_2$ be contained respectively in the
		principal left and right ideals of the Hecke
		algebra generated by $C'_{\mxml{J_*}}
			=q^{-\ell(\mxml{J_*})/2}\ds\sum_{x\in W_{J_*}}T_x$. 
	\cref{lem:bijections} and \cref{prop:hecke_mult} below
		imply that $\ds\sum_{x\in W}P_xT_x$ is in the principal
		right ideal of the Hecke algebra generated by
		$\ds\sum_{y\in W_J}T_y$ if and only if $J\subseteq\rdes
		P$, and hence that our Schur composition laws do agree with those of
		Williamson's Schur algebroid.
\end{rmrk}

\begin{lem}
\label{lem:bijections}
	Given subsets $J,J'\subseteq S$ we have that
		$W_J(q)=W_J^{J'}(q)\cdot W_{J\cap J'}(q)$, and when $J'$ commutes
		with $J\setminus J'$, we also have that 
		$W_J^{J'}(q)=W_{J\setminus J'}(q)$ and 
		$\ds\ds\sum_{x\in W_J}q^{\ell(x)}T_x
			=\qy*{\sum_{y\in W_{J\setminus J'}}q^{\ell(y)}T_y}
				\cdot\qy*{\sum_{z\in W_{J'}}q^{\ell(z)}T_z}$.
	Moreover, $\ds\sum_{x\in
		W}P_xT_x=\qy*{\sum_{x\in\rcos JW}P_xT_x}\cdot\qy*{\sum_{y\in
		W_J}T_y}$ when $J\subseteq\rdes P$.
	\begin{proof}
		Since $(W_J,J)$ is a Coxeter subsystem, we have a bijection
			$W_J^{J'}\times W_{J\cap J'}\cong W_J$ sending
			$(x,y)\mapsto xy$ such that $\ell(xy)=\ell(x)+\ell(y)$ and
			$T_xT_y=T_{xy}$, and when $J'\subseteq S$ commutes with
			$J\setminus J'$ also such that $W_J^{J'}=W_{J\setminus J'}$
			(see \cref{lem:coset-properties}(\ref{lem:coset-unique})).
		The equalities $W_J(q)=W_J^{J'}(q)\cdot W_{J\cap J'}(q)$ and 
		$\ds\sum_{x\in W_J}q^{\ell(x)}T_x
			=\qy*{\sum_{y\in W_{J\setminus J'}}q^{\ell(y)}T_y}
				\cdot\qy*{\sum_{z\in W_{J'}}q^{\ell(z)}T_z}$ are
			simply generating function versions of these bijections.
		Finally, if $J\subseteq\rdes P$, then
			$\ds\sum_xP_xT_x
				=\sum_{x\in W^J}\sum_{y\in W_J}P_{xy}T_{xy}
				=\sum_{x\in W^J}\sum_{y\in W_J}P_xT_xT_y
				=\sum_{x\in W^J}P_xT_x\sum_{y\in W_J}T_y$ is in the principal
				right ideal generated by $\ds\sum_{y\in W_J}T_y$.
	\end{proof}
\end{lem}

The explicit formulas of \cref{prop:hecke_mult} below are
	up to a uniform degree shift those of
	\cite[Proposition~2.3.3]{Williamson} that describe the action of
	$C'_{\mxml J}$ in the Schur algebroid in terms of an extension of the
	$\set{H_x}$-basis of the Hecke algebra to a basis for the Schur
	algebroid.
Because we need more precise results regarding involution-invariance and
	degree of the leading terms, we have chosen to prove them directly
	instead of introducing Williamson's basis and calculating the degree
	shift.

Recall from \cref{defn:right-parabolic-product} that for a family
	$(P_x)_{x\in W}\in\ZZ[q^{\frac12},q^{-\frac12}]$ and subsets
	$J_*,J\subseteq S$, we have defined the family
	$\rparp P{J_*}J$ by 
		$\rparp P{J_*}J_x
			=\ds\sum_{y\in\rcos[J]{J_*}W}q^{\ell(y)}P_{\rcos Jxy}$.
\begin{prop}
\label{prop:hecke_mult}
	Given a finitary subset $J\subseteq S$, we have that
		$\ds\qy*{\sum_xP_xT_x}
				\cdot\qy*{\sum_{z\in W_J}T_z}
			=W_{J_*}(q)\sum_x\rparp P{J_*}J_x\cdot T_x
			=\sum_x\rparp P\emptyset J\cdot T_x$
		for any $J_*\subseteq\rdes P\cap J$
	and $\ds\qy*{\sum_xp_xH_x}
			\cdot\qy*{\sum_{z\in W_J}T_z}
		=\ds W_{J_*}(q)\sum_x\sum_{y\in W_J^{J_*}}
			q^{(\ell(y)+\ell(x_J))/2}p_{x^Jy}H_x$.
	\begin{proof}
		$s\in J$ implies $\ds T_s\cdot\qy*{\sum_{z\in W_J}T_z}
			=T_s\cdot\qy*{\sum_{z\in\lcos[J]{\set s}W}(T_e+T_s)T_z}
			=\sum_{z\in\lcos[J]{\set s}W}(T_s+T_s^2)T_z
			=\sum_{z\in\lcos[J]{\set s}W}(T_s+(q-1)T_s+qT_e)T_z
			=\sum_{z\in\lcos[J]{\set s}W}(qT_s+qT_e)T_z
			=q\sum_{z\in W_J}T_z$.
		Hence $\ds T_y\cdot\qy*{\sum_{z\in
			W_J}T_z}=q^{\ell(y)}\sum_{z\in W_J}$ for any $y\in W_J$.
		Consequently, $\ds\qy*{\sum_{x\in W}P_xT_x}
				\cdot\qy*{\sum_{z\in W_J}T_z}
			=\sum_{x\in\rcos JW}\sum_{y\in W_J}P_{xy}T_xT_y\qy*{\sum_{z\in W_J}T_z}
			=\sum_{x\in\rcos JW}\sum_{y\in W_J}P_{xy}T_x q^{\ell(y)}\qy*{\sum_{z\in W_J}T_z}
			\break=\sum_{x\in W}\sum_{y\in W_J}P_{\rcos Jxy}T_x q^{\ell(y)}$.
		Since $J_*\subseteq\rdes P$ we have
			$P_{\rcos Jxy}=P_{\rcos Jx\rcos{J_*}y}$, and since $J_*\subseteq J$
			we have $W_J\cong \rcos[J]{J_*}W\times W_{J_*}$, and thus this is
			also
			$\ds\sum_{x\in W}\sum_{y\in\rcos[J]{J_*} W}
				P_{\rcos Jxy}T_x q^{\ell(y)}\sum_{z\in W_{J_*}}q^{\ell(z)}
			\break=W_{J_*}(q)\sum_{x\in W}\sum_{y\in\rcos[J]{J_*}W}q^{\ell(y)}
				P_{\rcos Jxy}T_x
			=W_{J_*}(q)\sum_x\rparp P{J_*}J_x\cdot T_x$, as desired.

		Since $H_x=q^{-\ell(x)/2}T_x$, we have $\ds
				\sum_x p_xH_x=\sum_x P_xT_x$ for $P_x=q^{-\ell(x)/2}p_x$ and
			$\ds W_{J_*}(q)\sum_{x\in\rcos JW}\sum_{y\in \rcos[J]{J_*}W}
				q^{\ell(y)}P_{\rcos Jxy}T_x
			=W_{J_*}(q)\sum_{x\in\rcos JW}\sum_{y\in \rcos[J]{J_*}W}
				q^{\ell(y)}q^{-\ell(\rcos Jxy)/2}p_{\rcos Jxy}q^{\ell(x)/2}H_x
			\break=W_{J_*}(q)\sum_{x\in\rcos JW}\sum_{y\in \rcos[J]{J_*}W}
				q^{(\ell(x_J)+\ell(y))/2}p_{\rcos Jxy}H_x$.
	\end{proof}
\end{prop}	

\cref{prop:invariant_hecke_mult} and \cref{cor:invariance-positive} are
	the key technical result of this section.
\begin{prop}
\label{prop:invariant_hecke_mult}
	Given a family of polynomials $(P_x)_{x\leq w}
			\in\ZZ[q^{\frac12},q^{-\frac12}]$ with $P_w\neq0$, a finitary
		subset $J\subseteq S$, and a subset $J_*\subseteq J\cap\rdes P$, if
		$\underline h^P=\ds\sum_xp_xH_x$, then 
		\[
			\underline h^{P^J}\circ_{J_*}\underline h^P
				=\ds\sum_x\sum_{y\in W_J^{J_*}}
					q^{(\ell(y)-\ell(\mxml J^{J_*})+\ell(x_J))/2}p_{x^Jy}H_x
		\]
	If also $P_x\in\NN[q^{\frac12},q^{-\frac12}]$ for all
		$x\in wW_J$, then $\qy*{\rparp P{J_*}J_x}_{x\leq\rcos Jw\mxml J}$ 
		with $0\neq\qy*{\rparp P{J_*}J}_{\rcos Jw\mxml J}
			\in\NN[q^{\frac12},q^{-\frac12}]$ and 
		$\underline h^{P^J}\circ_{J_*}\underline h^P
				=q^{d/2}\underline h^{\rparp P{J_*}J}$,
		where
		\begin{multline*}
			d=\deg\qy*{\rparp P{J_*}J}_{\rcos Jw\mxml J}(q)
			-\deg\qy*{\rparp P{J_*}J}_{\rcos Jw\mxml J}(q^{-1})\\
			-\qy*{\deg P_w(q)-\deg P_w(q^{-1})+\ell(\rcos[J]{J_*}w)}
		\end{multline*}
		does not depend on $J_*$, and has the property that, if 
		$\underline h^P$ is involution-invariant, then $d=0$ is equivalent
		to $\underline h^{\rparp P{J'_*}J}$ being involution-invariant for
		some, equivalently every, subset $J'_*\subseteq J\cap\rdes P$.

	Furthermore, if $P_x\in\NN[q^{\frac12}]$ for all $x\in wW_J$ and
		$P_w(0)\neq0\neq P_{\rcos Jw}(0)$, then
		$\rparp P{J_*}J_{\rcos Jw\mxml J}\in\NN[q^{\frac12}]$ with
		$\rparp P{J_*}J_{\rcos Jw\mxml J}(0)\neq0$, and moreover $0=d=
			(\deg\rparp P{J_*}J_{\rcos Jw\mxml J}(q)-\deg P_w(q))
				-\ell(\rcos[J]{J_*}w)$ holds if and only if 
		$\deg P_x(q)\leq\deg P_w(q)+\ell(w)-\ell(x)$ holds for all 
		$x\in wW_J$, which is implied by $\underline h^P$ being
		involution-invariant.
	\begin{proof}
		\cref{defn:schur}, \cref{exmp:r=1},
			and \cref{prop:hecke_mult} give the first claim:
			\begin{align*}
				\underline h^{P^J}\circ_{J_*}\underline h^P
					&=\qy*{q^{-\ell(\mxml{J_*})/2} W_{J_*}(q)}^{-1}
						{\underline h^P\cdot \underline h^{P^J}}\\
					&=q^{\len{\mxml{J_*}}/2}W_{J_*}(q)^{-1}
						\qy*{\sum_xp_xH_x}\cdot q^{-\len{\mxml J}/2}
							\qy*{\sum_{y\in W_J}T_y}\\
					&=q^{-\ell(\mxml J^{J_*})/2}\sum_x\sum_{y\in W_J^{J_*}}
						q^{(\ell(y)+\ell(x_J))/2}p_{x^Jy}H_x
			\end{align*}

		Next, since $P_x=0$ unless $x\leq w$, we have that $\rparp P{J_*}J_x
			=\ds\sum_{y\in\rcos[J]{J_*}W}q^{\ell(y)}P_{\rcos Jxy}$
			is nonzero only if $\rcos Jxy\leq w$ for some $y\in\rcos[J]{J_*}W$, or
			equivalently by \cref{prop:bruhat-coset} only if $\rcos
			Jx\leq\rcos Jw$, which is equivalent by
			\cref{prop:finitary} to $x\leq \rcos Jw\mxml J$.
		Now $(P_x)_{x\leq w}\in\NN[q^{\frac12},q^{-\frac12}]$ implies
			$(\rparp P{J_*}J_x)_{x\in W}\in\NN[q^{\frac12},q^{-\frac12}]$
			and also $\rparp P{J_*}J_{\rcos Jw\mxml J}\neq 0$ since it contains
			the summand $q^{\ell(\rcos[J]{J_*}w)}P_{\rcos Jw\rcos[J]{J_*}w}
			=q^{\ell(\rcos[J]{J_*}w)}P_w\neq0$, where the latter equality
			holds because $J_*\subseteq\rdes P$.
		In particular, the normalized Hecke algebra element
			$\underline h^{\rparp P{J_*}J}$ is well-defined and
			since $H_x=q^{-\ell(x)/2}T_x$ and $p_{\rcos Jxy}
				=q^{-(\deg P_w(q)-\deg P_w(q^{-1})+\ell(w)-\ell(\rcos Jxy))/2}
					P_{\rcos Jxy}$ by \cref{defn:normalized}, 
			we have $\underline h^{P^J}\circ_{J_*}\underline h^P
				=\ds\sum_x\sum_{y\in\rcos[J]{J_*}W}
					q^{d_{x,y}/2}P_{\rcos Jxy}T_x$
				where
			\begin{align*}
				d_{x,y}
				&=-\len{\rcos{J_*}{\mxml J}}+\ell(y)+\ell(x_J)\\
				&\qquad
					-\qy*{\deg P_w(q)-\deg P_w(q^{-1})+\ell(w)-\ell(\rcos Jxy)}
					-\ell(x)\\
				&=-\ell(\rcos Jw\mxml J^{J_*}w_{J_*})+2\ell(y)
					-\qy*{\deg P_w(q)-\deg P_w(q^{-1})}-\ell(\rcos[J]{J_*}w)\\
				&=-\ell(\rcos Jw\mxml J)-\qy*{\deg P_w(q)-\deg P_w(q^{-1})
					+\ell(\rcos[J]{J_*}w)}+2\ell(y)
			\end{align*}
			with $\rcos{J_*}{\mxml J}w_{J_*}=\rcos{J_*}{\mxml J}\mxml{J_*}=\mxml J$
				because of \cref{prop:finitary} and
				$J_*\subseteq J\cap\rdes P\subseteq J\cap\rdes w$ holding by
				\cref{lem:pre-descent}.
		Thus $\underline h^{P^J} \circ_{J_*}\underline h^P
				=q^{d/2}\underline h^{\rparp P{J_*}J}$
			for $d=\ell(\rcos Jw\mxml J)
				-\qy*{\deg\rparp P{J_*}J_{\rcos Jw\mxml J}(q)
				-\deg\rparp P{J_*}J_{\rcos Jw\mxml J}(q^{-1})}
				+d_{x,y}-2\ell(y)$ which simplifies to the desired value.

		To see that $d$ does not depend on the subset $J_*\subseteq
			J\cap\rdes P$, note that on the one hand by
			\cref{prop:hecke_mult} we have
			$\rparp P{J_*}J_{\rcos Jw\mxml J}=W_{J_*}(q)\rparp P\emptyset J$ 
			with $\deg W_{J_*}(q)=\ell(\mxml{J_*})$ and 
			$\deg W_{J_*}(q^{-1})=0$, while on the other hand
			\cref{lem:pre-descent} gives
				$J_*\subseteq\rdes P\subseteq\rdes w$ 
			so that $w_{J_*}=\mxml{J_*}$ by \cref{prop:finitary} and
			thus $\ell(w_J)=\ell(\rcos[J]{J_*}w)+\ell(\mxml{J_*})$.
		Then $\deg\rparp P{J_*}J_{\rcos Jw\mxml J}(q)
			-\deg\rparp P{J_*}J_{\rcos Jw\mxml J}(q^{-1})-\ell(\rcos[J]{J_*}w)
			=\deg\rparp P\emptyset J_{\rcos Jw\mxml J}(q)-\ell(\mxml{J_*})
			-\deg\rparp P\emptyset J_{\rcos Jw\mxml J}(q^{-1})-\ell(\rcos[J]{J_*}w)
			=\deg\rparp P\emptyset J_{\rcos Jw\mxml J}(q)
			-\deg\rparp P\emptyset J_{\rcos Jw\mxml J}(q^{-1})-\ell(w_J)$
			does not depend on the subset $J_*\subseteq J\cap\rdes P$, and
			hence neither does $d$.

		Next, suppose $\underline h^P$ is involution-invariant. 
		Then so is 
			$\qy*{q^{-\ell(\mxml{J_*})/2}W_{J_*}(q)}^{-1}\underline h^P
				\cdot\underline h^{P^J}
					=\underline h^{P^J}\circ_{J_*}\underline h^P
					=q^{d/2}\underline h^{\rparp P{J_*}J}$ because
			$q^{-\ell(\mxml{J_*})/2}W_{J_*}(q)$ is involution-invariant by
			\cref{lem:invariant_Poincare} and $\underline h^{P^J}$ is
			involution-invariant by \cref{exmp:r=1}.
		Thus, $d=0$ implies $\underline h^{\rparp P{J_*}J}$ is
			involution-invariant and conversely, if $\underline h^{\rparp
			P{J_*}J}$ involution-invariant, then 
			$q^{d/2}\underline h^{\rparp P{J_*}J}
				=\inv{q^{d/2}\underline h^{\rparp P{J_*}J}}
				=q^{-d/2}\inv{\underline h^{\rparp P{J_*}J}}
				=q^{-d/2}\underline h^{\rparp P{J_*}J}$, which implies
			$q^{d/2}=q^{-d/2}$ because $\underline h^{\rparp P{J_*}J}\neq 0$,
			and thus $d=0$ as desired.

		For the supplemental claim, suppose that 
			$P_x\in\NN[q^{\frac12}]$ for all $x\in wW_J$ and 
			$P_{\rcos Jw}(0)\neq0$, which evidently imply 
			$\rparp P{J_*}J_{\rcos Jw\mxml J}
				=\ds\sum_{y\in\rcos[J]{J_*}W}q^{\ell(y)}P_{\rcos Jxy} 
					\in\NN[q^{\frac12}]$ and
			$\rparp P{J_*}J_{\rcos Jw\mxml J}(0)=P_{\rcos Jw}(0)\neq0$. 
		Consequently, for $J_*=\emptyset$ we have
			$\deg\qy{\rparp P\emptyset J}_{\rcos Jw\mxml J}(q)
				=\max\set{\ell(y)+\deg P_{\rcos Jwy}(q):y\in W_J}
				=\max\set{\ell(x)-\ell(\rcos Jw)+\deg P_x(q):x\in wW_J}$ (recall
				from \cref{rmrk:bounded-has-pw-nonzero} that $\deg
				P_x(q)=-\infty$ if and only if $P_x=0$).

		Now, since also $P_w\in\NN[q^{\frac12}]$ with $P_w(0)\neq0$, we have
			that $d=\deg\qy{\rparp P\emptyset J}_{\rcos Jw\mxml J}(q)
				-\qy*{\deg P_w(q)+\ell(w_J)}
			=\max\set{\ell(x)-\ell(\rcos Jw)+\deg P_x(q):x\in wW_J}
				-(\deg P_w(q)+\ell(w_J))
			=\max\set{\deg P_x(q)-(\deg P_w(q)+\ell(w)-\ell(x)):x\in wW_J}$.
		Since for $x=w$ we have $0\neq\deg P_x(q)
			=\deg P_w(q)+\ell(w)-\ell(x)$, it follows that 
			$\deg P_x(q)\leq\deg P_w(q)+\ell(w)-\ell(x)$ for all $x\in wW$
			is equivalent to $0=d$, as desired.

		Finally, suppose again that $\underline h^P$ is
			involution-invariant, which by
			\cref{exmp:normalized-invariance-equation} implies that
			$\ds q^{\deg P_w(q)+\ell(w)-\ell(x)}\inv{P_x}
				=\sum_yR_{x,y}P_y$ for each $x$. 
		Fixing $x\in wW_J$, i.e.\ such that $\rcos Jx=\rcos Jw$, we have on
			the one hand that $R_{x,y}\neq0$ only if $x\leq y$ by
			\cref{prop:R-polynomials}(\ref{item:R-constant}), and
			on the other hand that $P_y\neq0$ only if $y\leq w$, whence by 
			\cref{prop:bruhat-coset} we have $R_{x,y}P_y\neq0$ only
			if $\rcos Jw=\rcos Jx\leq \rcos Jy\leq\rcos Jw$, i.e.\ only if
			$y\in wW_J$.
		But then since $R_{x,y}\in\ZZ[q]$ by
			\cref{prop:R-polynomials}(\ref{item:R-constant}) and
			$P_y\in\ZZ[q^{\frac12}]$ for $y\in wW_J$ by assumption,
			it follows that $q^{\deg P_w(q)+\ell(w)-\ell(x)}
				\inv{P_x}\in\ZZ[q^{\frac12}]$, whence
			$0\leq\deg P_w(q)+\ell(w)-\ell(x)-\deg\inv{P_x}(q^{-1})
				=\deg P_w(q)+\ell(w)-\ell(x)-\deg P_x(q)$, i.e.\
			$\deg P_x(q)\leq \deg P_w(q)+\ell(w)-\ell(x)$, as desired.
	\end{proof}
\end{prop}

\begin{cor}
\label{cor:invariance-positive}
	Given a finitary subset $J\subseteq S$ and a positive family
		$(P_x)_{x\leq w}\in\NN[q]$, for	any subset
		$J_*\subseteq\rdes P\cap J$ we have that
		$(\rparp P{J_*}J_x)_{x\leq\rcos Jw\mxml J}$ is positive,
		and that if $\underline h^P$ is involution-invariant, then 
		$\deg\rparp P{J_*}J_{\rcos Jw\mxml J}(q)
			=\deg P_w(q)+\ell(\rcos[J]{J_*}w)$ and
		$\underline h^{\rparp P{J_*}J}
			=\underline h^{P^J}\circ_{J_*}\underline h^P$ is
			involution-invariant.
	Furthermore, still assuming that $\underline h^P$ is
		involution-invariant, we have that 
		$\qy*{\rparp P{J_*}J}_{\rcos Jw\mxml J}=P_w$ if and only if
		$\ell(\rcos[J]{J_*}w)=0$, which is equivalent to $w_J=\mxml{J_*}$,
		implies $J_*=J\cap\rdes w$, and is implied by 
		$\deg\qy*{\rparp P{J_*}J}_{\rcos Jw\mxml J}(q)\leq0$.
	\begin{proof}
		Since $(P_x)_{x\leq w}\in\NN[q]$ with $P_w\neq0$, then
			$\rparp P{J_*}J_x
				=\ds\sum_{y\in\rcos[J]{J_*}W}q^{\ell(y)}P_{\rcos Jxy}
				\in\NN[q]$ has $\rparp P{J_*}J_x(0)=P_{\rcos Jx}(0)$
			so that $(\rparp P{J_*}J)_{x\leq\rcos Jw\mxml J}$ is indeed
			positive if $(P_x)_{x\leq w}$ is positive. 
		That $\underline h^P$ involution-invariant implies 
			$\underline h^{\rparp P{J_*}J}$ is involution-invariant and that
			$\deg\rparp P{J_*}J_{\rcos Jw\mxml J}(q)
				=\deg P_w(q)+\ell(\rcos[J]{J_*}w)$ follows immediately from the
				supplemental claim of \cref{prop:invariant_hecke_mult}.
		Consequently in that case $\len{\rcos[J]{J_*}w}$ is necessary for 
			$\qy*{\rparp P{J_*}J}_{\rcos Jw\mxml J}=P_w$, and since
			$(P_x)_{x\leq w}$ positive means $\deg P_w\geq 0$, it is also
			implied by $\deg\qy*{\qy*{\rparp P{J_*}J}_{\rcos Jw\mxml J}(q)}
			\allowbreak\leq 0$.
			
		To see that it is sufficient, note that if $\len{\rcos[J]{J_*}w}
			=0$, i.e.\ $\rcos[J]{J_*}w=e$, then for $y\in\rcos[J]{J_*}W$ we
			have $P_{\rcos Jwy}\neq 0$ only if $\rcos Jwy\leq w$, which
			by \cref{prop:bruhat-coset} is equivalent to
			$y\leq\rcos[J]{J_*}w=e$, i.e.\ to $y=e$.
		Thus, we have if $\len{\rcos[J]{J_*}w}=0$ then
			$\qy*{\rparp P{J_*}J}_{\rcos Jw\mxml J}
				=\ds\sum_{y\in\rcos[J]{J_*}W}q^{\ell(y)}P_{\rcos Jwy}
				=P_{\rcos Jw}$
			where since $e=\rcos[J]{J_*}w$ and $w_{J_*}\in W_{J_*}\subseteq
				W_{\rdes P}$ we have $P_{\rcos Jw}=P_w$ as desired.

		Finally, since $\rdes P\subseteq\rdes w$ by
			\cref{lem:pre-descent}, we have that $J_*\subseteq\rdes
			P\subseteq\rdes w$ and hence $w_{J_*}=\mxml{J_*}$ by
			\cref{prop:finitary}.
		Thus, $\ell(\rcos[J]{J_*}w)=0$, i.e.\ $\rcos[J]{J_*}w=e$, is
			equivalent to $w_J=w_{J_*}=\mxml{J_*}$. 
		But since $\rref{w_J}=\rref w\cap W_J$ and
			$\rref{\mxml{J_*}}=T\cap W_{J_*}$ by definition of $\mxml{J_*}$, we
			have by \cref{lem:injective-cocycle} that $\rcos[J]{J_*}w=e$,
			is equivalent to $\rref w\cap W_J=T\cap W_{J_*}$, and hence by
			\cref{lem:subsystem} that it implies $J_*=S\cap T\cap
			W_{J_*}=S\cap\rref w\cap W_J=\rdes w\cap J$.
	\end{proof}
\end{cor}

\begin{cor}
\label{cor:quotient}
	Given a factorization with overlaps $J_1,\dots,J_r$, the normalized
		Hecke algebra element $\underline h^{P^J}$ is in the closure of
		$\set{C'_{\mxml J}: J\subseteq S$ is finitary$}$ under the Schur
		composition laws, so is involution-invariant and furthermore
		satisfies
		\[
			C'_{\mxml{J_1}}\cdots C'_{\mxml{J_r}}
				=\ds\qy*{\prod_{1\leq j<t\leq r}
					q^{-\len{\mxml{\overlap Jjt}}/2} W_{\overlap Jjt}(q)}
						\underline h^{P^J}
		\]
	In particular, the family of defect polynomials is self-dual. 
	More generally given $k$ such that $0\leq k\leq r$, if we set
	$\overlap J{\leq k}{>k}=\ds\bigcup_{i\leq k<j}\overlap Jik$ 
		then for any subset
		$J_*\subseteq\overlap J{\leq k}{>k}$ we have
			$\underline h^{P^{\downslice J{k+1}}}
				\circ_{J_*}\underline h^{P^{\upslice Jk}}
			=q^{-\len{\rcos{J_*}{\mxml{\overlap J{\leq k}{>k}}}}/2}
				\rcos[\overlap J{\leq k}{>k}]{J_*}W(q)\underline h^{P^J}$
		so that $\underline h^{P^J}
			=\underline h^{P^{\downslice J{k+1}}}
				\circ_{\overlap J{\leq k}{>k}}\underline h^{P^{\upslice Jk}}$.
	\begin{proof}
		\cref{cor:recursive-defect-polys} implies
			$h^{P^{\upslice J{k+1}}}
				=\rparp{P^{\upslice Jk}}{\overlap J*{k+1}}{J_{k+1}}$
			with $\overlap J*{k+1}\subseteq\rdups Jk\cap J_{k+1}
				\subseteq\rdes{P^{\upslice Jk}}\cap J_{k+1}$ and 
				$(P^{\upslice Jk})_{x\leq w(\upslice Jk)}\in\NN[q]$ positive,
			while \cref{cor:length+defects<=meetings} gives
			$\deg P^{\upslice Jk}_x
				\leq\deg P^{\upslice Jk}_{w(\upslice Jk)}
					+\ell(w(\upslice Jk))-\ell(x)$ for all $x\in W$,
			so by \cref{prop:invariant_hecke_mult}
			$\underline h^{P^{J_{k+1}}}
					\circ_{\overlap J*{k+1}}\underline h^{P^{\upslice Jk}}
				=\underline h^{P^{\upslice J{k+1}}}$ for all $k$.

		Then induction on the length $r$ of the factorization gives by
			definition of the Schur composition laws
			(\cref{defn:schur}) that
			\begin{align*}
				\underline h^{P^J}
				&=\prod_{k=1}^r\qy*{q^{-\len{\mxml{\overlap J*k}}/2}
						W_{\overlap J*k}(q)}^{-1}
					\underline h^{P^{J_1}}\cdots \underline h^{P^{J_r}}
			\end{align*}
		Since $\overlap J*k=\ds\bigudj_{i=1}^{k-1}\overlap Jik$ 
			is the union of commuting sets, we have by
			\cref{lem:component-decompose}) and \cref{lem:bijections}
			that
			\[
				\prod_{k=1}^r\qy*{q^{-\len{\mxml{\overlap J*k}}/2}
					W_{\overlap J*k}(q)}
					=\prod_{1\leq j<t\leq r}q^{-\len{\mxml{\overlap Jjt}}/2}
						W_{\overlap Jjt}(q)
			\]
			giving the main claim.

		For the supplemental claim, since for any $j$ and $t$ such that 
			$1\leq j<t\leq k$ we have
			$\overlap{\qy*{\upslice Jk}}jt=\overlap Jjt$ and for any $j$ and
			$t$ such that $k+1\leq j<t\leq r$ we have
			$\overlap{\qy*{\downslice J{k+1}}}jt=\overlap Jjt$, the main claim
			gives us
			\begin{align*}
				&\prod_{1\leq j<t\leq r}q^{-\len{
					\mxml{\overlap Jjt}}/2}W_{\overlap Jjt}(q)
					\underline h^{P^J}\\
				&\qquad= C'_{\mxml{J_1}}\cdots C'_{\mxml{J_k}}
					C'_{\mxml{J_{k+1}}}\cdots C'_{\mxml{J_r}}\\
				&\qquad= \prod_{1\leq j<t\leq k}q^{-\len{
					\mxml{\overlap Jjt}}/2} W_{\overlap Jjt}(q)
					\underline h^{P^{\upslice Jk}}
					\prod_{k+1\leq j<t\leq r}q^{-\len{\mxml{ \overlap Jjt}}/2}
					W_{\overlap Jjt}(q) \underline h^{P^{\downslice J{k+1}}}
			\end{align*}
			and hence
			$\underline h^{P^{\upslice Jk}}\cdot
				\underline h^{P^{\downslice J{k+1}}}
				=\ds\prod_{\substack{1\leq j\leq k\\k+1\leq t\leq r}}
				q^{-\len{\mxml{\overlap Jjt}}/2}
				W_{\overlap Jjt}(q) \underline h^{P^J}$.
		Since by definition 
			$\ds\overlap J{\leq k}{>k}=\bigudj_{\substack{1\leq j\leq
			k\\k+1\leq t\leq r}} \overlap Jjt$ is a union of disjoint
			commuting subsets,
			\cref{lem:coset-properties}(\ref{lem:component-decompose}) and
			\cref{lem:bijections} give the supplemental claim.
	\end{proof}
\end{cor}

\begin{exmp}
	A factorization $J_1,\dots,J_r$ with Gelfand--MacPherson overlaps is
		the same data as a path
		$\emptyset\xrightarrow{C'_{\mxml{J_1}}}\overlap J12
			\xrightarrow{C'_{\mxml{J_2}}}\cdots
			\xrightarrow{C'_{\mxml{J_{r-1}}}}\overlap J{r-2}{r-1}
			\xrightarrow{C'_{\mxml{J_r}}}\emptyset$ of
		morphisms in Williamson's Schur algebroid.
	The composite morphism of such a path in the Schur algebroid is
		precisely its normalized Hecke algebra element $\underline h^{P^J}$.

	Accordingly, the normalized Hecke algebra element $\underline h^{P^J}$
		of a factorization with overlaps $J$ that is \emph{not} a
		Gelfand--MacPherson factorization is not necessarily a
		composite in the Schur algebroid of a path of morphisms because the
		repeated application of our Schur composition laws
		may at any stage reinterpret the resulting Hecke algebra
		element as representing a morphism with a different domain or
		codomain than the one naturally associated to it by composition in
		the Schur algebroid.
\end{exmp}

\section{Combinatorially interpreting Kazhdan--Lusztig polynomials}

\label{sec:deodhar}
We review the definition of Kazhdan--Lusztig basis elements, describe
	the generic algorithm for computing them, and use the results of
	Williamson \cite{Williamson} to show that Soergel's conjecture implies
	our generalization of Deodhar's theorem: the families of defect
	polynomials of \emph{admissible factorizations} of $w$ have
	restrictions giving a combinatorial interpretation of the coefficients
	of Kazhdan--Lusztig basis elements.
Along the way we prove the key fact that the properies of a factorization
    being admissible or tight are absolute, and also obtain a new
    derivation of the M\"obius function of the Bruhat order.

We also generalize Billey and Warrington's argument \cite[Proof of
	Theorem~1]{Billey--Warrington} to show that in the case of Weyl groups
	a factorization is tight if and only if it induces a particular
	factorization of the Poincar\'e polynomial of full intersection
	cohomology of the corresponding Schubert variety, and also realizes
	the polynomial as a generating function.

\begin{lem}
\label{lem:bounded}
	Given $x\in W$ and a family
		$(p_y)_{x<y}\in\ZZ[q^{\frac12},q^{-\frac12}]$, there is
		at most one $p_x\in\ZZ[q^{\frac12},q^{-\frac12}]$ such that 
		$\deg p_x(q)<0$ and such that $(p_y)_{x\leq y}$ extends to a family
		$(p_z)_{z\in W}$ such that the Hecke algebra element
		$h=\ds\sum_z p_zH_z$ is involution-invariant.

	In particular, given $x\leq w$ and
		$(P_y)_{x<y\leq w}\in\ZZ[q^{\frac12},q^{-\frac12}]$ with $P_w\neq0$,
		there is at most one $P_x\in\ZZ[q^{\frac12},q^{-\frac12}]$
		such that $\deg P_x(q)<\frac12\qy{\deg P_w(q)-\deg P_w(q^{-1})
			+\ell(w)-\ell(x)}$ and such that $(P_x)_{x\leq y\leq w}$ extends
		to a family $(P_z)_{z\leq w}\in\ZZ[q^{\frac12},q^{-\frac12}]$ with
		$\underline h^P$ involution-invariant.
	\begin{proof}
		\cref{prop:R-polynomials}(\ref{item:invariance})
			gives us that $\inv h=h$ only if
			$\inv{p_x}=\ds\sum_yq^{(\ell(x)-\ell(y))/2}R_{x,y}p_y$,
			or equivalently $p_x-\inv{p_x}
				=\ds-\sum_{x<y}q^{(\ell(x)-\ell(y))/2}R_{x,y}p_y$. 
		If $\deg p_x(q)<0$, then $p_x$ and $\inv{p_x}$ are linearly
			independent, giving the main claim.
		The supplemental claim follows from the main claim and the definition
			of $\underline h^P$.
	\end{proof}
\end{lem}

\begin{rmrk}
	The special case of the supplemental claim of \cref{lem:bounded}
		where $1=P_w$ is \cite[2.2.a]{Kazhdan--Lusztig}.
\end{rmrk}

\begin{defn}
\label{defn:bounded}
	Given a non-zero Hecke algebra element $h=\ds\sum_x p_xH_x$,
		we define the \define{set $U(h)$ of unbounded indices} to be
		$U(h)=\set{x\in W:\deg p_x(q)\geq0$ and there exists $y$ such that
		$y>x$ and $p_y\neq 0}$.
	We say that $h$ is \define{bounded}
		if it is non-zero, involution-invariant, and such that 
		$U(h)=\emptyset$.

	We say that  a family $(P_x)_{x\leq
		w}\in\ZZ[q^{\frac12},q^{-\frac12}]$ is itself
		\define{bounded} if $P_w\neq0$ and the normalized Hecke algebra
		$\underline h^P$ is bounded, or equivalently if $P_w\neq0$,
		$\underline h^P$ is involution-invariant, and
		$\deg P_x(q)<\frac12\qy{\deg P_w(q)-\deg P_w(q^{-1})
			+\ell(w)-\ell(x)}$ for all $x$ such that $x\neq w$
	(recall from \cref{rmrk:bounded-has-pw-nonzero} that $P_x=0=p_x$
		if and only if $\deg P_x(q)=-\infty=\deg p_x(q)$).

	We say that a family $(P_x)_{x\leq w}\in\ZZ[q^{\frac12},q^{-\frac12}]$
		is \define{monic} if $P_w=1$. 

	For each $w\in W$, the family $(P_{x,w})_{x\leq w}$ of
		\define{Kazhdan--Lusztig polynomials} is the (necessarily unique
		by \cref{lem:bounded}) monic bounded family, and the
		bounded normalized Hecke algebra element $\underline h^{P_{-,w}}$ is
		the \define{Kazhdan--Lusztig basis element} $C'_w$.
\end{defn}

\begin{rmrk}
\label{rmrk:bounded-existence}
	We introduce the notion of \emph{bounded} Hecke algebra elements
		because our generalization of Billey and Warrington's criterion
		for a factorization $J$ being tight amounts to requiring that
		the involution-invariant Hecke algebra element 
		$\underline h^{P^J}$  is bounded (see \cref{defn:tight},
		\cref{prop:bounded-monic}).
	
	The existence of the Kazhdan--Lusztig polynomials (i.e.\ a monic
		bounded family $(P_{x,w})_{x\leq w}$ for each $w\in W$) will
		follow as a special case of \cref{prop:dominance-dcc}
		below, which more generally shows for every non-zero
		$h=\ds\sum_x p_xH_x$ the existence of a (necessarily unique
		by \cref{lem:bounded}) bounded Hecke algebra element
		$h_B=\ds\sum_xp'_xH_x$ such that $w$ is maximal such that
		$p_w\neq0$ if and only if it is maximal such that $p'_w\neq0$, and
		in which case $p_w=p'_w$.
\end{rmrk}

\begin{exmp}
\label{exmp:finitary-C'w}
	For any finitary $J\subseteq S$, we have that 
		$C'_{\mxml J}=\underline h^{P^J}$ since $\underline h^{P^J}$ is
		involution-invariant, monic, and bounded by \cref{exmp:r=1}.
\end{exmp}

\begin{exmp}
\label{exmp:monic-singleton}
	The family of defect polynomials $(P^J_x)_{x\leq w}$ of a
		factorization $J$ induced by a reduced word $s_{i_1}\cdots
		s_{i_r}=w$ not only has an involution-invariant normalized Hecke
		algebra element $\underline h^{P^J}$ by
		\cref{cor:quotient}, but is also monic by
		\cref{exmp:admissible-singleton}.
	More generally, admissible factorizations are by definition those 
	    whose family of defect polynomials is monic.
\end{exmp}

\begin{prop}
\label{prop:bounded-rparp}
	Given a positive family $(P_x)_{x\leq w}\in\NN[q]$ such that 
		$\underline h^P$ is involution-invariant, for any finitary subset
		$J\subseteq S$ and any subset $J_*\subseteq\rdes P\cap J$ the
		family $(\rparp P {J_*}J_x)_{x\leq\rcos Jw\mxml J}$ defined in
		\cref{defn:right-parabolic-product} is bounded if and only
		if $J_*=\rdes w\cap J$, $\ell(\rcos[J]{J_*}w)=0$, and 
		$\ell(\rcos[J]{J_*}x)<\deg P_w(q)+\ell(w)-\ell(x)-2\deg P_x(q)$
		for all $x$ such that $x\not\in wW_J$, and 
		$0=\deg P_w(q)$ when $\emptyset\neq J$, in which case we also have
		$\ell(\rcos[J]{J_*}x)<\deg P_w(q)+\ell(w)-\ell(x)-2\deg P_x(q)$ for
		all $x\in wW_J\setminus\set w$.
	In particular, $(P_x)_{x\leq w}$ is unbounded only if 
		$(\rparp P{J_*}J_x)_{x\leq\rcos Jw\mxml J}$ is unbounded.
	\begin{proof}
		If $J=\emptyset$, then $\rcos[J]{J_*}W=\set e$, so that 
			$\rparp P{J_*}J_x
				=\ds\sum_{y\in\rcos[J]{J_*}W}q^{\ell(y)}P_{\rcos Jxy}
				=P_x$ and $\rcos[J]{J_*}x=e$ for all $x$.
		Then $0=\ell(\rcos[J]{J_*}w)=\ell(\rcos[J]{J_*}x)
				<\deg P_w(q)+\ell(w)-\ell(x)-2\deg P_x(q)$ holding for all $x$
			such that $x\neq w$ is the definition of $(P_x)_{x\leq w}
				=(\rparp P{J_*}J_x)_{x\leq\rcos Jw\mxml J}$ being bounded.

		Supposing now that $\emptyset\neq J$, we have by definition of
			$\rparp P{J_*}J_x
				=\ds\sum_{y\in\rcos[J]{J_*}W}q^{\ell(y)}P_{\rcos Jxy}
				\in\NN[q]$
			that $J\subseteq\rdes{\rparp P{J_*}J}$, and hence by
			\cref{lem:admissible-double-max} that
			$\deg P_x(q)
				-\frac12\qy{\deg\rparp P{J_*}J_{\rcos Jw\mxml J}(q)
					-\deg\rparp P{J_*}J_w(q^{-1})+\ell(\rcos Jw\mxml J)-\ell(x)}<0$
			for $x\in\rcos Jw\mxml JW_J\setminus\set{\rcos Jw\mxml J}
				 =wW_J\setminus\set{\rcos Jw\mxml J}$ if and only if
			$\deg\qy*{\qy*{\rparp P{J_*}J}_{\rcos Jw\mxml J}(q)}\leq 0$,
			which by \cref{cor:invariance-positive} is equivalent to
			$J_*=\rdes w\cap J$, $\ell(\rcos[J]{J_*}w)=0$, $w_J=\mxml{J_*}$, and
			$\deg P_w\leq 0$.

		Thus it remains to show that 
			$\ell(\rcos[J]{J_*}x)<\deg P_w(q)+\ell(w)-\ell(x)-2\deg P_x(q)$
			holds for all for all $x\not\in wW_J$ if and only if for all
			$x\not\in wW_J$ we have
			$\deg\rparp P{J_*}J_x(q)
				-\frac12\qy*{\deg\rparp P{J_*}J_{\rcos Jw\mxml J}(q)
					-\deg\rparp P{J_*}J_{\rcos Jw\mxml J}(q^{-1})
						+\ell(\rcos Jw\mxml J)-\ell(x)}<0$.

		Now $\underline h^P$ involution-invariant implies by
			\cref{cor:invariance-positive} that
			$\deg\rparp P{J_*}J_{\rcos Jw\mxml J}(q)
					-\deg\rparp P{J_*}J_{\rcos Jw\mxml J}(q^{-1})
				=\deg P_w(q)+\ell(\rcos[J]{J_*}w)$.
		Hence for $x\not\in W_J$, since $\rparp P{J_*}J_x$ by definition does
			not depend on $x_J$ while $\ell(x)$ achieves its maximum when
			$x_J=\mxml J$, it follows that the boundedness inequality holds for
			all $x\not\in wW_J$ if and only if 
			$2\deg\rparp P{J_*}J_x
				<\deg P_w(q)+\ell(\rcos[J]{J_*}w)+\ell(\rcos Jw\mxml J)-\ell(\rcos Jx\mxml J)
				=\deg P_w(q)+\ell(\rcos{J_*}w)-\ell(\rcos Jx)
				=\deg P_w(q)+\ell(w)-\ell(\rcos Jx\mxml{J_*})$  holds for all
				$x\not\in wW_J$.

		But positivity of $(P_x)_{x\leq w}$ and $J_*\subseteq\rdes P$
			also implies that $\deg\rparp P{J_*}J_z
				=\max\set{\ell(\rcos[J]{J_*}x)+\deg P_x(q):x\in zW_J}$, so
			that the latter inequality holding for all $x\not\in wW_J$ is
			equivalent to 
			$\ell(\rcos[J]{J_*}x)
				<\deg P_w(q)+\ell(w)-\ell(x)-2\deg P_x$ holding for all
			$x\not\in wW_J$, as desired.

	The supplemental claim follows because $(P_x)_{x\leq w}$ is
			unbounded precisely when $\deg P_w(q)+\ell(w)-\ell(x)-2\deg P_x(q)
				\leq0\leq\ell(\rcos[J]{J_*}x)$ for some $x$ such that $x\neq w$.
	\end{proof}
\end{prop}

\begin{rmrk}
\label{rmrk:bottom-up-tight-alg}
	\cref{prop:bounded-rparp} suggests an
		algorithm for enumerating all tight factorizations since a
		factorization $J_1,\dots,J_r$ satisfies $P^J=\rparp{P^{\upslice
		J{r-1}}}{\overlap J*r}{J_r}$ and hence is tight if and only
		if $\upslice J{r-1}$ is tight and $w(\upslice J{r-1})$ and $J_r$
		satisfy the conditions of \cref{prop:bounded-rparp}.
\end{rmrk}

\begin{cor}
\label{cor:factors}
	Any admissible factorization is absolutely admissible and any tight
		factorization is absolutely tight.
	\begin{proof}
		Since by \cref{cor:spread_defect}, the family of
			defect polynomials $P^J$ depends only on the induced element of
			the Cartier--Foata monoid, it follows that being admissible or
			tight also depends only on the element of the Cartier--Foata
			monoid.
		Moreover, \cref{cor:invariance-positive}, respectively
			\cref{prop:bounded-rparp}, imply that if $J$ is admissible,
			respectively tight, then so is $\upslice Jt$.
		Since defect polynomials are self-dual by \cref{cor:quotient}, it
			follows that $\slice Jkt$ is also admissible, respectively tight.
		Thus, the properties of being admissible or tight are absolute.
	\end{proof}
\end{cor}

Our generalization of Deodhar's combinatorial interpretation of
	Kazhdan--Lusztig polynomials, as well as of his and Scott and Sprowl's
	algorithms for computing them, relies on the following generalization
	to Hecke algebra elements of his dominance relation on subsets of
	masks introduced in \cite{Deodhar}.

\begin{defn}
\label{defn:dominance}
	Given a collection $\set{M_y}_{y\in W}$ of (non-empty) sets of monic
		families of polynomials
		$M_y\subseteq\set{(Q_x)_{x\leq y}\in\ZZ[q^{\frac12},q^{-\frac12}]:
		Q_y=1,\underline h^Q$ is involution-invariant$}$, we define the
		relation of \define{$M$-dominance} on Hecke algebra elements by the
		generating relations $h'\doml_M h$ if $h-h'=L\cdot\underline h^Q$
		for some non-zero $L\in\ZZ[q^{\frac12}+q^{-\frac12}]$ and $Q\in M_y$
		where $y\in U(h)$ is maximal in the Bruhat order on $U(h)$ and such
		that $y\not\in U(h')$.

	Once we know by \cref{cor:bounded-algorithm} that all
		Kazhdan--Lusztig basis elements exist,
		we will refer to the $M$-dominance relation for
		$M_y=\set{(P_{x,y})_{x\leq y}}$ simply as the \define{KL-dominance}
		relation, and we will simply write $h'\doml_{\mathrm{KL}} h$ instead
		of $h'\doml_M h$.
\end{defn}

\begin{prop}
\label{prop:dominance-dcc}
	Given a collection $\set{M_y}_{y\in W}$ as in
		\cref{defn:dominance}, its associated $M$-dominance
			relation has the following properties:
	\begin{enumerate}
		\item $\ds\sum_x P'_xT_x\doml_M\sum_x P_xT_x$ implies
			$P_w=P'_w$ for every $w$ maximal in the Bruhat order such that
			$P_w\neq0$ or $P'_w\neq0$;
		\item $h'\doml_Mh$ implies $\ordideal{U(h')}
			\subsetneq\ordideal{U(h)}$ where $\ordideal{U(h)}$ is the
			lower set generated by $U(h)$ in the Bruhat order;
		\item\label{prop:dominance-down}
		If $h=\ds\sum_x p_xH_x$, $y\in U(h)$ is maximal
			in the Bruhat order on $U(h)$, and $(Q_x)_{x\leq y}\in M_y$, then
			$h-L_y\cdot \underline h^Q\doml_Mh$ if and only if
			$L_y=p_y^{\geq 0}+\inv{p_y^{>0}}$ where $p_y^{\geq 0}$
			(respectively $p_y^{>0}$) denotes the sum of the terms in $p_y$
			with non-negative (respectively positive) degree. 
	\end{enumerate}
	In particular, the $M$-dominance relation generates a partial order
		such that if a Hecke algebra element $h=\ds\sum_xp_xH_x$
		is larger than $h'=\ds\sum_xp'_xH_x$ in this order,
		then $h=h'+\ds\sum_y L_yh_y$
		where $h_y=\underline h^Q$ for some $(Q_x)_{x\leq y}\in M_y$,
		$L_y=0$ unless $y\in\ordideal{U(h)}\setminus\ordideal{U(h')}$, and
		where for some $y\in U(h)$ maximal in the Bruhat order on $U(h)$ we
		have
		$L_y=p^{\geq0}_y+\inv{p^{>0}}_y\in\ZZ[q^{\frac12}+q^{-\frac12}]$.
	Furthermore, the generated partial order admits from any
		involution-invariant Hecke algebra element
		$h=\ds\sum_xp_xH_x$ a descending chain to the unique
		bounded Hecke algebra element $h_B=\ds\sum_xp'_xH_x$
		such that $w$ is maximal with $p'_w\neq 0$ if and only
		if it is maximal such that $p_w\neq 0$, in which case $p_w=p'_w$.
	\begin{proof}
		To verify (1) and (2) it will suffice to consider only the
			generating relations of the $M$-dominance relation. 
		To that end, suppose that $\ds\sum_xP'_xT_x=\sum_x p'_xH_x
				=h'\doml_M h=\sum_x p_xH_x=\sum_xP_xT_x$ with
			$h-h'=L\cdot\underline h^Q$, $L\in\ZZ[q^{\frac12}+q^{-\frac12}]$,
			and $Q\in M_y$ for some $y\in U(h)$ maximal in the Bruhat order on
			$U(h)$ and such that $y\not\in U(h')$.
		Both proofs will use the fact that if $\underline h^Q
			=\ds\sum_x q_xH_x$, so that 
			$p'_x=p_x-Lq_x$ for every $x$, then from $(Q_x)_{x\leq y}$ we have
			$p'_x=p_x$ unless $x\leq y$.
		
		(1) On the one hand, since by definition $y\in U(h)$ implies there
			exists $w\in W$ such that $y<w$ and $p_w\neq0$, we have that every
			$w\in W$ maximal such that $P_w\neq 0$ satisfies $w\not\leq y$,
			and hence $p'_w=p_w\neq 0$ and the desired $P'_w=P_w$. 

		On the other hand, if $p'_x\neq 0$ then either $p_x=p'_x$ or $x\leq
			y$, and there exists $w\geq x$ maximal such that $p_w\neq 0$ and
			hence such that $p'_w=p_w$ by the above.
		Consequently,	$w$ maximal in the Bruhat order such that $p'_w\neq 0$
			implies that $w$ is also maximal such that $p_w\neq 0$ and hence
			satisfies $P_w=P'_w$.

		(2) We have that $0\leq\deg p'_x(q)$ is equivalent to
			$0\leq\deg p_x(q)$ if $x\not\leq y$ and hence that
			$U(h)\setminus\ordideal{\set y}
				=U(h')\setminus\ordideal{\set y}$.
		But $y\in U(h)$ is maximal in the Bruhat order of $U(h)$ and
			satisfies $y\not\in U(h')$ by assumption, so that
			$\ordideal{U(h')}
				\subseteq\ordideal{U(h)}\setminus\set y
				\subsetneq\ordideal{U(h)}$. 

		(3) Since $y\in U(h)$ we have that
			$h-L_y\cdot\underline h^Q\doml_M h$ for $(Q_{x,y})_{x\leq y}
				\in M_y$ and $L_y\in\ZZ[q^{\frac12}+q^{-\frac12}]$
			if and only if $y\not\in U(h-L_y\cdot\underline h^Q)$, which
			since $y\in U(h)$ implies is equivalent to
			$\deg(p_y-L_y)<0$.
		But this is equivalent to $L_y^{\geq 0}=p_y^{\geq 0}$, which
				since $L_y\in\ZZ[q^{\frac12}+q^{-\frac12}]$
				is equivalent to $L_y=p_y^{\geq 0}+\inv{p_y^{>0}}$.
		
		For the supplemental claim, we have by (2) that the transitive
			closure of the $M$-dominance relation is a partial order,
			and furthermore since any lower order ideal $\ordideal{U(h)}$ is
			the finite set of Coxeter elements expressible as subwords of fixed
			reduced words for the finitely many elements of $U(h)$, this partial
			order is well-founded.
		Thus minimal elements exist, and (3) shows that the minimal elements
			have $U(h)=\emptyset$, and thus those comparable to an
			involution-invariant element are bounded Hecke algebra elements.
		(1) shows that all comparable elements have the same $w$ maximal with
			$P_w\neq 0$ and the same $P_w$ for such $w$.
		Uniqueness of the bounded Hecke algebra element $h_B$ below
			an involution-invariant Hecke algebra element $h$ thus follows from
			(1) and \cref{lem:bounded}.
		Now, the expression $h=h'+\ds\sum_y L_yh_y$ for
			$h_y=\underline h^Q$ for some $(Q_{x,y})_{x\leq y}\in M_y$ follows
			from the definition of the dominance order, and that
			$L_y=0$ unless $y\in\ordideal{U(h)}\setminus\ordideal{U(h')}$
			follows from (2).
		That there is some $y\in U(h)$ maximal in the Bruhat order with
			$L_y=p_y^{\geq 0}+\inv{p_y^{>0}}\in\ZZ[q^{\frac12}+q^{-\frac12}]$
			follows by (3) and (2), since after descending in the $M$-dominance
			order with the choice $y\in U(h)$, every subsequent descent
			in the order must involve only $y'$ such that $y'\not\geq y$.
	\end{proof}
\end{prop}

\begin{cor}
\label{cor:bounded-algorithm}
	For every $w\in W$, the Kazhdan--Lusztig basis element $C'_w$ exists.
	More precisely, for every pair of collections $\set{I_w}_{w\in W}$ and
		$\set{M_y}_{y\in Y}$ of non-empty
		sets of monic families of polynomials $I_x,M_x
			\subseteq\set{(P_w)_{x\leq w}\in\ZZ[q^{\frac12},q^{-\frac12}]
				:P_w=1,\underline h^P$ is
		involution-invariant$}$, the family of Kazhdan--Lusztig polynomials
		$(P_{x,w})_{x\leq w}$ can be computed by making an effective choice
		$P\in I_w$, and then descending from $\underline h^P$ along
		the $M$-dominance relation by making effective choices of
		$(Q_x)_{x\leq y}\in M_y$.
	\begin{proof}
		Observe that because \cref{exmp:monic-singleton} implies each
			factorization induced by a reduced word has a monic family of
			defect polynomials with involution-invariant normalized Hecke
			algebra element, the families of sets $\set{I_w}_{w\in W}
				=\set{M_w}_{w\in W}=\set{(P^J_x)_{x\leq w}:J$ is a factorization
				of $w$ such that $P^J_w=1}$ are non-empty for each $w\in W$.
		The result then follows from \cref{prop:dominance-dcc}.
	\end{proof}
\end{cor}

\begin{exmp}
\label{exmp:standard-existence}
	The existence proof of the Kazhdan--Lusztig basis elements given in
		\cite[Theorem~1.1]{Kazhdan--Lusztig} amounts to recursively
		computing the Kazhdan--Lusztig polynomials $P_{x,w}$ using 
		$I_w=\set{\rparp{(P_{x,ws})_{x\leq ws}}\emptyset{\set s}:s\in\rdes
		w}$ and descending along the KL-dominance relation.
	Unraveling the definitions, one can check using the explicit formula
		of \cref{prop:invariant_hecke_mult} that this algorithm
		amounts to the standard formula {\cite[2.3b]{Kazhdan--Lusztig}}
		$C'_{ws}=\ds C'_w\cdots C'_s-\sum_{x\in W^{\set s}}
			\mu(xs,w)C'_{xs}$, where $\mu(xs,w)$ is the coefficient of
		$q^{-\frac12}$ in $p_{xs,w}$ for 
		$C'_w=\ds\sum_xp_{x,w}H_x$ (equivalently the coefficient
		of $q^{(\ell(w)-\ell(xs)-1)/2}$ of the Kazhdan--Lusztig polynomial
		$P_{xs,w}$).
\end{exmp}

\begin{exmp}
\label{exmp:scott-sprowl}
	Scott and Sprowl's iterative algorithm for computing Kazhdan--Lusztig
		basis elements amounts to choosing $I_w=M_w=\set{P^J:J$ is a
		singleton factorization induced by a reduced word $s_{i_1}\cdots
		s_{i_r}=w}$, as in our proof of
		\cref{cor:bounded-algorithm}.
	We give our proposed modification to their algorithm in 
		\cref{exmp:new-algo}.
\end{exmp}

\begin{exmp}
	Deodhar's recursive algorithm (or rather, the formulation he alludes
		to in \cite[Remark~4.15]{Deodhar}) amounts to choosing
		$I_w=\set{(P^J_x)_{x\leq w}:J$ is a singleton factorization induced
		by a reduced word $s_{i_1}\cdots s_{i_r}=w}$ and descending along
		the KL-dominance relation.
	Note that this algorithm does \emph{not} depend on the non-negativity
		conjecture; rather the non-negativity conjecture (and in our more
		general setting, Soergel's conjecture) is used to extract a
		combinatorial interpretation of the Kazhdan--Lusztig polynomials
		from this algorithm, which we do below.
\end{exmp}

\begin{defn}
\label{defn:admissible}
	Given a factorization with overlaps $J_1,\dots,J_r$ of $w$, for each
		set $A$ of (equivalence classes of) masks on $J$ we define the
		family of \define{$A$-defect polynomials} $(P^J_{A,x})_{x\in W}$ by
		$P^J_{A,x}=\ds\sum_{[\sigma]\to_Jx\in A}q^{d_R(\sigma)}$.
	We say that a set $A$ of (equivalence classes of) masks on the
		factorization $J$ is \define{admissible} if $(P^J_{A,x})_{x\leq w}$
		is monic and $\underline h^{P^J_A}$ is involution-invariant,
		and say it is \define{KL-dominated} if $\underline h^{P^J_A}$ is
		either equal to or below $\underline h^{P^J}$ in the partial order
		generated by the KL-dominance relation.
\end{defn}

\begin{exmp}
	A factorization $J$ being admissible, i.e.\ $(P^J_x)_{x\leq w}$, 
	    being monic, is equivalent to the set of all (equivalence 
	    classes of) masks being admissible.
\end{exmp}

\begin{rmrk}
\label{rmrk:admissible-masks}
	By \cref{prop:dominance-dcc}(1), if a factorization with
		overlaps $J$ is admissible, then every KL-dominated set of
		(equivalence classes of) masks is admissible, and conversely
		if there exists an admissible KL-dominated set $A$ of
		(equivalence classes of) masks on $J$, then $J$ is itself
		admissible.
\end{rmrk}

\begin{rmrk}
	Deodhar calls \cite[Definition~4.2]{Deodhar}
		a set of masks on a factorization $J_1,\dots,J_r$
		induced by a reduced word $s_{i_1}\cdots s_{i_r}=w$
		\emph{admissible} if $(P^J_{A,x})_{x\leq w}$ is a monic family,
		$\underline h^{P^J_A}$ is involution-invariant, and
		$s_{i_r}\in\rdes{P^J_A}$.
	Since the equivalence relation on masks is trivial for factorizations
		induced by reduced words (see \cref{exmp:monic-singleton}),
		it follows that admissible sets in Deodhar's sense are also
		admissible in our sense.

	Furthermore, \cref{cor:invariance} below will show that any
		admissible (in our sense) KL-dominated set $A$ satisfies
		$\rdes{P^J}\subseteq\rdes{P^J_A}$, and hence 
		$s_{i_r}\in\rdes J\subseteq\rdes{P^J}\subseteq\rdes{P^J_A}$ in the
		case of a factorization associated to a reduced word.
	Thus, ours and Deodhar's notions of admissible KL-dominated
		sets coincide for factorizations associated to reduced words, which
		is the hypothesis for Deodhar's main results
		\cite[Proposition~4.7---Theorem~4.12]{Deodhar}.
\end{rmrk}

Before we prove our main theorem we generalize
	\cite[Proposition~4.3]{Deodhar} which asserts positivity of
	$(P^J_{A,x})_{x\in W}$ for any admissible set of masks $A$ on
	factorization $J$ induced by a reduced word.
\begin{prop}
\label{prop:admissible_subset_positive}
	Given a factorization with overlaps $J_1,\dots,J_r$ and a subset $A$
		of (equivalence classes of) masks on $J$, if $P^J_{A,w(J)}(0)=1$, 
			$\deg P^J_{A,w(J)}(q)=\deg P^J_{w(J)}(q)$, and
		$\underline h^{P^J_A}$ is involution-invariant, then
		$P^J_{A,x}(0)=1$ whenever $x\leq w(J)$, i.e.\ 
		$(P^J_{A,x})_{x\leq w(J)}$ is positive.
	In particular, $(P^J_{A,x})_{x\leq w(J)}$ is positive for $A$ a 
		KL-dominated admissible set of (equivalence classes of) masks.
	\begin{proof}
		Since $P^J_x(0)=1$ when $x\leq w$ by
			\cref{cor:recursive-defect-polys}, it will suffice to
			show by induction on $\ell(w(J))-\ell(x)$ that
			$P^J_{A,x}(0)=P^J_x(0)$ when $x\leq w$, where the base case
			$\ell(w(J))-\ell(x)=0$ follows because $P^J_{A,w(J)}=1$ by
			assumption.
		For the inductive step, we have that $\underline h^{P^J_A}$
			involution-invariant and $(P^J_{A,x})_{x\leq w(J)}\in\NN[q]$ imply
			by \cref{exmp:normalized-invariance-equation} the
			invariance equation
			$q^{\deg P^J_{A,w(J)}(q)+\ell(w)-\ell(x)}\cdot\inv{P^J_{A,x}}
				=\ds\sum_{y} R_{x,y} P^J_{A,y}$.
		Since \cref{cor:length+defects<=meetings}
			implies that when $x<w$ we have $\deg P^J_{A,x}(q)\leq\deg P^J_x(q)
				<\deg P^J_{w(J)}(q)+\ell(w)-\ell(x)
				=\deg P^J_{A,w(J)}(q)+\ell(w)-\ell(x)$,
			the left-hand side of the invariance equation has no constant
			term, and thus the right-hand side must also have no constant
			term, which implies by
			\cref{prop:R-polynomials}(\ref{item:R-constant})
			that $\ds\sum_{x\leq y\leq w}(-1)^{\ell(y)-\ell(x)}
				P^J_{A,y}(0)=0$.
		But since \cref{cor:quotient}
			implies that $\underline h^{P^J}$ is involution-invariant, 
			we also have
			$\ds\sum_{x\leq y\leq w}(-1)^{\ell(y)-\ell(x)}P^J_y(0)
				=0$, and so since $x<y\leq w(J)$ implies that
			$\ell(w(J))-\ell(y)<\ell(w(J))-\ell(x)$, the inductive
			hypothesis yields $P^J_{A,y}(0)=P^J_y(0)$ for such $y$ and hence
			that $P^J_{A,x}(0)=P^J_x(0)=1$, as desired.
	\end{proof}
\end{prop}

\begin{rmrk}
\label{rmrk:mobius}
	Implicit in the proof of
		\cref{prop:admissible_subset_positive} is a demonstration
		that the M\"obius function \cite[Section~3.7]{Stanley} of the Bruhat
		order is given by $\mu(x,y)=(-1)^{\ell(y)-\ell(x)}$, i.e.\ that
		$\ds\sum_{x\leq y\leq w}(-1)^{\ell(y)-\ell(x)}=0$.
	This derivation is actually quite short and follows from Deodhar's
		original framework since it relies only on
		\cref{prop:R-polynomials}(\ref{item:invariance}),
		which is essentially the definition of $R$-polynomials,
		\cref{prop:admissible_subset_positive} above, and
		\cite[Proposition~2.3(iii)]{Deodhar},
		\cite[Proposition~2.3(v)]{Deodhar}, and
		\cite[Proposition~3.5]{Deodhar}, which are the easy-to-prove
		special cases of \cref{cor:unique-reps},
		\cref{cor:length+defects<=meetings}, and
		\cref{cor:quotient} for a singleton factorization induced
		by a reduced word.
	Structurally, this derivation is analogous to (but more
		complicated than) the one given by Stembridge in
		\cite{stembridge_short_2007}, which uses the significantly simpler
		\emph{$0$-Hecke algebra} obtained by setting $q=0$ in the
		$\set{T_x}$-basis presentation of the Hecke algebra.
\end{rmrk}

\begin{lem}
\label{lem:kl-dom-expansion}
	Given KL-dominated sets of (equivalence classes of) masks $A$ and $B$
		on a factorization with overlaps $J_1,\dots,J_r$ of $w\in W$ such
		that $\underline h^{P^J_B}$ is below $\underline h^{P^J_A}$ in the
		partial order generated by the KL-dominance relation, then there are
		expansions 
		$\underline h^{P^J_A}=L_{A,w}C'_w+\ds\sum_{y<w}L_{A,y}C'_y$
		such that
		$\underline h^{P^J_B}=L_{B,w}C'_w+\ds\sum_{y<w}L_{B,y}C'_y$
		with $L_{A,w}=q^{-\deg P^J_w/2}P^J_w=L_{B,w}$, and
		$L_{A,y}\in\ZZ[q^{\frac12}+q^{-\frac12}]$ and 
		$L_{B,y}\in\set{L_{A,y},0}$ for all $y<w$.
	\begin{proof}
		Because $A$ and $B$ are KL-dominated,
			\cref{prop:dominance-dcc} implies that
		$(P^J_{A,x})_{x\leq w}$ and $(P^J_{B,x})_{x\leq w}$ satisfy
		$P^J_{A,w}=P^J_w=P^J_{B,w}$ and that the claimed expansions
		exist and satisfy $L_{A,w}=q^{-\deg P^J_w/2}P^J_w
				=L_{B,w}$, $L_{A,y}=0$ unless 
			$y\in\ordideal{U\qy*{\underline h^{P^J_A}}}$, and $L_{B,y}=0$
			unless $y\in\ordideal{U\qy*{\underline h^{P^J_B}}}$.

		Because $B$ is below $A$ in the partial order generated by the 
			KL-dominance relation, \cref{prop:dominance-dcc} also
			implies that $\underline h^{P^J_A}=\underline h^{P^J_B}
				+\ds\sum_{y<w}L_{\Delta,y}C'_y$ where
				$L_{\Delta,y}=L_{A,y}-L_{B,y}\neq0$ only if 
				$y\in\ordideal{U\qy*{\underline h^{P^J_A}}}
					\setminus\ordideal{U\qy*{\underline h^{P^J_B}}}$.
		But since $L_{B,y}\neq 0$ only if $y\in\ordideal{U\qy*{
			\underline h^{P^J_B}}}$, it follows that $L_{B,y}
				\in\set{L_{A,y},0}$ for all $y$ such that $y<w$ as desired.
	\end{proof}
\end{lem}

\begin{thm}[{\cite[Theorem 2,3]{Williamson}}]
\label{thm:schur}
	Soergel's conjecture (proven in \cite[Theorem 1.1]{Elias--Williamson})
		implies not only that the Kazhdan--Lusztig polynomials satisfy
		$P_{x,w}\in\NN[q^{\frac12},q^{-\frac12}]$, but also that the
		$\NN[q^{\frac12}+q^{-\frac12}]$-submodule of the Hecke algebra
		generated by the Kazhdan--Lusztig basis is closed under the Schur
		composition laws.
\end{thm}

\begin{thm}
\label{thm:1}
	Given a factorization with overlaps $J_1,\dots,J_r$ of $w\in W$, 
		any minimal (with respect to inclusion) admissible KL-dominated set
		$A$ of (equivalence classes of) masks satisfies
		$P^J_{A,x}=P_{x,w}$, i.e.\ is such that $\underline h^{P^J_A}=C'_w$.
	\begin{proof}
		If there exists an admissible KL-dominated set of (equivalence
			classes of) masks, the factorization $J$ of $w$ must itself be
			admissible by \cref{rmrk:admissible-masks}, so that for any
			KL-dominated set of (equivalence classes of) masks we have by
			\cref{lem:kl-dom-expansion} that
			$\underline h^{P^J_A}=C'_w+\ds\sum_{y<w}L_yC'_y$ where
				$L_y\in\ZZ[q^{\frac12}+q^{-\frac12}]$ with $L_y\neq 0$ only if
			$y\in\ordideal{U\qy*{\underline h^{P^J_A}}}$.
		Now if $U\qy*{\underline h^{P^J_A}}=\emptyset$ then
			$\underline h^{P^J_A}=C'_{w}$, so it will suffice to show that if
			$U\qy*{\underline h^{P^J_A}}\neq\emptyset$ then $A$ is not minimal
			admissible KL-dominated.
		Indeed, if $U\qy*{\underline h^{P^J_A}}\neq\emptyset$, then by
			\cref{prop:dominance-dcc} there is $y$ maximal in
			$U\qy*{\underline h^{P^J_A}}$ such that $L_y$ satisfies the condition
			of \cref{prop:dominance-dcc}(\ref{prop:dominance-down}),
			i.e.\ such that $\underline h^{P^J_A}-L_yC'_y\doml_{\mathrm{KL}}
				\underline h^{P^J_A}$, where
			$\underline h^{P^J_A}-L_yC'_y=\ds C'_w+\sum_{y'\neq y}L_yC'_y$.

		Furthermore, since \cref{cor:quotient} implies
			$\underline h^{P^J}$ is in the closure of $\set{C'_y:y\in W}$
			under the Schur composition laws, and since the set of all
			(equivalence classes of) masks on $J$ is KL-dominated,
			\cref{thm:schur} implies that $\underline h^{P^J}
				=\ds\sum_{y\leq w}L_{J,y}C'_y$
			with $L_{J,y}\in\NN[q^{\frac12}+q^{-\frac12}]$ for all $y\in W$.
		Since $\underline h^{P^J_A}$ is below $\underline h^{P^J}$ in the
			partial order generated by the KL-dominance relation,
			\cref{lem:kl-dom-expansion} implies that 
			$L_y\in\set{L_{J,y},0}$ for all $y$, so that
			$L_y\in\NN[q^{\frac12}+q^{-\frac12}]$ as well.

		Thus, if we write $(Q_x)_{x\leq w}$ for the family such that
			$\underline h^Q=\underline h^{P^J_A}-L_yC'_y
				\ds=C'_w+\sum_{y'\neq y}L_{y'}C'_{y'}$,
			we have that $Q_x=\ds P_{x,w}+\sum_{y'\neq y}L_{y'}P_{x,y'}
				\in\NN[q^{\frac12},q^{-\frac12}]$ for all $x$ by
			\cref{thm:schur}, where
			$\NN[q]\ni P^J_{A,x}=Q_x+q^{-(\ell(w)-\ell(y))/2}P_{x,y}$
			with $Q_x,P_{x,y}\in\NN[q^{\frac12},q^{-\frac12}]$, so that
			in fact $Q_x,P^J_{A,x}-Q_x\in\NN[q]$ for every $x\in W$.

		Then since $P^J_{A,x}=\ds\sum_{A\ni[\sigma]\to_Jx}q^{d_R(\sigma)}$
			and $Q_x,P^J_{A,x}-Q_x\in\NN[q]$, there is a subset of masks
			$B\subseteq A$ such that $P^J_{B,x}=Q_x$ for every $x\in W$,
			from which we have $\underline h^{P^J_B}=\underline h^Q
				\doml_{\mathrm{KL}}\underline h^{P^J_A}$.
		In particular, $B$ is KL-dominated, and by
			\cref{rmrk:admissible-masks} also admissible, whence $A$ is
			not minimal admissible KL-dominated as desired.
	\end{proof}
\end{thm}

\begin{cor}
\label{prop:tight-sum-weighted-defect}
	A factorization $J_1,\dots,J_r$ is tight if and only if we have that
		$\ds\sum_x q^{\ell(x)}P_{x,w(J)}(q)
			=\sum_x q^{\ell(x)} P^J_x
			=\sum_x\sum_{[\sigma]\to_Jx}q^{cd_R(\sigma)}$ where
		$cd_R(\sigma)$ is the ``either right-crossing or right-defect''
			statistic given by $cd_R(\sigma)=\ds
				\sum_{t\in T}\abs{\rcross\sigma t+\rdef\sigma t}$.
	\begin{proof}
		That $\ds\sum_x q^{\ell(x)}P^J_x
			=\sum_x\sum_{[\sigma]\to_Jx} q^{cd_R(\sigma)}$ follows
			from the fact that $q^{\ell(x)}P^J_x=\ds
				\sum_{[\sigma]\to_Jx} q^{\ell(x)+d_R(\sigma)}$ by definition,
			and $\ell(x)+d_R(\sigma)=\ds\sum_{t\in T}
				\abs{\rcross\sigma t+\rdef\sigma t}=cd_R(\sigma)$
			by \cref{cor:length+defects<=meetings}.

		If $J$ is tight, then by \cref{prop:bounded-monic}
			we have $P^J_x=P_{x,w(J)}$ so that
			$\ds\sum_xq^{\ell(x)}P^J_x=\sum_xq^{\ell(x)}P_{x,w(J)}$.
		For the converse, since $(P^J_x)_{x\leq w}\in\NN[q]$ we have
			$\deg\ds\sum_xq^{\ell(x)}P^J_x
				\geq\ell(w(J))+\deg P^J_{w(J)}$,
			and since $P_{w(J),w(J)}=1$,
			$\deg P_{x,w(J)}(q)<\dfrac12\qy*{\ell(w(J))-\ell(x)}$ for
			$x<w(J)$, and $P_{x,w(J)}=0$ otherwise, we have that
			$\ds\deg\sum_xq^{\ell(x)}P_{x,w(J)}(q)\leq\ell(w(J))$.
		Thus, $\ds\sum_xq^{\ell(x)}P^J_x
			=\sum_xq^{\ell(x)}P_{x,w(J)}$ only if $\deg P^J_w\leq 0$,
			or equivalently, only if $J$ is admissible.
		But for an admissible factorization, by \cref{thm:1} there
			exists a subset $A$ of (equivalence classes of) masks such that
			$P^J_{A,x}(q) = P_{x,w(J)}(q)$ for every $x\in W$, so that
			$\ds\sum_xq^{\ell(x)}P_{x,w(J)}(q)
				=\sum_xq^{\ell(x)}P^J_{A,x}(q)$.
		Since $P^J_{A,x}(q)\in\NN[q]$ and $P^J_x(q)-P^J_{A,x}(q)\in\NN[q]$,
			it follows that $\ds\sum_x q^{\ell(x)}P_{x,w(J)}(q)
			=\sum_xq^{\ell(x)}P^J_{A,x}=\sum_x q^{\ell(x)} P^J_x$ holds only if
			$P^J_x(q)=P^J_{A,x}(q)$ for all $x$, which by definition of $A$ is
			equivalent to $P^J_x(q)=P_{x,w(J)}(q)$ for all $x$, i.e.\ to $J$
			being a tight factorization. 
	\end{proof}
\end{cor}

\begin{lem}
\label{lem:sum-weighted-defect}
	Given a family $(P_x)_{x\in W}\in\ZZ[q^{\frac12},q^{-\frac12}]$, a
		finitary subset $J\subseteq S$ and a subset $J_*\subseteq\rdes P\cap
		J$, we have $\ds\sum_xq^{\ell(x)}\qy{\rparp P{J_*}J}_x(q)
			=\rcos[J]{J_*}W(q)\sum_xq^{\ell(x)}P_x(q)$. 
	In particular, for any factorization with overlaps $J_1,\dots,J_r$,
		its family of defect polynomials satisfies 
		$\ds\sum_x q^{\ell(x)}P^J_x(q)
			=\prod_{j=1}^r W_{J_j}^{\overlap J*j}(q)
			=\prod_{j=1}^r W_{J_j}(q)\qy*{\prod_{1\leq i<k\leq r}
				W_{\overlap Jik}(q)}^{-1}$.
	\begin{proof}
		The main claim follows from a straightforward computation:
	\begin{align*}
			\sum_x q^{\ell(x)}\qy{\rparp P{J_*}J}_x(q)
				&=\sum_x q^{\ell(x)}\sum_{y\in W_J^{J_*}}q^{\ell(y)}P_{x^Jy}(q)\\
				&=\sum_{x\in W^J}\sum_{y'\in W_J^{J_*}}\sum_{z\in
					W_{J_*}}q^{\ell(x)}q^{\ell(y')}q^{\ell(z)}\sum_{y\in
					W_J^{J_*}}q^{\ell(y)}P_{xy}(q)\\
				&=\sum_{y'\in W_J^{J_*}}q^{\ell(y')}
					\sum_{x\in W^J}\sum_{y\in W_J^{J_*}}\sum_{z\in W_{J_*}}
						q^{\ell(xyz)}P_{xyz}(q)\\
				&=W_J^{J_*}(q)
					\sum_x q^{\ell(x)}P_x(q).
		\end{align*}
		To verify the supplemental claim,  we show by induction that it
			holds for the slices $\upslice Jk$ whenever $k$ is such that 
			$1\leq k\leq r$.
		The base case $k=1$ is immediate since the desired equality reduces
			to $\ds\sum_{x\in W_{J_1}}q^{\ell(x)}=W_{J_1}(q)$.
		For the inductive step, we apply the main claim to the fact due to
			\cref{cor:recursive-defect-polys} that $P^{\upslice Jk}
				=\rparp{P^{\upslice J{k-1}}}{\overlap J*k}{J_k}$ to obtain 
			$\ds\sum_x q^{\ell(x)}P^{\upslice Jk}_x
				=\rcos[J_k]{\overlap J*k}W(q)
				\cdot\prod_{j=1}^{r-1}\rcos[J_j]{\overlap J*j}W(q)
			=\prod_{j=1}^k\rcos[J_j]{\overlap J*j}W(q)$.
		Since $\overlap J*j
			=\ds\bigudj_{1\leq i<k}\overlap Jik$ with $\overlap
			Jik$ commuting by definition, \cref{lem:bijections} implies we
			also have $\ds\prod_{j=1}^r\rcos[J_j]{\overlap J*j}W(q)
				=\prod_{j=1}^r W_{J_j}(q)
				\qy*{\prod_{1\leq i<k\leq r} W_{J_i^k}(q)}^{-1}$. 
	\end{proof}
\end{lem}

\begin{cor}
\label{cor:schubert}
	A factorization with overlaps $J_1,\dots,J_r$ of an element $w$ of a
		Weyl group is tight if and only if the Poincar\'e polynomial
		$\ds\sum_i(\dim\mathrm{IH}^{2i}(X_w))q^i$ of the
		intersection cohomology of the Schubert variety $X_w$ factors as
		$\ds\prod_{j=1}^rW_{J_j}^{\overlap J*j}(q)=
			\prod_{j=1}^rW_{J_j}(q)\qy*{\prod_{1\leq i<k\leq r}
				W_{J_i^k}(q)}^{-1}$, or
		equivalently, is given by 
	$\ds\sum_x\sum_{[\sigma]\to_Jx}q^{cd_R(\sigma)}$.
	In particular, in the case of factorizations with Gelfand--MacPherson
		overlaps, such factorizations of the Poincar\'e polynomial
		correspond to small Gelfand--MacPherson resolutions.
	\begin{proof}
	\cite[Corollary~4.9]{Kazhdan--Lusztig80}
		asserts that $\ds\sum_i(\dim\mathrm{IH}^{2i}(X_w))q^i
			=\sum_x q^{\ell(x)}P_{x,w}(q)$ while
		Lemma~\ref{lem:sum-weighted-defect} asserts
		$\ds\sum_x^{\ell(x)}q^{\ell(x)}P^J_x(q)
			=\prod_{j=1}^r W_{J_j}(q)\qy*{\prod_{1\leq i<k\leq r}
				W_{\overlap Jik}(q)}^{-1}$.
	\end{proof}
\end{cor}

Finally, we show that our notion of admissible sets of (equivalence
	classes of) masks coincides with Deodhar's in the case of KL-dominated
	admissible sets.

\begin{lem}
\label{lem:invariance}
	Given a Hecke algebra element $h=\ds\sum_xP_xT_x
			=\sum _xp_xH_x$, we have for every $s\in\rdes P$ and $y\in U(h)$
		maximal in the Bruhat order on $U(h)$ that either $s\in\rdes y$ or
		$w=ys$ is maximal in the Bruhat order such that $P_w\neq 0$ and
		satisfies $\deg p_w(q)\geq\frac12$.

	Furthermore, if $\deg p_w\leq0$ for every $w$ maximal such
		that $P_w\neq 0$, then given $s\in\rdes P$ and $\set{M_y}_{y\in W}$
		as in \cref{defn:dominance} such that $s\in\rdes y$
		implies $s\in\rdes Q$ for every $Q\in M_y$,
		we have that $h'=\ds\sum_x P'_xT_x$ is below $h$
		in the partial order generated by the $M$-dominance relation only
		if $s\in\rdes{P'}$.
	\begin{proof}
		Let $p_x=q^{\ell(x)/2}P_x$ so that $h=\ds\sum_x p_xH_x$,
			let $y\in U(h)$ be maximal in the Bruhat order on $U(h)$, and let
			$s\in\rdes P\setminus\rdes y$.
		Then $ys>y$, so $ys\not\in U(h)$ and $\ell(ys)-\ell(y)=1$, and hence
			$\deg p_{ys}(q)=\deg q^{\ell(ys)/2}P_{ys}(q)
				=\deg q^{\ell(ys)/2}P_y(q)=\frac12\qy*{\ell(ys)-\ell(y)}
					+\deg p_y(q)\geq\frac12+0=\frac12$.
		Since $ys\not\in U(h)$ we have by definition that $w=ys$ is maximal
			in the Bruhat order such that $P_w\neq 0$.
		
		For the supplemental claim, since \cref{prop:dominance-dcc}
			implies that $P_w=P'_w$ for every $w$ maximal with $P_w\neq 0$,
			it suffices to check the implication for the generating relations of
			$\doml_M$ by induction.
		If $h'=h-L_y\underline h^{Q}$ for some $y\in U(h)$ maximal in the
			Bruhat order on $U(h)$ and $Q\in M_y$, then by the main claim we have
			that $s\in\rdes y$ and so by assumption on $M_y$ we also have
			$s\in \rdes Q$. 
		Thus for each $x\in W$ we have 
			$P'_{xs}=P_{xs}-L_y\cdot q^{-\ell(y)/2}Q_{xs}
				=P_x-L_y\cdot q^{-\ell(y)/2}Q_x=P'_x$ and thus $s\in\rdes{P'}$
			as desired.
	\end{proof}
\end{lem}

\begin{prop}[{\cite[Equation~2.3.g]{Kazhdan--Lusztig}}]
\label{prop:invariant_kl}
	For every $w\in W$, the family of Kazhdan--Lusztig polynomials
		$(P_{x,w})_{x\leq w}$ satisfies $\rdes{P_{-,w}}=\rdes w$.
	\begin{proof}
		That $\rdes{P_{-,w}}\subseteq \rdes w$ follows from
			\cref{lem:pre-descent}, so it suffices to show that
			$s\in \rdes w$ implies $s\in \rdes{P_{-,w}}$.

		Given $s\in \rdes w$, define the sets $\set{M_y}_{y\in W}$ by
			$M_y=\set{(P^J_x)_{x\leq y}:J$ is induced by a reduced word
				$s_{i_1}\cdots s_{i_r}=y$ such that $s_{i_r}=s$ if
				$s\in\rdes y}$, and let $P\in M_w$ so that $s\in\rdes P$.
		Then since $P_w=1$ by \cref{exmp:monic-singleton},
			$s$ and $\set{M_y}_{y\in W}$ satisfy the conditions of
			the supplemental claim of \cref{lem:invariance} and thus
			since \cref{prop:dominance-dcc} implies
			$C'_w=\underline h^{P_{-,w}}$ is below $\underline h^P$ in the
			partial order generated by the $M$-dominance relations, it follows
			from \cref{lem:invariance} that $s\in \rdes{P_{-,w}}$ as
			desired.
	\end{proof}
\end{prop}

\begin{cor}
\label{cor:invariance}
	For any family $(P_x)_{x\in W}\in\ZZ[q^{\frac12},q^{-\frac12}]$, 
		if $\ds\sum_xP_xT_x
			=\sum_yL_yC'_y$, then $\rdes P=\bigcap\set{\rdes y:L_y\neq0}$.
	In particular, for any KL-dominated set $A$ of 
		(equivalence classes of masks) on a factorization $J$ we have 
		$\rdes{P^J}\subseteq\rdes{P^J_A}$.
	\begin{proof}
		We argue by induction on $\abs{\set{y\in W:L_y\neq0}}$ that
			$s\in\rdes P$ if and only if $s\in\bigcap\set{\rdes y:L_y\neq0}$.
		The base case of $\abs{\set{y\in W:L_y\neq0}}=0$ is immediate since
			$(L_y)_{y\in W}=0$ implies $(P_x)_{x\in W}=0$, so that $\rdes P=S
				=\bigcap\set{\rdes{y}:L_y\neq0}$.

		For the inductive step, let $y$ be maximal in the Bruhat order
			such that $L_y\neq 0$, and let $(Q_x)_{x\in W}$ be such
			that $\ds\sum_x Q_xT_x = -L_yC'_y + \sum_x P_xT_x
				= \sum_{y'\neq y} L_{y'}C'_{y'}$. 
		Then by the inductive hypothesis we have $\rdes Q =\ds
			\bigcap_{y'\neq y:L_{y'}\neq 0}\rdes{y'}$ and we need to show that
			$\rdes P = \rdes Q \cap \rdes y
			= \rdes Q \cap \rdes{(P_{x,y})_{x\leq y}}$ where the
			last equality holds by \cref{prop:invariant_kl}.

		For each $x\in W$ we have $P_x
			=\ds \sum_z q^{-\ell(z)/2}L_z P_{x,z}
			= Q_x + q^{-\ell(y)/2}L_yP_{x,y}$.
		Thus, $s\in \rdes Q \cap \rdes{(P_{x,y})_{x\leq y}}$
			implies $P_{xs}
				=Q_{xs} + q^{-\ell(y)/2}L_yP_{xs,y}
				=Q_x + q^{-\ell(y)/2}L_yP_{x,y}=P_x$
			and hence $\rdes P\supseteq\rdes Q\cap\rdes{(P_{x,y})_{x\leq y}}$.
		For the reverse inclusion, note first that since $y$ is maximal
			in the Bruhat order with $L_y\neq 0$, we also have that $y$ is
			maximal in the Bruhat order such that $P_y\neq 0$, and thus by
			\cref{lem:pre-descent} and \cref{prop:invariant_kl}
			we have that $\rdes P\subseteq\rdes y=\rdes{(P_{-,y}}$.
		Then $s\in\rdes P$ implies that for each $x\in W$ we have
			$Q_{xs}=P_{xs}-q^{\ell(y)/2}L_yP_{xs,y}
				=P_x-q^{\ell(y)/2}L_yP_{x,y}=Q_x$ so that
			$\rdes P \subseteq \rdes Q\cap\rdes{P_{-,y}}$, and
			hence $\rdes P=\rdes Q\cap\rdes y$ as desired.

		The supplemental claim then follows immediately from
			\cref{lem:kl-dom-expansion}, which implies that if
			$\underline h^{P^J}=\ds \sum_yL_yC'_y$ and 
			$\underline h^{P^J_A}=\ds \sum_yL'_yC'_y$ then
			$L'_y\neq0$ only if $L_y\neq0$.
	\end{proof}
\end{cor}

\section{Schur composites of Kazhdan--Lusztig basis elements}

\label{sec:schur}
We derive the connection of our notion of tight factorization with
	Skandera's question: that a Kazhdan--Lusztig basis element $C'_w$
	belongs to the closure of $\set{C'_{\mxml J}:J$ is finitary$}$ under
	the Schur composition laws if and only if $w$ has a tight
	factorization.
Our argument will be inductive with an inductive step requiring a more
	precise version of the result, which we record as \cref{thm:2}.
To that end, we introduce the notions of one factorization being a
	\emph{contraction} of another, of \emph{minimal factorizations},  and
	of a factorization being \emph{right-descent}, \emph{left-descent},
	and \emph{bidescent}, and show that every tight factorization has a
	tight right-descent contraction, and by duality also a tight
	left-descent contraction, and ultimately a tight absolutely bidescent
	contraction.

We begin by using Soergel's conjecture and Williamson's results on
	singular Soergel bimodules to show that Kazhdan--Lusztig basis
	elements are closed under ``decomposition'' under the Schur
	composition laws.
\begin{prop}
\label{prop:boundedness}
	Given Hecke algebra elements 
		$\ds h_1=\sum_xP_xT_x=\ds\sum_xL_xC'_x$ and 
		$\ds h_2=\sum_yQ_yT_y=\ds\sum_yM_yC'_y$ where
		$(L_x)_{x\in W},(M_y)_{y\in W}\in\NN[q^{\frac12},q^{-\frac12}]$,
		for any $J_*\subseteq\rdes P\cap\ldes Q$ their Schur composition
		$h_2\circ_{J_*}h_1
			=\qy*{q^{-\ell(\mxml{J_*})/2}W_{J_*}(q)}^{-1} h_1\cdot h_2$ is a
		Kazhdan--Lusztig basis element only if both $h_1$ and $h_2$ are 
		Kazhdan--Lusztig basis elements.
	\begin{proof}
		By \cref{prop:invariant_kl} and
			\cref{cor:invariance}, we have that
			$\rdes P=\bigcap\set{\rdes x:L_x\neq0}
				=\bigcap\set{\rdes{(P_{z,x})_{z\in W}}:L_x\neq0}$ and
			$\ldes Q=\bigcap\set{\ldes y:M_y\neq0}
				=\bigcap\set{\ldes{(Q_{z,y})_{z\in W}}:M_y\neq0}$, so that
			$J_*\subseteq\rdes{(P_{z,x})_{z\in W}}\cap\ldes{(Q_{z,y})_{z\in W}}$
			whenever $L_x\neq0\neq M_y$.
		Consequently, 
		\begin{align*}
			\qy*{q^{-\ell(\mxml{J_*})/2}W_{J_*}(q)}^{-1} h_1\cdot h_2
					=&\sum_{x,y}L_xM_y\cdot\qy*{q^{-\ell(\mxml{J_*})/2}
						W_{J_*}(q)}^{-1}C'_x\cdot C'_y\\
					=&\sum_{x,y}L_xM_y\cdot C'_y\circ_{J_*}C'_x
					=\sum_{x,y,z}L_xM_yN_{x,y,z}C'_z
		\end{align*}
			where by \cref{thm:schur} we have 
				$(N_{x,y,z})_{z\in W}\in\NN[q^{\frac12}+q^{-\frac12}]$ such that
				$C'_y\circ_{J_*} C'_x=\ds\sum_zN_{x,y,z}C'_z$,
				and hence $\ds\sum_zN_{x,y,z}(1)\geq1$ because 
				$C'_y\circ_{J_*} C'_x\neq0$.
		Now, $L_x,M_y,N_{x,y,z}\in\NN[q^{\frac12},q^{-\frac12}]$ imply
			$L_xM_yN_{x,y,z}\neq0$ only if $L_x(1),M_y(1),N_{x,y,z}(1)\geq1$.
		Thus $\ds\sum_{x,y,z} L_xM_y N_{x,y,z}C'_z=C'_{z'}$
			for $z'\in W$ only if $1=\ds
				\sum_{x,y,z}L_x(1)M_y(1)N_{x,y,z}(1)\allowbreak
				\geq\abs{\set{x:L_x\neq 0}}
					\cdot\abs{\set{y:L_y\neq 0}}$. 
		Hence, $L_x\neq0\neq M_y$ for unique $x$ and $y$, in which case
			$C'_{z'}=L_xM_y\ds\sum_z N_{x,y,z}C'_z$, and so
			$L_x=1=M_y$.
	\end{proof}
\end{prop}
Next, we identify a hypothesis on factorizations that guarantees that a
	Schur composite of the normalized Hecke algebra elements of their
	families of defect polynomials is itself a normalized Hecke algebra
	element $\underline h^{P^J}$ for some factorization~$J$. 
\begin{lem}
\label{lem:defect-ind}
	Let $J_1,\dots,J_r$ and $J'_1,\dots,J'_{r'}$ be factorizations such
		that $\rdes J=\rdes{P^J}$ and $\ldes{J'}=\ldes{P^{J'}}$.
	Then given a subset $J_*\subseteq\rdes J\cap\ldes{J'}$ we have
	$\underline h^{P^{J'}}\circ_{J_*}\underline h^{P^J}
		=\underline h^P$ for a monic family $(P_x)_{x\leq w}$
		only if $\underline h^P=\underline h^{P^{J,J'}}$ where
		$J,J'$ is the factorization $J_1,\dots,J_r,J'_1,\dots,J'_{r'}$,
		in which case $J,J'$ is an admissible factorization (of $w$).
	\begin{proof}
		By \cref{cor:quotient} we have
			\[
				\underline h^{P^{J'}}\circ_{J_*} \underline h^{P^J}
					=q^{-\len{\rcos{J_*}{\mxml{\rdes J\cap\ldes{J'}}}}/2}
					\rcos[\rdes J\cap\ldes{J'}]{J_*}W(q)\underline h^{P^{J,J'}}
			\]
			so that $w$ is maximal in the Bruhat order such that
			$P^{J,J'}_w\neq 0$ and thus $J,J'$ is a factorization of $w$.
		Then by definition of normalized element, we have
			\[q^{-\len{\rcos{J_*}{\mxml{\rdes J\cap\ldes{J'}}}}/2}
				\rcos[\rdes J\cap\ldes{J'}]{J_*}W(q)
				q^{-\deg P^{J,J'}_w(q)/2}
				P^{J,J'}_w=P_w=1\]
			where since $(P^{J,J'}_x)_{x\leq w}\in\NN[q]$ is positive
			and $\rcos[\rdes J\cap\ldes{J'}]{J_*}W(q) \in\NN[q]$ has
			constant term $1$, the left-hand side has a term
			$q^{-\len{\rcos{J_*}{\mxml{\rdes J\cap\ldes{J'}}}}/2
				-\deg P^{J,J'}_w/2}$ while the right-side is $1$,
			implying that $\rcos{J_*}{\mxml{\rdes J\cap\ldes{J'}}}=e$,
			i.e.\ $\rdes J\cap\ldes{J'}\subseteq J_*$ and
			$q^{-\len{\rcos{J_*}{\mxml{\rdes J\cap\ldes{J'}}}}/2}
				\rcos[\rdes J\cap\ldes{J'}]{J_*}W(q)=1$ so that
				$\underline h^{P^{J'}}\circ_{J_*} \underline h^{P^J}
					=\underline h^{P^{J,J'}}$, and
			$\deg P^{J,J'}_w=0$ so that $J,J'$ is admissible
			as desired.
	\end{proof}
\end{lem}

If every tight factorization $J_1,\dots,J_r$ satisfied $\rdes
	J=\rdes{P^J}$ and $\ldes J=\ldes{P^J}$, then combining
	\cref{prop:boundedness} with \cref{lem:defect-ind}
	would complete the inductive step of our inductive argument. 
This is not the case, however, as can be seen by the type $A$
	factorization $J_1=\set{s_1,s_2},J_2=\set{s_1}$ which has
	$(P^J_x)_{x\leq w(J)}=(P^{J_1}_x)_{x\leq \mxml{J_1}}
		=(P_{x,\mxml{J_1}})_{x\leq \mxml{J_1}}$ with
	$\rdes{P^J}=\rdes{P^{J_1}}=\rdes{J_1}=\set{s_1,s_2}
		\supsetneq\set{s_1}=\rdes{J}$,
	but we will show that any tight factorization $J$ can be
	transformed to another tight factorization $J'$ of $w(J)$ such
	that $\rdes{J'}=\rdes{P^{J'}}=\rdes{P^J}$.
To that end, we introduce the \emph{contraction} partial order on
	factorizations. 

\begin{defn}
\label{defn:contraction}
	We say that a factorization $J'_1,\dots,J'_r$ is a
		\define{right-contraction} $J_1,\dots,J_r$ if each $J'_k$ is a union
		of connected components of $J_k$ such that $J_k\setminus
		J'_k\subseteq \rdes{\upslice{J'}{k-1}}$.
	Similarly, we say that a factorization $J'_1,\dots,J'_r$ is a 
		\define{left-contraction} of $J_1,\dots,J_r$ if each $J'_k$ is a
		union of connected components of $J_k$ such that $J_k\setminus
		J'_k\subseteq \ldes{\downslice{J'}{k+1}}$.
	We recursively define a \define{contraction} of a factorization $J$ to
		be either $J$ itself or a left- or right-contraction of a
		contraction of $J$.
	We say that a factorization is \define{minimal} if it has no
		contractions other than itself.
\end{defn}

\begin{rmrk}
\label{rmrk:contraction-absorption}
	In the case of type $A$, contraction can be visualized as a kind of
		\emph{absorption} in which an interval $[a,b]$ corresponding to a
		connected component of some $J_k$ is contained entirely in an
		overlap interval $\overlap Jik$ or $\overlap Jki$, and can be
		``absorbed'' into $J_i$.
	For example, in \cref{fig:overlap}, the interval $[7,8]$ in $J_1$
		can be absorbed in to $J_2$, and correspondingly
		$J_1\setminus\set{s_7},J_2,J_3$ is a (left-)contraction of
		$J_1,J_2,J_3$.
\end{rmrk}

In order to strengthen our ultimate conclusion, we now show the property
	of being a minimal factorization is absolute.
\begin{lem}
\label{lem:slice-contraction}
	Given a factorization $J_1,\dots,J_r$, if $J'_i,\dots,J'_j$ is a
		contraction of $\slice Jij$, then
		$J'=J_1,\dots,J_{i-1},J'_i,\dots,J'_j,J_{j+1},\dots,J_r$ is a
		contraction of $J$.
	\begin{proof}
		By duality and the definition of contraction, it suffices to show
			that $J'_i,\dots,J'_j$ being a right-contraction of $\slice Jij$
			implies $J'$ is a right-contraction of $J_1,\dots,J_r$.
		In that case we have for $k\not\in[i,j]$ that $J_k\setminus J'_k
			=\emptyset\subseteq\rdes{\upslice{J'}{k-1}}$,
			while for $k\in[i,j]$ we have by assumption that
			$J_k\setminus J'_k\subseteq \rdes{\slice{J'}i{k-1}}
				=\ds\bigudj_{j=i}^{k-1} \rdeslev j{\downslice {J'}i}
				=\bigudj_{j=i}^{k-1} \rdeslev j{J'}
				\subseteq \bigudj_{j=1}^{k-1} \rdeslev j{J'}
				=\rdes{\upslice {J'}{k-1}}$.
		Thus, $J'$ is a right-contraction of $J$ as desired.
	\end{proof}
\end{lem}

\begin{prop}
\label{prop:contraction-cartier--foata}
	Given a factorization $J_1,\dots,J_r$ for which $J_k$ and $J_{k+1}$
		commute, then $J'_1,\dots,J'_r$ is a contraction of $J_1,\dots,J_r$
		if and only if
		$J'_1,\dots,J'_{k-1},J'_k\udj J'_{k+1},J'_{k+2},\dots,J'_r$
		is a contraction of
		$J_1,\dots,J_{k-1},J_k\udj J_{k+1},J_{k+2},\dots,J_r$.
	\begin{proof}
		It will suffice to show for a factorization $J_1,\dots,J_r$
			such that $J_k$ and $J_{k+1}$ commute,(i.e.\ such that 
			$J_k\cap\nbd{J_{k+1}}=\emptyset=\nbd{J_k}\cap J_{k+1}$, that
			$J'$ is a right-contraction of $J$ if and
			only if $J'_1,\dots,J'_k\udj J'_{k+1},\dots,J'_r$ is a
			right-contraction of $J_1,\dots,J_k\udj J_{k+1},\dots,J_r$.
		
		Given a pair of inclusions $J'_k\subseteq J_k$ and
			$J'_{k+1}\subseteq J_{k+1}$, we have by definition that
			$\nbd{J'_k\cup J'_{k+1}}=\nbd{J'_k}\cup\nbd{J'_{k+1}}$,
			$\nbd{J'_k}\subseteq\nbd{J_k}$, and
			$\nbd{J'_{k+1}}\subseteq\nbd{J_{k+1}}$.
		If $J_k$ and $J_{k+1}$ commute, then 
		$\nbd{J'_k\udj J'_{k+1}}\cap\qy{J_k\udj J_{k+1}}
			=\qy*{\nbd{J'_k}\cup\nbd{J'_{k+1}}}\cap\qy*{J_k\udj J_{k+1}}
			=\qy*{\nbd{J'_k}\cap J_k}\udj\qy*{\nbd{J'_{k+1}}\cap J_{k+1}}$.
		It follows that, $J'_k\udj J'_{k+1}$ is a union of connected
			components of $J_k\udj J_{k+1}$, i.e.\ that $J'_k\udj J'_{k+1}
				=\nbd{J'_k\udj J'_{k+1}}\cap\qy*{J_k\udj J_{k+1}}$, if and only
			if $J'_k=\nbd{J'_k}\cap J_k$ and 
			$J'_{k+1}=J'_{k+1}\cap\nbd{J_{k+1}}$, i.e.\ if and only if $J'_k$
			and $J'_{k+1}$ are unions of connected components of $J_k$ and
			$J_{k+1}$ respectively. 

		Furthermore, since $J_k$ and $J_{k+1}$ commuting imply $J'_k$ and
			$J'_{k+1}$ commute, we have that
			$\rdes{J'_1,\dots,J'_k,\dots,J'_j}=\rdes{J'_1,\dots,J'_k\udj
			J'_{k+1},\dots,J'_j}$ for any $j$ such that $j\not\in[k,k+1]$, so
			that the right-contraction conditions at index $j\not\in[k,k+1]$
			coincide for these two factorizations. 

		To verify the right-contraction conditions coincide 
			in the case where $j\in[k,k+1]$, note first that
			$\rdups{J'}k=J'_k\udj\qy*{\rdups{J'}{k-1}}\setminus\nbd{J'_k}$,
			so that $J_k$ and $J_{k+1}$ commuting imply
			$J_{k+1}\setminus J'_{k+1}\subseteq\rdups{J'}k$ if and only if
			$J_{k+1}\setminus J'_{k+1}\subseteq\rdups{J'}{k-1}$.
		Then $J_{k+1}\setminus J'_{k+1}\subseteq\rdups{J'}k$
			and $J_k\setminus J'_k\subseteq\rdups{J'}{k-1}$ if and only if
			$\rdups{J'}{k-1}\supseteq\qy{J_{k}\setminus J'_k}\udj\qy{J_{k+1}
					\setminus J'_{k+1}}
				=\qy{J_k\udj J_{k+1}}\setminus\qy{J'_k\udj J'_{k+1}}$,
			as desired.
	\end{proof}
\end{prop}

\begin{cor}
\label{cor:absolutely-minimal}
	Every minimal contraction of a factorization $J_1,\dots,J_r$ is
		absolutely minimal. 
	\begin{proof}
		Given a minimal factorization $J_1,\dots,J_r$, 
			\cref{prop:contraction-cartier--foata} implies any equivalent
			factorization is also minimal, and \cref{lem:slice-contraction}
			implies every slice of every equivalent factorization, i.e.\ every
			factor, is minimal.
	\end{proof}
\end{cor}

We now continue with the application of the partial order of
	contractions to transforming a tight factorization $J$ into $J'$ such
	that $\rdes{J'}=\rdes{P^{J'}}=\rdes{P^J}$.
	
\begin{lem}
\label{lem:contraction-increases-descents}
	If $J'_1,\dots,J'_r$ is a right-contraction of $J_1,\dots,J_r$,
		then $\rdes{\upslice Jk}\subseteq \rdes{\upslice {J'}k}$ for every
		$k$ such that $1\leq k\leq r$.
	\begin{proof}
		We proceed by induction on the length $r$ of the factorization
			$J_1,\dots,J_r$.
		The base case $r=1$ is immediate since $J'$ being a
			right-contraction of $J$ means by definition that
			$J_1\setminus J'_1 \subseteq\rdes{\upslice{J'}0}=\emptyset$ 
			and hence that $J_1=J'_1$.
		For the inductive step, $J'$ being a right-contraction of $J$
			implies that that $\upslice{J'}{r-1}$ is a right-contraction of
			$\upslice J{r-1}$ and so the inductive hypothesis implies
			$\rdes{\upslice Jk}\subseteq\rdes{\upslice {J'}k}$ for every
			$k$ such that $1\leq k\leq r-1$.
		It remains to show that $\rdes J\subseteq \rdes{J'}$.

		On the one hand, since, $J'_r\subseteq J_r$ implies
			$\nbd{J'_r}\subseteq\nbd{J_r}$ and since $\rdes{\upslice J{r-1}}
				\subseteq\rdes{\upslice{J'}{r-1}}$, 
			we have $\rdes J=J_r\udj\qy*{\rdes{\upslice J{r-1}}
					\setminus \nbd{J_r}}
				\subseteq J_r\cup\qy*{\rdes{\upslice{J'}{r-1}}
					\setminus\nbd{J'_r}}$.

		On the other hand, $J'$ being a right-contraction of $J$ by
			definition implies that $J_r\setminus J'_r
				\subseteq\rdes{\upslice{J'}{r-1}}$ and that $J'_r$ is a union of
			connected components of $J_r$, i.e.\ and that $\nbd{J'_r}\cap J_r
				=J'_r$, which is equivalent to 
			$\qy{J_r\setminus J'_r}\cap\nbd{J'_r}=\emptyset$. 
		Thus $J_r\setminus J'_r
			\subseteq\rdes{\upslice{J'}{r-1}}\setminus\nbd{J'_r}$ and hence
		 $J_r\cup\qy*{\rdes{\upslice{J'}{r-1}}\setminus\nbd{J'_r}}
				=\qy*{J'_r\udj\qy{J_r\setminus J'_r}}
					\cup\qy*{\rdes{\upslice{J'}{r-1}}\setminus\nbd{J'_r}}
				=J'_r\udj\qy*{\rdes{\upslice{J'}{r-1}}\setminus\nbd{J'_r}}
				=\rdes{J'}$. 

		Thus, $\rdes J
			\subseteq J_r\udj\qy*{\rdes{\upslice J{r-1}}\setminus\nbd{J_r}}
			\subseteq\rdes{J'}$, completing the induction.
	\end{proof}
\end{lem}

\begin{lem}
\label{lem:contraction}
	Given a family $(P_x)_{x\leq w}\in\ZZ[q^{\frac12},q^{-\frac12}]$, a
		finitary subset $J\subseteq S$, and a subset $J_*\subseteq
		J\cap\rdes{P}$, suppose that $J'\subseteq J$ is a union of
		connected components of $J$ such that $J\setminus
		J'\subseteq\rdes{P}$.
	Then
		$\rparp P {J_*} J=\rcos[J\setminus J']{J_*}{W}(q)
			\rcos[J'_*]{J_*}{W}(q) \rparp P {J'_*} {J'}$ for any $J'_*$
			such that $J_*\cap J'\subseteq J'_*\subseteq \rdes{P}\cap J'$. 
	\begin{proof}
		Since $J'$ is a  union of connected components of
			$J$, by \cref{lem:bijections} gives 
			$\ds\sum_{y\in W_J}q^{\ell(y)}T_y
				=\qy*{\sum_{y\in W_{J\setminus J'}}q^{\ell(y)}T_y}
					\cdot\qy*{\sum_{z\in W_{J'}}q^{\ell(z)}T_z}$.
		Since every element of $J_*\setminus J'$ commutes
			with every element of $J_*\cap J'$, \cref{lem:bijections}
			also gives $W_{J_*\setminus J'}(q)W_{J_*\cap J'}(q)=W_{J_*}(q)$.

		\cref{lem:coset-properties}(\ref{lem:coset-unique})
			gives us bijections $W_{J\setminus J'}\cong W_{J\setminus
			J'}^{J_*}\times W_{J_*\cap (J\setminus J')}
			\cong W_{J\setminus J'}^{J_*}\times W_{J_*\setminus J'}$
			and $W_{J'_*}\cong W_{J'_*}^{J_*}\times W_{J'_*\cap J_*}$, which
			by \cref{lem:bijections} yield equalities
			$W_{J\setminus J'}(q)=W_{J\setminus J'}^{J_*}(q)
			W_{J_*\setminus J'}(q)$ and $W_{J'_*}(q)
				=W_{J'_*}^{J_*}(q)W_{J'_*\cap J_*}(q)$.
		Since $J'_*\cap J_* = J_*\cap J'$, we have equalities
			$W_{J'_*\cap J_*}= W_{J_*\cap J'}$ and 
			$W_{J'_*\cap J_*}(q)=W_{J_*\cap J'}(q)$.

		By \cref{prop:hecke_mult} we then have
			\begin{align*}
				W_{J_*\setminus J'}(q)W_{J_*\cap J'}(q)
				&\sum_x\rparp P{J_*}J_x(q)T_x\\
					&=W_{J_*}(q)\ds\sum_x\rparp P{J_*}J_x(q)T_x\\
					&=\qy*{\sum_xP_xT_x}\cdot\qy*{\sum_{y\in W_J}q^{\ell(y)}T_y}\\
					&=\qy*{\sum_xP_xT_x}\cdot\qy*{\sum_{y\in W_{J\setminus J'}}
						q^{\ell(y)}T_y}\cdot\qy*{\sum_{z\in W_{J'}}q^{\ell(z)}T_z}\\
					&=W_{J\setminus J'}(q) \qy*{\sum_xP_xT_x}\cdot
						\qy*{\sum_{z\in W_{J'}}q^{\ell(z)}T_z}\\
					&=W_{J\setminus J'}(q)W_{J'_*}(q)
						\sum_x\rparp P{J'_*}{J'}_xT_x\\
					&=\rcos[J\setminus J']{J_*}W(q) W_{J_*\setminus J'}(q) 
						\rcos[J'_*]{J_*}W(q)W_{J_*\cap J'}(q)
						\sum_x\rparp P{J'_*}{J'}_xT_x.
			\end{align*}
		Cancelling common factors then gives the desired
			$\rparp P {J_*}J_x=\rcos[J\setminus J']{J_*}{W}(q)
			\rcos[J'_*]{J_*}{W}(q)\rparp P {J'_*}{J'}_x$ for all $x$.
	\end{proof}
\end{lem}

\begin{prop}
\label{prop:contract}
	For any contraction $J'$ of a factorization $J$, we have
		$P^J=p(q)P^{J'}$ where $p(q)\in\NN[q]$ is a product of Poincar\'e
		polynomials of parabolic subgroups and of sets of minimal coset
		representatives.
	In particular, for any contraction $J'$ of $J$ we have that
		$w(J')=w(J)$, and that $J$ is admissible (respectively tight)
		only if $J'$ is admissible (respectively tight).
	\begin{proof}
		We argue by induction on the length $r$ of the factorization 
			the claim for any \emph{right-}contraction $J'$ of $J$.
		The corresponding claim for left-contractions then follows by
			duality, and the main claim for general contractions then
			follows by induction on the number of left-contractions
			and right-contractions.

		The base case $r=1$ is trivial since always $J_1=J'_1$ and $\rdes
			J=\rdes{J'}$.
		For the inductive step, \cref{cor:recursive-defect-polys}
			asserts that $P^J=\rparp{P^{\upslice J{r-1}}}{\rdups J{r-1}\cap
			J_r}{J_r}$.

		Since $J'$ being a right-contraction of $J$ implies by definition
			that $\upslice{J'}{r-1}$ is a right-contraction of $\upslice
			J{r-1}$, the inductive hypothesis gives us 
			$P^{\upslice J{r-1}}=p(q)P^{\upslice{J'}{r-1}}$, and hence
				$\rdes{\upslice{J'}{r-1}}\subseteq\rdes{P^{\upslice{J'}{r-1}}}
				=\rdes{P^{\upslice J{r-1}}}$ by \cref{lem:pre-descent} and
				also $\rparp{P^{\upslice J{r-1}}}{\rdups J{r-1}\cap J_r}{J_r}
				=p(q)\rparp{P^{\upslice{J'}{r-1}}}{\rdups J{r-1}\cap J_r}{J_r}$.

		Since $J'$ being a right-contraction of $J$ also implies by
			definition that $J'_r$ is a union of connected components of $J_r$
			such that $J_r\setminus J'_r\subseteq \rdes{\upslice{J'}{r-1}}
				\subseteq \rdes{P^{\upslice J{r-1}}}$,
			and since \cref{lem:contraction-increases-descents} implies
			that $\rdes{\upslice J{r-1}}\cap J'_r\subseteq
				\rdes{\upslice{J'}{r-1}}\cap J'_r$, we can apply
				\cref{lem:contraction} to conclude that
			\begin{align*}
			P^J
				&=p(q)\rparp{P^{\upslice{J'}{r-1}}}{\rdups J{r-1}\cap J_r}{J_r}\\
				&=\rcos[J_r\setminus J'_r]{\rdups J{r-1}\cap J_r}W(q)
					\rcos[\rdups{J'}{r-1}\cap J'_r]{\rdups J{r-1}\cap J_r}W(q)
					p(q)\rparp{P^{\upslice{J'}{r-1}}}{\rdups{J'}{r-1}\cap J'_r}
						{J'_r}\\
				&=\rcos[J_r\setminus J'_r]{\rdups J{r-1}\cap J_r}W(q)
					\rcos[\rdups{J'}{r-1}\cap J'_r]{\rdups J{r-1}\cap J_r}W(q)
				p(q)P^{J'}
			\end{align*}
			as desired.
		The supplemental claim then follows because $P^J_x=0$ if and only if
			$P^{J'}_x = 0$, and furthermore $P^J_w=1$ implies $P^{J'}=P^J$.
	\end{proof}
\end{prop}

To complete the inductive step of our argument, it remains to show that
	every tight factorization $J$ has a contraction $J'$ such that
	$\rdes{J'}=\rdes{P^{J'}}$ (which by duality would imply it also has a
	contraction $J''$ such that $\ldes{J''}=\ldes{P^{J''}}$).
We will actually prove what turns out to be a more general fact: that a
	factorization $J$ such that $\rdes{P^{\upslice Jk}}
			=\rdes{w(\upslice Jk)}$ for all $k$ has a contraction that is
		\emph{right-descent}, which, as
		\cref{prop:rdescent-inductive-step} below shows, is a
		stronger but more tractable version of the property $\rdes J
			=\rdes{P^J}$.

\begin{defn}
\label{defn:descent}
	Given a factorization $J_1,\dots,J_r$, we denote by $\crd{J}_k$ the
		union of all connected components of $J_k$ that are also connected
		components of $\rdes{w(\upslice Jk)}$ and by $\cld J_k$ the union of
		all connected components of $J_k$ that are also connected components
		of $\ldes{w(\downslice Jk)}$.
	We say that $J_1,\dots,J_r$ is a \define{right-descent factorization}
		if $J_k=\crd J_k$ for all $k$, a \define{left-descent
		factorization} if $J_k=\cld J_k$ for all $k$. 
	We say that a factorization $J_1,\dots,J_r$ is a \define{bidescent
		factorization} if it is both a right-descent and left-descent, i.e.\
		if $J_k=\crd J_k=\cld J_k$ for all $k$. 
\end{defn}

\begin{exmp}
\label{exmp:descent-type-A}
	In type $A$, an interval $[a,b]$ corresponding to a connected
		component of $J_k$ fails to be a connected component of
		$w=w(\upslice Jk)$ if $w(a-1)>w(a)$ or $w(b)>w(b+1)$, i.e.\ if
		$\set{(a-1,a),(b,b+1)}\cap\rdes{w(\upslice Jk)}\setminus J_k
			\neq\emptyset$.
\end{exmp}

\begin{prop}
\label{prop:rdescent-inductive-step}
	Given a factorization $J_1,\dots,J_r$ and $k$ such that 
		$1\leq k\leq r$ and $\rdes{\upslice J{k-1}}
			=\rdes{w(\upslice J{k-1})}$, then $\rdes{\upslice Jk}
			=J_k\udj\qy{\rdes{w(\upslice Jk)}\setminus\nbd{J_k}}$, 
		and furthermore $\rdes{\upslice Jk}=\rdes{w(\upslice Jk)}$ if and
		only if $\crd J_k=J_k$.
	In particular, a factorization $J_1,\dots, J_r$ is right-descent
		if and only if the inclusions $\rdes{\upslice
			Jk}\subseteq\rdes{P^{\upslice Jk}} \subseteq\rdes{w(\upslice Jk)}$
			are equalities for all $k$.
	\begin{proof}
		By definition and the hypothesis, we have
			$\rdes{\upslice Jk}
				= J_k\udj\qy{\rdes{\upslice J{k-1}}\setminus\nbd{J_k}}
				= J_k\udj\qy{\rdes{w\qy*{\upslice J{k-1}}}
					\setminus\nbd{J_k}}
				= J_k\udj\qy{\rdes{w\qy*{\qy*{\upslice J{k-1}}}\cup J_k}
					\setminus\nbd{J_k}}$. 

		\cref{prop:T_L} and \cref{lem:pre-descent} 
			give containments 
			$J_k\udj\qy{\rdes{w\qy*{\qy*{\upslice J{k-1}}}\cup J_k}
					\setminus\nbd{J_k}} 
				\supseteq J_k\udj\qy*{\rdes{\rcos{J_k}{w\qy*{
					\upslice J{k-1}}}\mxml{J_k}} \setminus\nbd{J_k}}
				= J_k\udj\qy{\rdes{w\qy*{\upslice Jk}} \setminus\nbd{J_k}}
				\supseteq J_k\udj \qy{\rdes{\upslice Jk}\setminus\nbd{J_k}}$.
		Then by definition 
			$J_k\udj \qy{\rdes{\upslice Jk}\setminus \nbd{J_k}}
				=J_k\udj\qy{\qy{J_k\udj\qy{\rdes{\upslice J{k-1}}
					\setminus\nbd{J_k}}}\setminus\nbd{J_k}}
				=J_k\udj\qy{\rdes{\upslice J{k-1}}\setminus\nbd{J_k}}
				=\rdes{\upslice Jk}$. 

		Thus we have $\rdes{\upslice Jk}
				\supseteq J_k\udj\qy*{\rdes{w(\upslice Jk)}\setminus\nbd{J_k}}
				\supseteq\rdes{\upslice Jk}$, i.e.\  the desired 
			$\rdes{\upslice Jk}
				=J_k\udj\qy*{\rdes{w(\upslice Jk)}\setminus\nbd{J_k}}$, from
			which follows that $\rdes{w(\upslice Jk)}=\rdes{\upslice Jk}$ if
			and only if $\rdes{w(\upslice Jk)}\cap\nbd{J_k}=J_k$, i.e.\ if and
			only if $J_k$ is a union of connected components of
			$\rdes{w(\upslice Jk)}$, i.e.\ if and only if $J_k=\crd{J}_k$.
				
		The supplemental claim then follows by induction on $k$,
			since \cref{lem:pre-descent} implies
			that $\rdups Jk \subseteq \rdes{P^{\upslice Jk}}
				\subseteq \rdes{w\qy*{\upslice Jk}}$,
			and by the main claim equality of these inclusions
			given that $\upslice J{k-1}$ is right-descent is equivalent
			to $J_k=\crd J_k$.
	\end{proof}
\end{prop}

The advantage of being right-descent over the property $\rdes J
	=\rdes{P^J}$ is that, as \cref{lem:contraction-descent} below
	shows, any factorization has a natural candidate for a right-descent
	contraction.

\begin{lem}
\label{lem:contraction-descent}
	For a factorization $J_1,\dots,J_r$, we have that
		$\crd{J}$ is a right-contraction of $J$ if and only if
		$J_k\setminus \crd{J}_k\subseteq\rdes{w(\upslice J{k-1})}$
		for every $k$ such that $1\leq k \leq r$, in which case $\crd{J}$ is
		a right-descent factorization.
	\begin{proof}
		We will show the main claim holds for all slices $\upslice Jk$ of
			the factorization $J_1,\dots,J_r$ by induction on $k$.
		To that end, note that $\crd{\qy*{\upslice Jk}}=\upslice{\crd J}k$
			by definition.

		The base case $k=1$ is immediate since $J_1=\crd J_1$.
		For the inductive step, we have that $\upslice{\crd J}{k-1}$ is a
			right-contraction of $\upslice J{k-1}$ and a right-descent
			factorization, and hence by
			\cref{prop:rdescent-inductive-step} and
			\cref{prop:contract} that 
			$\rdes{\upslice{\crd J}{k-1}}=\rdes{w(\upslice{\crd J}{k-1})}
				=\rdes{w(\upslice J{k-1})}$.
		But then $\upslice{\crd J}k$ being a right-contraction of $\upslice
			Jk$ is equivalent to $J_k\setminus\crd J_k
				\subseteq\rdes{\upslice{\crd J}{k-1}}
				=\rdes{w(\upslice{\crd J}{k-1})}=\rdes{w(\upslice J{k-1})}$, in
			which case $w(\upslice{\crd J}k)=w(\upslice Jk)$ by
			\cref{prop:contract} and hence 
			$\crd{\qy*{\upslice{\crd J}k}}_k
				=\crd J_k$ so that $\upslice{\crd J}k$ is right-descent 
			as desired.
	\end{proof}
\end{lem}

\begin{rmrk}
	Note that $\crd J$ can be right-descent without $J_k\setminus
		\crd J_k\subseteq\rdes{w\qy*{\upslice J{k-1}}}$, i.e.\
		without being a right-contraction of $J$.
	For example, the type $A$ factorization of $w=4\ 3\ 2\ 1$ (in one line
		notation) given by $J_1,J_2,J_3$ corresponding to the intervals
		$[2,4],[1,2],[2,4]$ has $\crd J_3=\emptyset$ even though the top
		$[2,4]$ cannot be absorbed into the bottom as it is not disjoint
		from the middle $[1,2]$.
\end{rmrk}

We now proceed with showing that tight factorizations, and more
	generally factorizations with the property that
	$\rdes{P^{\upslice Jk}} =\rdes{w(\upslice Jk)}$ for all $k$,
	satisfy the property of \cref{lem:contraction-descent},
	i.e.\ are such that $\crd J$ is a right-descent right-contraction
	of $J$.

\begin{lem}
\label{lem:necessary}
	Given a finitary $J\subseteq S$, a positive family $(P_x)_{x\leq w}
		\in\NN[q]$ with $\underline h^P$ is involution-invariant,
		and a subset $J_*\subseteq\rdes P\cap J$,
		each $z\in W_{\rdes P\cap\rdes{\rparp P{J_*}J}}$ satisfies
		$\len{\rcos[J]{J_*}w} \geq\len{\rcos[J]{J_*}{(wz)}}$.
	\begin{proof}
		Since $J\subseteq\rdes{\rparp P{J_*}J}$
			implies $W_J\subseteq W_{\rdes{\rparp P{J_*}J}}$, it follows that
			$z\in\rdes{\rparp P{J_*}J}$ implies $\mxml J^{-1}w_Jz\in
				W_{\rdes{\rparp P{J_*}J}}$ so that $\ds
				\qy*{\rparp P{J_*}J}_{\rcos Jw\mxml J}= \qy*{\rparp P{J_*}J}_{wz}
				=\sum_{y\in\rcos[J]{J_*}W}q^{\ell(y)}P_{\rcos J{(wz)}y}$. 
		Since $(P_x)_{x\leq w}$ is positive, choosing
			$y=\rcos[J]{J_*}{\qy{wz}}$ gives
			$\deg{\qy*{\rparp P{J_*}J}_{\rcos Jw\mxml J}}
				\geq\len{\rcos[J]{J_*}{(wz)}}+\deg P_{\rcos J{\qy{wz}}
					\rcos[J]{J_*}{\qy{wz}}}
			=\len{\rcos[J]{J_*}{(wz)}}+\deg P_w$,
			where the last equality holds because $\rcos[J_*]{}{\qy{wz}}z^{-1}
				\in W_{J_*}W_{\rdes P}\subseteq W_{\rdes P}$.
		But since $\underline h^P$ is invariant and $(P_x)_{x\leq w}$
			is positive, \cref{cor:invariance-positive} implies that
			$\deg{\qy*{\rparp P{J_*}J}_{\rcos Jw\mxml J}}
				=\len{\rcos[J]{J_*}w}+\deg P_w$, so that
				$\len{\rcos[J]{J_*}w}+\deg P_w
					\geq\len{\rcos[J]{J_*}{(wz)}}+\deg P_w$, and hence
			$\len{\rcos[J]{J_*}w}\geq\len{\rcos[J]{J_*}{(wz)}}$ as
			desired.
	\end{proof}
\end{lem}

\begin{prop}
\label{prop:invariant_mult_connected}
	Suppose $(P_x)_{x\leq w}\in\NN[q]$ is a positive family with
		$\underline h^P$ involution-invariant.
	Then for any finitary subset $J\subseteq S$ and subset
		$J_*\subseteq\rdes P\cap J$, if there exists a connected component
		$J'$ of $J$ that is neither a connected component of $\rdes{\rparp
		P{J_*}J}$ nor contained in $\rdes w$, then $\rdes
		w\not\subseteq\rdes P$.
	\begin{proof}
		By definition $J'$ being a connected component of $J$ means that
			$J'$ is connected such that $\nbd{J'}\cap J=J'$, and hence $J'$
			failing to be a connected component of $\rdes{\rparp P{J_*}J}
				\supseteq J \supseteq J'$ means that there exists
			$s\in\nbd{J'}\cap\rdes{\rparp P{J_*}J}\setminus J'
				\subseteq\nbd{J'}\cap\rdes{\rcos Jw\mxml J}\setminus J'
				\subseteq\nbd{J'}\cap\qy*{\rdes w\cup J}\setminus J'
				=\nbd{J'}\cap\rdes w\setminus J'$,
				where the first containment holds by
				\cref{lem:pre-descent}, and the second holds by
				\cref{prop:T_L}.

		Furthermore, because $J'$ is connected but not contained in $\rdes
			w$, \cref{lem:wz-w-inequality} implies that corresponding to 
			$s\in\rdes w\cap\nbd{J'}\setminus J'$ there exists
			$z\in W_{\qy*{J'\cap\rdes w}\cup\set s}$ such that
			$\len{\rcos[J]{\rdes w}{(wz)}}-\len{\rcos[J]{\rdes w}w}>0$.
		Thus $\rdes w\subseteq\rdes P$ would imply
			$\rdes w\cap J\subseteq\rdes P\cap J$ and
			$z\in W_{\qy*{J'\cap\rdes w}\cup\set s}\subseteq
				W_{\rdes P\cap \rdes{\rparp P{\rdes w\cap J}J}}$
			but
			$\len{\rcos[J]{\rdes w\cap J}{(wz)}}-\len{\rcos[J]{\rdes w\cap J}w}
				>0$, contradicting \cref{lem:necessary}.
		Consequently, we must have $\rdes w\not\subseteq\rdes P$ as desired.
	\end{proof}
\end{prop}

\begin{cor}
\label{cor:invariant-rdes-contract}
	If a factorization $J_1,\dots,J_r$ satisfies $\rdes{P^{\upslice Jk}}
		=\rdes{w(\upslice Jk)}$ for all $k$ (e.g.\ if $J$ is tight),
		then $\crd J$ is a right-descent right-contraction of $J$.
	\begin{proof}
		We prove for each $k$ that $J_k\setminus \crd J_k
			\subseteq \rdes{w(\upslice J{k-1})}$, which is equivalent to
			the claim by \cref{lem:contraction-descent}.
		By \cref{cor:quotient} $\underline h^{P^{\upslice J{k-1}}}$
			is involution-invariant, and by
			\cref{lem:normalized-factorization} it is positive, so that
			by \cref{prop:invariant_mult_connected} we have that
			every connected component of $J_k$ which is not contained in
			$\rdes{P^{\upslice J{k-1}}}=\rdes{w(\upslice J{k-1})}$ is a
			connected component of $\rdes{
				\rparp{P^{\upslice J{k-1}}}{\rdups J{k-1}\cap J_k}{J_k}}
				=\rdes{P^{\upslice Jk}}=\rdes{w(\upslice Jk)}$ and so is
				contained in $\crd J_k$ by definition.

		To show that $J$ being tight implies the hypothesis, note that by
			\cref{prop:boundedness} and
			\cref{cor:recursive-defect-polys}
			$J$ tight implies $\upslice Jk$ tight for every $k$ such that
			$1\leq k\leq r$, and \cref{prop:invariant_kl}
			implies that in this case $\rdes{P^{\upslice Jk}}
				=\rdes{(P_{x,w(\upslice Jk)})_{x\leq w(\upslice Jk)}}
				=\rdes{w(\upslice Jk)}$.
	\end{proof}
\end{cor}

\begin{thm}
\label{thm:2}
	A factorization $J_1,\dots,J_r$ has the property that
		some iterated composition of Kazhdan--Lusztig basis elements
		$C'_{\mxml{J_1}},\dots, C'_{\mxml{J_r}}$ (in that order) using the Schur
		composition laws is a Kazhdan--Lusztig basis element $C'_w$ only if
		$J$ is a factorization of $w$ with a tight contraction that is
		admissible absolutely bidescent.
	\begin{proof}
		We argue by induction on the length $r$ of the factorization.
		The base case of length $r=1$ is immediate since any length $1$
			factorization is minimal by definition of the notion of
			contraction and is tight because $\underline
			h^{P^{J_1}}=C'_{\mxml{J_1}}$.

		For the inductive step, suppose that for some $k$ with $1\leq k<r$
			a Kazhdan--Lusztig basis element $C'_w$ is realized as a composite
			$\underline h^{P'}\circ_{J_*}\underline h^P$, where $\underline h^P$
			is an iterated composite of $C'_{\mxml{J_1}},\dots,C'_{\mxml{J_k}}$ and
			$\underline h^{P'}$ is an iterated composite of $C'_{\mxml{J_{k+1}}},
			\dots,C'_{\mxml{J_r}}$.
		By \cref{prop:boundedness}, both $\underline h^P$ and
			$\underline h^{P'}$ are Kazhdan--Lusztig basis elements, so by
			the inductive hypothesis there exist a pair of a tight
			contractions
			$J'_1,\dots,J'_k$ and $J'_{k+1},\dots,J'_r$ of $J_1,\dots,J_k$ and
			$J_{k+1},\dots,J_r$ with the property that
			$\underline h^P=\underline h^{P^{\upslice{J'}k}}$
			and $\underline h^{P'}=\underline h^{P^{\downslice{J'}{k+1}}}$.
		Since \cref{prop:contract} implies that taking a
			contraction of a tight factorization preserves the defect
			polynomials, we may assume by
			\cref{cor:invariant-rdes-contract} and its dual that
			$\upslice{J'}k$ and $\downslice{J'}{k+1}$ respectively right- and
			left-descent contractions of $\upslice Jk$ and $\downslice J{k+1}$.
		Then \cref{lem:defect-ind} implies that
			$C'_w=\underline h^{P'}\circ_{J_*}\underline h^P
				=\underline h^{P^{\downslice{J'}{k+1}}}
					\circ_{J_*}\underline h^{P^{\upslice{J'}k}}
				=\underline h^{P^{J'}}$, so that $J'$ is a tight factorization of
			$w$, and by \cref{lem:slice-contraction} a contraction of $J$.

		Finally, \cref{prop:contract} implies minimal contractions of tight
			factorization are tight, so we may assume $J$ is a minimal tight
			factorization.
		But if $J$ is a tight factorization, then by
			\cref{cor:invariant-rdes-contract} $\crd J$ is a right-descent
			right-contraction of $J$, so if $J$ is also minimal then $J=\crd
			J$ by definition, so that $J$ is right-descent.
		By duality, tight minimal factorizations are also left-descent and
			hence bidescent.
		Since $J$ tight implies $J$ is admissible by
			\cref{prop:bounded-monic} and absolutely tight by
			\cref{cor:factors}, and since $J$ minimal implies it is
			absolutely minimal by \cref{cor:absolutely-minimal}, it follows
			that $J$ is also admissible absolutely bidescent.
	\end{proof}
\end{thm}

\begin{rmrk}
\label{rmrk:simply-laced}
	A Coxeter system is \emph{simply-laced} if
		$s_is_js_i=\mxml{\set{s_i,s_j}}$ for any pair of non-commuting
		generators $s_i,s_j\in S$, i.e.\ if all short braids are maximal
		elements of parabolic subgroups.
	Thus, for such Coxeter systems short braids are neither right- nor
		left-descent, and hence the property of being absolutely bidescent
		is a generalization of short braid avoidance.

	This is not the case for other Coxeter systems, where short braid
		avoidance must be additionally imposed.
	It would be interesting to determine a common generalization of short
		braid avoidance and being (admissible) absolutely  bidescent for
		those Coxeter systems.
\end{rmrk}

\section{Admissible right-descent and absolutely bidescent factorizations}

\label{sec:bidescent}
We present our results toward constructing tight factorizations in
	arbitrary Coxeter system.
We first describe a generic algorithm for computing admissible
	right-descent factorizations, which we suggest be used as the basis
	for a modification of Deodhar's and Scott and Sprowl's algorithms for
	computing Kazhdan--Lusztig basis elements.

We then give examples indicating that being right-descent is not in
	general an absolute property.
Finally, we prove that nevertheless every Coxeter element has an
	absolute bidescent factorization and then indicate the difficulty in
	verifying that every Coxeter element has an \emph{admissible}
	absolutely bidescent factorization.
At the end, we introduce the property of \emph{strong right-descent}
	factorizations, which is easily seen to be absolute and plays a key
	role in our subsequent results for type $A$.

We begin by characterizing admissible factorizations among 
	right-descent ones.
\begin{lem}
\label{lem:admissible-right-descent}
	A factorization $J_1,\dots,J_r$ is admissible only if for each $k$
		there exists a finitary subset $\overlap J*k\subseteq J_k$ such that
		$w(\upslice J{k-1})_{J_k}=\mxml{\overlap J*k}$, in which case
		$\overlap J*k=\rdes{\upslice J{k-1}}\cap J_k
			=\rdes{P^{\upslice J{k-1}}}\cap J_k
			=\rdes{w(\upslice J{k-1})}\cap J_k$.
	Conversely for factorizations such that $\rdups J{k-1}\cap J_k
				=\rdes{w(\upslice J{k-1})}\cap J_k$ for all $k$ (e.g., by
		\cref{prop:rdescent-inductive-step}, for right-descent
		factorizations), the existence of $\overlap J*k$ such that 
		$w(\upslice J{k-1})_{J_k}=\mxml{J^k_*}$ is sufficient for the
		factorization to be admissible.
	\begin{proof}
		By \cref{cor:factors} if $J$ is admissible then so are
			$\upslice J{k-1}$ and $\upslice Jk$ for every $k$ such that
			$1\leq k\leq r$, so that since
			\cref{cor:recursive-defect-polys} asserts that
			$P^{\upslice Jk}
				=\rparp{P^{\upslice J{k-1}}}{\rdups J{k-1}\cap J_k}{J_k}$
			while \cref{cor:quotient} guarantees that $\underline
			h^{P^{\upslice J{k-1}}}$ is involution-invariant,
			\cref{cor:invariance-positive} implies that $\upslice Jk$
			is admissible, i.e.\ that $P^{\upslice Jk}_{w(\upslice Jk)}=1$, 
			only if $\overlap J*k=\rdups J{k-1}\cap J_k$ and 
			$w(\upslice J{k-1})_{J_k}=\mxml{\overlap J*k}$, with the latter
			implying $\overlap J*k=S\cap T\cap W_{\overlap J*k}
				=S\cap\rref{\mxml{\overlap J*k}}
				=S\cap\rref{(w(\upslice J{k-1})_{J_k}}
				=\rdes{w(\upslice J{k-1})}\cap J_k$ by
			\cref{lem:subsystem} and \cref{prop:finitary}.
		Since $\rdups J{k-1}\subseteq \rdes{P^{\upslice J{k-1}}}
			\subseteq\rdes{w(\upslice J{k-1})}$ by
			\cref{lem:pre-descent}, the main claim follows.

		For the converse, we argue by induction on the length $r$ of the
			factorization.
		The base case of $r=1$ is immediate since such factorizations are
			admissible and satisfy $w(\upslice J0)=e$.
		For the inductive step, assume the claims hold for factorizations of
			length less than $r$ and consider a factorization $J_1,\dots,J_r$.
		By the inductive hypothesis, the existence of appropriate 
			$\overlap J*k$ for every $k$ such that $1\leq k\leq r-1$ implies
			that $\upslice J{r-1}$ is admissible.
		Now if $\rdes{\upslice J{r-1}}\cap J_r
			=\rdes{w(\upslice J{r-1})}\cap J_r$, then because $w(\upslice
			J{r-1})_{J_r}=\mxml{\overlap J*r}$ implies as above that 
			$\overlap J*r=\rdups J{r-1}\cap J_r$,
			it follows by \cref{cor:recursive-defect-polys}
			that $P^J=\rparp{P^{\upslice J{r-1}}}{\overlap J*r}{J_r}$,
			which since $\underline h^{P^J}$ is invariant by
			\cref{cor:quotient} implies by \cref{cor:invariance-positive}
			that $P^J_{w(J)}=P^{\upslice J{r-1}}_{w(\upslice J{r-1})}=1$ so
			that $J$ is admissible as desired.
	\end{proof}
\end{lem}

\begin{exmp}
\label{exmp:adm-rdes-typea}
	In the case of Coxeter systems of type $A$, recall that $\rcos Jw$ is
		obtained from $w$ by sorting in increasing order the sequences of
		outputs $w(a),\dots,w(b)$ for each interval $[a,b]$ corresponding to
		a connected component of $J$.
	For $J_*\subseteq J$ it then follows that $\rcos Jw\mxml{J_*}$ is
		obtained by further sorting for each interval $[a,b]$ in decreasing
		order each sequence of outputs within intervals $[a_i,b_i]$ where
		$[a=a_1,b_1],[b_1+1=a_2,b_2],\dots,[b_{m-1}+1=a_m,b_m=b]$ is the
		decomposition of the interval $[a,b]$ arising from the interval
		decomposition corresponding to $J_*\subseteq S$.

	Consequently, $w_J=\mxml{J_*}$ for some $J_*$ if and only if for any
		pair of positions $p<q\in[a,b]$, where $[a,b]$ is an interval
		corresponding to a connected component of $J$, we have
		$w(p)<w(q)$ if and only if $p\in[a_i,b_i]$ and $q\in[a_j,b_j]$ with
		$i<j$.
	Moreover, if $wv=\rcos Jww_J$, then \cref{prop:T_L} implies that 
		for $p<q\in[a,b]$ we have $w(p)<w(q)$ equivalent to $v(p)>v(q)$.
\end{exmp}

\begin{exmp}
\label{exmp:skandera}
	A factorization $J_1,\dots,J_r$ satisfies the property that for each
		$k$ there exists $J^k_*\subseteq J$ such that $w(\upslice
		J{k-1})_{J_k}=\mxml{J^k_*}$ if and only if it is a right-contraction
		of the factorization $J'_1,\dots,J'_{2r+1}$ given by $J'_{2k+1}=J_k$
		and $J'_{2k}=J^k_*$, and which has the further property that
		$\mxml{J'}=\mxml{J'_1}\cdots \mxml{J'_{2r+1}}$.
	As mentioned in \cref{rmrk:skandera}, Skandera's zig-zag
		factorizations are instances of such factorizations. 
\end{exmp}

\begin{prop}
\label{prop:rigidity}
	Given finitary connected subsets $J,J'\subseteq S$ and a Coxeter
		element $w$ such that $w\neq w^J\mxml J=w^{J'}\mxml{J'}$, if
		$J\not\subseteq J'$, then $w_J\neq \mxml{J_*}$ for any
		subset $J_*\subseteq J$.

	In particular, given a sequence $w_1,\dots,w_r$ of distinct Coxeter
		elements, there is at most one admissible right-descent factorization
		$J_1,\dots,J_r$ such that each $J_k$ is connected and
		$w(\slice{J}{}{k})=w_k$ for each $1\leq k\leq r$.
	\begin{proof}
		Since $w_J=\mxml{J_*}$ only if $J_*=\rdes w\cap J$ by
			\cref{lem:subsystem}, in which case $\rcos[J]{\rdes w}w
				=\rcos[J]{J_*}w=e$,
			it will suffice to show that 
			$\len{\rcos[J]{\rdes w}w}>0$, or equivalently by
			\cref{lem:length-coset-rep}, to construct 
			$t\in\rref w\cap W_J\setminus W_{\rdes w}$.

		Let $v\in W$ be such that $\rcos Jw\mxml J=wv=\rcos{J'}w\mxml{J'}$.
		Then \cref{prop:T_L} implies $J\cup J'\subseteq\rdes{wv}
			\subseteq\rdes w\cup\qy*{J\cap J'}$, so that
			$J\setminus\rdes w=\qy*{J\cap J'}\setminus\rdes w
			=J'\setminus\rdes w\neq\emptyset$ because $\rcos Jw\mxml J\neq w
				\neq\rcos Jw\mxml{J'}$, and $\emptyset\neq J+J'
				=\qy*{J\cup J'}\setminus\qy*{J\cap J'}\subseteq\rdes w$.
		Furthermore, since $J\cup J'\subseteq\rdes{wv}$ implies $J\cup J'$
			is finitary by \cref{prop:finitary} and so
			by \cref{prop:finitary} that $W_{J\cup J'}$ is finite,
			by the classification of finite Coxeter systems \cite{Coxeter35},
			or more directly \cite[Lemme~2]{Tits}, 
			we have that $W_{J\cup J'}$ has an acyclic Coxeter--Dynkin
			diagram, and thus the non-empty intersection $J\cap J'$
			of the connected subsets $J,J'\subseteq S\cap W_{J\cup J'}$
			is itself connected.

		Thus, since $J\cap J'\subsetneq J$ is connected, there exists
			$s\in J\cap\nbd{J\cap J'}\setminus\qy{J\cap J'}
				\subseteq J\setminus J'\subseteq\rdes w$, and in
			particular such that
			$s\in\rdes w\cap J \cap\nbd{J\cap J'}\setminus(J\cap J')$.
		Now, since $J\cap J'$ is connected and
			$\emptyset\neq J\cap J'\setminus\rdes w$, i.e.\ since 
			$J\cap J'$ is not contained in $\rdes w$,
			\cref{lem:connected} implies there exist 
			$s'\in J\cap J'\setminus\rdes w\subseteq J$ and 
			$z\in W_{\qy{\rdes w\cap J\cap J'}\cup\set s}\subseteq W_J$ such
			that $s\in\rref{z^{-1}s'z}$.
		Then $t=z^{-1}s'z\in T\cap W_J$ satisfies
			$J\cap J'\not\ni s\in\rref t$, and hence also 
			$t\not\in W_{J\cap J'} \supseteq \lref v$. 
		By \cref{prop:T_L} it follows that 
			$t\in T\cap W_J\setminus\lref v=\rref w\cap W_J$,
			and since $s'\not\in\rdes w$ but $z\in W_{\rdes w}$, we have
			$t=z^{-1}s'z\not\in W_{\rdes w}$, so that
			$t\in \rref w\cap W_J\setminus W_{\rdes w}$. 

		The supplemental claim follows from immediately from
			\cref{lem:admissible-right-descent} and the main claim.
	\end{proof}
\end{prop}

\begin{defn}
\label{defn:admissible_descent_algo}
	The \define{generic admissible right-descent factorization procedure}
		takes as input a Coxeter element $w\in W$ and produces a
		factorization of $w$ as follows:
		\begin{algorithmic}[1]
				\State Initialize $x \gets w$ and an empty sequence of subsets
					of $S$.
				\While{$x\neq e$, or equivalently $\rdes{x}\neq\emptyset$,}
					\parState{Choose a union $J$ of connected components of $\rdes
						x$ and add it to the beginning of the
						sequence.}%
					\State Choose $v=\mxml{J_*}$ for some subset $J_*\subseteq J$
						and assign $x \gets x^Jv$.%
				\EndWhile
				\State \Return the sequence of subsets.
		\end{algorithmic}
\end{defn}

\begin{exmp}
	In the case of type $A$, choosing a union $J$ of connected components
		of $\rdes w$ for a permutation $w$ amounts to choosing right-descent
		intervals $[a,b]$ of the permutation $w$.
	Choosing $J_*\subseteq J$ then amounts to choosing a partition
		$[a=a_1,b_1],[b_1+1=a_2,b_2],\dots,[b_{m-1}+1=a_m,b_m=b]$ of each of
		the chosen right-descent intervals. 
	Finally, replacing $w$ with $\rcos Jw\mxml{J_*}$ amounts to sorting the
		sequences of outputs $w(a),\dots,w(b)$ so that the outputs within
		each interval $[a_k,b_k]$ of the partition are in decreasing order
		but the outputs within any $[a_i,b_i]$ are in order with the outputs
		within any $[a_j,b_j]$ with $i<j$.
\end{exmp}

We call this procedure ``generic'' because it does not 
	restrict the choices of $J$ and $J_*\subseteq J$.
Combined with appropriate restrictions  and
	\cref{cor:bounded-algorithm}, our generic procedure yields in
		\cref{exmp:new-algo} our proposed modification of Scott and
		Sprowl's algorithm for computing Kazhdan--Lusztig polynomials of
		arbitrary Coxeter elements, and in \cref{prop:monotone-algorithm} of
		$4231$-$45312$-avoiding permutations.

\begin{lem}
\label{lem:factorization_procedure}
	Any specialization of the generic admissible right-descent
		factorization procedure has the following properties:
	\begin{enumerate}
		\item if $v$ is chosen to be the maximal element $\mxml J$ of $W_J$,
			then $\rcos Jxv=x$ and otherwise $\ell(\rcos Jxv)<\ell(x)$; in
			particular the procedure terminates after at most $\ell(w)$
			choices of $v\neq \mxml J$;
		\item any output of the procedure is an admissible right-descent
			factorization of~$w$.
	\end{enumerate}
	Furthermore the generic procedure itself has the properties that:
	\begin{enumerate}
	\setcounter{enumi}{2}
		\item it is always possible to choose a non-empty $J\subseteq\rdes
			w$ and $J_*$ such that $\ell(\rcos Jx\mxml{J_*})<\ell(x)$;
		\item every admissible right-descent factorization $J_1,\dots,J_r$
			of $w$ such that $J_k\neq\emptyset$ for each $k$ such that
			$1\leq k\leq r$ is a possible output of the generic procedure.
	\end{enumerate}
	\begin{proof}
		(1) $J\subseteq\rdes x$ if and only if $x=\rcos Jx\mxml J$ by
			\cref{prop:finitary}.
			Furthermore, for $v\in W_J$ we have $\ell(\rcos Jxv)
				=\ell(x^J)+\ell(v)$ by
				\cref{lem:coset-unique})
				and $\ell(v)<\ell(\mxml J)$ if $W_J\ni v\neq \mxml J$ by
				\cref{prop:finitary}.
			It follows that $\ell(\rcos Jxv)\leq\ell(x)$ with equality only if
				$v=\mxml J$.

		(2) We argue by induction on the length $r$ of the constructed
			factorization $J_1,\dots,J_r$.
		The base case of $r=1$ is immediate since by construction that
			happens if and only if $e=\rcos Jw\mxml{J_*}$ where $J$ is a union
			of connected components of $\rdes w$, in which case
			$J_*=\emptyset$ and $w=\rcos Jww_J=w_J=\mxml J$ by
			\cref{prop:finitary}.

		For the inductive step, we have that $\upslice J{r-1}$ is an
			admissible right-descent factorization of $\rcos{J_r}wv$ where
			$J_r$ is a union of connected components of $\rdes w$ and
			$v=\mxml{J_*}$ for some $J_*\subseteq J_r$.
		It follows from \cref{prop:finitary} that
			that $w=\rcos{J_r}w\mxml{J_r}=\rcos{J_r}{(\rcos{J_r}wv)}\mxml{J_r}$,
			i.e.\ that $J$ is a factorization of $w$. 
		Since $\upslice J{r-1}$ is a right-descent factorization and $J_r$
			is a union of connected components of $\rdes{w(J)}$, we have by
			definition that that $J$ is also a right-descent factorization.
		Finally, it is admissible by
			\cref{lem:admissible-right-descent}.

		(3) Since $J$ is a non-empty union of connected components of $x$,
			choosing $J_*=\emptyset$ gives $\ell(x^J\mxml{J_*})
				=\ell(x^J)<\ell(x^J)+\ell(\mxml J)=\ell(x^J\mxml J)$ by 
			\cref{lem:coset-unique}).

		(4) Any admissible right-descent factorization $J_1,\dots,J_r$ with
			all $J_k\neq\emptyset$ is output by the procedure by successively
			choosing the subset $J_k$ and the element
			$v=w(\upslice J{k-1})_{J_k}$ for each $k$ such that $1\leq k\leq r$
			starting with $k=r$, where $v=\mxml{J^*_k}$ by
			\cref{lem:admissible-right-descent}.
	\end{proof}
\end{lem}

\begin{exmp}
\label{exmp:new-algo}
	Our proposed generalization of Scott and Sprowl's algorithm for
		computing Kazhdan--Lusztig basis elements is to apply 
		\cref{cor:bounded-algorithm} with $I_x=M_x
		=\set{P^J:J_1,\dots,J_r$ is an admissible right-descent
		factorization of $w$ with each $J_k$ connected$}$, where each choice
		of a factorization $J$ is to be accomplished using the procedure of
		\cref{defn:admissible_descent_algo}, and the family
		$(P^J)_{x\in W}$ is to be computed using the formulas of
		\cref{cor:recursive-defect-polys} (in the
		$\set{T_x}$-basis) or \cref{prop:invariant_hecke_mult}
		(in the $\set{H_x}$-basis).

	We restrict to choices of connected $J_k$ in order to minimize the
		possibility of obtaining a right-descent factorization that is not
		absolutely right-descent (see
		\cref{exmp:spread-rearrangement_fail} and
		\cref{prop:connected-factorizations} below).

	Compared to Scott and Sprowl's original algorithm, our proposal uses
		shorter monomials that are more complex to multiply out.
	It would be interesting to determine whether this makes our algorithm
		faster in practice.

	Finally, note that Deodhar's recursive algorithm can likewise be
		modified so that it sill uses $M_y=\set{(P_{x,y})_{x\in W}}$ but
		initializes using $I_x =\set{P^J:J_1,\dots,J_r$ is an admissible
		right-descent factorization of $w$ with each $J_k$ connected$}$
		instead of $I_x=\set{P^J:J_1,\dots,J_r$ is a singleton factorization
		of $x$ induced by a reduced word$}$.
\end{exmp}

We now derive a simpler criterion for whether a factorization is
	absolutely right-descent, which we then use to show that every Coxeter
	element has an absolutely bidescent factorization.

\begin{exmp}
\label{exmp:spread-rearrangement_fail}
	Recall that \cref{cor:factors} asserts that a factorization is
		admissible only if it is absolutely admissible, and tight only if it
		is absolutely tight.
	It is \emph{not} true that a factorization is right-descent only if it
		is absolutely right-descent.
	In other words, it is not necessarily the case that every other
		factorization inducing the same element of the Cartier--Foata monoid
		as a right-descent factorization is also a right-descent
		factorization.
	As a specific counter-example, $\set{s_1,s_3}$ has interval
		decomposition $[1,2],[3,4]$ of $ws_1s_3=4231$ for $w=2413$,
		but $\set{s_1}$ is not a right-descent interval $[1,2]$ of
		$w\set{s_1}=4213$, nor is
		$\set{s_3}$ a right-descent interval $[3,4]$ of $ws_3=2431$. 

	Another counter-example is that even though $\set{s_1}$ is a
		right-descent interval of $ws_1=4123$ for $w=1423$ and $\set{s_3}$
		is a right-descent interval $[3,4]$ of $ws_1s_3=4132$, it is not the
		case that $\set{s_3}$ is a right-descent interval $[3,4]$ of
		$ws_3=1432$.
\end{exmp}

\begin{lem}
\label{lem:rdes-slice}
	A factorization $J_1,\dots,J_r$ is right descent only if every slice
		$\slice Jik$ is right-descent.
	\begin{proof}
		Given a slice $\slice Jik$, by self-duality of factorizations there
			exists $v\in W$ with $w(\upslice Jj)=vw(\slice Jij)$ with
			$\ell(w(\upslice Jj)=\ell(v)+\ell(w(\slice Jij))$. 
		Thus, $\rref{w(\upslice Jj)}
				=w(\slice Jij)^{-1}\rref vw(\slice Jij)\udj\rref{w(\slice Jij)}
				\supseteq \rref{w(\slice Jij)}$ and hence
		$J_j\subseteq \rdes{\slice Jij}\subseteq\rdes{w(\slice Jij)}
			=\rdes{w(\upslice{\slice Jik}j}$ by \cref{lem:pre-descent}.
		In particular, $J_j$ a connected component of 
			$\rdes{w(\upslice Jj)}$ implies $J_j$ is a connected
			component of $\rdes{w(\slice Jij)}$, and hence $J$ being
			right-descent implies $\slice Jik$ is right-descent. 
	\end{proof}
\end{lem}

\begin{lem}
\label{lem:rdes-union}
	Suppose $J_1,\dots,J_r$ is a factorization and $k$ such that
		$1\leq k<r$ and $J_k$ commutes with $J_{k+1}$.
	Then $J_1,\dots,J_r$ is right-descent only if the factorization
		$J'$ given by
		$J_1,\dots,J_{k-1},J_k\udj J_{k+1},J_{k+2},\dots,J_r$ is
		right-descent.
	\begin{proof}
		Note that for $i<k$ we have $J_i=J'_i$ and $w(\upslice Ji)
			=w(\upslice{J'}i)$ so that $\crd J_i=\crd{J'}_i$ and hence 
			$J_i=\crd J_i$ implies $J'_i=\crd{J'}_i$.
		Similarly, for $i>k+1$ we have $J_i=J'_{i-1}$ and
			$w(\upslice Ji)=w(\upslice{J'}{i-1})$, so again $\crd
			J_i=\crd{J'}_{i-1}$ and hence
			$\crd J_i=J_i$ implies $\crd{J'}_{i-1}=J'_{i-1}$.

		It remains to show that $J_k=\crd J_k$ and $J_{k+1}=\crd J_{k+1}$
			imply $J'_k=\crd{J'}_k$. 
		Since $w(\upslice{J'}k)
			=w(\upslice J{k+1})$ by \cref{cor:spread_defect}, this is
			by definition equivalent to showing that 
			$\rdes{\upslice Jk}\cap\nbd{J_k}=J_k$ and 
			$\rdes{\upslice J{k+1}}\cap\nbd{J_{k+1}}=J_{k+1}$ 
			imply $\rdes{w(\upslice J{k+1})}\cap\nbd{J_k\udj J_{k+1}}
				=J_k\udj J_{k+1}$.

		To that end, first compute 
		$\rdes{w(\upslice J{k+1})}\cap\nbd{J_k\udj J_{k+1}}
			=\qy*{\rdes{w(\upslice J{k+1})}\cap\nbd{J_k}}
				\cup\qy*{\rdes{w(\upslice J{k+1})}\cap\nbd{J_{k+1}}}
			=\qy*{\rdes{w(\upslice J{k+1})}\cap\nbd{J_k}}\cup J_{k+1}$.
		Then taking $v\in W$ such that $w(\upslice J{k+1})=w(\upslice Jk)v$,
			\cref{prop:T_L} gives us 
			$\rdes{w(\upslice J{k+1})}\cap\nbd{J_k}
				=\qy*{\rdes{w(\upslice Jk)}\cap\rdes{w(\upslice J{k+1})}
					\udj\lref v} \cap \nbd{J_k}
				=\nbd{J_k}\cap\rdes{w(\upslice Jk)}\cap\rdes{w(\upslice J{k+1})}
				=J_k\cap\rdes{w(\upslice J{k+1})}$
				because $v\in W_{J_{k+1}}$ implies $\lref v\cap\nbd{J_k}
				\subseteq W_{J_{k+1}}\cap W_{\nbd{J_k}}=W_{\emptyset}$. 

		Since $\rdes{w(\upslice J{k+1})}\supseteq\rdups J{k+1}
			=J_{k+1}\udj\qy{\rdups Jk\setminus \nbd{J_{k+1}}}\supseteq J_k$ 
			(\cref{lem:pre-descent}), and since $J_k$ and $J_{k+1}$
			commute, it follows that $J_k\cap\rdes{w(\upslice J{k+1})}=J_k$
			and hence that $\rdes{w(\upslice J{k+1})}\cap\nbd{J_k\udj J_{k+1}}
				=J_k\udj J_{k+1}$, as desired.
	\end{proof}
\end{lem}

\begin{prop}
\label{prop:connected-factorizations}
	A factorization $J_1,\dots,J_r$ is absolutely right-descent if and
		only if every factorization $J'_1,\dots,J'_{r'}$ with all $J'_k$
		connected that determines the same element of the Cartier--Foata
		monoid as $J$ is right-descent.
	\begin{proof}
		The only if direction is immediate by definition of a property
			holding absolutely.
		For the reverse direction, suppose that every factorization
			$J'_1,\dots,J'_{r'}$ with all $J'_k$ connected that determines the
			same element of the Cartier--Foata monoid as $J$ is right-descent,
			and let $J''_1,\dots,J''_k$ be a factor of $J$.
		By definition of factor, there is a factorization
			$J'_1,\dots,J'_{r'}$ with all $J'_i$ connected that determines the
			same element of the Cartier--Foata monoid as $J$ and that is
			obtained from a factorization
			$J'_1,\dots,J'_i,J''_1,\dots,J''_k,J'_{j+1},\dots,J'_r$ by
			replacing each $J''_m$ with a sequence of its connected 
			components. 
		Then by \cref{lem:rdes-slice} we have that
			$\slice{J'}ij$ is right-descent, and then repeated applications
			of \cref{lem:rdes-union} imply that $J''_1,\dots,J''_k$ is
			right-descent as desired.
	\end{proof}
\end{prop}

We now show that any factorization of an element $w$ can be modified
	into a (not necessarily admissible) absolutely bidescent
	factorization of $w$.
\begin{defn}
	We say that a factorization $J'_1,\dots,J'_{r'}$ is an
		\define{expansion} of a factorization $J_1,\dots,J_r$ if
		$w(J')=w(J)$, all $J'_k$ are connected, and there is a factorization
		$J''_1,\dots,J''_{r'}$ inducing the same element of the
		Cartier--Foata monoid as $J$ such that $J''_k\subseteq J'_k$ for all
		$k$.
	We say the expansion is \define{proper} if there exists $k$ such that
		either $J_k$ is not connected or $J''_k\subsetneq J'_k$.
\end{defn}
\begin{rmrk}
\label{rmrk:expansion-transitivity}
	The relation of $J'$ being an expansion of $J$ is transitive, since
		an expansion can only decrease the set of generating relations
		of the Cartier--Foata monoid.
\end{rmrk}

\begin{prop}
\label{prop:pre-order}
	Given a factorization $J_1,\dots,J_r$ of a Coxeter element $w$,
		replacing a connected component of $J_k$ with the connected
		component of $\rdes{w(\slice Jik)}$ or $\ldes{w(\slice Jki)}$
		containing it still gives a factorization of $w$.
	\begin{proof}
		For $J'_k\supseteq J_k$ the connected component of
			$\rdes{w(\slice Jik)}$ containing $J_k$, it suffices by definition
			of factorization to show that $w(\upslice Jk)
				=\rcos{J'_k}{w(\upslice J{k-1})}\mxml{J'_k}$.
		Since $J_k\subseteq J'_k$ implies
			$\rcos{J'_k}{w(\upslice J{k-1})}
				=\rcos{J'_k}{\qy*{\rcos{J_k}{w(\upslice J{k-1})}}}
				=\rcos{J'_k}{\qy*{\rcos{J_k}{w(\upslice J{k-1})}\mxml{J_k}}}
				=\rcos{J'_k}{w(\upslice Jk)}$,
			by \cref{prop:finitary} this is equivalent to showing
			$J'_k\subseteq\rdes{w(\upslice Jk)}$, so by definition of
			$J'_k$ it is enough to show that
			$\rdes{w(\slice Jik)}\subseteq\rdes{w(\upslice Jk)}$.
		For this, self-duality of factorizations implies
			the existence of an element $v\in W$ such that $w(\upslice Jk)=
				vw(\slice Jik)$ and $\ell(w(\upslice Jk))
					=\ell(v)+\ell(w(\slice Jik))$, so that in fact
			$\rref{w(\slice Jik)}\subseteq\rref{w(\upslice Jk)}$.
			
		The dual statement follows from self-duality of the notion of
			factorization.
	\end{proof}
\end{prop}

\begin{cor}
\label{cor:absolutely-bidescent}
	Every factorization has a maximal expansion (i.e.\ an expansion with
		no proper expansions), and such maximal expansions are absolutely
		bidescent.
	In particular, every $w\in W$ has an absolutely bidescent
		factorization.
	\begin{proof}
		By \cref{prop:connected-factorizations}, a
			factorization $J$ fails to be absolutely bidescent if and only
			if some factorization $J'$ with all $J'_k$ connected that induces
			the same element of the Cartier--Foata monoid as $J$ is such that
			$J'_k$ is not a connected component of $\rdes{w(\upslice Jk)}$ or 
			$\ldes{w(\downslice Jk)}$ for some $k$, in which case $J$ has a
			proper expansion by \cref{prop:pre-order}.
		Thus any maximal expansion of a factorization is absolutely
			bidescent.
		
		To see that maximal expansions exists, note first that every
			factorization $J$ of $w$ has a (possibly non-proper) expansion with
			all $J_k$ connected, since by definition replacing each $J_k$ with a
			sequence of its connected components yields a factorization inducing
			the same element of the Cartier--Foata monoid as $J$.
		Furthermore, each expansion $J'$ of a factorization
			$J_1,\dots,J_r$ with all $J_k$ connected, has the same length
			$r$ and satisfies $\ds\sum_k\abs{J_k}\leq\sum_k\abs{J'_k}$,
			with strict inequality if the expansion is proper.
		Since \cref{prop:defect-ideals} implies for every
			$k$ and $s\in J_k$ that $s\leq w(J)$ (by considering the mask
			$\sigma_k=s$ and $\sigma_i=e$ for $i\neq k$), we have that
			$\bigcup_{1\leq k\leq r}J_k\subseteq\set{s\in S:s\leq w(J)}$,
			where the right-hand side has size at most $\ell(w)$ since
			each $s\in S$ with $s\leq w(J)$ must appear in some fixed
			reduced word of $w(J)$ by \cref{lem:bruhat}.
		Thus, $\ds\sum_{k=1}^r\abs{J_k}\leq\ell(w(J))\cdot r$ has
			a finite upper bound, and hence a maximal factorization of $w$ must
			exist, as desired.
	\end{proof}
\end{cor}

Unfortunately, \cref{prop:rigidity} implies that the simple
	method for expanding factorizations of
	\cref{prop:pre-order} \emph{cannot} produce an admissible
	factorization when applied to an admissible factorization.
We are thus left with the following question.
\begin{ques}
\label{ques}
	Does every Coxeter element $w\in W$ have a factorization without
		proper expansions that is absolutely admissible bidescent?
\end{ques}

We end this section by giving a stronger variant of the notion of
	right-descent which we show holds absolutely if it holds, but such
	that arbitrary $w\in W$ need not have such factorizations.
This notion will be crucial for our type $A$ results in
	\cref{sec:4231-45312}.

\begin{defn}
\label{defn:strong-rdes}
	A factorization $J_1,\dots,J_r$ is \define{strong right-descent} if
		for every $k$ such that $1\leq k\leq r$, every connected component
		$J'$ of $J_k$ is a connected component of $\rdes{w(\upslice Jk)v}$
		for every $v\in W_{S\setminus\nbd{J'}}$. 
\end{defn}

\begin{rmrk}
	By choosing $v=e$, any strong right-descent factorization is
		right-descent.
\end{rmrk}

\begin{exmp}
\label{exmp:strong-descent-type-A}
	In type $A$, an interval $[a,b]$ corresponding to a connected
		component $J'$ of $J_k$ is a connected component of
		$\rdes{w(\upslice Jk)v}$ for every $v\in W_{S\setminus\nbd{J'}}$ if
		and only if $w(a)>w(x)$ when $x<a$ and $w(b)<w(x)$ when $b<x$.
\end{exmp}

\begin{lem}
\label{lem:strong-rdes-slice}
	$J_1,\dots,J_r$ is strong right-descent only if $\slice Jik$ is strong
		right-descent for every $1\leq i\leq k\leq r$.
	\begin{proof}
		That $J$ strong right-descent implies $\upslice Jk$ strong
			right-descent is immediate from the definition, so we need only
			show that $\upslice Jk$ strong right-descent implies $\slice Jik$
			strong right-descent.
		Let $J'$ be a connected component of $J_j$ for $i\leq j\leq k$, and
			let $v\in W_{S\setminus\nbd{J'}}$.
		Note that $J'$ is a connected component of
			$\rdes{w(\slice Jij)v}$ if and only if
			$J'=\nbd{J'}\cap\rref{w(\slice Jij)v}
				=\nbd{J'}\cap\qy*{v^{-1}\rref{w(\slice Jij)}v+\rref v}
				=\nbd{J'}\cap v^{-1}\rref{w(\slice Jij)}v$,
			where the last equality holds because 
			$\rref v\subseteq W_{S\setminus \nbd{J'}}$ by \cref{lem:subsystem}.

		By self-duality of factorizations (\cref{cor:quotient}), there
			exists $x\in W$ with $w(\upslice Jj)=xw(\slice Jij)$ and
			$\ell(w(\upslice Jj))=\ell(x)+\ell(w(\slice Jij))$, whence 
			$\rref{w(\slice Jij)}\subseteq\rref{w(\upslice Jj)}$ by
			\cref{lem:coset-properties}(\ref{lem:coset-unique}).
		So $\nbd{J'}\cap v^{-1}\rref{w(\slice Jij)}v
				=\nbd{J'}\cap v^{-1}\qy*{\rref{w(\upslice Jj)}
					\cap\rref{w(\slice Jij)}}v
				\allowbreak=\nbd{J'}\cap v^{-1}\rref{w(\upslice Jj)}v
					\cap v^{-1}\rref{w(\slice Jij)}v
				=J'\cap v^{-1}\rref{w(\slice Jij)}v=J'$,
			where the second to last equality holds because $J'$ is a
			connected component of $\rdes{w(\slice J1jv)}$, and
			the last equality holds because $J'\subseteq J_j\subseteq
			\rref{w(\slice Jij)}$ and $v\in S_{w\setminus\nbd{J'}}$.
	\end{proof}
\end{lem}

\begin{lem}
\label{lem:strong-rdes-cf}
	Suppose $J_1,\dots,J_r$ is a factorization and $k$ such that
		$1\leq k<r$ and $J_k$ commutes with $J_{k+1}$.
	Then $J_1,\dots,J_r$ is strong right-descent if and only if the
		factorization $J'$ given by
		$J_1,\dots,J_{k-1},J_k\udj J_{k+1},J_{k+2},\dots,J_r$ is
		strong right-descent.
	\begin{proof}
		Note that if $i<k$ we have $J_i=J'_i$ and $w(\upslice Ji)
			=w(\upslice{J'}i)$, and for $i>k+1$ we have $J_i=J'_{i-1}$
			and $w(\upslice Ji)=w(\upslice{J'}{i-1})$, so we need only
			consider the connected components of $J_k$, $J_{k+1}$, and
			$J'_k$.
		Furthermore, since $w(\upslice{J'}k)=w(\upslice J{k+1})$, we
			have that every connected component $J'$ of $J_{k+1}$ satisfies
			the strong right-descent property for $\upslice{J'}k$ if and only
			if it does for $\upslice J{k+1}$.

		Thus, we need only show that for $J'$ a connected component $J'$ of
			$J_k$ we have $\nbd{J'}\cap\rref{w(\upslice Jk)v}=J'$
			for every $v\in W_{S\setminus\nbd{J'}}$ if and only if
			$\nbd{J'}\cap\rref{w(\upslice{J'}k)}
				=\nbd{J'}\cap\rref{w(\upslice J{k+1})v}=J'$ for every
			$v\in W_{S\setminus\nbd{J'}}$.
		But $w(\upslice Jk)W_{S\setminus\nbd{J'}}=w(\upslice J{k+1})
			W_{S\setminus\nbd{J'}}$ because $w(\upslice J{k+1})W_{J_{k+1}}
				=w(\upslice Jk)W_{J_{k+1}}$ and $W_{J_{k+1}}\subseteq
			W_{S\setminus\nbd{J_k}}\subseteq W_{S\setminus\nbd{J'}}$
			since $J_k$ and $J_{k+1}$ commute by assumption.
	\end{proof}
\end{lem}

\begin{cor}
\label{cor:strong-rdes-absolute}
	A factorization $J_1,\dots,J_r$ is strong right-descent if and only
		if it is absolutely strong right-descent.
\end{cor}

\section{Tight and monotone factorizations of \pdfmath{4231}-avoiding permutations}
\label{sec:4231-45312}
In this section, we begin our analysis of the special case of Coxeter
	systems of type $A$ and show that a factorization $J_1,\dots,J_r$ of a
	\emph{$4231$-avoiding} permutation $w$ is tight only if it is
	absolutely \emph{$4231$-$45312$-avoiding} and has a contraction that
	is what we call a \emph{monotone} factorization.
The definition of monotone factorizations in this section, as well as
	its analysis in subsequent sections, will be in terms of a notion of
	\emph{permutation-directedness} of positions, which will also allow us
	to characterize \emph{$45312$-avoiding} permutations among the
	$4231$-avoiding ones.

The reasoning in this section is largely self-contained without
	relying on the details of the previous sections, but requires our
	main result \cref{thm:2} to justify restricting our attention to
	absolutely bidescent factorizations, and also requires the 
	characterization of admissibility in \cref{exmp:adm-rdes-typea}. 

We begin with \cref{lem:unbounded}, which will allow us to deduce
	that tight factorizations of $4231$-avoiding permutations are
	absolutely $4231$-avoiding and have a monotone contraction.
	
\begin{lem}
\label{lem:unbounded}
	Given a permutation $w$, a finitary subset $J\subseteq S$, and
		$v\in W_J$ such that $wv=\rcos Jw\mxml J$, suppose that for some
		interval $[a,b]$ corresponding to a connected component of $J$
		there exist positions $p_4,p_1,p_3$ satisfying
		$v(p_4)<a\leq v(p_1)<v(p_3)\leq b$ and $wv(p_1)<wv(p_3)<wv(p_4)$,
		and such that $wv(r)<wv(p_1)$ for all $r$ for which $v(p_4)<v(r)<a$.
	Then for any positive family of polynomials $(P_x)_{x\leq w}
		\in\NN[q]$ such that $\underline h^P$ is involution-invariant and
		for any subset $J_*\subseteq\rdes P\cap J$, the family 
		$(\rparp P {J_*}J_x)_{x\leq\rcos Jw\mxml J}$ defined in
		\cref{defn:right-parabolic-product} is unbounded.
	\begin{figure}[htbp]
		\centering
		\begin{tikzpicture}
	\node (vp4) at (0,0) {$v(p_4)$};
	\node (vr) at (1.5,0) {$v(r)$};
	\node (a) at (3,0) {$a\vphantom{v(p_4)}$};
	\node (vp1) at (4,0) {$v(p_1)$};
	\node (vp3) at (7,0) {$v(p_3)$};
	\node (b) at (8,0) {$b$};

	\draw[decoration={brace,mirror},decorate] (vp1.north) -- node[yshift=2ex] {$wv(r)<wv(p_1)$} (vr.north);
	\draw[decoration={brace,mirror},decorate] (vp1.south) --node[yshift=-2ex] {$wv(p_3)>wv(p_1)$} (vp3.south);
	\draw[decoration={brace,mirror},decorate] (vp4.south) --node[yshift=-2ex] {$wv(p_4)>wv(p_3)$} (vp1.south);
\end{tikzpicture}
		\caption{Setup for \cref{lem:unbounded}}
		\label{fig:lem-unbounded-setup}
	\end{figure}
	\begin{proof}
		\cref{prop:bounded-rparp} implies it suffices to show that in the
			case where $J_*=\rdes w\cap J$, $\ell(\rcos[J]{J_*}w)=0$, and
			$w_J=\mxml{J_*}$, there exists $x$ such that $x<w$ in the Bruhat
			order and $\ell(\rcos[J]{J_*}x)\geq\ell(w)-\ell(x)$.

		Since $[a,b]$ corresponds to a connected component of $J$ and we
			have $w_J=w_0(J_*)$, then as in \cref{exmp:adm-rdes-typea} there
			exists a partition
			$[a=a_1,b_1],[b_1+1=a_2,b_2],\dots,[b_{m-1}+1=a_m,b_m=b]$ of
			$[a,b]$ such that each $[a_i,b_i]$ with $a_i<b_i$ corresponds to
			a connected component of $J_*=\rdes w\cap J$, and such that
			for $p<q\in[a,b]$ we have $w(p)<w(q)$ if and only if
			$p\in[a_i,b_i]$ $q\in[a_j,b_j]$ where $i<j$. 

		In particular, $a\leq v(p_1)<v(p_3)\leq b$ and $wv(p_1)<wv(p_3)$,
			imply there exist $i,j$ such that $v(p_1)\in[a_i,b_i]$,
			$v(p_3)\in[a_j,b_j]$ and $i<j$, which in turn imply $wv(p_1)\leq
			w(a_i)$ and $wv(p_3)\leq w(a_j)$.
		Thus we have $v(p_4)<a\leq a_i\leq v(p_1)\leq b_i<a_j\leq v(p_3)
			\leq b_j\leq b$, and 
			$wv(p_1)\leq w(a_i)<w(b_j)\leq wv(p_3)<wv(p_4)$.

		Evidently, the transposition $t=(v(p_4),a_i)$ satisfies
			$t\in\rref w$ since $v(p_4)<a_i$ but $w(a_i)<wv(p_4)$ (so that
			$wt<w$ in the Bruhat order (\cref{defn:bruhat})), and we claim
			$\ell(\rcos[J]{J_*}{(wt)})\geq 1$ and $1=\len w-\len{wt}$.
		\cref{fig:lem-unbounded-argument} illustrates the argument below.
		\begin{figure}[htbp]
			\centering
			\begin{tikzpicture}
	\node (vp4) at (0,0) {$v(p_4)$};
	\node (vr) at (2.5,0) {$v(r)$};
	\node (ai) at (5,0) {$a_i\vphantom{v(p_4)}$};
	\node (aiexpl) at (5,-1.2) {$wv(r)<wv(p_1)\leq w(a_i)<wv(p_3)<wv(p_4)$};
	\node (vp1) at (6,0) {$v(p_1)$};
	\node (aj) at (8,0) {$a_j$};
	\node (vp3) at (11,0) {$v(p_3)$};

	\draw[decoration={brace,mirror},decorate] (vp4.south) -- node[yshift=-2ex] {$T_R(w)\cap T_R(t)$} (vr.south west);
	\draw[decoration={brace,mirror},decorate] (vr.south east) -- node[yshift=-2ex] {$T_R(t)\setminus T_R(w)$} (ai.south west);
	
	\draw[decoration={brace},decorate] (vp4.north) --node[yshift=2ex] {$t=(v(p_4),a_i)\in T_R(w)\cap T_R(t)$} (ai.north west);
	\draw[decoration={brace},decorate] (ai.north east) --node[yshift=2ex] {$t'=(a_i,v(p_3))\in T_R(wt)\cap W_J\setminus W_{J_*}$} (vp3.north);
	\draw[decoration={brace},decorate] (0,1.2) --node[yshift=2ex] {$tt't=(v(p_4),v(p_3))\in T_R(w)\setminus T_R(t)$} (11,1.2);
	\draw (vp4.140) -- (0,1.2);
	\draw (vp3.40) -- (11,1.2);
\end{tikzpicture}
			\caption{Argument for \cref{lem:unbounded}}
			\label{fig:lem-unbounded-argument}
		\end{figure}

		First, using the identification (\cref{exmp:type-A-tr}) of $\rref t$
			with position-based inversions of the transposition
			$(v(p_4),a_i)$, we have $\rref t=\set{(v(p_4),a_i)}
				\udj\set{(v(p_4),v(r)):v(p_4)<v(r)<a_i}
				\udj\set{(v(r),a_i):v(p_4)<v(r)<a_i}$. 
		In particular, $t'=(a_i,v(p_3))\not\in\rref w$ but 
			$t'\in W_J\setminus W_{J_*}$ because $a\leq a_i<v(p_3)<b$,
			$v(p_3)\not\in[a_i,b_i]$, and $wv(p_3)<w(a_i)$.
		Thus $t,tt't=(v(p_4),v(p_3))\in\rref w$, so 
			$t'\in(t\rref wt+\rref t)\cap W_J\setminus
				W_{J_*}=\rref{wt}\cap W_J\setminus W_{J_*}$ and hence
			$\ell(\rcos[J]{J_*}{(wt)})
				=\abs{\rref{wt}\cap W_J\setminus W_{J_*}}
				\geq\abs{\set{t'}}=1$.

		Second, by \cref{lem:self-dual} we have 
			$\ell(wt)=\abs{\rref w+\lref t}=\abs{\rref w+\rref t}$, so that
			$\ell(w)-\ell(wt)=\abs{\rref t\cap\rref w}
				-\abs{\rref t\setminus\rref w}$.
		Now, for any $r$ with $v(p_4)<v(r)<a\leq a_i$ we have by assumption
			that $wv(r)<wv(p_1)\leq w(a_i)<wv(p_3)<wv(p_4)$, so
			$(v(p_4),v(r))\in\rref w\not\ni(v(r),a_i)$.
		But also for any $r$ with $v(p_4)<a\leq v(r)<a_i$ we have
			$v(r)\in[a_k,b_k]$ with $k<i$, whence $wv(r)<w(a_i)<wv(p_4)$, and
			so again $(v(p_4),v(r))\in\rref w\not\ni(v(r),a_i)$.
		Thus,
			$\ell(w)-\ell(wt)=\abs{\set t}
				+\abs{\set{(v(p_4),v(r)):v(p_4)<v(r)<a_i}}
				-\abs{\set{(v(r),a_i):v(p_4)<v(r)<a_i}}
				=1
				+\abs{\set{r:v(p_4)<v(r)<a_i}}
				-\abs{\set{r:v(p_4)<v(r)<a_i}}
				=1$.
	\end{proof}
\end{lem}

Next, we introduce the key notions of \emph{permutation-directedness} of
	positions. 
\begin{defn}
	A permutation $w$ is said to \define{realize the pattern $4231$} in 
		positions $p_4<p_2<p_3<p_1$ if $w(p_1)<w(p_2)<w(p_3)<w(p_4)$. 
	If a permutation $w$ does not realize $4231$ in any quadruple of
		positions $p_4<p_2<p_3<p_1$, it is said that it \define{avoids
		$4231$}.

	Given a permutation $w$ and a position $p$, we define $\lambda_w(p)$
		to be the position such that $w(\lambda_w(p))
			=\max\set{w(\lambda):\lambda\leq p}$ and say that $p$ is
		\define{left-$w$-directed} if $\lambda<\lambda_w(p)$ implies
		$w(\lambda)<w(p)$ and \define{strong left-$w$-directed} if
		$\lambda_w(p)=p$.
	Similarly, we define $\rho_w(p)$ to be the position such that
		$w(\rho_w(p))=\min\set{w(\rho):p\leq\rho}$ and say that $p$ is
		\define{right-$w$-directed} if $\rho_w(p)<\rho$ implies
		$w(p)<w(\rho)$ and \define{strong right-$w$-directed} if
		$p=\rho_w(p)$.
	We say that $p$ is \define{$w$-directed} if it is left- or
		right-$w$-directed, that it is \define{$w$-uncrossed} if it is both
		strong right- and strong left-$w$-directed, and that it is
		\define{properly $w$-directed} if
		it is $w$-directed but not $w$-uncrossed.

	We say that an interval $[a,b]$ with $w(a)>\cdots>w(b)$ a
		decreasing sequence is a \define{strong right-descent interval} of
		$w$ if $a$ is properly strong left-$w$-directed and $b$ is
		properly strong right-$w$-directed, i.e.\ if $a<b$, and
		$w(\lambda)<w(a)$ and $w(b)<w(\rho)$ for all $\lambda,\rho$ such
		that $\lambda<a$ and $b<\rho$.
	We say that the interval $[a,b]$ is a \define{right-monotone} interval
		of $w$ if it is strong right-descent and every position $p\in[a,b]$
		is $w$-directed.

	We say that a factorization $J_1,\dots,J_r$ is 
		\define{strong right-descent}, resp.\ \define{right-monotone}, if
		each $J_k$ is a union of connected components
		corresponding to strong right-descent, resp.\ right-monotone,
		intervals of $w(\upslice Jk)$. 
	We define a factorization being \define{strong left-descent},
		respectively \define{left-monotone}, if its dual $J_r,\dots,J_1$
		is strong right-descent, respectively, left-monotone.
	We say the factorization is \define{strong bidescent}, resp.\
		\define{monotone}, if it is both strong right-descent and strong
		left-descent, respectively both right- and left-monotone.
\end{defn}

\begin{rmrk}
\label{rmrk:strong-descent-interval}
	Note that as in \cref{exmp:strong-descent-type-A},
		$[a,b]$ is a strong right-descent interval of $w$ if
		and only if $[a,b]$ is a connected component of
		$\rdes{wv}$ for every $v\in W_{S\setminus\nbd{J}}$ where
		$J$ is the parabolic subset corresponding to $[a,b]$.
	In particular, the above definition of strong right-descent
		factorizations coincides with that of \cref{defn:strong-rdes} and
		so a factorization is strong-right descent only if it is absolutely
		strong right-descent (\cref{cor:strong-rdes-absolute}). 
\end{rmrk}

We note here some basic properties of these notions that follow
	immediately from the definitions.
\begin{lem}
\label{lem:directed-basic}
	Given a permutation $w$ and a position $p$, we have that:
	\begin{enumerate}
		\item $p$ is strong left-$w$-directed (respectively strong
			right-$w$-directed) only if it is left-$w$-directed (respectively
			right-$w$-directed)
		\item $\lambda_w(\lambda_w(p))=\lambda_w(p)\leq p
			\leq \rho_w(p)=\rho_w(\rho_w(p))$, in particular $\lambda_w(p)$
			and $\rho_w(p)$ are strong left-$w$-directed and strong
			right-$w$-directed respectively.
		\item $\lambda_w$ and $\rho_w$ are constant on intervals $[a,b]$
			with $w(a)>\cdots>w(b)$, so $p\in[a,b]$ is left-$w$-directed with
			$\lambda_w(p)=a$ if and only if $w(x)<w(p)$ for all $x<a$, and
			right-$w$-directed with $\rho_w(p)=b$ if and only if $w(p)<w(x)$
			for all $x>b$.
	\end{enumerate}
	\begin{proof}
		Immediate from the definitions.
	\end{proof}
\end{lem}

\begin{exmp}
	For a permutation $\mxml J$ that sorts in decreasing order the
		positions within disjoint intervals $[a_i,b_i]$ with $a_i<b_i$
		(corresponding to connected components of $J$), we have that every
		position is both left- and right-$\mxml J$-directed but only the
		positions within the intervals $[a_i,b_i]$ are properly
		$\mxml J$-directed (that is, every position not in some $[a_i,b_i]$
		is $\mxml J$-uncrossed).
	Furthermore, the only strong left-$\mxml J$-directed positions are the
		left endpoints $a_i$ of the intervals, and the only strong
		right-$\mxml J$-directed positions are the
		right endpoints $b_i$ of the intervals, i.e.\ all the intervals
		$[a_i,b_i]$ are right-monotone intervals of $\mxml J$.
\end{exmp}

\begin{exmp}
	The interval $[2,4]$ is a strong right-descent interval of the
		permutation $45312$, but is not right-monotone because $3$ is
		neither left- nor right-$45312$-directed.
\end{exmp}

\begin{exmp}
	\label{exmp:rdes-int-directed}
	For a strong right-descent interval $[a,b]$ of a permutation $w$,
		there are (necessarily unique) positions $r,\ell\in[a,b]$ such that
		$p\in[a,b]$ is left- (resp.\ right-) $w$-directed if and only if
		$p\in[a,\ell]$ (resp.\ $p\in[r,b]$).
	This property and the name ``permutation directedness'' will be
		further explored in \cref{prop:monotone-structure}.
\end{exmp}

\begin{rmrk}
\label{rmrk:type-A-symmetry}
	In this section, we will make extensive use of the symmetry of
		the group $A_n$ sending $s_i$ to $s_{n+1-i}$, equivalently the
		string $i$ to the string $n+1-i$.
	For example, we will prove results only for the operator $\lambda$
		but also apply the symmetric results for $\rho$.
\end{rmrk}

We now use \cref{lem:unbounded} to relate tight factorizations of
	$4231$-avoiding permutations with monotonicity and to show that tight
	factorizations of a $4231$-avoiding permutaiton are absolutely
	$4231$-avoiding.
\begin{prop}
\label{prop:4231-unbounded}
	Given a positive family $(P_x)_{x\leq w}\in\NN[q]$ and a 
		subset $J\subseteq S$ such that $\rcos Jw\mxml J$ avoids $4231$,
		for any subset $J_*\subseteq\rdes P\cap J$ the family $(\rparp
		P{J_*}J_x)_{x\leq\rcos Jw\mxml J}$ is bounded only if each interval
		$[a,b]$ corresponding to a connected component $J'$ of $J$ not
		contained in $\rdes w$ is a right-monotone interval for $\rcos
		Jw\mxml J$ such that $\ell(\rcos[J']{\rdes w}w)=0$.
	\begin{proof}
		If $J=\emptyset$ the condition is vacuous, so we need only consider
			the case $\emptyset\neq J$, for which
			\cref{prop:bounded-rparp} yields that
			$(\rparp P{J_*}J_x)_{x\leq\rcos Jw\mxml J}$ is bounded only if
			$\ell(\rcos[J]{J_*}w)=0$ which by \cref{cor:invariance-positive} is
			equivalent to $J_*=\rdes w\cap J$ and $w_J=\mxml{J_*}$, and
			thus requires 
			$\len{\rcos[J']{\rdes w}w}\leq\len{\rcos[J]{\rdes w}w}
				=\len{\rcos[J]{J_*}w}=0$.

		The fact that $J'\not\subseteq\rdes w$ implies we have
			a partition $[a=a_1,b_1],[b_1+1=a_2,b_2],\dots,[a_m,b_m=b]$
			such that $m>1$ and each $[a_k,b_k]$ with $a_k<b_k$ corresponds to
			a connected component of $J'\cap\rdes w$. 
		Furthermore, as in \cref{exmp:adm-rdes-typea} we have for
			positions $a\leq q<r\leq b$ that $w(q)<w(r)$ is equivalent to
			$q\in[a_i,b_i]$ and $r\in[a_j,b_j]$ for some $i<j$, and also to
			$v(q)>v(r)$ where $wv=\rcos Jww_J$.

		Since $m>1$, each $p\in[a,b]$ satisfies $p\in\set{p_3,p_1}$ 
			where $p_3\in[a_i,b_i]$ and $p_1\in[a_j,b_j]$ for some $i<j$,
			i.e.\ where $w(p_1)<w(p_3)$ and $v(p_1)>v(p_3)$.
		We claim that $p_3$ is left-$wv$-directed with $\lambda_{wv}(p_3)
			=a$, and $p_1$ is right-$wv$-directed with $\rho_{wv}(p_1)=b$.
		But $J'\subseteq\rdes{wv}$ implies
			$wv(a)\geq wv(p_3)>\cdots>wv(p_1)\geq\cdots\geq wv(b)$, and so
			$p_3$ being left-$wv$-directed with $\lambda_{wv}(p_3)=a$ is
			equivalent $x<a$ (i.e.\ $v(x)<a$ since $v\in W_J$) implying
			$wv(x)<wv(p_3)$.

		Assuming for the sake of contradiction $p_3$ is not
			left-$wv$-directed, i.e.\ that there exists $v(x)<a$ with
			$wv(x)>wv(p_3)$, let $v(p_4)<a$ be maximal such that
			$wv(p_4)>wv(p_3)$.
		Then every $r$ such that $v(p_4)<v(r)<a$ has $wv(r)<wv(p_3)<wv(p_4)$
			by maximality of $v(p_4)$, so that $(v(p_4),v(r))\in\rref{w}$
			and thus $(p_4,r)=v^{-1}(v(p_4),v(r))v\in v^{-1}\rref wv\subseteq
				v^{-1}\rref wv\udj\rref v=\rref{wv}$ by \cref{prop:T_L}, so that
				$p_4<r$.
		Thus $p_4<r<p_3<p_1$ with $wv(p_4)>wv(p_3)>\max\set{wv(r),wv(p_1)}$,
			so since $wv$ avoids $4231$ and so does not realize $4231$
			in positions $p_4<r<p_3<p_1$, we have $wv(r)<wv(p_1)$ for all
			such $r$.
		But then \cref{lem:unbounded} implies that $\rparp P{J_*}J$ is
			unbounded, a contradiction, and so $p_3$ is left-$wv$-directed
			with $\lambda_{wv}(p_3)=a$ as desired.

		That $p_1$ is right-$wv$-directed with $\rho_{wv}(p_1)=b$ follows
			by symmetry (c.f.\ \cref{rmrk:type-A-symmetry}).
		Finally, since $a$ can be chosen as $p_3$ and $b$ as $p_1$,
			we get that $[a,b]$ is also a strong right-descent interval of
			$wv$, and hence right-$wv$-monotone as desired.
	\end{proof}
\end{prop}

\begin{prop}
\label{prop:4231-recurse}
	Given a permutation $w$ that realizes $4231$ and a subset 
		$J\subseteq S$, if $\rcos Jw\mxml J$ avoids $4231$, then for
		every positive family $(P_x)_{x\leq w}\in\NN[q]$ such that
		$\underline h^P$ is involution-invariant and subset
		$J_*\subseteq\rdes P\cap J$, the family
		$(\rparp P{J_*}J_x)_{x\leq \rcos Jw\mxml J}$ is unbounded, and
		moreover satisfies $\deg\qy*{\rparp P{J_*}J}_{\rcos Jw\mxml J}=0$
		only if some interval $[a,b]$ corresponding to a connected component
		of $J$ not contained in $\rdes w$ is not
		right-$\rcos Jw\mxml J$-monotone.
	\begin{proof}
		Suppose $w$ realizes $4231$ and $\rcos Jw\mxml J$ avoids $4231$, so
			that in particular $J\neq\emptyset$.
		By \cref{cor:invariance-positive},
			$\deg\qy*{\rparp P{J_*}J}_{\rcos Jw\mxml J}\leq 0$ only if
			$J_*=\rdes w\cap J$ and $w_J=\mxml{J_*}$, i.e.\ $\rref w\cap W_J
				=T\cap W_{J_*}$,
			which by \cref{prop:bounded-rparp} is also necessary
			for $\rparp P{J_*}J$ to be bounded.
		Thus, by \cref{prop:4231-unbounded}
			it suffices to show that in this case some interval $[a,b]$
			corresponding to a connected component of $J$ is not
			right-$wv$-monotone where $wv=\rcos Jw\mxml J$.

		Let $v(p_4)<v(p_2)<v(p_3)<v(p_1)$ be positions in which $w$ realizes
			$4231$, i.e.\ such that $wv(p_4)>wv(p_3)>wv(p_2)>wv(p_1$). 
		Then $\lambda_w(v(p_2))<v(p_2)$ with
			$w(\lambda_w(v(p_2))\geq wv(p_4)>wv(p_3)$ and
			$v(p_3)<w(\rho_w(v(p_2)))$ with $wv(p_3)>wv(p_1)\geq
			w(\rho_w(v(p_2)))$, so we may assume without loss of generality
			that $v(p_4)=\lambda_w(v(p_2))$ and $v(p_1)=\rho_w(v(p_3))$.

		Moreover, $(p_4,p_2)=v^{-1}(v(p_4),v(p_2))v\in v^{-1}\rref wv
			\subseteq v^{-1}\rref wv\udj\rref v=\rref{wv}$ by \cref{prop:T_L},
			and so we have $p_4<p_2$, and similarly $p_3<p_1$.
		But because $wv$ avoids $4231$, we cannot have $p_4<p_2<p_3<p_1$, so
			we must have $p_3<p_2$, and hence $(p_2,p_3)\in\lref v  
				=T\cap W_J\setminus W_{J_*}$ by \cref{prop:T_L}.
		Thus, there must exist an interval $[a,b]$
			corresponding to a connected component $J'$ of $J$ such that
			$v(p_2),v(p_3)\in[a,b]$, and in particular such $J'$ is not
			contained in $\rdes w$.

		In particular, we have a partition
			$[a=a_1,b_1],[b_1+1=a_2,b_2],\dots,[a_m,b_m=b]$
			such that $m>1$ and each $[a_k,b_k]$ with $a_k<b_k$ corresponds to
			a connected component of $J'\cap\rdes w$. 
		Furthermore, as in \cref{exmp:adm-rdes-typea} we have for
			positions $a\leq q<r\leq b$ that $w(q)<w(r)$ is equivalent to
			$q\in[a_i,b_i]$ and $r\in[a_j,b_j]$ for some $i<j$. 

		Thus $v(p_2)\in[a_i,b_i]$ and $v(p_3)\in[a_j,b_j]$ with $i<j$, and
			so $a\leq v(p_4)<v(p_2)$ would require
			$v(p_4)\in[a_i,b_i]\cap[a_j,b_j]=\emptyset$.
		Consequently, $v(p_4)<a$ and hence $p_4<a$ because $[a,b]$
			corresponds to a connected component of $J$ where $v\in W_J$.
		Since $wv(p_4)=w(v(p_4))>w(v(p_2))=wv(p_2)$, it follows that 
			$p_2\in[a,b]$ is not left-$wv$-directed with
			$\lambda_{wv}(p_2)=a$. 
		Similarly, we also have $v(p_1)>b$, $p_1<b$, and 
			$wv(p_1)=w(v(p_1))<w(v(p_2))=wv(p_2)$ so that
			$p_2$ is not right-$wv$-directed with $\rho_{wv}(p_2)=b$.
		Thus, the interval $[a,b]$ is not right-$wv$-monotone as desired.
	\end{proof}
\end{prop}

\begin{cor}
\label{cor:4231-monotone-slice}
	A factorization $J_1,\dots,J_r$ of a $4231$-avoiding permutation $w$
		is tight only if for each $k$ such that $1\leq k\leq r$,
		$w(\upslice Jk)$ is $4231$-avoiding and every connected component
		of $J_k$ not contained in $\rdes{w(\upslice J{k-1})}$ is
		right-$w(\upslice Jk)$-monotone.
	In particular, every tight factorization of a $4231$-avoiding
		permutation is absolutely $4231$-avoiding with a monotone
		contraction.
	\begin{proof}
		Since tight factorizations are absolutely tight by
			\cref{cor:factors}, \cref{prop:4231-recurse}
			implies that if $w(\upslice Jk)$ realizes $4231$ then $J$ is tight
			only if it also realizes 4231, giving the first claim.
		The second claim then follows by \cref{prop:4231-unbounded}.
		For the supplemental claim, we have by \cref{cor:factors,thm:2}
			that a minimal contraction $J$ of a tight factorization of a
			$4231$-avoiding permutation is absolutely bidescent, and so has
			$\rdes{w(\upslice J{k-1})}=\rdes{\upslice J{k-1}}\not\supseteq
			J'_k$ for each $k$ and connected component $J'_k$ of $J_k$
			by \cref{prop:rdescent-inductive-step} and definition of
			contraction, so that the main claim implies that such a
			minimal contraction is monotone as desired. 
	\end{proof}
\end{cor}

The remainder of the proof that tight factorizations of
	$4231$-avoiding permutations must also absolutely avoid $45312$ does
	not require additional arguments about boundedness, and instead
	only deals with the notions of permutation-directedness, strong
	(right-)descent intervals, and (right-)monotone intervals.

\begin{defn}
	A permutation $w$ is said to \define{realize the pattern $45312$} in
		positions $p_4<p_5<p_3<p_1<p_2$ if
		$w(p_5)>w(p_4)>w(p_3)>w(p_2)>w(p_1)$, and
		\define{avoids $45312$} if it does not realize $45312$.
\end{defn}

$45312$-avoidance is closely related to $w$-directedness, and for
	$4231$-avoiding permutations is characterized in terms of it.

\begin{lem}
\label{lem:45312-directed}
	A position $p$ is not $w$-directed if and only if there exist positions
		$p_4<\lambda_w(p)$ and $p_2>\rho_w(p)$ such that $w$ realizes
		$45312$ in positions $p_4<\lambda_w(p)<p<\rho_w(p)<p_2$.
	In particular, if either $w$ avoids $4231$ or $[\lambda_w(p),\rho_w(p)]$
		is a strong right-descent interval of $w$, then $w$ realizes $45312$
		in some positions $p_4<p_5<p<p_1<p_2$ if and only if $p$ is not
		$w$-directed.
	\begin{proof}
		The main claim is immediate from the definition of directedness.
		For the supplemental claims, we show that if $w$ realizes
			$45312$ in positions $p_4<p_5<p<p_1<p_2$, then
			$p_4<\lambda_w(p)$ and $\rho_w(p)<p_2$, so that $w$ also realizes
			$45312$ in positions $p_4<\lambda_w(p)<p<\rho_w(p)<p_2$.
		Indeed, we have by definition that $w(\lambda_w(p))\geq w(p_5)>w(p_4)
			>w(p)$ and $p_4<p_5<p$, so that if $\lambda_w(p)\leq p_4$ then
			$w(\lambda_w(p))>w(p_5)>w(p_4)>w(p)$ and
			$\lambda_w(p)<p_4<p_5,p$ so that $[\lambda_w(p),p]$ is not a
			right-descent interval and $w$ realizes $4231$ in positions
			$\lambda_w(p)<p_4<p_5<p$.
		Thus, $p_4<\lambda_w(p)<p<\rho_w(p)<p_2$ with $w(\lambda_w(p))\geq
			w(p_5)>w(p_4)>w(p)>w(p_2)>w(p_1)\geq w(\rho_w(p))$.
	\end{proof}
\end{lem}

\begin{cor}
\label{cor:strong-rdes-monotone-iff-45312}
	A strong right-descent factorization $J_1,\dots,J_r$ is right-monotone
		if and only if $w(\upslice Jk)$ avoids $45312$ for every $k$ such
		that $1\leq k\leq r$.
	\begin{proof}
		By induction on $k$, it suffices to prove that if $\upslice J{k-1}$
			is right-monotone and $w(\upslice Jj)$ avoids $45312$ for
			$j\leq k-1$, we have that $\upslice Jk$ is right-monotone if and
			only if $w(\upslice Jk)$ avoids $45312$.

		Now if $w(\upslice Jk)$ realizes $45312$ in positions
			$p_4<p_5<p_3<p_1<p_2$, then since $w(\upslice J{k-1})$ avoids
			$45312$ we must have that $p_3\in[a,b]$ for an interval
			$[a,b]$ corresponding to a connected component $J'$ of $J_k$
			not contained in $\rdes{w(\upslice J{k-1})}$.
		But then since $\upslice Jk$ is strong right-descent, we have that
			$[a=\lambda_w(p_3),b=\rho_w(p_3)]$ is a strong right-descent
			interval of $w(\upslice Jk)$, so that by
			\cref{lem:45312-directed} $p_3$ is not $w$-directed and
			$\upslice Jk$ is not right-monotone.

		Conversely, if $\upslice Jk$ is not right-monotone, then by
			definition some position $p$ is not $w(\upslice Jk)$-directed,
			so that $w(\upslice Jk)$ realizes $45312$ by
			\cref{lem:45312-directed}.
	\end{proof}
\end{cor}

A key element of our proofs is understanding how directedness interacts
	with the map $w\mapsto wv=\rcos Jw\mxml J$ for a finitary subset
	$J\subseteq S$.
\begin{lem}
\label{lem:directed-up}
	Given a permutation $w$ and a finitary subset $J\subseteq S$, let 
		$v\in W$ be such that $wv=\rcos Jw\mxml J$.
	Then for every position $p$, $wv(\lambda_{wv}(p))
		\geq w(\lambda_w(v(p)))$ with equality if and only if
		$v(\lambda_{wv}(p))\leq v(p)$, in which case $p$ is
		left-$wv$-directed if $v(p)$ is left-$w$-directed.
	Furthermore, equality fails to hold only if $p,v(p)\in[a,b]$ for some
		interval $[a,b]$ corresponding to a connected component of $J$ not
		contained in $\rdes w$ such that $a$ is strong left-$wv$-directed.
	\begin{proof}
		If $x$ is such that $v(x)\leq v(p)$ and $wv(x)\geq wv(p)$, then
			$\rref w\ni(v(x),v(p))=v(x,p)v^{-1}$ so that
			$(x,p)\in v^{-1}\rref wv\subseteq\rref{wv}$, and thus
			$x\leq p$.
		In particular, we have that $\max\set{wv(x):x\leq p}
			\geq\max\set{wv(x):v(x)\leq v(p)}=\max\set{w(x):x\leq v(p)}$,
			so that $wv(\lambda_{wv}(p))\geq w(\lambda_w(v(p)))$.

		Now if $wv(\lambda_{wv}(p))=w(\lambda_w(v(p)))$ then
			$v(\lambda_{wv}(p))=\lambda_w(v(p))\leq v(p)$, and conversely
			$v(\lambda_{wv}(p))\leq v(p)$ implies $w(\lambda_w(v(p)))
			\geq wv(\lambda_{wv}(p))$, which by the main inequality implies
			$wv(\lambda_{wv}(p))=w(\lambda_w(v(p)))$ as desired.

		Note that $x<\lambda_{wv}(p)$ implies $wv(x)<wv(\lambda_{wv}(p))$ by
			definition, so that $(x,\lambda_{wv}(p))\not\in\rref{wv}
				=v^{-1}\rref wv\udj\rref v$, and thus $v(x)<v(\lambda_{wv}(p))$.
		But now, if equality holds and $v(p)$ is left-$w$-directed, then
			$x<\lambda_{wv}(p)$ implies $v(x)<v(\lambda_{wv}(p))\leq v(p)$
			so that $wv(x)=w(v(x))<w(v(p))=wv(p)$ and thus $p$ is
			left-$wv$-directed.

		Finally, if equality doesn't hold, then $v(\lambda_{wv}(p))>v(p)$
			and $\lambda_{wv}(p)\leq p$ so that $(\lambda_{wv}(p),p)\in\rref v
			\subseteq W_J$, and thus there is an interval $[a,b]$ corresponding
			to a connected component of $J$ such that $\lambda_{wv}(p),p,v(p)
				\in[a,b]$.
		Since by definition $wv(a)>\cdots>wv(b)$ we thus have that
			$\lambda_{wv}(p)\in[a,b]$ implies $\lambda_{wv}(p)=a$ so that
			$a$ is strong left-$wv$-directed, and also $wv(a)
				=wv(\lambda_{wv}(p))>w(\lambda_w(v(p)))$ implies that 
			$[a,b]$ is not contained in $\rdes w$.
	\end{proof}
\end{lem}

\begin{cor}
\label{cor:directed-induct}
	Given a permutation $w$, a finitary subset $J\subseteq S$, let
		$v\in W$ be such that $wv=\rcos Jw\mxml J$.
	Then a position $p$ is such that $v(p)$ is $w$-directed but
		$p$ is not $wv$-directed only if $p\in[a,b]$ for a strong
		right-descent interval of $wv$ corresponding to a connected
		component of $J$ not contained in $\rdes w$ which is not
		right-$wv$-monotone.
	\begin{proof}
		This follows immediately from
			\cref{lem:directed-up} and its horizontal-symmetric dual.
	\end{proof}
\end{cor}

\begin{cor}
\label{cor:4231-45312-monotone-slice}
	A factorization $J_1,\dots,J_r$ of a $4231$-avoiding permutation is
		tight only if for each $k$ such that $1\leq k\leq r$,
		$w(\upslice Jk)$ is $4231$-$45312$-avoiding and every connected
		component of $J_k$ not contained in $\rdes{w(\upslice J{k-1})}$
		is right-$w(\upslice Jk)$-monotone.
	\begin{proof}
		By \cref{cor:4231-monotone-slice} we need only show
			that $w(\upslice Jk)$ is $45312$-avoiding, which by
			\cref{lem:45312-directed} is equivalent to showing that
			each position $p$ is $w(\upslice Jk)$-directed.
		But \cref{cor:directed-induct} implies this inductively
			since every connected component of $J_k$ not contained in
			$\rdes{w(\upslice J{k-1})}$ is right-$w(\upslice Jk)$-monotone.
	\end{proof}
\end{cor}

\begin{thm}
\label{thm:4231-45312-avoid}
	A factorization $J_1,\dots,J_r$ of a $4231$-avoiding permutation $w$
		is tight only if it is absolutely $4231$-$45312$-avoiding and has a
		monotone contraction.
	\begin{proof}
		Because tight factorizations are absolutely tight and
			the generating relations of the Cartier--Foata monoid preserve
			$w(J)$, the main result follows from
			\cref{cor:4231-45312-monotone-slice} and its dual,
			noting that $4231$ and $45312$ are self-dual.
		That such factorizations have a monotone contraction follows from
			\cref{cor:4231-monotone-slice}.
	\end{proof}
\end{thm}

\section{Strong bidescent factorizations of \pdfmath{45312}-avoiding permutations}
\label{sec:45312}
In this section, we characterize strong right-descent and strong
	bidescent factorizations in terms of their \emph{heap of connected
	components}, and then use this characterization to show that monotone
	factorizations are exactly those strong bidescent factorizations that
	satisfy the equivalent conditions of being right-monotone,
	being left-monotone, or $w$ being $45312$-avoiding. 
We also show that that a factorization is monotone only if it is
	absolutely monotone.

\begin{defn}
\label{defn:heap}
	Given a factorization $J_1,\dots,J_r$ we define its \define{heap of
		connected components} to be the disjoint union of the family of
		connected components of each $J_1,\dots,J_r$.
	We denote by $[i_f,j_f]$ the interval corresponding to the connected
		component $f\in H(J)$ and by $\lvl f$ the index such that $f$ is a
		connected component of $J_{\lvl f}$.
	Furthermore, we denote by juxtaposition the partially-defined function
		$H(J)\times\set{+\vec u_1,+\vec u_2,-\vec u_1,-\vec u_2}
			\dashrightarrow H(J)$ given by
		\begin{enumerate}
			\item $f+\vec u_1$ is the unique (if it exists) element $h\in
				H(J)$ such that $\lvl f<\lvl h$ and $j_f\in[i_h,j_h]$ but
				$j_f\not\in[i_g,j_g]$ for any $g\in H(J)$ such that
				$\lvl f<\lvl g<\lvl h$;
			\item $h-\vec u_1$ is the unique (if it exists) element $f\in
				H(J)$ such that $i_h\in[i_f,j_f]$ but $i_h\not\in[i_g,j_g]$ for
				any $g\in H(J)$ such that $\lvl f<\lvl g<\lvl h$;
			\item $f+\vec u_2$ is the unique (if it exists) element $h\in
				H(J)$ such that $\lvl f<\lvl h$ and $i_f\in[i_h,j_h]$ but
				$i_f\not\in[i_g,j_g]$ for any $g\in H(J)$ such that
				$\lvl f<\lvl g<\lvl h$;
			\item $h-\vec u_2$ is the unique (if it exists) element $f\in
				H(J)$ such that $j_h\in[i_f,j_f]$ but $j_h\not\in[i_g,j_g]$ for
				any $g\in H(J)$ such that $\lvl f<\lvl g<\lvl h$.
		\end{enumerate}
	The partially-defined function $H(J)\times\set{+\vec u_1,+\vec u_2,
		-\vec u_1,-\vec u_2}\dashrightarrow H(J)$ then induces a unique 
		\define{partially defined monoid action} 
		$H(J)\times M_4\dashrightarrow H(J)$ such that
		$(f+\vec v)\pm\vec u_i=f+(\vec v\pm\vec u_i)$ where $M_4$ is the
		monoid freely generated by
		$\set{+\vec u_1,+\vec u_2,-\vec u_1,-\vec u_2}$, i.e.\ whose
		elements are finite sequences of elements of 
		$\set{+\vec u_1,+\vec u_2,-\vec u_1,-\vec u_2}$.
	For convenience, for any $m\in\ZZ_{\geq0}$ we abbreviate the elements 
		$\underbrace{+\vec u_i\cdots+\vec u_i}_m$ and 
		$\underbrace{-\vec u_i\cdots-\vec u_i}_m$ of $M_4$ 
		as $+m\vec u_i$ and $-m\vec u_i$.
\end{defn}

\begin{figure}[htbp]
	\centering
	\begin{tikzpicture}[yscale=0.7]

	\draw (0,0) -- node[below] {$e$} (8,0);
	\draw (-1,2.5) -- node[above] {$f$} (3,2.5);
	\draw (-2,5) -- node[above] {$g$} (5,5);
	\draw (6.5,5) -- node[above] {$h$} (10,5);

	\draw[dashed,-Latex] (0,0.2) --node[left] {$+u_2$} (0,2.35);
	\draw[dashed,-Latex] (-1,2.65) --node[left] {$+u_2$} (-1,4.8);
	\draw[dashed,-Latex] (3,2.65) --node[left] {$+u_1$} (3,4.8);
	\draw[dashed,-Latex] (3,2.35) --node[left] {$-u_2$} (3,0.2);
	\draw[dashed,-Latex] (5,4.8) --node[left] {$-u_2$} (5,0.2);
	\draw[dashed,-Latex] (6.5,4.8) --node[left] {$-u_1$} (6.5,0.2);
	\draw[dashed,-Latex] (8,0.2) --node[left] {$+u_1$} (8,4.8);

\end{tikzpicture}
	\caption{Visualization of the partially-defined monoid action.}
	\label{fig:monoid-action}
\end{figure}
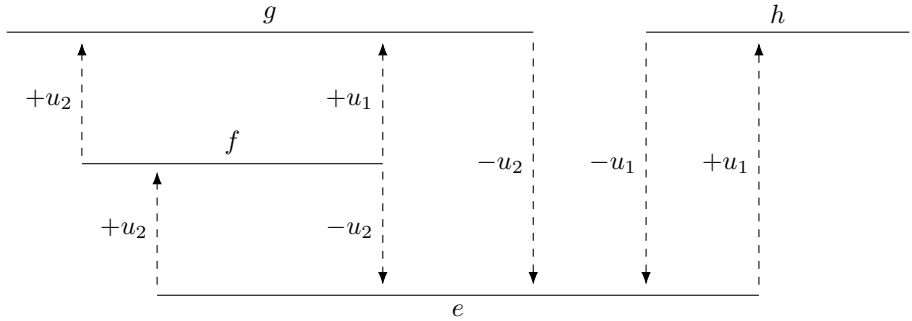

\begin{rmrk}
	Our heap of connected components is a more highly structured instance
		of a Viennot heap \cite{Viennot}.
\end{rmrk}
\begin{lem}
\label{lem:strong-right-descent}
	Given a permutation $w$ and a finitary subset $J\subseteq S$,
		let $v\in W$ be such that $wv=\rcos Jw\mxml J$.
	Then for $[a,b]$ an interval corresponding to a connected component
		of $J$, the following are equivalent: (1) $a$ is strong
		left-$wv$-directed; (2) $\lambda_w(b)=v(a)$; (3)
		$v(a)$ is strong left-$w$-directed; (4) $\lambda_w(p)\in[a,b]$
		for some $p$.
	\begin{proof}
		Note first that $v\in W_J$, so that $x<a$ if and only if $v(x)<a$
			and $x\in[a,b]$ if and only if $v(x)\in[a,b]$ for every $x$.

		(1) implies (2): if $\lambda_{wv}(a)=a$, then for every
			$v(x)<b$, we have either $x<a$, in which case $w(v(x))=wv(x)
				<wv(a)=w(v(a))$ where the inequality is since $a$ is strong
			left-$wv$-directed, or $a\leq x\leq b$, in which case
			$w(v(x))=wv(x)\leq wv(a)=w(v(a))$ where the inequality is since
			$J\subseteq\rdes{wv}$.

		(2) implies (3) by \cref{lem:directed-basic}, and (3)
			implies (4) since $v(a)\in[a,b]$ and thus $\lambda_w(v(a))
				=v(a)\in[a,b]$.

		(4) implies (1): If $\lambda_w(p)\in[a,b]$, then for all $x<a$ we
			have $v(x)<a\leq\lambda_w(p)$ so $wv(x)<w(\lambda_w(p))
				=wv(v^{-1}(\lambda_w(p)))\leq wv(a)$, where the last inequality
			is because $v^{-1}(\lambda_w(p))\in[a,b]$ and $J\subseteq\rdes{wv}$.
	\end{proof}
\end{lem}

We begin with the following proposition linking the property of being
	strong-descent to the monoid action.

\begin{prop}
\label{prop:staircase}
	Given a factorization $J_1,\dots,J_r$ and an integer $0\leq j\leq r$,
		the following are equivalent:
		\begin{enumerate}
			\item $i_g$ is a strong left-$w(\upslice J{\lvl g})$-directed
				position for each $g\in H(J)$ such that $\lvl g\leq j$;
			\item if $f=g-\vec u_2$ for $f,g\in H(J)$ with $\lvl g\leq j$,
				then $g=f+m\vec u_2$ for some $m\in\ZZ_{>0}$;
			\item for all $f,g\in H(J)$ with $\lvl f<\lvl g\leq j$,
				$\wlvl f(i_f)=\wlvl g(i_g)$ if and only if $g=f+m\vec u_2$
				for some $m\in\ZZ_{>0}$.
		\end{enumerate}
	\begin{figure}[htbp]
		\centering
		\begin{tikzpicture}[yscale=0.7]

	\draw (0,0) -- node[below] {$f$} (7,0);
	\draw (-4,6) -- node[above] {$g$} (5,6);

	\draw (-1, 2) -- (3, 2);
	\draw (-3, 4) -- (2, 4);

	\draw[dashed,-Latex] (5,5.8) --node[left] {$-u_2$} (5,0.2);

	\draw[dashed,-Latex] (0,0.2) --node[left] {$+u_2$} (0,1.8);
	\draw[dashed,-Latex] (-1,2.2) --node[left] {$+u_2$} (-1,3.8);
	\draw[dashed,-Latex] (-3,4.2) --node[left] {$+u_2$} (-3,5.8);

	\node at (3.6, 3.1) {$\Longleftarrow$};

\end{tikzpicture}
		\caption{Illustration of \cref{prop:staircase}(2)}
		\label{fig:staircase}
	\end{figure}
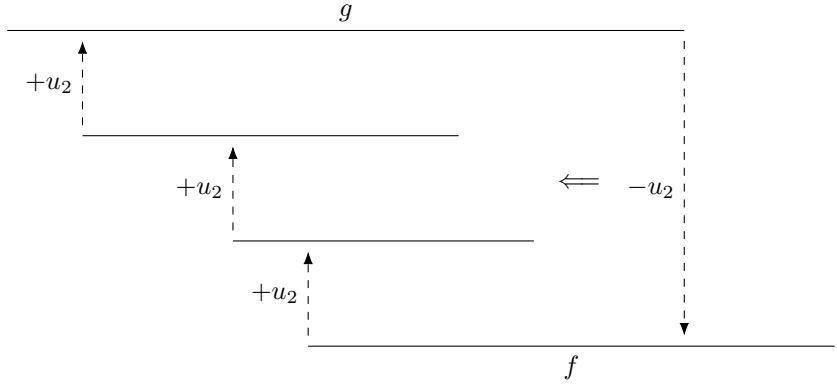
	\begin{proof}
		We argue by induction on $j$, with the base case $j=0$ immediate
			since all the conditions are vacuous in this case.

		Now, assume that (1-3) hold for $j-1$, and let $g\in H(J)$ be such
			that $\lvl g\leq j$.

		If $f=g-\vec u_2$ does not exist, then $j_g$ is
			$w(\upslice Ji)$-uncrossed for every $i<\lvl g$, and so
			$j_g$ is strong left-$w(\upslice J{\lvl g-1})$-directed and
			\cref{lem:strong-right-descent} implies $i_g$ is strong
			left-$w(\upslice J{\lvl g})$-directed with $\wlvl g(i_g)
			=w(\upslice J{\lvl g-1})(j_g)=j_g$.
		Furthermore, $j_g$ being $w(\upslice J{\lvl g-1})$-uncrossed
			also implies that $j_g\neq \wlvl f(i_f)$ for every $f\in H(J)$
			with $\lvl f<\lvl g$.

		We have to show that $i_g$ is strong left-$w(\upslice J{\lvl g})$-%
			directed if and only if $f=g-\vec u_2$ implies $g=f+m\vec u_2$
			for some $m\in\ZZ_{\geq 0}$, in which case $w(\upslice J{\lvl g})
			(i_g)=w(\upslice J{\lvl f})(i_f)$, and otherwise $w(\upslice J{\lvl g})
			=j_g$.

		Otherwise, if $f\in H(J)$ such that $f=g-\vec u_2$ does exist,
			then $\lvl f<j$, and so the inductive hypothesis implies that $i_f$
			is a strong left-$w(\upslice J{\lvl f})$-directed position, and
			so in particular $i_f=\lambda_{w(\upslice J{\lvl f})}(j_g)$ since
			$j_g\in[i_f,j_f]$ for $[i_f,j_f]$ a descent interval.
		Since $j_g\not\in [i_e,j_e]$ for all $e\in H(J)$ with $\lvl f
			<\lvl e<\lvl g$, \cref{lem:directed-up} implies that
			$\wlvl g(i_g)\leq \wlvl g(\lambda_{\wlvl g}(j_g))
				=\wlvl f(i_f)$, and so
			$i_g$ is strong left-$\wlvl g$-directed if and only if
			$\wlvl f(i_f)=\wlvl g(i_g)$.

		Now, let $h$ be the element of $H(J)$ with maximal $\lvl h$ such that
			$\lvl h<\lvl g$ and $\wlvl f(i_f)=\wlvl h(i_h)$, so that 
			by the inductive hypothesis $i_h$ is strong left-$\wlvl h$-directed
			and $h=f+m_h\vec u_2$ for some $m_h\in\ZZ_{\geq 0}$ (where
			the $m_h=0$ case is equivalent to $h=f$).
		It thus remains to show that $\wlvl h(i_h)=\wlvl g(i_g)$
			(equivalently $\wlvl g(i_g)\geq \wlvl h(i_h)$ by the previous
			paragraph) if and only if there exists $m\in\ZZ_{>0}$ with
			$g=h+m\vec u_2$.
		By definition of $h$, we have that either $h+\vec u_2$ does not
			exist or has $\lvl{h+\vec u_2}\geq\lvl g$, so $g=h+m\vec u_2$
			for some $m\in\ZZ_{>0}$ if and only if $g=h+\vec u_2$, which
			is again by definition of $h$ equivalent to $i_h\in[i_g,j_g]$.
		But by \cref{lem:directed-up} $i_h
			=\lambda_{(\upslice J{\lvl g-1})}(j_g)$ with
			$w(\upslice J{\lvl g-1})(i_h)=\wlvl h(i_h)$, so that
			$\wlvl g(i_g)\geq \wlvl h(i_h)$ if and only if $i_h\in[i_g,j_g]$
			as desired.
		\end{proof}
\end{prop}

\begin{cor}
\label{cor:strong-rdes-char}
	For a factorization $J_1,\dots,J_r$ the following are equivalent:
	\begin{enumerate}
		\item $J$ is strong right-descent;
		\item $f=g-\vec u_i$ for $f,g\in H(J)$ implies there exists
			$m\in\ZZ_{>0}$ such that $f+m\vec u_i=g$;
		\item for $f,g\in H(J)$ with $\lvl f<\lvl g$,
			$\wlvl f(i_f)=\wlvl g(i_g)$ if and only if there exists
			$m\in\ZZ_{>0}$ such that $f+m\vec u_2=g$, and $\wlvl f(j_f)
			=\wlvl g(j_g)$ if and only if there exists $m\in\ZZ_{>0}$ such
			that $f+m\vec u_1=g$.
	\end{enumerate}
	\begin{proof}
		Follows from \cref{prop:staircase} and its horizontal-%
			symmetric dual (in the sense of \cref{rmrk:type-A-symmetry})
			with $j=r$.
	\end{proof}
\end{cor}

We now get the following characterization of strong bidescent
	factorizations in terms of invertibility of the monoid action.
\begin{cor}
\label{cor:strong-bidescent-char}
	For a factorization $J_1,\dots,J_r$, the following are equivalent.
	\begin{enumerate}
		\item $J$ is strong bidescent;
		\item a pair $f,g\in H(J)$ satisfy $f=g\pm\vec u_i$ only if
			there exists $m\in\ZZ_{>0}$ such that $g=f\mp m\vec u_i$;
		\item a pair $f,g\in H(J)$ satisfy $f=g\pm\vec u_i$ if and only
			if $g=f\mp\vec u_i$;
		\item a pair $f,g\in H(J)$ satisfy $\wlvl f(i_f)=\wlvl g(i_g)$
			if and only if there exists $m\in\ZZ$ such that $f=g+m\vec u_2$,
			and satisfy $\wlvl f(j_f)=\wlvl g(j_g)$ if and only if there exists
			$m\in\ZZ$ such that $f=g+m\vec u_1$.
	\end{enumerate}
	\begin{proof}
		\cref{cor:strong-rdes-char} and its dual imply the
			equivalence of (1) and (2), and (4).
		That (3) implies (2) is immediate, so it remains to show that
			(2) implies (3).
		Let $f=g\pm\vec u_i$, so that there exists $m\in\ZZ_{>0}$ with
			$g=f\mp m\vec u_i=f\mp(m-1)\vec u_i\mp \vec u_i$.
		Then since for $h=f\mp(m-1)\vec u_i$ we have $g=h\mp\vec u_i$,
			there exists $m'\in\ZZ_{>0}$ with $h=g\pm m'\vec u_i
			=g\pm \vec u_i \pm(m'-1)\vec u_i=f\pm (m'-1)\vec u_i$,
			and thus $f\mp(m-1)\vec u_i=h=f\pm (m'-1)\vec u_i$.
		But the first equality implies $\lvl f\leq \lvl h$ with equality
			if and only if $m=1$, and the second implies $\lvl h\leq \lvl f$
			with equality if and only if $m'=1$, so that $\lvl h=\lvl f$
			and $m=1=m'$.
	\end{proof}
\end{cor}

Using the above characterization, we show that strong bidescent
	factorizations that avoid $45312$ do so absolutely, and so
	we get a characterization of strong bidescent factorizations
	which are monotone.
\begin{prop}
\label{prop:strong-bides-45312-absolute}
	Let $J_1,\dots,J_r$ be a strong bidescent factorization such that
		$w(J)$ avoids $45312$.
	Then for every $k$ such that $1\leq k\leq r$ we have that
		$w(\upslice Jk)$ avoids $45312$.
	\begin{proof}
		Assume for the sake of contradiction that there is some minimal
			$1\leq k\leq r$ with $w(\upslice Jk)$ realizing $45312$
			in positions $p_4<p_5<p_3<p_1<p_2$.
		By minimality of $k$, there is $e\in H(J)$ with $\lvl e=k$
			and $p_3\in[i_e,j_e]$, and this interval $[i_e,j_e]$ contains
			at least one of $p_5$ and $p_1$ (as otherwise $w(\upslice J{k-1})$
			would realize $45312$ with the same strings).
		Since $J$ is strong right-descent, we know that $[i_e,j_e]$ is a
			strong right-descent interval of $w(\upslice Jk)$, and thus
			$p_4<i_e=\lambda_{w(\upslice Jk)}(p_3)\leq p_5<p_3
				<p_1\leq\rho_{w(\upslice Jk)}(p_3)=j_e<p_2$.
		Furthermore, we can take $p_4$ to be $\lambda_{w(\upslice Jk)}(
			i_e-1)$ and $p_2$ to be $\rho_{w(\upslice Jk)}(j_e+1)$,
			which by \cref{lem:directed-up} implies that there are
			components $f,g\in H(J)$ with $i_f=p_4$ and $j_g=p_2$,
			$\lvl f, \lvl g\leq\lvl e$, and $p_4\not\in [i_{f'},j_{f'}]$
			for $f'\in H(J)$ with $\lvl f\leq \lvl{f'}\leq\lvl e$, and
			similarly for $g$.
		But now, by \cref{cor:strong-bidescent-char}, since
			$J$ is strong bidescent we have
			that the maximal $m$ such that $f+m\vec u_2$ exists has
			$w(i_{f+m\vec u_2})=\wlvl{f+m\vec u_2}(i_{f+m\vec u_2})
				=\wlvl f(i_f)=w(\upslice Jk)(p_4)$, so that
			$i_{f+m\vec u_2}$ is strong left-$w$-directed, and thus
			$w^{-1}(w(\upslice Jk)(p_4))<w^{-1}(w(\upslice Jk)(p_5))$.
		Similarly, we get $w^{-1}(w(\upslice Jk)(p_1))
			<w^{-1}(w(\upslice Jk)(p_2))$, so that $w$ realizes
			$45312$ with the same strings as $w(\upslice Jk)$.
	\end{proof}
\end{prop}

\begin{cor}
\label{cor:monotone-45312}
	For a strong bidescent factorization $J_1,\dots,J_r$ of $w$, the
		following are equivalent:
		\begin{enumerate}
			\item $w$ avoids $45312$;
			\item $J$ is absolutely 45312-avoiding, i.e.\ for every
				factor $J'$ of $J$, $w(J')$ avoids $45312$.
			\item $J$ is absolutely monotone;
			\item $J$ is monotone;
			\item $J$ is right-monotone;
			\item $J$ is left-monotone;
		\end{enumerate}
	\begin{proof}
		\cref{cor:strong-rdes-absolute} implies that $J$ strong
			bidescent implies that it is absolutely strong bidescent, so
			\cref{prop:strong-bides-45312-absolute} and its
			dual give the implication (1) implies (2).
		\cref{cor:strong-rdes-monotone-iff-45312} and its
			dual then give the implications (2) implies (3), (5) implies
			(1), and (6) implies (1).
		Finally, the implications (3) implies (4) and (4) implies (5) and
			(6) are by definition.
	\end{proof}
\end{cor}

\section{Rigidity properties of strong bidescent and monotone factorizations}
\label{sec:rigid}
In this section, we show that a permutation $w$ has a monotone
	factorization if and only if it is $4231$-$45312$-avoiding, in which
	case the factorization is admissible and any two minimal monotone
	factorizations are equivalent. 
In particular, we conclude that $4231$-$45312$-avoiding permutations
	have at most one minimal tight factorization up to equivalence.
Moreover, we give an effective procedure for determining a minimal
	monotone factorization of a $4231$-$45312$-avoiding permutation,
	thereby obtaining an embedding of $4231$-$45312$-avoiding
	permutations in the Cartier--Foata monoid generalizing that of
	$321$-avoiding permutations (i.e.\ of short-braid avoiding elements).

The key to the above results are the surprisingly strong
	rigidity properties satisfied by strong bidescent and monotone
	factorizations, which we now describe.

\begin{cor}
\label{cor:strong-bides-no-intermediate}
	Given a strong bidescent factorization $J_1,\dots,J_r$ and
		$e,f\in H(J)$ with $f=e+\vec u_1$ (respectively
		$f=e+\vec u_2$) is such that every $g\in H(J)$
		with $\lvl e<\lvl g<\lvl f$ has $i_g>j_e$ or
		$j_g<i_f$ (respectively $i_g>j_f$ or $j_g<i_e$), so that
		for every $k$ with $\lvl e\leq k<\lvl f$ and $p\in
		[i_f,j_e]$ (respectively $p\in[i_e,j_f]$) we have
		$w(\upslice Jk)(p)=\wlvl e(p)$.
	\begin{proof}
		Assume for the sake of contradiction
			there exists $g\in H(J)\setminus\set{e,f}$ with minimal
			$\lvl g$ such that $i_g\leq j_e$ and $j_g\geq i_f$
			(respectively $i_g\leq j_f$ and $j_g\geq i_e$).
		We consider two cases.
		If $i_g\leq i_f$ (respectively $j_g\geq j_f$)
			then $i_g\leq i_f\leq j_g$ and $f-\vec u_1\neq e$ (respectively
			$j_g\geq j_f\geq i_g$ and $f-\vec u_2\neq e$),
			contradicting \cref{cor:strong-bidescent-char}.
		Similarly, if $i_g>i_f$ (respectively $j_g<j_f$) then by minimality
			of $\lvl g$ we have $g-\vec u_1=e$ (respectively $g-\vec u_2=e$),
			so that by \cref{cor:strong-bidescent-char} we have
			$e+\vec u_1=g$ (respectively $e+\vec u_2=g$), contradicting
			the definition of $f$.
	\end{proof}
\end{cor}

\begin{cor}
\label{cor:strong-bides-notmin}
	A strong bidescent factorization $J_1,\dots,J_r$ fails to be minimal
		if and only if there are $e,f\in H(J)$ such that $e+\vec u_1
		=f=e+\vec u_2$, in which case $i_e=i_f$, $j_e=j_f$, and
		every $g\in H(J)\setminus\set{e,f}$
		with $\lvl e\leq\lvl g\leq\lvl f$ has $j_g<i_e$ or
		$i_g>j_e$, so that for all $k$ with $\lvl e \leq k\leq\lvl f$
		and all $p\in[i_e,j_e]$ it holds that $\wlvl f(p)
			=w(\upslice Jk)(p)=\wlvl e(p)$.
\end{cor}
\begin{proof}
	The ``in which case'' claims are immediate from
		\cref{cor:strong-bidescent-char,cor:strong-bides-no-intermediate}.
	For the if direction, we have that if $e+\vec u_1=f=e+\vec u_2$,
		then the ``in which case'' claims hold and so
		$J_1,\dots,J_{\lvl f-1},J_{\lvl f}\setminus [i_f,j_f]$
		is a contraction of $\upslice J{\lvl f}$, and thus
		by \cref{cor:absolutely-minimal} we have that $J$ is not minimal.

	For the only if direction, suppose there is a right-contraction
		$J'$ of $J$ and a component $f\in H(J)\setminus H(J')$ with
		$\lvl f$ minimal (the case of a left-contraction is analogous
		by choosing $e\in H(J)\setminus H(J')$ with $\lvl e$ maximal).
	Then since $[i_f,j_f]$ is contained in $\rdes{\upslice J{\lvl f-1}}$,
		which is a union of disjoint commuting connected subsets, we have
		that also $[i_f,j_f]$ is contained in some connected component
		of the overlap subset $\overlap Jk{\lvl f}$ for some $k<\lvl f$, and
		thus contained in $[i_e,j_e]$ for some $e\in H(J)$ with
		$\lvl e=k$.
	Then $f-\vec u_1=e=f-\vec u_2$ and the claims follows from
		\cref{cor:strong-bidescent-char,cor:strong-bides-no-intermediate}.
\end{proof}

\begin{rmrk}
	The above proposition implies that non-minimal strong bidescent
		factorizations are simply strong bidescent factorizations
		with redundant repeated components, so that to understand all
		strong bidescent factorizations it suffices to understand only
		the minimal ones.
\end{rmrk}

We first note that minimal strong bidescent factorizations have a
	simple structure in which strings go directly from $c\in H(J)$ to
	$e\in H(J)$ only for very restricted $c$.
\begin{cor}
\label{cor:strong-bides-structure}
	Let $J_1,\dots,J_r$ be a minimal strong bidescent factorization and
		$e\in H(J)$.
	Then if there exist $p\in[i_e,j_e]$ and some $c\in H(J)$ with maximal
		$\lvl c$ among all $c\in H(J)$ with $\lvl c<\lvl e$ and
		$p\in[i_c,j_c]$, then either (1) $c=e-\vec u_1$, (2) $c=e-\vec u_2$, or
		(3) $f=e-\vec u_1$ and $g=e-\vec u_2$ exist, and $c=f-\vec u_2
		=g-\vec u_1$.
	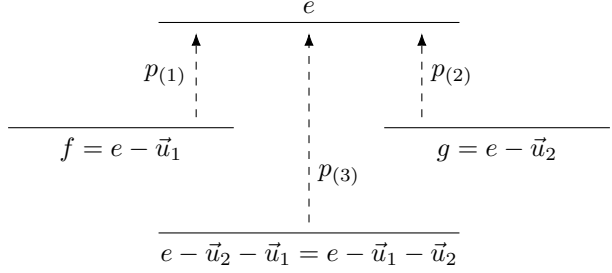
\begin{figure}[htbp]
		\centering
		\begin{tikzpicture}[yscale=0.7]

	\draw (2,4) -- node[above] {$e$} (6,4);

	\draw (2,0) -- node[below] {$e-\vec u_2-\vec u_1=e-\vec u_1-\vec u_2$} (6,0);

	\draw (0,2) -- node[below] {$f=e-\vec u_1$} (3,2);
	\draw (5,2) -- node[below] {$g=e-\vec u_2$} (8,2);

	\draw[dashed,-Latex] (2.5,2.2) -- node[left] {$p_{(1)}$} (2.5,3.8);
	\draw[dashed,-Latex] (4,0.2) -- node[right] {$p_{(3)}$} (4,2) -- (4,3.8);
	\draw[dashed,-Latex] (5.5,2.2) -- node[right] {$p_{(2)}$} (5.5,3.8);
\end{tikzpicture}
		\caption{Possible values of $c\in H(J)$ in \cref{cor:strong-bides-structure}}
		\label{fig:strong-bides-structure}
	\end{figure}
	\begin{proof}
		Let $c$ and $p$ be as in the claim, and note that
			$p\in[i_e,j_e]\cap[i_c,j_c]$ implies that
			$i_c\leq j_e$ and $j_c\geq i_e$.

		We claim that if $f=e-\vec u_1$ does not exist or $f$ exists
			but $\lvl c\geq \lvl f$, then $c=e-\vec u_2$.
		Indeed, in either case, then for every $m\in\ZZ_{\geq 0}$
			such that $c+m\vec u_2$ exists and $\lvl c\leq \lvl{c+m\vec u_2}
			<\lvl e$, we have that $i_e<i_{c+m\vec u_2}\leq i_c\leq j_e$ (as
			otherwise $e-\vec u_1$ would exist and not be $f$).
		Thus, for maximal such $m$ we have that $(c+m\vec u_2)+\vec u_2
			=e$, and so by \cref{cor:strong-bidescent-char} we have
			$c+m\vec u_2=e-\vec u_2$ which implies $i_{c+m\vec u_2}\leq i_c
			\leq j_e \leq j_{c+m\vec u_2}$, and in particular that
			$[i_{c+m\vec u_2},j_{c+m\vec u_2}]\cap [i_e,j_e]
				\supseteq [i_c,j_c]\cap [i_e,j_e]$, so that the existence of
			$p$ for which $c$ is maximal with $\lvl c<\lvl e$ and
			$p\in[i_e,j_e]\cap[i_c,j_c]$ implies $m=0$ and $c=e-\vec u_2$.

		An analogous argument implies that if $g=e-\vec u_2$ does not exist
			or $\lvl c\geq \lvl g$, then $c=e-\vec u_1$, so the only case
			that remains is that $f=e-\vec u_1$ and $g=e-\vec u_2$ exist, and
			$\lvl c<\min\set{\lvl f,\lvl g}$.
		Then since $\lvl c$ is maximal such that $p\in[i_c,j_c]$
			and $\lvl c<\lvl e$, we have that $p\in(j_f,i_g)$.
		Now, let $d\in H(J)$ have maximal $\lvl d$ among all
			$d\in H(J)\setminus\set{f,g}$ with $\lvl d<\lvl e$ such that there
			exists some $q\in (j_f,i_g) \cap[i_d,j_d]$ (note that $c$ has
			maximal $\lvl c$ for the particular $p\in(j_f,i_g)$, but $d$ is
			maximal over the set of all $q\in(j_f,i_g)$).
		Since $d\not\in\set{f,g}$, by the previous argument we have
			that $\min\set{\lvl f,\lvl g}>\lvl d\geq\lvl c$.
		Furthermore, by maximality of $\lvl d$ that
			since $d\neq f=e-\vec u_1$ we have $e\neq d+\vec u_1$
			and $i_d\leq j_f$, and analogously since
			$d\neq g=e-\vec u_2$ we have $e\neq d+\vec u_2$
			and $j_d\geq i_g$.
		Then since $\lvl c\leq \lvl d$ and $p\in (j_f,i_g)
			\subseteq [i_d,j_d]$, we have that $c=d$, and
			it remains to show that $d=f-\vec u_2=g-\vec u_1$.

		Since $\lvl d<\lvl f$ and $j_f\in[i_d,j_d]$ we have that
			$f-\vec u_2$ exists, and furthermore by minimality of
			$J$ and \cref{cor:strong-bidescent-char,cor:strong-bides-notmin}
			we have $j_f+1\in [i_{f-\vec u_2},j_{f-\vec u_2}]$, from
			which by maximality of $\lvl d$ and the fact that
			$j_f+1\in[i_d,j_d]$ implies $d=f-\vec u_2$ as desired.
		That $d=g-\vec u_1$ also holds by symmetry.
	\end{proof}
\end{cor}

Next, we show that monotone factorizations have even greater rigidity
	than minimal strong bidescent factorizations.
\begin{prop}
\label{prop:monotone-structure}
	Let $J_1,\dots,J_r$ be a minimal monotone factorization,
		$e\in H(J)$, and $v\in W$ such that $\wlvl e=\wlvl[-1]ev$.
	Then there are (necessarily unique) positions
		$r_e,l_e\in[i_e,j_e]$ such that $p\in[i_e,j_e]$ is left-
		(resp.\ right-) $\wlvl e$-directed if and only if $p\in[i_e,l_e]$
		(resp.\ $p\in[r_e,j_e]$) as in \cref{exmp:rdes-int-directed},
		which satisfy:
		\begin{enumerate}
			\item if $f=e-\vec u_1$ exists, then $j_e-l_e=j_f-i_e+1$, every
				$p\in(l_e,j_e]$ has $v(p)\in[i_e,j_f]$, and $r_f\leq i_e$, that
				is, the set of strings (in the right-canonical mask
				$\omega_R(J)$) in $e$ not left-$\wlvl e$-directed is equal to
				the set of strings coming from $f$, each of which was in a
				position that was right-$\wlvl f$-directed.  If $e-\vec u_1$
				does not exist, then $j_e-l_e=0$, i.e.\ $l_e=j_e$, so that every
				string in $e$ is left-$\wlvl e$-directed.
			\item if $g=e-\vec u_2$ exists, then $r_e-i_e=j_e-i_g+1$, every
				$p\in[i_e,r_e)$ has $v(p)\in[i_g,j_e]$, and $l_g\geq j_e$, that
				is, the set of strings (in the right-canonical mask
				$\omega_R(J)$) in $e$ not right-$\wlvl e$-directed is equal to
				the set of strings coming from $g$, each of which was in a
				position that was left-$\wlvl g$-directed.  If $e-\vec u_2$ does
				not exist, then $r_e-i_e=0$, i.e.\ $r_e=i_e$, so that every
				string in $e$ is right-$\wlvl e$-directed.
			\item if $p\in[i_e,j_e]$ is such that $p\not\in[i_g,j_g]$
				and $p\not\in[i_f,j_f]$, then $p$ is $\wlvl[-1]e$-uncrossed.
		\end{enumerate}
	\begin{proof}
		We first prove (3).
		\begin{figure}[htbp]
			\centering
			\begin{tikzpicture}[yscale=0.7]

	\draw (2,4) -- node[above] {$e$} (6,4);

	\draw (2,0) -- node[below] {$f-\vec u_2=c=e-\vec u_1$} (6,0);

	\draw (0,2) -- node[below] {$f=e-\vec u_1$} (3,2);
	\draw (-1,3) -- node[above] {$f'$} (1,3);

	\draw (5,2) -- node[below] {$g=e-\vec u_2$} (8,2);
	\draw (7,4) -- node[above] {$g'$} (9,4);

	\node at (-0.2, 2) {$\mathbf{4}$};
	\node at (-1.2, 3) {$\mathbf{4}$};
	\node at (1.8, 0) {$\mathbf{4}$};

	\node at (4.8, 2) {$\mathbf{5}$};
	\node at (1.8, 4) {$\mathbf{5}$};

	\node at (4, -0.9) {$\mathbf 3$};
	\node at (4, 4.8) {$\mathbf 3$};

	\node at (3.2, 2) {$\mathbf{1}$};
	\node at (6.2, 4) {$\mathbf{1}$};

	\node at (8.2, 2) {$\mathbf{2}$};
	\node at (9.2, 4) {$\mathbf{2}$};
	\node at (6.2, 0) {$\mathbf{2}$};

	\draw[dashed,-Latex] (4,0.2) -- node[right] {$p$} (4,2) -- (4,3.8);
\end{tikzpicture}
			\caption{Argument of (3) in \cref{prop:monotone-structure}}
			\label{fig:monotone-structure-3}
		\end{figure}
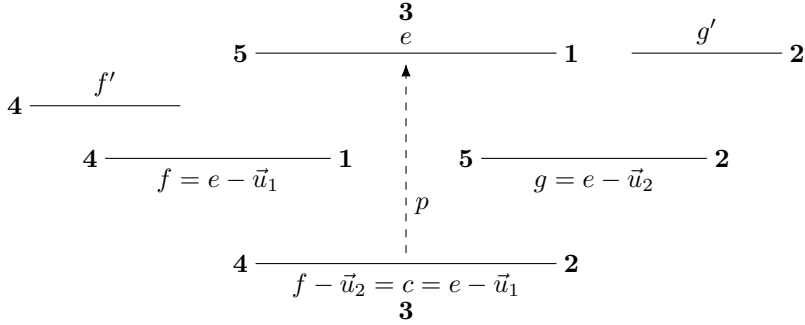
		Assume for the sake of contradiction there is $p\in[i_e,j_e]$
			with $p\not\in[i_f,j_f]$ and $p\not\in[i_g,j_g]$ that is not
			$\wlvl[-1]e$-uncrossed.
		Then by \cref{cor:strong-bides-structure} $f$ and $g$
			exist and there is a $c\in H(J)$ with $p\in[i_c,j_c]$,
			$\lvl c$ maximal among all such $c$ with $\lvl c<\lvl e$,
			and this $c$ satisfies $f=c+\vec u_2$ and $g=c+\vec u_1$.
		Letting $f'=f+m_f\vec u_2$ for $m_f\in\ZZ_{\geq 0}$ maximal such that
			$\lvl{f'}\leq \lvl e$, we have by
			\cref{cor:strong-bidescent-char,cor:strong-bides-no-intermediate}
			that $i_{f'}\leq i_f<i_e$ and
			$\wlvl {f'}(i_{f'})=\wlvl f(i_f)=\wlvl c(i_c)>\wlvl c(p)
				=\wlvl e(v^{-1}(p))$.
		Similarly letting $g'=g+m_g\vec u_1$ for $m_g\in\ZZ_{\geq 0}$
			maximal such that $\lvl{g'}\leq \lvl e$, we have that
			$j_{g'}\geq j_g>j_e$ and
			$\wlvl{g'}(j_{g'})=\wlvl g(j_g)=\wlvl c(j_c)<\wlvl c(p)
				=\wlvl e(v^{-1}(p))$,
			and so since $[i_e,j_e]$ is a strong right-descent interval of
			$\wlvl e$, we have by \cref{lem:directed-up} that $\wlvl e$
			realizes $45312$ in positions $i_{f'}<i_e<v^{-1}(p)<j_e<j_{g'}$,
			contradicting \cref{cor:monotone-45312}.
		This argument is illustrated in \cref{fig:monotone-structure-3}.

		It remains to prove (1), since (2) follows from (1) by symmetry.
		We consider two cases depending on whether $f=e-\vec u_1$ exists.
		If $e-\vec u_1$ does not exist, then there is no
			$c\in H(J)$ with $\lvl c<\lvl e$ and $i_c\leq i_e-1<i_e\leq j_c$,
			so that every $p\geq i_e$ has $\wlvl e(p)\geq i_e$ and thus
			every $p\in[i_e,j_e]$ is left-$\wlvl e$-directed, e.g.\ $l_e=j_e$
			as desired.

		Now suppose $f=e-\vec u_1$ exists.
		We first prove $r_f\leq i_e$, which we illustrate in
		\cref{fig:monotone-structure-directed}.
		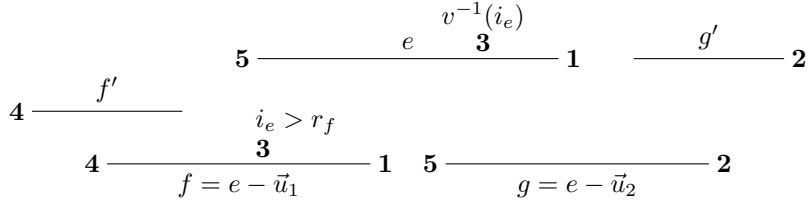
\begin{figure}[htbp]
			\centering
			\begin{tikzpicture}[yscale=0.7]

	\draw (2,4) -- node[above] {$e$} (6,4);

	\draw (0,2) -- node[below] {$f=e-\vec u_1$} (3.5,2);
	\draw (-1,3) -- node[above] {$f'$} (1,3);

	\draw (4.5,2) -- node[below] {$g=e-\vec u_2$} (8,2);
	\draw (7,4) -- node[above] {$g'$} (9,4);

	\node at (-0.2, 2) {$\mathbf{4}$};
	\node at (-1.2, 3) {$\mathbf{4}$};

	\node at (4.3, 2) {$\mathbf{5}$};
	\node at (1.8, 4) {$\mathbf{5}$};

	\node[anchor=west] at (1.85, 2.3) {$\mathbf 3$};
	\node[anchor=west] at (1.85, 2.8) {$i_e>r_f$};

	\node at (5, 4.3) {$\mathbf 3$};
	\node at (5, 4.8) {$v^{-1}(i_e)$};

	\node at (3.7, 2) {$\mathbf{1}$};
	\node at (6.2, 4) {$\mathbf{1}$};

	\node at (8.2, 2) {$\mathbf{2}$};
	\node at (9.2, 4) {$\mathbf{2}$};

\end{tikzpicture}
			\caption{Argument of $r_f\leq i_e$ in \cref{prop:monotone-structure}}
			\label{fig:monotone-structure-directed}
		\end{figure}
		Assume for the sake of contradiction that $r_f>i_e$, so that
			$i_e$ is not right-$\wlvl f$-directed, that is, there exists
			$q>j_f$ with $\wlvl f(i_e)>\wlvl f(q)$.
		Then we must have that $j_f+1$ is not $\wlvl f$-uncrossed, so
			that \cref{cor:strong-bidescent-char} implies $g=e-\vec u_1$
			exists, which since $J$ is strong bidescent implies
			$\rho_{\wlvl g}(j_g)=j_g$ and thus $\wlvl f(i_f)>\wlvl g(i_g)$.
		But then letting $f'$ and $g'$ as in the proof of (3), we have
			(again as in the proof of (3)) that $\wlvl e$ realizes $45312$
			in positions $i_{f'}<i_e<v^{-1}(i_e)<j_e<j_{g'}$, contradicting
			\cref{cor:monotone-45312}.

		We now prove the remaining claims.
		Since $i_e$ is strong left $\wlvl e$-directed, $p\in[i_e,j_e]$ is
			left-$\wlvl e$-directed if and only if $\wlvl e(p)>\wlvl
			e(\lambda_{\wlvl e}(i_e-1))$.
		Then since $J$ is minimal, we have $i_e-1\in [i_f,j_f]$,
			and furthermore by \cref{cor:strong-bides-no-intermediate}
			no $c\in H(J)$ with $\lvl f<\lvl c\leq\lvl e$ has
			$[i_e-1,i_e]\subseteq [i_c,j_c]$, and thus
			\cref{lem:directed-up} implies $\wlvl e(\lambda_{\wlvl e}(i_e-1))
				=\wlvl f(\lambda_{\wlvl f}(i_e-1))=\wlvl f(i_f)$.
		Thus, a position $p\in[i_e,j_e]$ is left-$\wlvl e$-directed if
			and only if $\wlvl e(p)=\wlvl[-1]e(v(p))<\wlvl f(i_f)$ where
			$v(p)\in[i_e,j_e]$.
		We conclude by considering cases on the range of $v(p)$.
		If (i) $v(p)\in[i_e,j_f]$, we have by
			\cref{cor:strong-bides-no-intermediate} that
			$\wlvl[-1]e(v(p))=\wlvl f(v(p))>\wlvl f(i_f)$.
		If (ii) $v(p)\in (j_f,i_g)$, we have by (3) of this
			proposition that $v(p)$ is $\wlvl[-1]e$-uncrossed so
			$\wlvl[-1]e(v(p))=v(p)>\wlvl f(i_f)$.
		If (iii) $v(p)\in[i_g,j_e]$, since we have shown (by symmetry) that
			$l_g\geq j_e$, we have that $v(p)$ is left-$\wlvl g$-directed, and
			thus $\wlvl[-1]e(v(p))=\wlvl g(v(p))>\wlvl f(i_f)$.
		Together, we get that that $p$ is
			left-$\wlvl e$-directed if and only if $\wlvl[-1]e(v(p))<\wlvl
			f(i_f)$ if and only if $v(p)\in[i_e,j_f]$ as desired.
	\end{proof}
\end{prop}

\begin{cor}
\label{cor:monotone-avoids-4231}
	If $J_1,\dots,J_r$ is a monotone factorization of $w\in W$, then
		$J$ is absolutely $4231$-avoiding.
	\begin{proof}
		Since monotone factorizations are absolutely monotone
			by \cref{cor:monotone-45312}, it suffices to prove
			that $w(J)$ avoids $4231$, which we do by induction on $r$,
			with the base case $r=0$ being immediate.
		For the inductive step, suppose $w(\upslice J{r-1})$
			avoids $4231$, and assume for the sake of contradiction
			that $w$ realizes $4231$ in positions
			$p_4<p_2<p_3<p_1$, where without loss of generality
			we have that $p_4=\lambda_w(p_4)$ and $p_1=\rho_w(p_1)$.
		Letting $v\in W$ be such that $w=w(\upslice J{r-1})v$,
			we have since $w(\upslice J{r-1})$ avoids $4231$ that
			either $v(p_4)>v(p_2)$ or $v(p_3)>v(p_1)$, assume
			without loss of generality (by symmetry) that it is the former.
		Thus, there is a connected component $e\in H(J)$
			with $\wlvl e=r$ and $i_e\leq v(p_2)<v(p_4)\leq j_e$,
			which since $J$ is strong right-descent implies
			also $j_e<p_3$.
		In particular, we have that $p_4$ and $p_2$ are not
			right-$w$-directed, so \cref{prop:monotone-structure} implies that
			$v(p_4),v(p_2)\in [i_g,j_e]$ for $g=e-\vec u_2$.
		However, $[i_g,j_e]$ is a right-descent interval of
			$\wlvl[-1]e=w(\upslice J{r-1})$ by
			\cref{cor:strong-bides-no-intermediate},
			contradicting the fact that $w(p_2)
				=w(\upslice J{r-1})(v(p_2))
				<w(\upslice J{r-1})(v(p_4))
				=w(p_4)$.
	\end{proof}
\end{cor}

\begin{cor}
\label{cor:monotone-admissible}
	If $J_1,\dots,J_r$ be a monotone factorization of $w\in W$, then
		$J$ is (absolutely) admissible.
	\begin{proof}
		Since admissible factorizations are absolutely admissible
			by \cref{cor:factors}, it suffices to show that
			$J$ is admissible.
		Suppose that two strings in the right-canonical mask
			$\omega_R(J)$ have a right-meeting at $e\in H(J)$
			and are in positions $p$ and $q$ after, and let $v\in W$ be
			such that $\wlvl e=\wlvl[-1]e v$, so that $v(p)<v(q)$
			are not both in the same connected component of
			$\rdes{\wlvl[-1]e}$ by definition of right-meeting.
		If either $v(p)$ or $v(q)$ are $\wlvl[-1]e$-uncrossed, then	
			the right-meeting is not a right-defect since the strings
			were in-order, so assume neither $v(p)$ nor $v(q)$
			are $\wlvl[-1]e$-uncrossed.
		Then since $v(p)$ and $v(q)$ are not in the same connected
			component of $\rdes{\wlvl[-1]e}$,
			\cref{prop:monotone-structure} implies that
			$v(p)\in[i_h,j_h]$ for $h=e-\vec u_1$, $v(q)\in[i_g,j_g]$
			for $g=e-\vec u_2$, and $v(p)$
			is right-$\wlvl[-1]e$-directed, whence
			$\wlvl[-1]e(v(p))<\wlvl[-1]e(v(q))$.
		In particular, $v(p)$ and $v(q)$ are in-order in
			$\wlvl[-1]e$, and so their right-meeting at $e$
			is not a right-defect.
		Since $e\in H(J)$ and the right-meeting were arbitrary,
			$\omega_R(J)$ has no right-defects and thus $J$ is
			admissible as desired.
	\end{proof}
\end{cor}

\begin{cor}
\label{cor:monotone-abs-all}
	Monotone factorizations are absolutely monotone,
		absolutely $4231$-$45312$-avoiding, and absolutely admissible.
\end{cor}

We now show that every $4231$-$45312$ avoiding permutation $w$
	has a unique (up to equivalence in the Cartier--Foata monoid)
	minimal monotone factorization.

\begin{lem}
\label{lem:strong-right-descent-pattern}
	Given a non-identity permutation $w$, the position $a=
		\max\set{p:\lambda_w(p)=p\neq\rho_w(p)}$ exists and is such that
		$[a,\rho_w(a)]$ is a strong right-descent interval of
		$w$ unless $w$ realizes $4231$ in positions $p_4<p_2<p_3=p_2+1<p_1$.
	In particular, every non-identity $4231$-avoiding permutation has a
		strong right-descent interval.
	\begin{proof}
		Note that if $q$ is $w$-uncrossed then for any $p\leq q\leq r$ we have
			$w(p)\leq w(q)\leq w(r)$ so that $(p,r)\not\in\rref w$, and thus
			$w$ not the identity implies the existence of an interval
			$[p,r]$ such that $p<r$ and no position $q\in[p,r]$ is
			$w$-uncrossed, and moreover this implies that 
			for a maximal such interval we have
			that $\set{w(q):q\in[p,r]}=\set{q:q\in[p,r]}$.
		In particular, $q=w^{-1}(r)\in[p,r]$ is strong left-$w$-directed
			but not $w$-uncrossed, i.e.\ $\lambda_w(q)=q\neq\rho_w(r)$, so
			that the position $a$ exists as desired.

		Now, let $b=\max\set{q\geq a:w(a)>w(a+1)>\dots>w(q)}$ be
			the maximal position such that $a\leq b$ and the sequence
			$w(a),\dots,w(b)$ is decreasing, so that $\rho_w(a)=\rho_w(b)$.
		Then $b$ fails to be strong right-$w$-directed if and only if
			there exists $p_1$ such that $b<p_1$ and $w(b)>w(p_1)$.
		Such a $p_1$ exists only if the position $b+1$ also exists,
			which by definition of $b$ implies that $w(b+1)>w(b)>w(p_1)$,
			so that $b+1<p_1$ and $\rho_w(b+1)\neq b+1$.
		By maximality of $a$, we have that $\lambda_w(b+1)\neq b+1$, so that
			since $[a,b]$ is a descent interval with $a$ strong
			left-$w$-directed we have that $\lambda_w(b+1)=a$.
		In particular, we get that $w$ realizes $4231$ in positions
			$a<b<b+1<p_1$, as desired.
	\end{proof}
\end{lem}

\begin{lem}
\label{lem:maximal-strong-rdes-intervals}
	If $J_1,\dots,J_r$ is a strong right-descent factorization of $w$,
		then $[a,b]$ is a right-descent interval of $w$ if and only if there
		exists $e\in H(J)$ such that $[a,b]\subseteq[i_e,j_e]$ and such that
		for every $f\in H(J)$ with $\lvl f>\lvl e$ it holds that
		$[a,b]\cap[i_f,j_f]=\emptyset$.
	In particular, $[a,b]$ is a strong right-descent interval of $w$
		if and only if there exists $e\in H(J)$ such that $e+\vec u_1$
		and $e+\vec u_2$ don't exist and $[i_e,j_e]=[a,b]$.
	\begin{proof}
		For the if direction, we have if $e$ is as in the lemma then
			for every $p\in[a,b]$ we have $w(p)=\wlvl e(p)$ so that
			$\wlvl e(a)>\cdots \wlvl e(b)$ implies the claim.

		Let $e\in H(J)$ have maximal $\lvl e$ such that
			$[a,b]\cap[i_e,j_e]\neq\emptyset$ (such an $e$ exists since $a$ is
			not $w$-uncrossed).
		We claim that $i_e\leq a\leq b\leq j_e$.
		If $i_e>a$ then $a\leq i_e-1<i_e\leq b$,
			and by maximality of $\lvl e$ we have $w(i_e-1)=\wlvl e(i_e-1)$ and
			$\wlvl e(i_e)=w(i_e)\geq w(b)$, but since $J$ is strong
			right-descent that $\wlvl e(i_e)>\wlvl e(i_e-1)$,
			contradiction.
		The case $j_e<b$ is analogous, completing the proof of the main
			claim.

		The in particular claim follows immediately from the main claim
			by \cref{lem:directed-up}.
	\end{proof}
\end{lem}

\Cref{prop:monotone-structure} is useful for proving uniqueness of
	minimal monotone factorizations, since it can be interpreted as
	saying for any component $e\in H(J_1,\dots,J_r)$ with $\lvl e=r$,
	one can ``read off'' $e-\vec u_1$ and $e-\vec u_2$ from $w(J)$.
The next lemma gives the complementary result showing that following
	this procedure results in an element of $w$ to which this procedure
	can again be applied.

\begin{lem}
\label{lem:4231-45312-ind}
	Let $J\subseteq S$ corresponding to an interval $[a,b]$ and
		$\rcos Jw\in\rcos JW$ be such that $[a,b]$ is a strong right-descent
		interval of the $4231$-$45312$-avoiding permutation $\rcos Jww_0(J)$.
	Then for $r,l\in[a,b]$ such that $p\in[a,b]$ is left- (resp.~right-)
		$\rcos Jww_0(J)$-directed if and only if $p\in[a,l]$
		(resp.~$p\in[r,b]$) as in \cref{exmp:rdes-int-directed}, and
		$J_*=[a,a+b-l)\cup (b-r+a,b]$, we have that $\rcos Jww_0(J_*)$ is
		$4231$-$45312$ avoiding.
	Furthermore, every position $p\in[a+b-l,b-r+a]$ is
		$\rcos Jww_0(J_*)$-uncrossed.
	\begin{proof}
		For convenience, write $w=\rcos Jww_0(J_*)$ and $v$ such that
			$wv=\rcos Jww_0(J)$.
		Since $wv$ is $4231$-avoiding and $wv=\rcos Jww_0(J)$,
			we have that $w$ realizes $4231$ only if
			there are positions $p_4<p_2<p_3<p_1$ realizing
			$4231$ in $w$ with $p_2,p_3\in[a,b]$.
		Since $[a,b]$ is a monotone interval of $\rcos Jww_0(J)$,
			this implies that each position $q\in[a,b]$ is such that
			either $w(q)>w(\lambda)$ for every $\lambda<a$ or
			$w(q)<w(\lambda)$ for every $\lambda>b$, and hence at least
			one of $p_4$ and $p_1$ must also be in $[a,b]$, assume
			without loss of generality (by symmetry) that it is $p_4$.
		Then since $w_J=w_0(J_*)$ we have as in \cref{exmp:adm-rdes-typea} 
			that $w(p_2)<w(p_3)$ implies that $p_2\not\in(b-r+a,b]$,
			so that $w(p_4)<w(p_2)$ implies $p_4,p_2\in[a,a+b-l)$, which
			implies $w(p_4)<w(p_3)$, contradiction.

		For the main claim, it remains to show that $w$ avoids $45312$.
		Since $wv$ avoids $45312$, we must have that if $w$
			realizes $45312$ in positions $p_4<p_5<p_3<p_1<p_2$,
			then ether $p_4,p_5\in[a,b]$ or $p_1,p_2\in[a,b]$, so
			assume without loss of generality (by symmetry) that
			$p_4,p_5\in[a,b]$.
		We claim $p_3\not\in[a,b]$: if $p_3\in[a,b]$, then
			since $w_J=w_0(J_*)$ we have as in \cref{exmp:adm-rdes-typea}
			that $p_4<p_5$ implies $p_5\not\in[a,a+b-l)$, but then
			$p_4<p_5<p_3$ implies $w(p_4)<w(p_3)$, contradiction.
		Then since $[a,b]$ is a strong right-descent interval
			of $wv$, we have that $b<p_3$ and $wv(v^{-1}(p_4))
				=w(p_4)>w(p_3)>w(b)$, and thus $v^{-1}(p_4)$ is
			not right-$wv$-directed.
		But then $v^{-1}(p_5)<v^{-1}(p_4)$ implies that
			$v^{-1}(p_5)$ is also not right-$wv$-directed, so
			by definition of $J_*$ we have that
			$p_5,p_4\in(b-r+a,b]$ and so $w(p_5)>w(p_4)$,
			contradiction.

		For the supplemental claim, we have for each $p\in[a+b-l,b-r+a]$
			that $v^{-1}(p)$ is both left- and right- $wv$-directed, so that
			$w(q)<w(p)<w(r)$ for every $q<a$ and $b<r$, and thus by
			\cref{exmp:adm-rdes-typea} we have that $p$ is $w$-uncrossed as
			desired.
	\end{proof}
\end{lem}

\begin{prop}
\label{prop:monotone-algorithm}
	For every $4231$-$45312$-avoiding permutation $w\in W$, there is a
		unique (up to equivalence in the Cartier--Foata monoid) minimal
		monotone factorization $J_1,\dots,J_r$ of $w$, and such a
		factorization is produced by the following specialization of the
		admissible right-descent factorization procedure.
		\begin{algorithmic}[1]
				\State Initialize $x \gets w$ and an empty sequence of subsets
					of $S$.
				\While{$x\neq e$, or equivalently $\rdes{x}\neq\emptyset$,}
					\State Let $C_1,\dots,C_m$ correspond to the strong right-descent
						components of $x$.
					\State Add $\bigcup_i C_i$ to the beginning of the
						sequence.
					\State For each $i$, let $J_{*,i}$ be as in
						\cref{lem:4231-45312-ind}, and let $J_*=\bigcup_i J_{*,i}$.
					\State $x\gets \rcos Jxw_0(J_*)$.
				\EndWhile
				\State \Return the sequence of subsets.
		\end{algorithmic}
	\begin{proof}
		We prove this by strong induction on $\ell(w)$, where the base case
			$\ell(w)=0$, i.e.~$w=e$, is immediate.

		For the inductive step, by \cref{lem:strong-right-descent-pattern}
			there exist $C_1,\dots,C_m$ as in the proposition, and letting
			$J_*$ as in the proposition we have by \cref{lem:4231-45312-ind}
			that $x=\rcos Jww_0(J_*)<w$ is $4231$-$45312$-avoiding.
		Thus, by the inductive hypothesis there is a unique (up to
			equivalence in the Cartier--Foata monoid) minimal monotone
			factorization of $x$, a representative $J_1,\dots,J_{r-1}$ of
			which is produced by the algorithm.

		We claim for $J_r=\bigcup_i C_i$ that $J=J_1,\dots,J_r$ is a
			minimal monotone factorization of $w$.
		Note that the claim of uniqueness will then follow from
			\cref{lem:maximal-strong-rdes-intervals,prop:monotone-structure},
			so since $x$ is $45312$-avoiding it suffices by
			\cref{cor:monotone-45312} to show that $J$ is strong bidescent.
		Since $J$ is strong right-descent by construction, we need only
			show that it is strong left-descent, for which by
			\cref{cor:strong-rdes-char} it is sufficient to show that that
			$e=f+\vec u_i$ implies $f=e-\vec u_i$ for $e,f\in H(J)$, where
			since $\upslice J{r-1}$ is strong bidescent by the inductive
			hypothesis this holds by \cref{cor:strong-bidescent-char} when
			$\lvl e,\lvl f\leq r-1$.
		Since $e=f+\vec u_i$ only if $\lvl f<\lvl e$, it remains to
			consider the case $r=\lvl e>\lvl f$.

		Now, let $v\in W$ be such that $xv=w$, and suppose $e=f+\vec u_1$
			(the case $e=f+\vec u_2$ is symmetric).
		Then since $j_f$ is strong right-$x$-directed by
			\cref{lem:directed-up}, we have $v^{-1}(j_f)$ is
			right-$w$-directed, and thus by \cref{lem:4231-45312-ind} since
			$j_f$ is not $x$-uncrossed we have $j_f\in[i_e,i_e+j_e-l_e)$, and
			thus $x(i_e)>\cdots>x(j_f)$.
		If $i_e=j_f$ then $f=e-\vec u_1$ is by definition, and if $i_e<j_f$
			then by \cref{lem:maximal-strong-rdes-intervals} since $j_f\not\in
			[i_{f'},j_{f'}]$ for every $f'\in H(\upslice J{r-1})$ with
			$\lvl{f'}>\lvl f$ we have $[i_e,j_f]\subseteq[i_f,j_f]$ and
			$e=f-\vec u_1$ as desired.
	\end{proof}
\end{prop}

\begin{rmrk}
    The resolutions of cograssmanian permutations
		obtained by Zelevinski\u{\i} are constructed by an algorithm
		similar to that of \cref{prop:monotone-algorithm}, and it
		would be interesting to know if the outputs of that algorithm
		exhaust all possible Gelfand--MacPherson resolutions, as they
		do for our algorithm in the case of $4231$-$45312$-avoiding
		permutations. 
\end{rmrk}

\section{Tight factorizations of
	\pdfmath{4231}-\pdfmath{45312}-\pdfmath{45123}-\pdfmath{34512}-avoiding
	permutations}
\label{sec:height2}
In this last section, we apply the combinatorial criterion to show that
	$4231$-$45312$-$45123$-$34512$-avoiding permutations have exactly one
	minimal tight factorization up to equivalence.

To do so, we develop a bit more of the structure theory of heaps of
	connected components of strong bidescent factorizations, so that we
	may relate the avoidance of patterns of the form $m_2\ m_2+1\
	\dots\ m_2+m_1\ 1\ 2\ \dots\ m_2$ by a permutation to restrictions on
	the complexity of the heaps of connected components of minimal strong
	bidescent factorizations.
We then use the fact that $45123$-$34512$-avoidance corresponds to
	particularly simple heaps for which we can give an inequality relating
	defects and bounces.
\begin{lem}
\label{lem:min-strong-bides-is-conv}
	Let $c\in H(J)$ be a connected component of the minimal strong
		bidescent factorization $J_1,\dots,J_r$.
	Then if $i\in\set{1,2}$ is such that $e=c+\vec u_i
		+\vec u_{3-i}$ exists, then $c+\vec u_{3-i}+\vec u_i$
		exists and is also equal to $e$.
	\begin{proof}
		We consider the case $i=1$, the remaining case follows by symmetry.
		Let $g=c+\vec u_1$, so that by \cref{cor:strong-bidescent-char} we
			have $c=g-\vec u_1$ and $c=e-\vec u_2$.
		Now, by minimality of $J$, we have that $\max\set{i_e,i_c}<i_g$, so
			let $h\in H(J)$ have maximal $\lvl h$ such that
			$\lvl c\leq\lvl h<\lvl e$ and $i_g-1\in[i_h,j_h]$.
		Then by \cref{cor:strong-bides-structure} we have either
			$f=e-\vec u_1$ exists and $h=f-\vec u_2=g-\vec u_1=c$ (in which
			case we are done), or $h=f=e-\vec u_1$ exists, which we assume
			in the remainder of the proof.
		Since $g=c+\vec u_1\neq h=f$ we have that $c\neq h=f$ by
			\cref{cor:strong-bidescent-char}, and since $\lvl h>\lvl c$,
			by \cref{cor:strong-bides-no-intermediate} we have that
			$j_f=i_g-1$ and $c=f-\vec u_2$, from which we conclude by
			\cref{cor:strong-bidescent-char} as desired.
	\end{proof}
\end{lem}
\begin{rmrk}
	The special case of \cref{lem:min-strong-bides-is-conv} for
		strong bidescent factorizations induced by reduced words, i.e.\ for
		factorizations induced by $321$-avoiding permutations, is the
		\emph{lateral convexity} lemma \cite[Lemma~1]{Billey--Warrington} of
		Billey and Warrington.
\end{rmrk}
\begin{defn}
	We define $-\colon M_4\to M_4$ to be the monoid anti-involution
		induced by the map $\pm\vec u_i\mapsto\mp\vec u_i$ on
		$\set{+\vec u_1,+\vec u_2,-\vec u_1,-\vec u_2}$.
	We define the map $\iota:M_4\to\ZZ(1,1)\oplus \ZZ(-1,1)$ induced
		by $\pm u_1\mapsto \pm (1,1)$ and $\pm u_2\mapsto \pm (-1,1)$.
	We define for each pair $c,e\in H(J)$ the \define{set of paths
		from $c$ to $e$} to be $L(c,e)=\set{\vec v\in M_4:c+\vec v=e}$.
\end{defn}

\begin{prop}
\label{prop:lattice-realization}
	Let $J_1,\dots,J_r$ be minimal strong bidescent factorization.
	Then there is an injection $\iota:H(J)\to\ZZ^2$ such that
		for every $c,e\in H(J)$ and $\vec v\in L(c,e)$ we have
		$\iota(e)=\iota(c)+\iota(\vec v)$, and furthermore such
		that $\iota(e)=\iota(c) +\qy[\big]{(-1)^{1-i},1}$ implies
		$e=c+\vec u_i$.
	\begin{proof}
		We first show the claim without the furthermore restriction.
		For this, t suffices (by a greedy construction) to show that
			for every $c,e\in H(J)$ and $\vec v_1,\vec v_2\in L(c,e)$ we have
			$\iota(\vec v_1)=\iota(\vec v_2)$.
		Since $J$ is strong bidescent, we have $-\vec v_2\in L(e,c)$
			by \cref{cor:strong-bidescent-char} so that $\vec v_1
				-\vec v_2\in L(c,c)$, so equivalently we need
			to show $\iota(\vec v)=(0,0)$ for every $0\neq
				\vec v\in L(c,c)$.

		We first consider the case where $\vec v$ does not contain both
			$\vec u_i$ and $-\vec u_i$ for any $i\in\set{1,2}$.
		Then since $e=c+\vec u_i$ implies $\lvl e>\lvl c$, we have that
			$\vec v$ consists only of elements of $\set{\vec u_1,-\vec u_2}$
			or only of $\set{\vec u_2, -\vec u_1}$.
		However, since $J$ is minimal, we have in the first case that
			$i_{c+\vec v}>i_c$, and in the second case that $i_{c+\vec v}<i_c$,
			so that $c+\vec v\neq c$.

		Now suppose $\vec v$ does not contain both $\vec u_i$ and $-\vec u_i$
			for some $i\in\set{1,2}$.
		Consider a decomposition $\vec v=\vec v_1+\vec v_2+\vec v_3$ such
			that $\vec v_2$ contains both of $\set{\pm\vec u_i}$ for some
			$i\in\set{1,2}$ and is minimal with this property.
		Then $\vec v_2=\pm\vec u_i$ and $\vec{v_2'}\mp\vec u_i$ where
			$\vec{v_2'}=m\vec u_{3-i}$ for some $m\in\ZZ$, so that by
			\cref{lem:min-strong-bides-is-conv} we have
			that $e=c+\vec v_1+\vec{v_2'}+\vec v_3$, i.e.~
			$\vec{v'}=\vec v_1+\vec{v_2'}+\vec v_3\in L(c,c)$
			has $\iota(\vec{v'})=\iota(\vec v)$ and is of
			shorter length.
		By iterating this reduction, we get the claim by reduction
			to the earlier case.

		We now show the existence of $\iota$ also satisfying the furthermore
			restriction.
		By \cref{cor:strong-bidescent-char}, the relation $\sim$ on $H(J)$
			defined by $c\sim e$ if $L(c,e)\neq\emptyset$ is an equivalence
			relation, and furthermore the equivalence classes $C_1,\dots,C_m$
			are themselves strong bidescent factorizations.
		The reduction in the previous paragraph shows that every injection
			$\iota_i:C_i\to\ZZ^2$ such that $\iota_i(c)
				=\iota_i(e)+\iota(\vec v)$ for $\vec v\in L(c,e)$ satisfies the
			furthermore restriction, from which the claim follows by taking
			$\iota$ as a disjoint union of the $\iota_i$.
	\end{proof}
\end{prop}
\begin{lem}
\label{lem:diamond-pattern}
	Let $J_1,\dots,J_r$ be a minimal strong bidescent factorization of
		a permutation $w\in W$.
	Then for every $m_1,m_2\in\ZZ_{>0}$ and $e\in H(J)$,
		the following are equivalent:
	\begin{enumerate}
		\item $w$ realizes the pattern $m_2+1,m_2+2,\cdots,m_2+m_1,
			1,2,\cdots,m_2$ in positions $p_{m_2+1}<\cdots<p_{m_2+m_1}
			<p_1<\cdots<p_{m_2}$ such that $e$ has minimal $\lvl e$ for which
			exist distinct $k,k'\in\set{1,2,\dots,m_1+m_2}$ for which
			$\wlvl e^{-1}(w(p_k)),\wlvl e^{-1}(w(p_{k'}))\in[i_e,j_e]$.
		\item $e+(m_1-1)\vec u_1+(m_2-1)\vec u_2$ exists.
		\item $e+(m_1-1)\vec u_1+(m_2-1)\vec u_2$ exists, and
			$w$ realizes the pattern $m_2+1,m_2+2,\cdots,m_2+m_1,
			1,2,\cdots,m_2$ in positions $p_{m_2+1}<\cdots<p_{m_2+m_1}
			<p_1<\cdots<p_{m_2}$ such that
			$w(p_{m_2+1+k})=\wlvl{e+k\vec u_1}(i_{e+k\vec u_1})$
			for every $k\in\set{0,\dots,m_1-1}$ and
			such that $w(p_{1+k})=\wlvl{e+k\vec u_2}(j_{e+k\vec u_2})$
			for every $k\in\set{0,\dots,m_2-1}$.
	\end{enumerate}
	In particular, $w$ realizes the pattern $m_2+1,m_2+2,\cdots,m_2+m_1,
		1,2,\cdots,m_2$ if and only if there exists $e\in H(J)$ such that
		$e+(m_1-1)\vec u_1+(m_2-1)\vec u_2$ exists.
	\begin{proof}
		Note (3) implies (1) by \cref{cor:strong-bidescent-char}.
		We prove the implications (1) implies (2) implies (3) by induction
			on $m_1+m_2$, where the base case $m_1+m_2=2$ (i.e.\ $m_1=1=m_2$)
			is immediate since then (3) holds for every $e\in H(J)$.

		For the inductive step, we first suppose that (1) holds.
		By minimality of $\lvl e$, we have that $\wlvl[-1]e^{-1}(p_1)
			<\cdots<\wlvl[-1]e^{-1}(p_{m_1+m-2})$, and thus since $w$
			realizes the given pattern in the given positions we have
			that $\set{k\leq i\leq m_1+m_2 : \wlvl e^{-1}(p_k)\in[i_e,j_e]}
			=\set{m_2+1,m_2}$.
		Now let $f\in H(J)$ have minimal $\lvl f>\lvl e$ such
			there exist distinct $k,k'\in\set{1,2,\dots,m_1+m_2}$ with
			that $\wlvl f^{-1}(w(p_k)),\wlvl f^{-1}(w(p_{k'}))
			\in[i_f,j_f]$.
		By an analogous argument, we have that
			$\set{1\leq i\leq m_1+m_2 : \wlvl f^{-1}(p_k)\in[i_f,j_f]}$
			is either $\set{m_2+2,m_2}$ or $\set{m_2+1,m_2-1}$,
			assume without loss of generality (by symmetry)
			that it is the former.
		Note that $w$ realizes the pattern $m_2+1,m_2+m_1-1,1,2,\dots,m_2$
			in positions $p_{m_2+2}<\cdots<p_{m_2+m_1}<p_1<\cdots<p_{m_2}$
			and $f$ realizes (1) for these positions, so by the inductive
			hypothesis we have that $f+(m_1-2)\vec u_1+(m_2-1)\vec u_2$
			exists.
		Since $\wlvl e^{-1}(w(p_{m_2}))\in[i_e,j_e]$ and
			$\wlvl f^{-1}(w(p_{m_2}))\in[i_f,j_f]$, we have by
			\cref{cor:strong-bides-structure,lem:min-strong-bides-is-conv}
			that there are $a,b\geq 0$ such that $f=e+a\vec u_1+b\vec u_2$.
		But since $w(p_{m_2+1})>w(p_{m_2})$, we have that
			$a=0$ would imply $\wlvl f^{-1}(w(p_{m_2+1})),
				\wlvl f^{-1}(w(p_{m_2+2})) \in[i_f,j_f]$,
			contradicting the fact that $w(p_{m_2+1})<w(p_{m_2+2})$, so that
			$a>0$ and thus by \cref{lem:min-strong-bides-is-conv} we have that
			$e+(m_1-1)\vec u_1+(m_2-1)\vec u_2$ exists as desired.

		Now suppose that (2) holds, and assume without loss of generality
			(by symmetry) that $m_1\geq m_2$, and so in particular $m_1\geq 2$.
		Then $f=e+\vec u_1$ is such that $f+((m_1-1)-1)\vec u_1
			+ (m_2-1)\vec u_2$ exists, so that by the inductive
			hypothesis we have that $w$ realizes the pattern
			$m_2+1,\dots,m_2+m_1-1,1,2,\dots,m_2$ in positions
			$p_{m_2+2}<\cdots<p_{m_2+m_1}<p_1<\cdots,p_{m_2}$
			(note the indexing)
			with $w(p_{m_2+2+k})=\wlvl{f+(k\vec u_1}(i_{f+k\vec u_1})$
			for every $k\in\set{1,\dots,m_1-1}$ and such that
			$w(p_k)=\wlvl{f+k\vec u_2}(j_{f+k\vec u_2})$ for
			every $k\in\set{1,\dots,m_2}$.
		Now, let $p_{m_2+1}=w^{-1}(\wlvl e(i_e))$.
		By \cref{cor:strong-bidescent-char} we have that
			$p_{m_2+1}<p_{m_2+2}$, and also that $\wlvl e(j_e)=\wlvl f(j_f)$
			so we have $w(p_{m_2+1})>w(p_{m_2})$.
		Finally, since $J$ is strong right-descent, we have by
			\cref{lem:directed-up} that $\wlvl f(i_f)<\wlvl e(i_e)$, and thus
			we have (3) as desired.
	\end{proof}
\end{lem}

In order to apply the combinatorial criterion, it will be convenient to
	overload the notation concerning defects, crossings, and bounces to
	refer to elements of the heap of connected components. 
\begin{defn}
	Let $J_1,\dots,J_r$ be a factorization of a permutation $w$.
	Then for a mask $\sigma\to_Jx$ and $e\in H(J)$, define
		$\rmeet\sigma e=\set{t\in T:\lvl e\in \rmeet\sigma t
			\land t^\sigma_{\lvl e-1}\in[i_e,j_e]}$,
		$\rdef\sigma e=\rmeet\sigma e\cap\set{t\in T:\lvl e\in
			\rdef\sigma t}$,
		$\rcross\sigma e=\rmeet\sigma e\cap\set{t\in T:\lvl e\in
			\rcross\sigma t}$,
		and $\rbounce\sigma e=\rmeet\sigma e\cap\set{t\in T:
			\lvl e\in\rbounce\sigma t}$.
	For a non-existent component $c$ of $H(J)$, we write
		$\rmeet\sigma c=\rdef\sigma c=\rcross\sigma c=\rbounce\sigma c
		=\emptyset$.

	Similarly, for a position $p$ and mask $\sigma\to_Jx$ define
		$p^\sigma_k=(\sigma_1\cdots\sigma_k)^{-1}(p)$ the position
		at which $\sigma$ places $p$ at level $k$.
	Finally, for $e\in H(J)$ and mask $\sigma\to_Jx$ write
		$P^\sigma_R(e)=\set{p:p^\sigma_{\lvl e}\in[i_e,j_e]}$ the
		set of right-labels of strings $\sigma$ sends through~$e$.
\end{defn}

\begin{lem}
\label{lem:defect-implies-square}
	Let $J_1,\dots,J_r$ be a minimal strong bidescent factorization.
	Then for every $\sigma\to_Jx$ and $e\in H(J)$, and $t\in T$, we have
		that $t\in \rdef\sigma e$ only if there are $m_1,m_2\in\ZZ_{>0}$
		such that $c=e-m_1\vec u_1-m_2\vec u_2$ exists and
		$t\in\rmeet\sigma c$.
	In particular, $\rdef\sigma e\neq\emptyset$ only if
		$e-\vec u_1-\vec u_2$ exists.
	\begin{proof}
		Suppose that there is $t\in\rdef\sigma e$.
		Then there exists $c\in H(J)$ with minimal $\lvl c$ such
			that $t\in\rmeet\sigma c$, and this $c$ has
			$\lvl c<\lvl e$.
		By \cref{lem:min-strong-bides-is-conv} there are
			$m_1,m_2\in\ZZ_{\geq 0}$ for which $e=c+m_1\vec u_2
				+m_2\vec u_2$.
		We claim that $m_i>0$ for $i\in\set{1,2}$: if $m_i=0$,
			then we must have $t^\sigma_{\lvl{c+k\vec u_{3-i}}-1}
				\in[i_{\lvl{c+k\vec u_{3-i}}},j_{\lvl{c+k\vec u_{3-i}}}]$
			for every $0\leq k\leq m_{3-i}$, and thus
			$t\not\in\rmeet\sigma e$, contradiction.
		The supplemental claim follows from
			\cref{lem:min-strong-bides-is-conv}.
	\end{proof}
\end{lem}

\begin{lem}
\label{lem:defect-square-prod-bounce}
	Let $J_1,\dots,J_r$ be a minimal monotone factorization,
		and $e\in H(J)$ be such that $c=e-\vec u_1-\vec u_2$
		does not exist or $c$ exists but $c-\vec u_1$ and
		$c-\vec u_2$ do not.
	Then for every $\sigma\to_Jx$, we have that 
		$\abs{\rdef\sigma e}^2 \leq\abs{\rbounce\sigma{e-\vec u_1}}
			\cdot\abs{\rbounce\sigma{e-\vec u_2}}$.
	\begin{proof}
		If $\rdef\sigma e=\emptyset$ there is nothing to prove, so
			assume $\rdef\sigma e\neq\emptyset$.
		Then by \cref{lem:defect-implies-square} we have that
			$c=e-\vec u_1-\vec u_2$ exists, and since by assumption that
			$c-\vec u_1$ and $c-\vec u_2$ do not exist we have by
			\cref{prop:monotone-structure}(3) that every $p\in [i_c,j_c]$ is
			$\wlvl[-1]c$-uncrossed.
		In particular, by \cref{lem:defect-implies-square}
			this implies that $\rdef\sigma e
				\subseteq\rcross\sigma c
				\subseteq\rmeet\sigma c$.

		Write $f=c+\vec u_2=e-\vec u_1$ and $g=c+\vec u_1=e-\vec u_2$,
			and define the sets
			$L_\sigma=P^\sigma_R(f)\cap P^\sigma_R(e)\cap P^\sigma_R(c)$
			and
			$R_\sigma=P^\sigma_R(g)\cap P^\sigma_R(e)\cap P^\sigma_R(c)$.
		Then \cref{lem:defect-implies-square,prop:monotone-structure}
			imply that $p<q$ are such that $(p,q)\in\rdef
				\sigma e$ only if $p\in L_\sigma$ and $q\in R_\sigma$, so that
			$\abs{\rdef\sigma e}\leq \abs{L_\sigma}\cdot\abs{R_\sigma}$.
		We now claim that
			$\abs{\rbounce\sigma f}\geq\abs{L_\sigma}^2$
			and $\abs{\rbounce\sigma g}\geq\abs{R_\sigma}^2$.
		We prove the claim about $f$, the claim about $g$ follows by
			symmetry.

		Write $L'_\sigma=P^\sigma_R(f)\setminus P^\sigma_R(e)
			\setminus P^\sigma_R(c)$.
		Then $L_\sigma\times L'_\sigma\subseteq\rbounce\sigma f$ since
			$c-\vec u_1$ does not exist, so it suffices to show that
			$\abs{L'_\sigma}\geq\abs{L_\sigma}$.
		We have
		\begin{align*}
			\abs{L'_\sigma}
			&=
			\abs{P^\sigma_R(f)\setminus P^\sigma_R(e)\setminus P^\sigma_R(c)}
			\\
			&=
			\abs{P^\sigma_R(f)\setminus P^\sigma_R(e)}
				-\abs*{\qy[\big]{P^\sigma_R(f)\setminus P^\sigma_R(e)}	
					\cap P^\sigma_R(c)}\\
			&=
			\abs{P^\sigma_R(f)\setminus P^\sigma_R(e)}
				-\abs{P^\sigma_R(f)\cap P^\sigma_R(c)}
				+\abs{P^\sigma_R(f)\cap P^\sigma_R(e)\cap P^\sigma_R(c)}\\
			&=
			\abs{P^{\omega_R(J)}_R(f)\setminus P^{\omega_R(J)}_R(e)}
				-\abs{P^{\omega_R(J)}_R(f)\cap P^{\omega_R(J)}_R(c)}
				+\abs{L_\sigma}\\
			&\geq \abs{L_\sigma}
		\end{align*}
		where the last equality is because the number of strings entering $f$
			and not going to $e$ (resp.~coming from $c$) does not depend
			on the mask $\sigma$, and the last inequality is because
			$P^{\omega_R(J)}_R(f)\cap P^{\omega_R(J)}_R(c)
				\subseteq P^{\omega_R(J)}_R(f)\setminus P^{\omega_R(J)}_R(e)$
			by \cref{prop:monotone-structure}(1).
	\end{proof}
\end{lem}

\begin{thm}
\label{thm:4}
	Let $w$ be a $4231$-$45312$-$45123$-$34512$-avoiding permutation.
	Then the unique minimal (up to equivalence in the Cartier--Foata
		monoid) monotone factorization $J_1,\dots,J_r$ of $w$ is tight.
	\begin{proof}
		Let $\sigma\to_J x<w$ be a mask on $J$ such that $d_R(\sigma)>0$.
		We claim that $d_R(\sigma)<\sum_{t\in T}\abs{\rbounce\sigma t
			\setminus\rdef\sigma t}$.
		We have by definition that $d_R(\sigma)
			=\sum_{e\in H(J)}\abs{\rdef\sigma e}$ and
			$\sum_{t\in T}\abs{\rbounce\sigma t
				\setminus\rdef\sigma t}
				=\sum_{e\in H(J)}\abs{\rbounce\sigma e
					\setminus\rdef\sigma e}$.

		By \cref{lem:diamond-pattern} and assumption on $w$ we have for
			every $e\in H(J)$ that $e-\vec u_1-\vec u_2-\vec u_i$ and
			$e-\vec u_i-\vec u_1-\vec u_2$ do not exist, so
			\cref{lem:defect-square-prod-bounce} and the AM--GM inequality
			imply
			\[
				\rdef\sigma e
					\leq\sqrt{\abs{\rbounce\sigma{e-\vec u_1}}
					\cdot\abs{\rbounce\sigma{e-\vec u_2}}}
					\leq\frac12\qy[\big]{\abs{\rbounce\sigma{e-\vec u_1}}
						+\abs{\rbounce\sigma{e-\vec u_2}}}
			\]
			for each $e\in H(J)$, and \cref{lem:defect-implies-square}
			implies that $\rbounce\sigma{e-\vec u_i}
				=\rbounce\sigma{e-\vec u_i}\setminus\rdef\sigma{e-\vec u_i}$.
		Note that by \cref{cor:strong-bidescent-char} we have
			for each $f\in H(J)$ and $i\in\set{1,2}$ that there is
			$e\in H(J)$ with $\rdef\sigma e\neq\emptyset$ and
			$f=e-\vec u_i$ if and only if $e=f+\vec u_i$ and
				$\rdef\sigma{f+\vec u_i}\neq\emptyset$.
		Thus, we can re-index a sum over all $e-\vec u_i$ to get
			\[
				d_R(\sigma)\leq\sum_{e\in H(J)}
					\abs{\rbounce\sigma e \setminus\rdef\sigma e}
					\cdot\frac12\abs*{
						\set{i\in\set{1,2}:\rdef\sigma{e+\vec u_i}
							\neq\emptyset}}\,.
			\]
		Now, let $e\in H(J)$ have minimal $i_e$ such that
			$\rdef\sigma e\neq\emptyset$.
		Then by \cref{lem:defect-square-prod-bounce} we
			have $\abs{\rbounce\sigma{e-\vec u_1}
				\setminus\rdef\sigma{e-\vec u_1}}>0$
			and $\abs*{\set{i\in\set{1,2}:\rdef\sigma{e-\vec u_1
				+\vec u_i} \neq\emptyset}}=1$, so
			\begin{multline*}
				\sum_{e\in H(J)}
					\abs{\rbounce\sigma e \setminus\rdef\sigma e}
					\cdot\frac12\abs*{
						\set{i\in\set{1,2}:\rdef\sigma{e+\vec u_i}
							\neq\emptyset}}
					\\<
				\sum_{e\in H(J)}
					\abs{\rbounce\sigma e \setminus\rdef\sigma e}
					\cdot\frac22
				=\sum_{t\in T}\abs{\rbounce\sigma t
					\setminus\rdef\sigma t}
			\end{multline*}
			as desired.
	\end{proof}
\end{thm}

\section*{Acknowledgements}

The authors would like to thank Pavlo Pylavskyy, who directed them
	toward this question at the conclusion of their 2011 Summer REU
	project at the University of Minnesota, who suggested that the
	combinatorial criterion used by Billey and Warrington for their result
	should have an extension to the setting suggested by Skandera, and who
	specifically suggested the statement of the combinatorial criterion
	(supplemental claim of \cref{thm:combi}) for type $A$.
A few of the ideas in this paper were subsequently worked out by
	the authors under Pylyavskyy's supervision during the 2012 Summer REU.
	
Portions of this work were written as part of the second-named author's
    doctoral dissertation \cite{V21}, and he is grateful to Dima Arinkin
    for advising him and in particular for encouraging him to develop 
    the geometric aspects of this work, as described in \cref{sec:gm}.

The authors would also like to thank Ben Elias for pointing out the
	quotients of products of Kazhdan--Lusztig basis elements analyzed in
	this work are instances of composites in the Schur algebroid.

\bibliographystyle{amsalpha}
\bibliography{kl-bibliography}
\end{document}